\documentclass[11pt]{report}

\usepackage[utf8]{inputenc}
\usepackage[T1]{fontenc}
\usepackage{lmodern}

\usepackage{suthesis-2e}
\usepackage{bigfoot} %
\usepackage{tocbibind}
\usepackage[shortlabels]{enumitem}

\usepackage{booktabs}
\usepackage{subcaption}
\usepackage[section]{placeins}

\usepackage[cache=.,frozencache,chapter]{minted}
\setminted{fontsize=\small}
\usepackage{tcolorbox}

\usepackage{amsmath,amssymb,amsthm,mathtools}
\usepackage{bbm}
\usepackage{./math-macros}

\usepackage{tikz}
\usetikzlibrary{calc,positioning}
\usetikzlibrary{shapes.geometric,arrows.meta,decorations.markings}
\usetikzlibrary{cd}
\tikzcdset{arrow style=tikz}
\newcommand*{\ucorner}{\rotatebox{-45}{$\ulcorner$}}

\usepackage{etoolbox}
\usetikzlibrary{external}
\AtBeginEnvironment{tikzcd}{\tikzexternaldisable}
\AtEndEnvironment{tikzcd}{\tikzexternalenable}

\usepackage{hyperref}
\usepackage[capitalize,noabbrev]{cleveref}
\hypersetup{bookmarksnumbered,bookmarksopen,hidelinks}
\usepackage{doi}

\usepackage[symbols,abbreviations,sort=none]{glossaries-extra}
\renewcommand{\glossarysection}[2][]{} %

\usepackage[maxbibnames=10,style=alphabetic]{biblatex}
\renewbibmacro{in:}{}

\theoremstyle{definition}
\newtheorem{theorem}{Theorem}[section]
\newtheorem{corollary}[theorem]{Corollary}
\newtheorem{lemma}[theorem]{Lemma}
\newtheorem{proposition}[theorem]{Proposition}
\newtheorem{definition}[theorem]{Definition}
\newtheorem{example}[theorem]{Example}
\newtheorem{remark}[theorem]{Remark}

\title{The Algebra and Machine Representation \\ of Statistical Models}
\author{Evan Patterson}

\titlegraphic{%
  \tikzexternaldisable%
  \input{wiring-diagrams/zoo-statistics/glm-n-sampling}%
  \tikzexternalenable}

\dept{Statistics}
\principaladviser{Emmanuel Cand{\`e}s}
\firstreader{John Chambers}
\secondreader{Trevor Hastie}
\thirdreader{Mark Musen}

\submitdate{June 2020}
\copyrightyear{2020}
\figurespagefalse
\tablespagefalse

\begin{document}

\newcommand*{\symbolentry}[3]{%
  \newglossaryentry{#1}{name={#2},description={#3},type={symbols}}}

\symbolentry{ObSet}{$\Ob{\cat{C}}$}{collection of all objects in category
$\cat{C}$}

\symbolentry{HomSet}{$\cat{C}(x,y)$}{collection of all morphisms $x \to y$ in
category $\cat{C}$}

\symbolentry{Simplex}{$\Delta^d$}{standard $d$-simplex in $\R^{d+1}$}

\symbolentry{Copy}{$\Copy_x$}{copying, or duplication, morphism $x \to x \otimes
x$}

\symbolentry{Delete}{$\Delete_x$}{deleting, or discarding, morphism $x \to I$}

\symbolentry{Expectation}{$\E(X)$}{expectation of random variable $X$}

\symbolentry{Field}{$\K$}{a generic field, usually $\R$}

\symbolentry{Model}{$\Model{\cat{T}}$}{category of models of theory $\cat{T}$}

\symbolentry{Naturals}{$\N$}{rig of natural numbers, $\{0, 1, \dots\}$}

\symbolentry{Normal}{$\Normal$}{family of univariate normal distributions}

\symbolentry{NormalAlt}{$\Normal_d$}{family of $d$-dimensional normal
distributions}

\symbolentry{ProbSpace}{$\ProbSpace(\sspace{X})$}{space of probability measures
on measurable space $\sspace{X}$}

\symbolentry{Rig}{$R$}{a generic commutative rig (``ring without negatives''),
usually $\R$ or $\R_+$}

\symbolentry{Reals}{$\R$}{field of real numbers}

\symbolentry{NonnegativeReals}{$\R_+$}{rig of nonnegative real numbers}

\symbolentry{PosReals}{$\PosR$}{multiplicative group of positive real numbers}

\symbolentry{SymMat}{$\SymMat^d$}{space of $d \times d$ symmetric real matrices}

\symbolentry{PSD}{$\PSD^d$}{cone of $d \times d$ positive semi-definite real
matrices}

\symbolentry{Braid}{$\Braid_{x,y}$}{braiding morphism $x \otimes y \to y \otimes
x$}

\symbolentry{Theory}{$\Theory{\cat{C}}$}{theory (a category) corresponding to
category $\cat{C}$}

\symbolentry{Variance}{$\Var(X)$}{variance of random variable $X$}

\symbolentry{Integers}{$\Z$}{ring of integers}

\newcommand*{\catentry}[3]{\newabbreviation{#1}{#2}{#3}}

\catentry{Ab}{$\Ab$}{abelian groups and group homomorphisms}
\catentry{Aff}{$\Aff_\K$}{affine spaces over $\K$ and affine maps}
\catentry{Bimon}{$(\CAT{C})\Bimon$}{(bicommutative) bimonoids and bimonoid
homomorphisms}
\catentry{Cart}{$\Cart$}{cartesian monoidal categories and cartesian functors}
\catentry{Cat}{$\Cat$}{categories and functors}
\catentry{CCC}{$\CCC$}{cartesian closed categories and cartesian closed
functors}
\catentry{Comon}{$(\CAT{C})\Comon$}{(commutative) comonoids and comonoid
homomorphisms}
\catentry{Cone}{$\Cone$}{conical spaces and conic-linear maps}
\catentry{Conv}{$\Conv$}{convex spaces and convex-linear maps}
\catentry{Graph}{$\Graph$}{graphs and graph homomorphisms}
\catentry{Group}{$\Group$}{groups and group homomorphisms}
\catentry{Lawvere}{$\Lawvere$}{Lawvere theories and strict cartesian functors
preserving object generators}
\catentry{Markov}{$\Markov$}{Polish measurable spaces and Markov kernels}
\catentry{Meas}{$\Meas$}{Polish measurable spaces and measurable maps}
\catentry{Mod}{$\Mod_R$}{$R$-modules and $R$-linear maps}
\catentry{Mon}{$(\CAT{C})\Mon$}{(commutative) monoids and monoid homomorphisms}
\catentry{MonCat}{$(\CAT{S})\MonCat$}{(symmetric) monoidal categories and strong
(symmetric) monoidal functors}
\catentry{PRO}{$\PRO(\CAT{P})$}{PRO(P)s and strict (symmetric) monoidal functors
preserving object generators}
\catentry{Poset}{$\Poset$}{partially ordered sets (posets) and monotone maps}
\catentry{Preord}{$\Preord$}{preordered sets (prosets) and monotone maps}
\catentry{Rel}{$\Rel$}{sets and binary relations}
\catentry{Set}{$\Set$}{sets and functions}
\catentry{Stat}{$\Stat$}{vector space subsets and Markov kernels}
\catentry{Top}{$\Top$}{topological spaces and continuous maps}
\catentry{Vect}{$\Vect_\K$}{vector spaces over $\K$ and linear maps}

\beforepreface

\prefacesection{Preface}

As the twin movements of open science and open source bring an ever greater
share of the scientific process into the digital realm, new opportunities arise
for the meta-scientific study of science itself, including of data science and
statistics. Future science will likely see machines play an active role in
processing, organizing, and perhaps even creating scientific knowledge. To make
this possible, large engineering efforts must be undertaken to transform
scientific artifacts into useful computational resources, and conceptual
advances must be made in the organization of scientific theories, models,
experiments, and data.

This dissertation takes steps toward digitizing and systematizing two major
artifacts of data science, statistical models and data analyses. Using tools
from algebra, particularly categorical logic, a precise analogy is drawn between
models in statistics and logic, enabling statistical models to be seen as models
of theories, in the logical sense. Statistical theories, being algebraic
structures, are amenable to machine representation and are equipped with
morphisms that formalize the relations between different statistical methods.
Turning from mathematics to engineering, a software system for creating machine
representations of data analyses, in the form of Python or R programs, is
designed and implemented. The representations aim to capture the semantics of
data analyses, independent of the programming language and libraries in which
they are implemented.

\subsection*{Guide to reading}

Data science and statistics, category theory and categorical logic, program
analysis and programming language theory all interact in this dissertation.
Individually these are large fields, but the group of people conversant in them
all cannot be very large. In the hope of expanding the readership, I have
exposited much that is already known in addition to my own contributions. So, if
the expert in some area is surprised to find me explaining something they know
very well, I hope that they will forgive me on the grounds that it is part of
the design.

After a general introduction in the first chapter, the dissertation proceeds in
two fairly independent directions. In
\cref{ch:algebra-statistics,ch:zoo-statistics}, I develop the algebra of
statistical theories, statistical models, and their morphisms. Knowledge of
probability theory and exposure to abstract algebra are assumed. Measure theory
is used, but more as an unifying language for probability than as a substantive
mathematical theory. The necessary background in category theory is presented in
\cref{ch:category-theory}. Readers familiar with category theory may omit this
chapter, referring to it as needed for the less standard definitions and
examples. An impressionistic view can be formed by reading the Introduction and
then perusing the many examples of \cref{ch:zoo-statistics}, with the
understanding that the mathematical underpinnings are supplied by
\cref{ch:category-theory,ch:algebra-statistics}.

In the second major part, I describe the design and implementation of a software
system for creating semantic representations of data science workflows.
\cref{ch:program-analysis}, on the analysis of data science code, assumes a
working knowledge of the programming languages Python and R. In
\cref{ch:semantic-enrichment}, I present the Data Science Ontology and a
procedure for the semantic enrichment of idealized computer programs. The formal
description of the ontology and procedure is category-theoretic, for which the
necessary background is again contained in \cref{ch:category-theory}. Although
they describe components of a unified software system,
\cref{ch:program-analysis,ch:semantic-enrichment} can be read independently of
each other. The concluding \cref{ch:conclusion} describes limitations of the
work, suggests directions for future work, and offers a general outlook on how
the structuralist approach to data analysis might transform the scientific
process.

Pictorially, the dependencies between the chapters are:

\begin{center}
  \begin{tikzcd}[math mode=false]
    \cref{ch:introduction} \arrow[r] \arrow[dr] \arrow[ddr]
      & \cref{ch:algebra-statistics} \arrow[r]
      & \cref{ch:zoo-statistics} \\
    \cref{ch:category-theory} \arrow[ur] \arrow[dr]
      & \cref{ch:program-analysis} \\
      & \cref{ch:semantic-enrichment}
  \end{tikzcd}
\end{center}

\subsection*{Acknowledgments}

Over the course of this project I have benefited from the advice and support of
many people. I am deeply grateful to my advisor, Emmanuel Cand{\`e}s, for
sharing his broad knowledge of statistics and applied mathematics and for
consistently supporting my research explorations, no matter where they led me.
Through his commitment to careful, intellectually serious science, he has been
and remains a great role model to me . I also thank the other members of my
thesis committee, namely John Chambers, Trevor Hastie, Mykel Kochenderfer, and
Mark Musen, for their advice. I am especially thankful to John Chambers for his
enthusiastic support of both parts of this project and for his helpful comments
on the manuscript.

The semantic modeling of data science code began during my time at IBM Research,
where I was fortunate to collaborate with Ioana Baldini, Aleksandra Mojsilovic,
and Kush Varshney. I thank them for their mentorship and for introducing me to a
line of research that I likely would not have otherwise pursued. This part of my
thesis has also benefited from conversations with Roxana Danger, Julian Dolby,
Yaoli Mao, Robert McBurney, Holie Schmidt, and Gustavo Stolovitzky.

The algebra of statistical theories and models grew out of interactions with the
field of applied category theory. I am grateful to John Baez and David Spivak
for encouraging me to explore this area at a time when I knew very little about
it. I thank James Fairbanks, Brendan Fong, Jade Master, and David Spivak for
helpful conversations related to this part of the thesis. I also thank
Arquimedes Canedo for the opportunity to work on applied category theory at
Siemens.

As a student in the Statistics Department at Stanford I benefited from
interactions with many excellent peers and colleagues. From my original student
cohort, I am especially grateful for time spent with my friends Paulo Orenstein
and Feng Ruan. From the Cand{\`e}s research group, I thank Stephen Bates, Qijia
Jiang, and Yaniv Romano for conversations and collaboration. I also thank the
department administrative staff, particularly Susie Ementon, for their patience
and willingness to help.

Finally, I wish to express my deep appreciation for the friends and family who
have supported me over a period of time that has not always been simple or easy.
Especially to my mother, Morgan; my brother, Wes; my uncle, Brian; and my
significant other, Julia, I am ever grateful for your unconditional love and
encouragement.

\prefacesection{List of Notation}

\subsection*{Commonly used symbols}

\printunsrtglossary[type=symbols,style=long]

\subsection*{Named categories}

\printunsrtglossary[type=abbreviations,style=long]

\afterpreface

\chapter{Introduction}
\label{ch:introduction}

In 1988, retired physicist and information scientist Don Swanson announced his
discovery of an intriguing possible connection between migraine headaches and
magnesium deficiency. Although Swanson had no formal medical training and had
conducted no experiments, his argument was found compelling enough by the
medical community to warrant further study. The argument follows a simple
schema. Due to the vastness of the scientific literature, and its tendency to
cluster along disciplinary lines, there surely exist causal connections between
factors $A$ and $C$ that are not explicitly documented in the literature, but
that are still \emph{implicitly} present in it through mediating factors $B$
with documented links between $A$ and $B$, and $B$ and $C$. In present case,
where $A$ is magnesium deficiency and $C$ is migraine headaches, Swanson found
little direct literature overlap but no less than eleven possible mediators $B$,
including calcium channel blockers, epilepsy, and platelet aggregation.
Subsequent confirmatory studies reported a statistically and practically
significant effect of magnesium in reducing migraine headaches.

This unorthodox research program was crucially aided by the online database
Medline, which provides bibliographic information for the life sciences and
biomedicine. The decades following Swanson's discovery have seen large growth in
the number and size of scientific publication databases, accompanied by the
creation of aggregators like Web of Science and Google Scholar and the adoption
of preprint servers like arXiv and bioRxiv. As the open science movement gathers
momentum, the prospect of a world where most or all scientific publications are
open and available online begins to resemble reality. Nevertheless,
meta-scientific research in the style of Swanson remains a niche activity that
has not dramatically accelerated scientific discovery. It is instructive to ask
why this is.

One reason is that a fishing expedition through the scientific literature may
easily yield false discoveries. Swanson mitigates the statistical problems of
multiple comparisons and spurious correlations by searching for evidence that,
once noticed, appears overwhelming. The migraine-magnesium link was supported by
not one but eleven distinct mediating correlations, reducing the chances of
finding a weak or nonexistent direct correlation or of confusing a purely
correlational relationship with a causal one. In effect, Swanson targets a
regime where the signal so overwhelms the noise that statistics is superfluous.
Considering the breadth of science and practical necessity for scientists to
specialize, it seems likely that many other potential discoveries are ``hiding
in plain sight.''

So, without diminishing the serious statistical problems that it brings, the
challenge facing meta-scientific research today is more basic. It is that any
systematic processing of the scientific literature involving logical inference,
no matter how trivial, is extremely difficult. Swanson conducted his searches
manually, aided only by the Medline search engine. This methodology is
inherently unscalable, and few people possess the patience and good scientific
judgment required to carry it out effectively, even in a limited domain. The
need for machine assistance is clear.

In order to create useful software to support meta-scientific research,
scientific knowledge must somehow be encoded on a computer. This is challenging
because the primary bearer of scientific knowledge---the scientific paper---is
fundamentally a human artifact, written by humans to be read by humans. Two
broad approaches to the computational representation of scientific knowledge,
neither mutually exclusive nor exhaustive, can be distinguished. The first is to
extract information algorithmically from scientific papers, textbooks, and other
natural language texts. This approach has the advantage of being highly
scalable, but, despite significant recent advances in natural language
processing, extracting precise logical information from unstructured text
continues to be unreliable. An alternative is to directly specify scientific
knowledge in machine-interpretable form. Knowledge engineering tends to yield
greater precision at the expense of greater human effort and, in some cases,
increased brittleness in unforeseen situations.

Any method to digitize scientific knowledge includes, at least implicitly, an
understanding about what kind of information is to be represented. In fact,
knowledge about science takes many forms, ranging from experimental designs and
collected data, to data processing and analysis, to mathematical and statistical
models, and finally to the larger scientific theories and paradigms in which
models are embedded. These common elements interact with each other in complex
ways. Moreover, they manifest differently across different fields of the natural
and social sciences.

In short, the program of digitizing scientific knowledge is complex and
multi-faceted. It has inspired and involved workers from many fields, each with
their own perspectives. Being a dissertation in statistics, the present work
restricts itself to that part of the scientific process concerning statistical
modeling and data analysis. Certain other fields of science, such as genetics,
proteomics, and biomedicine, have seen large, organized efforts to record their
accumulated findings in machine-interpretable form. In statistics and machine
learning, there has been comparatively little such work, despite statistics
serving an important auxiliary role in scientific inference across the sciences.
As a result, an essential element of the scientific process remains largely
impervious to introspection by machines.

A good computational representation of data analysis would enable data
scientists to communicate, collaborate, and extend existing work more
efficiently. Today, data analysis is communicated almost entirely through
written descriptions in the methods and results sections of scientific reports.
While providing valuable opportunity to develop a narrative and motivate
analytical decisions, these descriptions often omit crucial details. Moreover,
statistical modeling has many degrees of freedom: it is rarely the case that a
single model presents itself as obviously superior to all others. A critical
reader may wish to test the selected model's assumptions or fit, extend the
model, or evaluate it against competing models. To do so, they must either
reproduce the original analysis from scratch, or else obtain and modify the
analysis source code. The former is often impossible, and the latter requires,
at minimum, that the code be understood well enough to use it. All of this would
be enormously simplified by a computational representation that allowed models
to be inspected and manipulated \emph{as data}.

At a larger scale, the ability to treat data analysis as data itself is clearly
essential to both meta-analysis and meta-learning. In meta-analysis, the aim is
to increase statistical precision and power by aggregating analyses from
independent studies; in meta-learning, it is to improve machine learning
algorithms by learning from previous experiments. Statistical analyses and
machine learning metadata, respectively, are thus treated as data themselves.
The inability to represent data analysis in a computationally useful way is a
major obstacle to automation in meta-analysis.

Even when restricted to its statistical aspects, the digitization of science
remains far too broad a project to be encompassed within a single work. This
dissertation investigates two specific topics in the representation of data
analysis. These topics share similar motivations, yet are quite distinct in
concept and execution. In part, the difference reflects that between the
relatively well-defined field of statistics and its younger and more amorphous
cousin, data science.

Although it has a long prehistory, statistics first established itself as an
independent field in the early twentieth century. Its central concepts are the
\emph{statistical model}, as a parametric family of probabilistic data
generating mechanisms, and \emph{statistical inference}, as the approximate
inversion of a statistical model to infer model parameters from observed data.
Statistical inference is often classified according to whether it is frequentist
or Bayesian and whether it concerns parameter estimation, hypothesis testing, or
prediction. While the methods of statistics have evolved to meet the
opportunities and challenges afforded by greater computational power and
availability of data, the conceptual core of the field has held fairly constant.

The discipline of data science emerged more recently, in response to several
developments. First, in the late 1980s and early 1990s, machine learning
coalesced out of computer science as a field sharing much of its subject matter
with statistics but having rather different culture and priorities, as well as
its own theoretical framework. So, at the very least, ``data science'' is an
umbrella term meant to encompass both statistics and machine learning. But data
science reflects a trend deeper than changing disciplinary boundaries, namely
that gains to computing and data acquisition technology have forced us to expand
our understanding of what data analysis is. No longer does data analysis consist
in fitting, more or less by hand, an analytically tractable model to a small
sample of data collected in a carefully designed experiment. Data is now often
gathered voraciously from heterogeneous sources under a policy of ``collect
first, ask questions later.'' Utilizing such data might require database queries
and integration, data cleaning and preprocessing, exploratory data analysis and
visualization, and large-scale and distributed computing. Statistical modeling
remains important, but it is now only one task among many. Data science embraces
all of these activities.

The two topics of this dissertation concern the representation of statistical
models and data science workflows, respectively. In the first case, the
representation is so far mathematical rather than computational, although the
mathematics is purposefully designed to admit a computer implementation. The
second topic has a stronger engineering focus and is accompanied by a prototype
implementation. Both are introduced in the sections below.

\section{The algebra of statistical models}
\label{sec:introduction-algebra-statistics}

In the practice of statistics, in contrast to much of the theory, it is rare to
posit a single fixed model, chosen before seeing the data and forbidden from
being revised afterwards. Instead, the statistician initially has in mind a
number of reasonable models and, after receiving the data, they will choose one
on the basis of simplicity, accuracy, goodness of fit, interpretability,
scientific plausibility, and other criteria. Applied statistics is thus a
complex process of exploration through a whole space of models.

The models under consideration are not selected at random but are related to
each other in meaningful ways. For example, one model might be related to
another by adding extra predictors, incorporating interactions or
nonlinearities, generalizing the sampling distribution to account for
overdispersion, adjusting the hyperparameters or priors, or altering the
internal components of the model, such as the link function in a generalized
linear model or the activation function in a neural network. At present, there
exists no formal or systematic language in which to express such relationships.
This dissertation aims to introduce one.

Although motivated by practical considerations, the problem of formalizing
relationships between statistical models has a philosophical background relevant
to the path taken here. Beginning with the logical positivists in the early
twentieth century, philosophers of science have tried to explicate how
scientific theories and models are related to theories and models in
mathematical logic. A bold step in this program was taken when Patrick Suppes
proposed that the word ``model'' has essentially the same meaning in science as
it does in logic \cite{suppes1961}. This does not imply that models have the
same \emph{uses} in logic and science, but rather that scientific models can be
reconstructed as models of logical theories in the sense of Tarski. For Suppes,
scientific theories are tested by experiments through a hierarchy of models,
descending from the general theory to models of experiments and models of data
\cite{suppes1966}. The latter can be regarded as statistical models.

Despite Suppes being, among other things, a statistician, later advocates of
what is now called the \emph{semantic view} of scientific theories have placed
little emphasis on statistics. Statisticians, for their part, generally eschew
philosophical questions about how statistical models relate to other kinds of
models. Consequently, the relationship between logical and statistical models
has never been properly clarified. So, in addition to its practical purpose,
this work has a conceptual or philosophical purpose in drawing an exact
connection between logical and statistical models.

The mathematics fitted for both purposes is algebra, particularly category
theory and categorical logic. Because any algebraic structure inherits from its
axiomatization a notion of structure-preserving map, or \emph{homomorphism},
reconstructing statistical models as algebraic structures would immediately
yield a notion of morphism between models. Morphisms would then formalize the
relationships between statistical models described earlier. But what kind of an
algebraic structure is a statistical model? Category theory provides the class
of compositional structures needed to answer this question. In fact, categorical
logic has already established a bridge between category theory and logic that
can serve as a template for the algebra of statistical models.

\begin{table}
  \centering
  \begin{tabular}{lll}
    \toprule
    Category theory & Mathematical logic & Statistics \\ \midrule
    Category & Theory & Statistical theory \\
    Functor & Model & Statistical model \\
    Natural transformation & Model homomorphism
                            & Morphism of statistical model \\
    \bottomrule
  \end{tabular}
  \caption{Dictionary between category theory, logic, and statistics}
  \label{tbl:category-logic-statistics}
\end{table}

\cref{tbl:category-logic-statistics} states a dictionary translating between
concepts of category theory, mathematical logic, and statistics. These
correspondences are now described more carefully, beginning with that between
category theory and logic.

Categorical logic originates with a formal analogy between natural
transformations in category theory and model homomorphisms in logic. Given two
functors $F,G: \cat{C} \to \Set$ from a category $\cat{C}$ into the category of
sets, the equations
\begin{equation*}
  \begin{tikzcd}
    Fx \arrow[r, "\alpha_x"] \arrow[d, "Ff"'] & Gx \arrow[d, "Gf"] \\
    Fy \arrow[r, "\alpha_y"] & Gy
  \end{tikzcd},
  \qquad \forall f: x \to y \text{ in } \cat{C},
\end{equation*}
constraining the components $(\alpha_x)_{x \in \cat{C}}$ of a natural
transformation $\alpha: F \to G$ are formally the same as the equations making
the functions $(\alpha_x)_{x \in \cat{C}}$ into a model homomorphism, when the
morphisms $f$ in $\cat{C}$ are interpreted as function symbols in a logical
theory and the functions $F(f)$ and $G(f)$ as set-theoretic models of $f$. In
his seminal PhD thesis \cite{lawvere1963}, William Lawvere transformed this
simple analogy into a deep connection between category theory and logic,
initiating a research program that continues today. According to the dictionary
of categorical logic, logical theories are small categories, usually with some
extra structure; models are functors out of these categories, preserving the
extra structure; and model homomorphisms are natural transformations. Choosing
the extra structure, say that of cartesian categories or elementary toposes,
amounts to choosing a logical system.

\begin{table}
  \centering
  \begin{tabular}{rll}
    \toprule
    & Logic & Categorical logic \\ \midrule
    Logical system & Syntax & Algebra (2-categories) \\
    Theories & Syntax & Algebra (small categories) \\
    Models & Semantics & Algebra (large categories) \\
    Model homomorphisms & Semantics & Algebra (large categories) \\
    \bottomrule
  \end{tabular}
  \caption{Syntax versus semantics in conventional and categorical logic. In
    categorical logic, everything is algebraic, with the formerly syntactical
    parts of logic happening in ``small'' structures and semantical parts
    happening in ``large'' ones. This distinction is stated precisely in
    \cref{ch:category-theory}.}
  \label{tbl:syntax-vs-semantics}
\end{table}

Categorical logic is an \emph{algebraization} of logic that obliterates the
traditional distinction between syntax and semantics
(\cref{tbl:syntax-vs-semantics}). As conventionally understood, logic is an
interplay between syntax and semantics. The logical system and theories
expressed in the system are \emph{syntactical}, made out of strings of symbols
manipulated according to certain rules, whereas models of theories and their
homomorphisms are \emph{semantical}, giving interpretation to the symbols as
ordinary mathematical objects. By contrast, in categorical logic, the logical
system, theories, models, and model homomorphisms are all algebraic structures
or morphisms thereof. These mathematical entities, like any others, can be
described in a formal language,\footnote{Formal languages for categorical
  structures include essentially algebraic theories \cite[\S 3.D]{adamek1994}
  and generalized algebraic theories \cite{cartmell1978,cartmell1986}.} as is
important for computer implementation, or in the natural language of written
mathematics, as in this text.

For applications to statistics, categorical logic has three advantages over
classical logic, all consequences of algebraization. The first is
\emph{functorial semantics}. Because, in categorical logic, a model of a theory
$\cat{C}$ is just a functor $\cat{C} \to \Set$, the target category of sets and
functions can easily be replaced by another large category $\cat{S}$. The
\emph{$\cat{S}$-valued models of $\cat{C}$}, or functors $\cat{C} \to \cat{S}$,
then offer an alternative interpretation of the theory. In this work, the target
category will be a category of sets in Euclidean space and Markov kernels
between them, thus interpreting the morphisms of the theory as probabilistic, or
randomized, functions. By comparison, in theoretical statistics, a statistical
model is classically defined as a Markov kernel $P: \Omega \to \sspace{X}$ from
a parameter space $\Omega$ to a sample space $\sspace{X}$, assigning to every
parameter $\theta \in \Omega$ a probability distribution $P_\theta$ on
$\sspace{X}$. Functorial semantics in a category of Markov kernels thus makes
immediate contact with statistics.

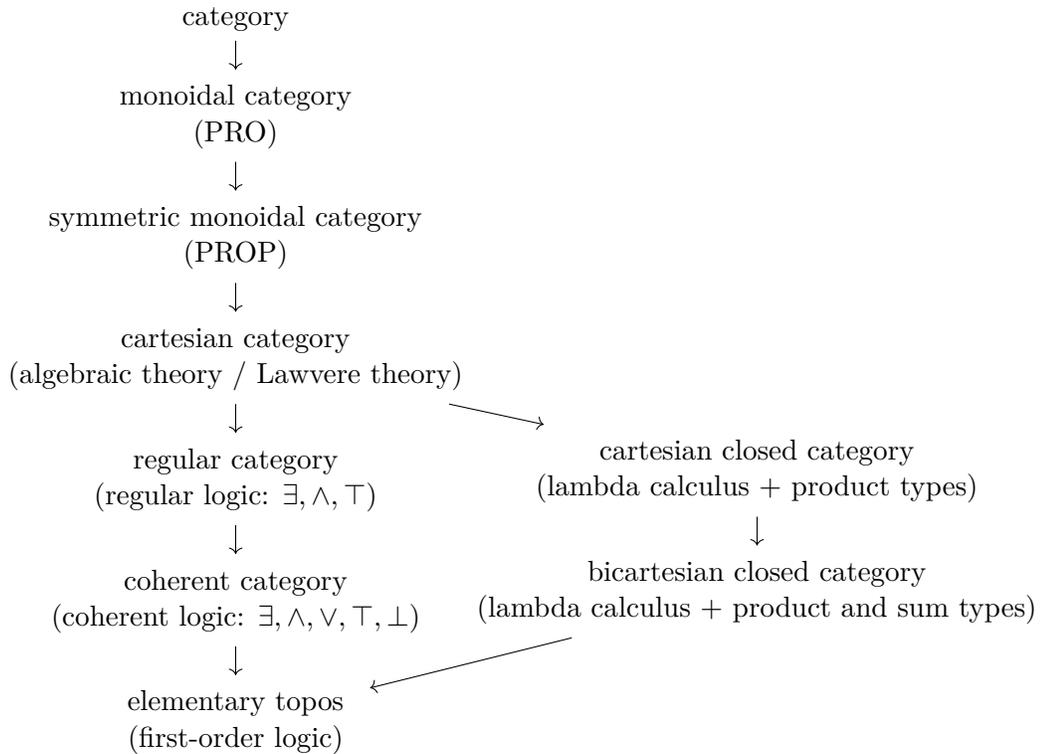
\begin{figure}
  \centering
  \begin{tikzpicture}[node distance=1em]
    \node (Cat) {category};
    \node[below=of Cat] (MonCat)
    {\begin{tabular}{c}
       monoidal category \\ (PRO)
     \end{tabular}};
    \node[below=of MonCat] (SMonCat)
    {\begin{tabular}{c}
       symmetric monoidal category \\ (PROP)
     \end{tabular}};
    \node[below=of SMonCat] (Cart)
    {\begin{tabular}{c}
       cartesian category \\ (algebraic theory / Lawvere theory)
     \end{tabular}};
    \node[below=of Cart] (RegCat)
    {\begin{tabular}{c}
       regular category \\ (regular logic: $\exists, \wedge, \top$)
     \end{tabular}};
    \node[below=of RegCat] (CohCat)
    {\begin{tabular}{c}
       coherent category \\
       (coherent logic: $\exists, \wedge, \vee, \top, \bot$)
     \end{tabular}};
    \node[below=of CohCat] (Topos)
    {\begin{tabular}{c}
      elementary topos \\ (first-order logic)
     \end{tabular}};
    \node[below right=of Cart] (CCC)
    {\begin{tabular}{c}
       cartesian closed category \\
       (lambda calculus + product types)
     \end{tabular}};
    \node[below=of CCC] (BCC)
    {\begin{tabular}{c}
       bicartesian closed category \\
       (lambda calculus + product and sum types)
     \end{tabular}};
    \draw[->] (Cat) -- (MonCat);
    \draw[->] (MonCat) -- (SMonCat);
    \draw[->] (SMonCat) -- (Cart);
    \draw[->] (Cart) -- (RegCat);
    \draw[->] (Cart) -- (CCC);
    \draw[->] (RegCat) -- (CohCat);
    \draw[->] (CohCat) -- (Topos);
    \draw[->] (CCC) -- (BCC);
    \draw[->] (BCC) -- (Topos);
  \end{tikzpicture}
  \caption{An incomplete family tree of categorical logic, including some
    fundamental logical systems in mathematics and computer science. The arrows
    are inclusions.}
  \label{fig:categorical-logics}
\end{figure}

Second, categorical logic is not a single logic but a whole family of logics
defined by the presence or absence of various algebraic gadgets, allowing great
flexibility in the construction of new, possibly unconventional logics. Since
the publication of Lawvere's thesis, categorical logic has subsumed many
important logical systems in mathematics and theoretical computer science, a few
of which are shown in the family tree of \cref{fig:categorical-logics}. Most of
the categorical structures extend the cartesian categories studied by Lawvere to
accommodate increasingly expressive subsystems of typed first-order logic or
other type theories.

\begin{figure}
  \centering
  \begin{tikzpicture}[node distance=1em]
    \node (SMonCat) {symmetric monoidal category};
    \node[below=of SMonCat] (SMonCatCComon)
    {\begin{tabular}{c}
       symmetric monoidal category \\ supplying commutative comonoids
     \end{tabular}};
    \node[below=of SMonCatCComon] (MarkovCat) {Markov category};
    \node[below right=of MarkovCat] (Cart) {cartesian category};
    \node[below left=of MarkovCat] (MarkovCatSupp)
    {\begin{tabular}{c}
       linear algebraic \\ Markov category
     \end{tabular}};
    \node[below=4.5em of MarkovCat] (CartSupp)
    {\begin{tabular}{c}
       linear algebraic \\ cartesian category
     \end{tabular}};
    \draw[->] (SMonCat) -- (SMonCatCComon);
    \draw[->] (SMonCatCComon) -- (MarkovCat);
    \draw[->] (MarkovCat) -- (MarkovCatSupp);
    \draw[->] (MarkovCat) -- (Cart);
    \draw[->] (MarkovCatSupp) -- (CartSupp);
    \draw[->] (Cart) -- (CartSupp);
  \end{tikzpicture}
  \caption{Systems of categorical logic relevant to probability and statistics}
  \label{fig:categorical-logics-stats}
\end{figure}
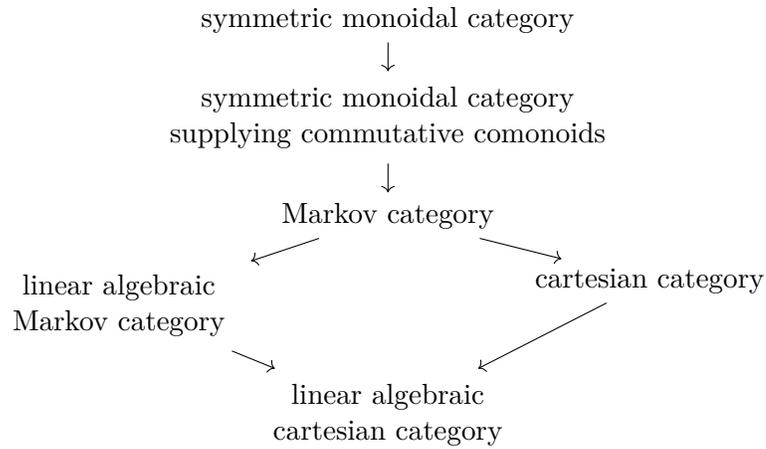

The algebra of statistical modeling belongs to a different branch of the family
tree, shown in \cref{fig:categorical-logics-stats}. The starting point,
\emph{Markov categories},\footnote{The complex provenance of the notion of
  Markov category is reviewed in the Notes to \cref{ch:algebra-statistics}.} are
intermediate between symmetric monoidal categories and cartesian categories. In
fact, a Markov category satisfies all the laws of a cartesian category, except
one: a morphism $M: \sspace{X} \to \sspace{Y}$ can be \emph{nondeterministic} in
that it fails to preserve the copying of data:
\begin{equation*}
  \input{wiring-diagrams/algebra-statistics/markov-kernel-mcopy-lhs} \neq
  \input{wiring-diagrams/algebra-statistics/markov-kernel-mcopy-rhs}.
\end{equation*}
In addition to this, almost all parametric statistical models require vector
space or similar structure in the parameter space and possibly also in the
sample space. Thus, in \emph{linear algebraic} Markov categories, objects may be
supplied with operations for taking linear, affine, conical, or convex
combinations, enabling them to be treated as vector, affine, conical, or convex
spaces. A surprisingly large amount of statistics can be formulated in this
setting. Aligning \cref{fig:categorical-logics,fig:categorical-logics-stats}
shows precisely how the algebra of statistical models is related to other
categorical structures and logical systems.

Lastly, returning to the original motivation, theories in categorical logic are
algebraic structures and thus admit morphisms between them. Morphisms of
theories, in turn, induce morphisms of the theories' categories of models. The
theories belonging to the algebra of statistics are called \emph{statistical
  theories}, and the ensuing morphisms between statistical theories, and between
their models and categories of models, all serve to formalize the relations that
exist between statistical models.

Statistical models are classically understood as models of phenomena or data, in
an informal sense of the word ``model,'' but not as models of a theory, in the
logical sense. This work separates statistical theories, which are small
categories typically presented by generators and relations, from models of
statistical theories, which are statistical models as usually understood. More
precisely, a \emph{statistical theory} is a small linear algebraic Markov
category $\cat{T}$ with a distinguished morphism $p: \theta \to x$, representing
the sampling distribution, and a \emph{statistical model} is a
structure-preserving functor $M: \cat{T} \to \Stat$, realizing the sampling
morphism as the Markov kernel $M(p)$. A \emph{model homomorphism} between two
models $M$ and $N$ of a statistical theory $\cat{T}$ is a monoidal natural
transformation $M \to N$. Having established these three fundamental
definitions, a \emph{theory morphism} $(\cat{T},p) \to (\cat{T},p')$ between two
different statistical theories is then a functor $F: \cat{T} \to \cat{T}'$
between the underlying categories that preserves the sampling morphism, taking
$p$ to $p'$ strictly or laxly. The functor $F$ also induces a \emph{model
  migration functor} $F^*$ in the opposite direction, taking models of
$\cat{T}'$ to models of $\cat{T}$.

This machinery is developed systematically in \cref{ch:algebra-statistics}. For
now, we simply emphasize that the notion of statistical theory gives rise to a
whole network of morphisms, all of which cease to exist when statistical models
are divorced from their theories.\footnote{To express this point differently, a
  Markov kernel $P: \Omega \to \sspace{X}$, viewed as a classical statistical
  model, is generally a model of many different statistical theories. What is a
  valid morphism of the model depends on which theory is adopted.} As an example
to be presented more fully in \cref{ch:zoo-statistics}, consider the standard
fact that a linear model is the special case of a generalized linear model where
the link function is the identity and the random component is the normal
distribution. Formalizing this, the theory of a linear model is a certain linear
algebraic Markov category $\cat{LM}_n$ together with a sampling morphism
$p_n: \beta \otimes \sigma^2 \to y^{\otimes n}$ of the form
\begin{equation*}
  \input{wiring-diagrams/zoo-statistics/lm-n-sampling}.
\end{equation*}
Similarly, the theory of a generalized linear model with a dispersion parameter
is another category $\cat{GLM}_n$ together with a sampling morphism
$p_n: \beta \otimes \phi \to y^{\otimes n}$ of form
\begin{equation*}
  \input{wiring-diagrams/zoo-statistics/glm-n-sampling},
\end{equation*}
where $h := g^{-1}: \eta \to \mu$ is the inverse of the link morphism
$g: \mu \to \eta$. A morphism of statistical theories
$F: (\cat{GLM}_n, p_n) \to (\cat{LM}_n, p_n)$ sends both parameter objects $\mu$
and $\eta$ in $\cat{GLM}_n$ to $\mu$ in $\cat{LM}_n$; sends the parameter object
$\phi$ to $\sigma^2$; sends the link morphism $g: \mu \to \eta$ to the identity
$1_\mu$,
\begin{equation*}
  F: \input{wiring-diagrams/zoo-statistics/glm-link}
  \mapsto  \input{wiring-diagrams/zoo-statistics/mu-id};
\end{equation*}
and preserves the other generators. When applied to linear models, the model
migration functor $F^*: \Model{\cat{LM}_n} \to \Model{\cat{GLM}_n}$ recovers the
generalized linear models with identity link and the normal family.

\section{Machine representation of data analyses}
\label{sec:introduction-semantic-enrichment}

While research papers have traditionally been, and to a large extent still are,
the primary information carriers of science, digital artifacts are becoming
increasingly valuable in disseminating and reusing scientific work. Chief among
digital artifacts are data and code. Large-scale observational studies and
high-throughput measurement devices now produce data of greater size and
complexity than ever before. Effective analysis of this data requires, in turn,
greater computational power and more sophisticated software. Modern scientific
datasets are unlikely to fit easily inside a data table in a scientific report,
and their analyses are unlikely to be faithfully captured by a few paragraphs of
written summary. Data and code thus take on new importance not just as
ancillaries to scientific research but as essential products of it.

The elevated status of data and code stems from both the need to verify existing
scientific work and the opportunity to produce new, derivative works. On the one
hand, as data and data analysis become more complex, reconstructing an analysis
becomes ever more difficult without access to the data behind it or the computer
programs implementing it. The movement for \emph{reproducible research} aims to
ensure that the analysis conducted in any study can be fully reproduced by
scientists besides the original authors. This requires at minimum that all
relevant data and code be preserved and made available. On the other hand,
collecting richer data sets increases the value of data reuse through creative
analyses not necessarily anticipated by the original data collectors. The
prospect of productive data reuse has motivated large-scale collection efforts,
of digital data and physical biosamples, in genetics, biomedicine, and other
fields. Likewise, publishing data analysis code as open source software allows
the whole community to modify and extend it, possibly in ways unforeseen by the
original authors.

Data and code, wherever they are available, provide new material for metascience
that goes beyond the traditional research paper. As data without accompanying
metadata is effectively meaningless, guidelines and standards have been proposed
for annotating published data. The FAIR Data Principles aim to simplify data
reuse, for machines as well as humans, by making datasets more ``FAIR'':
findable, accessible, interoperable, and reusable \cite{wilkinson2016}.
Beginning with early standards like the Minimal Information About a Microarray
Experiment \cite{brazma2001}, ``minimal information'' standards have been
defined for many kinds of experimental data \cite{mcquilton2016,sansone2019}.
Generic tools to assist the acquisition of metadata and improve its quality are
also being developed \cite{musen2015}.

Source code performing data analysis can now be found on general-purpose
platforms for open source software, such as GitHub, as well as on more
specialized, data-scientific platforms, such as Kaggle, Code Ocean, and DREAM
Challenges \cite{stolovitzky2007,saez-rodriguez2016}. Although code sharing,
like data sharing, is far from being a universally adopted practice, it should
be expected to grow along with the broader open science movement. Eventually it
is likely to become a norm or requirement across the sciences.\footnote{In
  fields involving human subjects, privacy protections place inherent
  limitations on data sharing and raise difficult questions about
  reproducibility. Code sharing faces no such limitations. It is hindered only
  by institutional and technological deficits.} But how to make effective
metascientific use of this new artifact is not yet clear.

The second part of this dissertation concerns the semantic analysis of data
science code, as a complement to the natural language processing of scientific
papers. It may appear that computer programs, unlike natural text, are already
interpretable by machines, but that is true only in the weakest of senses.
Computer programs must be sufficiently detailed and precise to be executed
without ambiguity by a computer. But a computer program, as understood by its
authors and its users, is more than a series of instructions for moving bits
around in memory. The program embodies generic concepts from its application
domain, which in the context of data science include loading and transforming
data, fitting and evaluating models, and making plots and other visualizations.
Such concepts are not transparent to the machine. They may not be transparent to
a human either, depending on how clearly the code is written and who is reading
it.

This difficulty is only compounded by the massive proliferation of programming
languages, frameworks, and libraries in the data science community, with
different traditions tending to favor certain tools over others. The Python
programming language has become very popular within machine learning, while the
R language is preferred by many statisticians \cite{rlang2020}. Both the Python
and R ecosystems are backed by lower-level code written in C, C++, and Fortran,
invisible to most users but necessary for certain purposes. The Julia language
has emerged recently to address this ``two-language problem'' \cite{julia2017}.
Outside these general-purpose programming languages, statistical models are also
built in probabilistic programming languages like Stan \cite{carpenter2017}.

The commonly used frameworks and libraries for data science are far too numerous
to summarize. One general observation is that a single programming language can
easily support multiple competing or complementary software ecosystems, written
in different styles or for different purposes but providing overlapping
functionality. In Python, the NumPy and SciPy packages define performant
multidimensional arrays and core routines for scientific computing, in a mainly
procedural style, while packages like Pandas, Scikit-learn, and Matplotlib build
on this foundation to provide data structures, algorithms, and visualization for
data science, in a more object-oriented style
\cite{scipy2020,mckinney2010,scikit-learn2011,hunter2007}. Within the R
community, the ``base R'' ecosystem comprising R's standard library and other
core packages is complemented by the newer ``tidyverse'' \cite{wickham2019}. The
latter set of packages makes extensive use of metaprogramming to create
domain-specific language for data analysis.

Consider the highly simplified data analysis shown in three variants in
\cref{lst:toy-kmeans-scipy,lst:toy-kmeans-sklearn,lst:toy-kmeans-r}. The first
is written in Python, using NumPy and SciPy; the second also in Python, but now
using Pandas and Scikit-learn; and the third in R, using R's standard library.
The first and third programs are written in a mostly functional style, whereas
the second is written in object-oriented style with mutating operations. Despite
differences in programming paradigm, language, and packages, all three programs
perform the same analysis: read the Iris dataset from a CSV file, drop the last
column (labeling the flower species), fit a $k$-means clustering model with
three clusters to the remaining columns, and return the cluster assignments and
centroids. The programs are thus \emph{semantically equivalent}, at least at a
certain level of abstraction.\footnote{The programs may not actually produce
  identical results, due to numerical error or different initializations of the
  iterative k-means algorithm. The question of what it means for two programs to
  be ``the same'' is more subtle than it may initially appear.}

\begin{listing}
  \begin{tcolorbox}
    \inputminted{python}{code/clustering-kmeans-scipy.py}
  \end{tcolorbox}
  \caption{$k$-means clustering in Python, using NumPy and SciPy}
  \label{lst:toy-kmeans-scipy}
\end{listing}

\begin{listing}
  \begin{tcolorbox}
    \inputminted{python}{code/clustering-kmeans-sklearn.py}
  \end{tcolorbox}
  \caption{$k$-means clustering in Python, using Pandas and Scikit-learn}
  \label{lst:toy-kmeans-sklearn}
\end{listing}

\begin{listing}
  \begin{tcolorbox}
    \inputminted{r}{code/clustering-kmeans-base-r.R}
  \end{tcolorbox}
  \caption{$k$-means clustering in R}
  \label{lst:toy-kmeans-r}
\end{listing}

This dissertation introduces a method and accompanying software system for
creating semantic representations of data science code. The representations are
machine-interpretable and, under certain assumptions, they precisely identify
semantic equivalence of code, irrespective of the programming paradigm,
language, or libraries used by the code.

On each of the simple programs in
\cref{lst:toy-kmeans-scipy,lst:toy-kmeans-sklearn,lst:toy-kmeans-r}, the system
produces the output shown in \cref{fig:toy-kmeans-semantic}. The boxes and wires
in the diagram represent function calls and objects, respectively, in an
idealized programming language. The labels on the boxes and wires refer to
concepts in an ontology about data science, while the absence of a label
indicates unknown semantics in a part of the program. Since it produces the same
result on all three programs, the system correctly identifies the semantic
equivalence of the programs, modulo the blank boxes. Of course, realistic data
analyses will never be perfectly semantically equivalent, just more or less
\emph{semantically similar}. This similarity will be reflected by overlap
between their semantic representations.

\begin{figure}
  \centering
  \input{wiring-diagrams/semantic-enrichment/clustering-kmeans}
  \caption{Semantic flow graph for
    \cref{lst:toy-kmeans-scipy,lst:toy-kmeans-sklearn,lst:toy-kmeans-r}}
  \label{fig:toy-kmeans-semantic}
\end{figure}

The semantic representations are constructed through a composite of two
high-level processes, which form the topics of
\cref{ch:program-analysis,ch:semantic-enrichment}. First, the source code is
subjected to a mixture of static and dynamic computer program analysis. The aim
of this analysis is to capture the program's data flow as a string diagram, the
\emph{raw flow graph}. The boxes in this diagram represent function calls,
method calls, or other computational units in the program and the wires
represent objects. The raw flow graphs constructed for
\cref{lst:toy-kmeans-scipy,lst:toy-kmeans-sklearn,lst:toy-kmeans-r} are shown in
\cref{fig:toy-kmeans-scipy,fig:toy-kmeans-sklearn,fig:toy-kmeans-r}. Although
the labels in the raw flow graph correspond closely to names in the original
program, a variety of transformations are needed to bring the program into
conformance with an idealized model of computation. For example, the mutating
method \verb|fit| in \cref{lst:toy-kmeans-sklearn} is reinterpreted in
\cref{fig:toy-kmeans-sklearn} as a function taking one \verb|KMeans| object as
input and returning another as output. The program analysis method for producing
raw flow graphs has been implemented for both Python and R programs, in the
respective languages.\footnote{Source code for the program analysis tools is
  available at \url{https://github.com/IBM/pyflowgraph} and
  \url{https://github.com/IBM/rflowgraph}.}

\begin{figure}
    \centering
    \input{wiring-diagrams/program-analysis/clustering-kmeans-scipy.raw}
    \caption{Raw flow graph for \cref{lst:toy-kmeans-scipy}}
    \label{fig:toy-kmeans-scipy}
    \vspace{0.25in}
\end{figure}
  
\begin{figure}
  \begin{minipage}{0.5\textwidth}
    \centering
    \input{wiring-diagrams/program-analysis/clustering-kmeans-sklearn.raw}
    \caption{Raw flow graph for \cref{lst:toy-kmeans-sklearn}}
    \label{fig:toy-kmeans-sklearn}
  \end{minipage}%
  \begin{minipage}{0.5\textwidth}
    \centering
    \input{wiring-diagrams/program-analysis/clustering-kmeans-base-r.raw}
    \caption{Raw flow graph for \cref{lst:toy-kmeans-r}}
    \label{fig:toy-kmeans-r}
  \end{minipage}
\end{figure}

The raw flow graph is expressed in terms of the programming language and
libraries of the original code. In the second major step, of \emph{semantic
  enrichment}, the raw flow graph is transformed into a representation that is
language and library independent. The resulting \emph{semantic flow graph} has
already been shown in \cref{fig:toy-kmeans-semantic}.

For semantic enrichment, two essential sources of information are supplied by
the \emph{Data Science Ontology}, a knowledge base about data
science.\footnote{The Data Science Ontology is browsable at
  \url{https://www.datascienceontology.org}, with the source files available at
  \url{https://github.com/IBM/datascienceontology}.} \emph{Concepts} in the
ontology, cataloging ideas from machine learning, statistics, and computing on
data, constitute the language in which semantic flow graphs are expressed. For
instance, the \textsf{k-means} type concept in \cref{fig:toy-kmeans-semantic}
refers to the k-means clustering model and the \textsf{fit} function concept
refers to the process of fitting an unsupervised model to data. The ontology
also contains \emph{annotations} mapping code from data science libraries onto
concepts. For example, an annotation identifies the \verb|KMeans| class in
Scikit-learn as an instance of the \textsf{k-means} concept. Annotations
determine the transformations made by the semantic enrichment algorithm.
Although agnostic to the programming language of its input data, semantic
enrichment is itself implemented in Julia.\footnote{Source code for semantic
  enrichment is available at \url{https://github.com/IBM/semanticflowgraph}. An
  important dependency is the package Catlab.jl, available at
  \url{https://github.com/epatters/Catlab.jl}.}

\begin{figure}
  \centering
  \includegraphics[width=\textwidth]{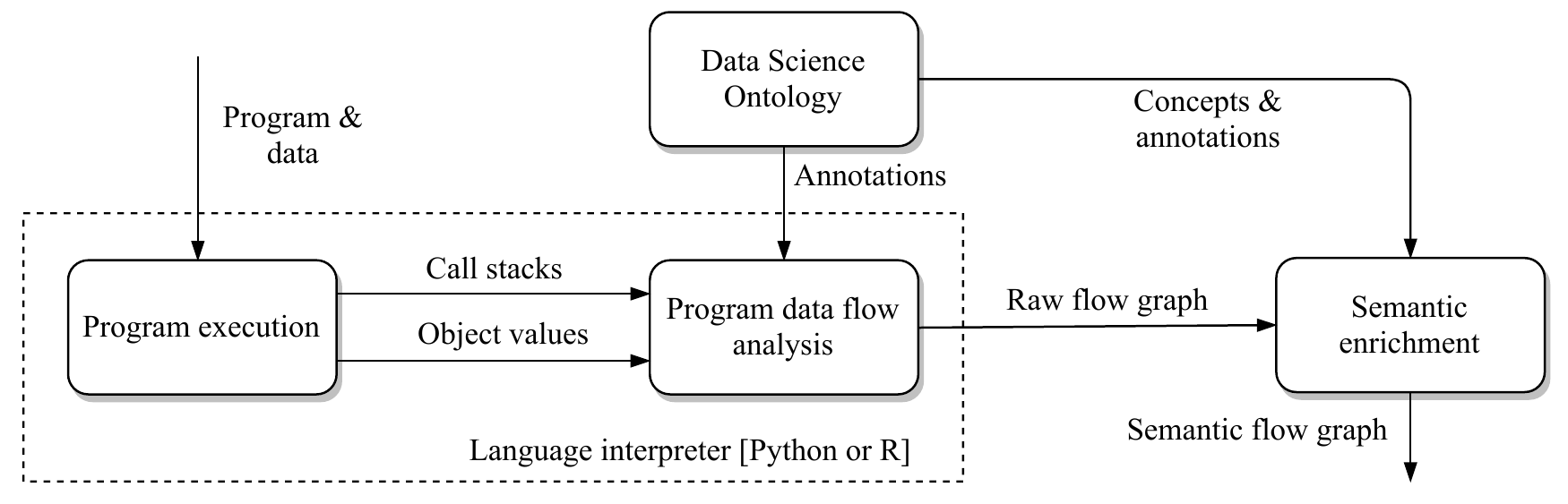}
  \caption{Architecture of software system for program analysis and semantic
    enrichment}
  \label{fig:system-architecture}
\end{figure}

The architecture of software system, spanning program analysis and semantic
enrichment, is summarized by \cref{fig:system-architecture}. The system is fully
automated, inasmuch as it requires nothing from the data analyst besides the
analysis itself and the ability to execute it. The method does, however,
indirectly depend on human input, since the Data Science Ontology is at present
entirely human-constructed. Thus, human effort is involved in the definition of
concepts and library annotations, which may then be reused for any analyses
involving those concepts or libraries. The implicit assumption behind this
division of labor is that data scientists use methods implemented in commonplace
software packages, like Pandas or Scikit-learn, rather than inventing their own.
This assumption is often, but of course not always, satisfied.

Our methodology is not universally applicable, but neither is it specific in all
respects to data science. It is designed for scripts and computational notebooks
written by data scientists, which tend to be shorter, more linear, and clearer
semantically than the large-scale codebases written by software engineers. It
would be wholly unsuited to the analysis of complex software systems like
compilers or web servers. However, it might fruitfully be extended to other
scientific domains with a computational focus, such as bioinformatics or
computational neuroscience, through integration with existing domain-specific
ontologies.

\section{Notes and references}

References related to the context and philosophy of this dissertation are
collected below. For technical references, see the Notes to subsequent chapters.

\paragraph{Open and networked science}

Swanson announced his discovery, by literature mining, of a possible connection
between migraine headaches and magnesium in a now classic paper
\cite{swanson1988}. The review \cite{swanson1990} describes the ``$ABC$
syllogism'' for implicit connections in the scientific literature and summarizes
the migraine-magnesium study, along with an earlier proposal linking Raynaud's
syndrome and fish oil \cite{swanson1986}. Clinical studies on magnesium
treatment for migraines, showing mostly positive results, are reviewed up to
1998 in \cite{mauskop1998}.

In \emph{Reinventing Discovery} \cite{nielsen2012}, Nielsen gives an inspiring
account of a future science that is more open, networked, data- and
machine-driven than today. Chapter 6 explores the consequences of digitizing of
scientific knowledge, telling the story of Don Swanson among others. The edited
collection \cite{hey2009} explores how the ``fourth
paradigm''\footnote{According to Jim Gray, as transcribed in the introduction to
  \cite{hey2009}, the first paradigm is \emph{empirical}, dating back to
  antiquity; the second paradigm is \emph{theoretical}, beginning with the
  revolution of physics in the 17th century; and the third paradigm is
  \emph{computational}, only a few decades old. Despite the implied progression,
  the paradigms are evidently not mutually exclusive.} of data-intensive science
will impact specific scientific fields, as well as scientific computing and
communication generally.

\paragraph{Statistics and data science}

Efron and Hastie trace the evolution of statistics and the outgrowth of data
science, from the 1950s to the present \cite{efron2016}. A recurring theme is
how statistical methodology has changed to take advantage of greater
computational power and meet the demands imposed by large-scale data collection.

All are agreed that data science has something to do with an algorithmic and
computational turn in data analysis, but one could be forgiven for thinking that
every data scientist has their own idea of what data science is or should be.
Perspectives on the relationship between statistics and data science, some
significantly predating the term ``data science,''\footnote{Tukey coined the
  term ``data analysis'' to encompass an activity larger than mathematical
  statistics and resembling what is now called ``data science''
  \cite{tukey1962}. Chambers later spoke of ``greater statistics,'' in contrast
  to the ``lesser statistics'' that is mainly confined to the mathematics of
  probability models \cite{chambers1993b}.} include
\cite{tukey1962,chambers1993b,breiman2001,donoho2017,carmichael2018}. In his
essay on ``50 years of data science'', Donoho identifies ``science about data
science'' as the last of six fundamental divisions\footnote{The other five
  divisions proposed by Dohono are data gathering, preparation, and exploration;
  data representation and transformation; computing with data; data modeling,
  both generative and predictive; and data visualization and presentation.}
within data science \cite{donoho2017}. This division encompasses the
meta-scientific study of data analysis workflows used by practitioners, as well
as the foundational work needed to make such metascience possible. The present
work belongs to the latter category. In our view, science about data science is
quite obviously the least developed of the six divisions and the farthest from
the mainstream of statistics.

\paragraph{Models in philosophy of science}

Already in the early twentieth century, with mathematical logic still in its
infancy, the logical positivists made the first attempts to reconstruct
scientific theories as logical ones. Present-day philosophers of science
emphasize scientific models at least as much as theories, viewing models not as
mere analogies or visual aides but as mathematical objects mediating between
theories and observations. The history of scientific models in the philosophy of
science, as well as contemporary views, are summarized by Bailer-Jones
\cite{bailer-jones2009}.

Models figure prominently in the \emph{semantic view} of scientific theories.
Although antecedents may be found in Evert Beth's semantic analysis of classical
and quantum mechanics, Patrick Suppes made the first general arguments for
identifying models in science with those in mathematical logic
\cite{suppes1961,suppes1966,suppes1967}. Later influential workers in this
tradition include van Fraassen \cite{vanfraassen1980,vanfraassen1987}, Sneed
\cite{sneed1979}, and Suppe \cite{suppe1977,suppe1989}. A useful summary and
synthesis of the semantic tradition is given by Ruttkamp \cite{ruttkamp2002}.
Despite technical differences between the accounts, proponents of the semantic
view generally agree that scientific models can be understood as models of
logical theories, even if the practical uses of scientific and logical models
differ considerably.

Another recurring idea within the semantic tradition is that scientific
theories, while they may be formulated linguistically, are better understood as
extra-linguistic entities whose primary purpose is to determine a class of
models. For Suppes, this amounts to working in set theory under the slogan that
``to axiomatize a theory is to define a set-theoretical predicate'' \cite[\S
2.3]{suppes2002}. Unfortunately, the set-theoretical approach lends itself
poorly to a self-contained formalization or computer implementation, as it
assumes that all the relevant mathematics from real analysis and probability
theory has been encoded in set theory. Using the modern tools of categorical
logic, it is a simple matter to incorporate analytical and probabilistic
elements into theories and models while retaining an easily formalized, purely
equational language for stating theories. In construing theories as algebraic
structures, categorical logic agrees with the semantic view that theories are
extra-linguistic entities while implementing this principle differently.

\chapter{Elements of category theory and categorical logic}
\chaptermark{Category theory and categorical logic}
\label{ch:category-theory}

The structures of abstract algebra are stylized mathematical theories of
commonplace physical or mathematical ideas. Groups are an algebraic theory of
symmetry, groupoids are an algebraic theory of local symmetry
\cite{vistoli2011}, and rings are an algebraic theory of number systems. From
this point of view, categories are an algebraic theory of
\emph{compositionality}: the possibility of taking mappings, or relations, or
processes, or some other kind of generalized arrow from one object to another,
and composing them to get new arrows.

Outside of pure mathematics, algebra is useful insofar as it successfully models
our informal human concepts or aspects of natural phenomena. The theory of
groups and group representations is useful to physicists because symmetry has
emerged as a fundamental organizing principle of physics \cite{gross1996}.
Category theory is potentially useful to a range of scientific and engineering
fields because many physical and computational processes are intrinsically
compositional. Chemists compose multi-step chemical processes out of basic
reactions, programmers develop elaborate software systems out of built-in
routines, and statisticians build complex statistical models out of simpler
ones. Category theory offers a common mathematical language to describe and
reason about such seemingly different activities.

A striking success of category theory has been the \emph{algebraization} of
logic. Mathematical logic is conventionally understood to be an interplay
between syntax and semantics \cite{girard1999}. Logical theories are syntactic
objects, made out of symbols arranged into lists or trees, whereas the models of
a theory, its semantics, are ordinary mathematical objects. Categorical logic
turns this distinction on its head by making every constituent of logic,
including the theories, into algebraic structures. In the dictionary of
categorical logic, theories become small categories, models of theories become
functors, and model homomorphisms become natural transformations. The result is
to introduce, or at least to make evident, the compositionality of logic.
Theories, like any other algebraic structure, have morphisms, and morphisms of
theories can be composed with models to obtain new models, of different
theories. Dually, morphisms of model categories can be composed with models,
yielding new models in different categories.

This chapter introduces a selection of essential ideas from category theory and
categorical logic and serves as background for the rest of the thesis. Apart
from \cref{sec:interacting-supplies}, its content is not original, except
possibly in presentation. No knowledge of category theory or mathematical logic
is assumed of the reader, although some exposure to abstract algebra will be
helpful. Suggestions for further reading are provided at the end of the chapter.

\section{Categories, functors, and natural transformations}
\label{sec:categories}

The three fundamental notions of category theory are categories themselves;
morphisms of categories, or functors; and morphisms of functors, or natural
transformations.

A category axiomatizes the composition of mappings, relations, or other directed
arrows between different objects. In plain words, a category is a directed
graph, possibly infinite, in which it is possible to compose two adjacent edges
to obtain another edge. The composition is subject to axioms of unitality and
associativity familiar from abstract algebra.

\begin{definition}[Category]
  A \emph{category} $\cat{C}$ consists of a collection $\Ob{\cat{C}}$ of
  \emph{objects}, written $x,y,z,\ldots \in \cat{C}$, and for each pair of
  objects $x$ and $y$, a collection $\cat{C}(x,y)$ of \emph{morphisms}, written
  $f,g,h,\ldots \in \cat{C}(x,y)$. A morphism $f \in \cat{C}(x,y)$ has
  \emph{domain} $x$ and \emph{codomain} $y$ and is denoted by $f: x \to y$.

  For each object $x$, there is an \emph{identity} morphism $1_x: x \to x$ and
  for any two morphisms $f: x \to y$ and $g: y \to z$, where the codomain of $f$
  is equal to the domain of $g$, there is a \emph{composite} morphism
  $f \cdot g: x \to z$,\footnote{Composition of morphisms is written in
    diagrammatic (left-to-right) order. Thus, the usual notation for function
    composition translates to $g \circ f := f \cdot g$.} or simply $fg$, subject
  to two axioms:
  \begin{enumerate}
  \item (Unitality) For any morphism $f: x \to y$, both $1_x \cdot f$ and
    $f \cdot 1_y$ are equal to $f$.
  \item (Associativity) For any composable morphisms
    $w \xrightarrow{f} x \xrightarrow{g} y \xrightarrow{h} z$, the composites
    $(f \cdot g) \cdot h$ and $f \cdot (g \cdot h)$ are equal, and are hence
    denoted simply by $f \cdot g \cdot h$.
  \end{enumerate}

  A morphism $f: x \to y$ in $\cat{C}$ is an \emph{isomorphism}, or is
  \emph{invertible}, if there exists a morphism $g: y \to x$ such that
  $f \cdot g = 1_x$ and $g \cdot f = 1_y$. In this case, the objects $x$ and $y$
  are \emph{isomorphic}, denoted $x \cong y$.
\end{definition}

By tradition, the first example of a category is the category $\Set$, having
sets as objects and functions between them as morphisms. Composition in $\Set$
is the usual composition of functions, and the identities $1_X: X \to X$ are the
usual identity functions $1_X: x \mapsto x$. The isomorphisms in $\Set$ are the
bijections.

Many objects in mathematics are defined as ``sets with extra structure.'' In
each case, there is a category whose objects are sets with that structure and
whose morphisms are functions preserving the structure. There is a category
$\Group$ of groups and group homomorphisms, a category $\Mon$ of monoids and
monoid homomorphisms, a category $\Graph$ of directed graphs and graph
homomorphisms, a category $\Top$ of topological spaces and continuous maps, a
category $\Meas$ of measurable spaces and measurable maps, and a category of
$\Poset$ of partially ordered sets (posets) and monotone maps, among countless
others.

All of these categories are \emph{large} in the sense that their collections of
objects are too large to be sets (famously, there is no set of all
sets).\footnote{They are, however, \emph{locally small}: categories $\cat{C}$
  for which all collections $\cat{C}(x,y)$ are sets, the \emph{hom-sets}. In
  this text, all categories will be locally small.} In a \emph{small} category,
the collections of all objects and of all morphisms are sets. Small categories
are algebraic structures in the vein of groups and rings. In fact, every group
$G$ furnishes an example of a small category: this category has a single object
$*$ and the morphisms $g: * \to *$ are the group elements $g \in G$. Turning
this around, a group could be defined as a category with one object and
invertible morphisms. Likewise, a monoid could be defined as a category with one
object. Groups and monoids regarded as small categories should not be confused
with the large categories $\Group$ and $\Mon$, in which groups and monoids are
\emph{objects}.

Categories with one object constitute a degenerate case in the taxonomy of
categories. At the other extreme are \emph{thin} categories, which have at most
one morphism between any two objects. A thin category is a \emph{preorder}: the
elements of the preorder are the objects $x,y,z,\ldots$ of the category, and
$x \leq y$ if and only if there exists a morphism $x \to y$. The preorder axioms
of reflexivity ($x \leq x$) and transitivity ($x \leq y$ and $y \leq z$ implies
$x \leq z$) are then equivalent to the category axioms of unitality and
associativity.

Every category has an opposite category obtained by turning around the arrows.
Formally, the \emph{opposite category} $\cat{C}^\opposite$ of a category
$\cat{C}$ has the same objects as $\cat{C}$ but has hom-sets
$\cat{C}^\opposite(x,y) := \cat{C}(y,x)$ for all $x,y \in \cat{C}$. Composition
in $\cat{C}^\opposite$ is the same as in $\cat{C}$, except that the order is
reversed. For example, if $X$ is a preorder, then $x \leq y$ in $X^\opposite$ if
and only if $y \leq x$ in $X$. In category theory as in order theory, duality is
powerful principle because every theorem that holds for arbitrary categories
automatically gives another, dual theorem by reversing the arrows.

Beyond the degenerate cases of monoids and preorders, the most common way to
construct small categories is to \emph{present} them by generators and relations
\cite[\S 5.4]{spivak2014}, \cite[\S II.8]{maclane1998}. This is perfectly
analogous to, and indeed generalizes, the presentation of groups by generators
and relations. In this text, categories defined by presentations will often
represent logical or statistical theories, where the objects of the category
correspond to the types of the theory and the morphisms to the operations.

\begin{example}[Theories of graphs]
  The \emph{theory of graphs}, specifically directed graphs,\footnote{Without
    further qualification, a ``graph'' will always be a directed multigraph,
    possibly with self-loops.} is the category freely generated by two objects
  and two parallel morphisms:
  \begin{equation*}
    \Theory{\Graph} := \left\langle
      \begin{tikzcd}
        E \rar[shift left=1]{\src}
        \rar[shift right=1][below]{\tgt}
        & V
      \end{tikzcd}
    \right\rangle.
  \end{equation*}
  The objects $V$ and $E$ represent vertices and edges, and the morphisms
  $\src, \tgt: E \to V$ represent the source and target vertices of an edge.

  The \emph{theory of symmetric graphs} is the category generated by an
  additional morphism $\inv: E \to E$, subject to three relations:
  \begin{equation*}
    \setlength{\jot}{0pt}
    \Theory{\SGraph} := \left\langle
      \begin{tikzcd}
        E \lar[loop left]{\inv}
        \rar[shift left=1]{\src}
        \rar[shift right=1][below]{\tgt}
        & V
      \end{tikzcd}
      \middle|\
      {\footnotesize
        \begin{aligned}
          \inv^2 = 1_E \\
          \inv \cdot \src = \tgt \\
          \inv \cdot \tgt = \src
        \end{aligned}}
    \right\rangle.
  \end{equation*}
  A symmetric graph is a directed graph with an involution on edges, matching
  every edge in the graph with an oppositely oriented edge \cite{boldi2002}. For
  most purposes, symmetric graphs are interchangeable with undirected graphs.
\end{example}

\begin{example}[Theory of discrete dynamical systems]
  The \emph{theory of discrete dynamical systems} is the category freely
  generated by one endomorphism:
  \begin{equation*}
    \Theory{\cat{DDS}} := \left\langle
      \begin{tikzcd}
        * \rar[loop right]{T}
      \end{tikzcd}
    \right\rangle.
  \end{equation*}
  This category is a monoid, the free monoid on one generator. It is isomorphic
  to $(\N,+,0)$, the additive monoid of natural numbers.
\end{example}

Any algebraic structure has a concomitant notion of homomorphism, or
structure-preserving map. For categories, that notion is a \emph{functor}.
Taking the pedestrian view of categories as directed graphs with composable
arrows, a functor between categories is just a graph homomorphism that preserves
the composition.

\begin{definition}[Functor]
  A \emph{functor} $F: \cat{C} \to \cat{D}$ between categories $\cat{C}$ and
  $\cat{D}$ consists of a map of objects $F: \Ob{\cat{C}} \to \Ob{\cat{D}}$ and
  for each pair of objects $x$ and $y$ in $\cat{C}$, a map of morphisms
  $F: \cat{C}(x,y) \to \cat{D}(Fx,Fy)$, satisfying the \emph{functorality
    axioms}:
  \begin{enumerate}
  \item For each object $x$ in $\cat{C}$, $F(1_x) = 1_{Fx}$.
  \item For any composable morphisms $x \xrightarrow{f} y \xrightarrow{g} z$ in
    $\cat{C}$, $F(f \cdot g) = Ff \cdot Fg$.
  \end{enumerate}
\end{definition}

Functors $F: \cat{C} \to \cat{D}$ and $G: \cat{D} \to \cat{E}$ compose by
composing the underlying maps on objects and morphisms, and every category
$\cat{C}$ has an identity functor $1_{\cat{C}}: \cat{C} \to \cat{C}$. Thus, the
small categories and functors between them form a category of their own: $\Cat$,
the \emph{category of small categories}.

A functor between groups, interpreted as categories with one object, is a group
homomorphism. Similarly, a functor between monoids is a monoid homomorphism. A
functor between posets or preorders, regarded as thin categories, is a monotone
map. The categories $\Group$, $\Mon$, $\Poset$, and $\Preord$ thus all belong to
$\Cat$ as subcategories.

Any category $\cat{C}$ of ``sets and functions with extra structure'' has a
\emph{forgetful functor} $U: \cat{C} \to \Set$ that discards this extra
structure. For example, the underlying set of a group defines a functor
$U: \Group \to \Set$. In general, a category $\cat{C}$ equipped with a faithful
(injective on hom-sets) functor $U: \cat{C} \to \Set$ is called a \emph{concrete
  category}.\footnote{Some authors refer to general categories as
  \emph{abstract} \cite{adamek2004}, but this usage is not entirely standard.}
Most of the large categories considered above, such as $\Group$, $\Mon$, $\Top$,
and $\Meas$, are concrete.

However, not all large categories are concrete, nor do categories that have sets
as objects necessarily have functions as morphisms. A typical example is the
category $\Rel$ of sets and binary relations, where the composite
$R \cdot S: X \to Z$ of relations $R: X \to Y$ and $S: Y \to Z$ is defined by
taking $(R \cdot S)(x,z)$ to be true if and only if there exists $y \in Y$ such
that $R(x,y)$ and $S(y,z)$. Composition of relations extends the usual
composition of functions, making $\Set$ into a subcategory of $\Rel$. Another
example, to be defined more carefully in \cref{ch:algebra-statistics}, is the
category $\Markov$ of measurable spaces and Markov kernels (informally,
``probabilistic functions'').

A forgetful functor $U: \Cat \to \Graph$ gives the underlying graph of a
category, discarding the composition and identity operations. In the other
direction, a functor $F: \Graph \to \Cat$ constructs the free category generated
by a graph. The free category $F(G)$ on a graph $G$ is characterized by its
universal property, but can be constructed explicitly by taking the vertices of
$G$ as objects and the directed paths in $G$ as morphisms. Composition is then
concatenation of paths, and the identity morphisms are empty paths (paths of
length zero). The theories of graphs and discrete dynamical systems, defined
above, are free categories. A category presented by generators and relations,
such as the theory of symmetric graphs, can be constructed explicitly by forming
the free category on the graph of generators and then quotienting out by the
relations. However, it is usually preferable to work with free categories and
presented categories abstractly, eliding the details of any particular
construction.\footnote{An exception is the computer implementation of
  categorical algebra, where the choice an appropriate concrete representation
  is essential.}

In the \emph{functorial semantics} of categorical logic, the \emph{models} of a
theory $\cat{C}$, a small category, are functors $X: \cat{C} \to \Set$ from
$\cat{C}$ into the category of sets.

\begin{example}[Graphs]
  A functor $G: \Theory{\Graph} \to \Set$ consists of a set $G(V)$ of vertices
  and a set $G(E)$ of edges, together with functions
  $G(\src), G(\tgt): G(E) \to G(V)$ assigning source and target vertices to each
  edge. Such a functor is simply a graph. Similarly, a functor
  $\Theory{\SGraph} \to \Set$ is a symmetric graph.
\end{example}

Group actions are set-valued functors of groups. Specifically, viewing a group
$G$ as a category on a single object $*$, a functor $X: G \to \Set$ consists of
a set $X := X(*)$ and functions $X(g): X \to X$, for $g \in G$, that preserve
the multiplication and identity of $G$. This is exactly an action of the group
$G$ or a $G$-set. The abstract group $G$ acts as a theory and the functor
$X: G \to \Set$, a transformation group, as a model of the theory. In
conventional notation, we conflate the function $X(g)$ with the group element
$g$ and write $g \cdot x := X(g)(x)$. Turning this around, a functor
$X: \cat{C} \to \Set$ is interpreted as an \emph{action} of the category
$\cat{C}$, or a \emph{$\cat{C}$-set} \cite{reyes2004}. When no confuse will
arise, we may conflate a morphism $f$ in $\cat{C}$ with its image $X(f)$ in
$\Set$.

Functorial semantics liberates logical semantics from the traditional setting of
sets and functions. For any category $\cat{S}$, typically large, a \emph{model}
of $\cat{C}$ in $\cat{S}$ is a functor $X: \cat{C} \to \cat{S}$. Taking
$\cat{S} = \Set$ recovers the usual set-valued semantics, but taking other
categories allows for models with extra, or even fundamentally different,
structure.

\begin{example}[Discrete dynamical systems]
  A model of the theory of discrete dynamical systems, or a functor
  $X: \Theory{\cat{DDS}} \to \Set$, is a \emph{discrete dynamical system}: a set
  $X := X(*)$ of states together with a state transition function $X \to X$. A
  model of the same theory, but in the category of Markov kernels, is a functor
  $X: \Theory{\cat{DDS}} \to \Markov$. Such a functor defines a \emph{Markov
    chain}: a measurable space $X := X(*)$ together with a probabilistic
  transition kernel $X \to X$. Simple though it may be, this observation is a
  point of departure for \cref{ch:algebra-statistics}, where statistical models
  are reinterpreted as models of theories in a category of Markov kernels.
\end{example}

\begin{example}[Group representations]
  Let $\Vect_\K$ be the category of vector spaces over a field $\K$, with linear
  transformations as morphisms. A functor $G \to \Vect_\K$ is a \emph{group
    representation}, the linear-algebraic variant of a group action.
\end{example}

Categories are distinguished from other common algebraic structures by having
morphisms between their morphisms. A morphism of functors is a natural
transformation.

\begin{definition}[Natural transformation]
  Let $F,G: \cat{C} \to \cat{D}$ be parallel functors between categories
  $\cat{C}$ and $\cat{D}$. A \emph{natural transformation} $\alpha: F \to G$
  between $F$ and $G$ consists of, for each object $x \in \cat{C}$, a morphism
  $\alpha_x: Fx \to Gx$ in $\cat{D}$, called the \emph{component} of $\alpha$ at
  $x$. Moreover, the components must satisfy the \emph{naturality axiom}: for
  each morphism $f: x \to y$ in $\cat{C}$, the square of morphisms
  \begin{equation*}
    \begin{tikzcd}
      Fx \arrow[r, "\alpha_x"] \arrow[d, "Ff"']
        & Gx \arrow[d, "Gf"] \\
      Fy \arrow[r, "\alpha_y"]
        & Gy
    \end{tikzcd}
  \end{equation*}
  in $\cat{D}$ commutes.

  A \emph{natural isomorphism} is a natural transformation $\alpha: F \to G$ in
  which every component $\alpha_x: Fx \to Gx$ is an isomorphism. Natural
  isomorphisms are also written $\alpha: F \cong G$.
\end{definition}

The arithmetic of sets is a classic source of natural isomorphisms. Let
$A \times B$ be the cartesian product of sets $A$ and $B$, and let $A+B$ be
their disjoint union. The distributive law does not strictly hold, i.e.,
$A \times (B + C)$ is not equal to $(A \times B) + (A \times C)$, but the law
does hold \emph{up to natural isomorphism}:
\begin{equation*}
  A \times (B + C) \cong (A \times B) + (A \times C).
\end{equation*}
The precise meaning of this statement is that the functors
$F,G: \Set \times \Set \times \Set \to \Set$, defined on objects by
$F(A,B,C) := A \times (B + C)$ and $G(A,B,C) := (A \times B) + (A \times C)$,
are related by a natural isomorphism $\alpha: F \to G$, whose components are the
bijections
\begin{equation*}
  \alpha_{A,B,C}: A \times (B+C) \to (A \times B) + (A \times C), \quad
    (a,(i,x)) \mapsto (i,(a,x)) ,
\end{equation*}
assuming the set-theoretic construction of $B+C$ as
$\{(0,b): b \in B\} \cup \{(1,c): c \in C\}$. In set arithmetic, the
associativity and commutativity laws of products and sums also hold only up to
natural isomorphism:
\begin{align*}
  (A \times B) \times C &\cong A \times (B \times C)
  \qquad\qquad
  A \times B \cong B \times A \\
  (A + B) + C &\cong A + (B + C)
  \qquad\qquad
  A + B \cong B + A.
\end{align*}
These and other natural isomorphisms are summarized by saying that the
arithmetic of sets is a ``categorification'' of the arithmetic of natural
numbers. In general, \emph{categorification} is a process of generalizing
algebraic laws by replacing equalities with natural isomorphisms
\cite{baez1998}.

A natural transformation $f \to k$ of group homomorphisms $f,k: G \to H$,
regarded as functors, is always an isomorphism and is determined by a single
element $h \in H$ that makes the homomorphism $k$ \emph{conjugate} to $f$,
meaning that $k(g) = h^{-1} \cdot f(g) \cdot h$ for all $g \in G$.
Geometrically, the homomorphisms $f$ and $k$ are equivalent up to a symmetry
transformation of the codomain. Turning from groups to group actions, a natural
transformation $X \to Y$ of $G$-sets $X$ and $Y$, regarded as functors
$G \to \Set$, is map of sets $\phi: X \to Y$ that preserves the action of $G$,
so that $\phi(g \cdot x) = g \cdot \phi(x)$ for all $g \in G$ and $x \in X$.
These morphisms of $G$-sets are called \emph{equivariant maps}. Similarly, a
natural transformation of group representations, regarded as functors
$G \to \Vect_\K$, is an \emph{equivariant map} or \emph{intertwining} of
representations.

Natural transformations compose in not one, but two, dimensions. To illustrate
this, draw a natural transformation $\alpha: F \to G$ between functors
$F,G: \cat{C} \to \cat{D}$ as the two-dimensional figure
\begin{equation*}
  \begin{tikzcd}
    \cat{C}
      \arrow[r, bend left=45, "F"{above}, ""{name=F, below}]
      \arrow[r, bend right=45, "G"{below}, ""{name=G, above}]
    & \cat{D}
      \arrow[Rightarrow, "\alpha", from=F, to=G]
  \end{tikzcd}.
\end{equation*}
Given three parallel functors $F,G,H: \cat{C} \to \cat{D}$, the \emph{vertical
  composite} $\alpha \cdot \beta: F \to H$ of natural transformations
$\alpha: F \to G$ and $\beta: G \to H$ is depicted as
\begin{equation*}
  \begin{tikzcd}[column sep=3em]
    \cat{C}
      \arrow[r, bend left=60, "F"{above}, ""{name=F, below}]
      \arrow[r, "G"{name=G, anchor=center, fill=white}]
      \arrow[r, bend right=60, "H"{below}, ""{name=H, above}]
    & \cat{D}
      \arrow[Rightarrow, from=F.center, to=G]
      \arrow[Rightarrow, from=G, to=H.center]
  \end{tikzcd}
  \qquad\leadsto\qquad
  \begin{tikzcd}[column sep=3em]
    \cat{C}
      \arrow[r, bend left=45, "F"{above}, ""{name=F, below}]
      \arrow[r, bend right=45, "H"{below}, ""{name=H, above}]
    & \cat{D}
      \arrow[Rightarrow, "\alpha \cdot \beta", from=F, to=H]
  \end{tikzcd}.
\end{equation*}
It is defined componentwise by setting
$(\alpha \cdot \beta)_x := \alpha_x \cdot \beta_x$ for all $x \in \cat{C}$. The
naturality of the vertical composite $\alpha \cdot \beta$ follows by pasting
together the naturality squares for $\alpha$ and $\beta$.

The second mode of composition, in the horizontal direction, is defined in terms
of a simpler operation called \emph{whiskering}, wherein a natural
transformation is composed with a functor. Given parallel functors
$G,H: \cat{C} \to \cat{D}$, the \emph{pre-whiskering} $F\beta$ of a natural
transformation $\beta: G \to H$ by a functor $F: \cat{B} \to \cat{C}$, is
depicted as
\begin{equation*}
  \begin{tikzcd}
    \cat{B}
      \arrow[r, "F"]
    & \cat{C}
      \arrow[r, bend left=45, "G"{above}, ""{name=G, below}]
      \arrow[r, bend right=45, "H"{below}, ""{name=H, above}]
    & \cat{D}
      \arrow[Rightarrow, "\beta", from=G, to=H]
  \end{tikzcd}
  \qquad\leadsto\qquad
  \begin{tikzcd}[column sep=3em]
    \cat{B}
      \arrow[r, bend left=45, "FG"{above}, ""{name=FG, below}]
      \arrow[r, bend right=45, "FH"{below}, ""{name=FH, above}]
    & \cat{D}
      \arrow[Rightarrow, "F \beta", from=FG, to=FH]
  \end{tikzcd}
\end{equation*}
and is defined by pre-composing along components, as
$(F\beta)_x := \beta_{Fx}: G(F(x)) \to H(F(x))$ for all $x \in \cat{B}$. Dually,
given functors $F,G: \cat{C} \to \cat{D}$, the \emph{post-whiskering} $\alpha H$
of a natural transformation $\alpha: F \to G$ by a functor
$H: \cat{D} \to \cat{E}$ is depicted as
\begin{equation*}
  \begin{tikzcd}
    \cat{C}
      \arrow[r, bend left=45, "F"{above}, ""{name=G, below}]
      \arrow[r, bend right=45, "G"{below}, ""{name=H, above}]
    & \cat{D}
      \arrow[Rightarrow, "\alpha", from=G, to=H]
      \arrow[r, "H"]
    & \cat{E}
  \end{tikzcd}
  \qquad\leadsto\qquad
  \begin{tikzcd}[column sep=3em]
    \cat{C}
      \arrow[r, bend left=45, "FH"{above}, ""{name=FH, below}]
      \arrow[r, bend right=45, "GH"{below}, ""{name=GH, above}]
    & \cat{E}
      \arrow[Rightarrow, "\alpha H", from=FH, to=GH]
  \end{tikzcd}
\end{equation*}
and is defined by post-composing as
$(\alpha H)_x := H(\alpha_x): H(F(x)) \to H(G(x))$ for all $x \in \cat{C}$.
Finally, given pairs of parallel functors $F,G: \cat{C} \to \cat{D}$ and
$H,K: \cat{D} \to \cat{E}$, the \emph{horizontal composite}
$\alpha * \beta: FH \to GK$ of natural transformations $\alpha: F \to G$ and
$\beta: H \to K$ is drawn as
\begin{equation*}
  \begin{tikzcd}
    \cat{C}
      \arrow[r, bend left=45, "F"{above}, ""{name=F, below}]
      \arrow[r, bend right=45, "G"{below}, ""{name=G, above}]
    & \cat{D}
      \arrow[Rightarrow, "\alpha", from=F, to=G]
      \arrow[r, bend left=45, "H"{above}, ""{name=H, below}]
      \arrow[r, bend right=45, "K"{below}, ""{name=K, above}]
    & \cat{E}
      \arrow[Rightarrow, "\beta", from=H, to=K]
  \end{tikzcd}
  \qquad\leadsto\qquad
  \begin{tikzcd}[column sep=3em]
    \cat{C}
      \arrow[r, bend left=45, "FH"{above}, ""{name=FH, below}]
      \arrow[r, bend right=45, "GK"{below}, ""{name=GK, above}]
    & \cat{E}
      \arrow[Rightarrow, "\alpha * \beta", from=FH, to=GK]
  \end{tikzcd}.
\end{equation*}
It is defined by setting, for each $x \in \cat{C}$, the component
$(\alpha * \beta)_x: H(F(x)) \to K(G(x))$ to be the common composite in the
commutative square
\begin{equation*}
  \begin{tikzcd}[row sep=2.5em, column sep=2.5em]
    (FH)x
      \arrow[r, "(F \beta)_x"]
      \arrow[d, "(\alpha H)_x", swap]
      \arrow[dr, dashed, "(\alpha * \beta)_x"{anchor=center, fill=white}]
    & (FK)x
     \arrow[d, "(\alpha K)_x"] \\
    (GH)x
      \arrow[r, "(G \beta)_x", swap]
    & (GK)x
  \end{tikzcd}.
\end{equation*}
Standard lemmas in category theory affirm that this square does indeed commute,
and that the transformations defined by pre-whiskering, post-whiskering, and
horizontal compositional are all natural.

Vertical and horizontal composition each introduce an additional categorical
structure into $\Cat$, the category of categories. For any two categories
$\cat{C}$ and $\cat{D}$, the \emph{functor category}
$[\cat{C},\cat{D}] := \Cat(\cat{C},\cat{D})$ has functors $\cat{C} \to \cat{D}$
as objects and natural transformations between them as morphisms, with
composition provided by vertical composition and identities by the identity
transformations $1_F: F \to F$, where $(1_F)_x := 1_{Fx}$. Also, there is a
second category, besides $\Cat$ itself, with the small categories as objects,
but now the morphisms $\cat{C} \to \cat{D}$ are the natural transformations
$\alpha: F \to G$ between any functors $F,G: \cat{C} \to \cat{D}$. Composition
is provided by horizontal composition and identities by the transformations
$1_{1_{\cat{C}}}: 1_{\cat{C}} \to 1_{\cat{C}}$, where
$1_{\cat{C}}: \cat{C} \to \cat{C}$ is the usual identity functor. Finally,
vertical and horizontal composition commute with each other according to the law
of \emph{middle four interchange} \cite[Lemma 1.7.7]{riehl2016}. The effect of
all this to make $\Cat$ into not just a category, but a two-dimensional
categorical structure known as a \emph{2-category}. It is beyond the scope of
this text to define a 2-category in generality.

If, in categorical logic, functors serve as models, then natural transformations
must be model homomorphisms. Continuing the example, let
$G,H: \Theory{\Graph} \to \Set$ be graphs. A natural transformation
$\phi: G \to H$ consists of a vertex map, $\phi_V: G(V) \to H(V)$, and an edge
map, $\phi_E: G(E) \to H(E)$, making the two naturality squares commute:
\begin{equation*}
  \begin{tikzcd}
    G(E)
      \arrow[r, "\phi_E"]
      \arrow[d, "\src", swap] &
    H(E)
      \arrow[d, "\src"] \\
    G(V)
      \arrow[r, "\phi_V"] &
    H(V)
  \end{tikzcd}
  \qquad\qquad
  \begin{tikzcd}
    G(E)
    \arrow[r, "\phi_E"]
    \arrow[d, "\tgt", swap] &
    H(E)
    \arrow[d, "\tgt"] \\
    G(V)
    \arrow[r, "\phi_V"] &
    H(V)
  \end{tikzcd}.
\end{equation*}
This condition says that the edge map preserves the source and target vertices.
Thus, a natural transformation of graphs is simply a graph homomorphism, and
there is an isomorphism of categories
\begin{equation*}
  \Graph \cong [\Theory{\Graph},\Set].
\end{equation*}
In general, under the dictionary of categorical logic, functor categories are
categories of models and model homomorphisms.

\section{Monoidal categories and their graphical language}
\label{sec:monoidal-categories}

Although small categories can be interpreted as logical theories, the logical
system they comprise is not very expressive. None of the classical theories of
abstract algebra, such as the theories of groups or monoids, can be presented as
categories, since the logic of categories admits only unary operations. Binary
operations and operations of other arity require categories with extra
structure. Monoidal categories are the minimal elaboration of categories
furnishing that extra structure.

In categorical logic and beyond, monoidal categories are an algebraic theory of
operations, mappings, or other arrows that can have multiple inputs or outputs.
Morphisms in monoidal categories are often depicted by \emph{string diagrams},
also known as \emph{wiring diagrams}. Equations between morphisms can then be
proved by manipulating the diagrams according to the rules of a graphical
calculus that is intuitive yet fully rigorous. This combination of mathematical
precision and easy interpretability has made monoidal categories into an
indispensable tool of applied category theory.

Loosely speaking, a monoidal category is a category in which both objects and
morphisms can be juxtaposed ``in parallel'' through an associative, functorial
binary operation, the monoidal product. Like ordinary categories, monoidal
categories come in the small and in the large, with large categories typically
having more than one interesting monoidal structure. The category $\Set$ has two
important monoidal products, the cartesian product $\times$ and the disjoint
union $+$. As noted previously, the associative laws for set products and sums,
such as $(A \times B) \times C \cong A \times (B \times C)$, hold only up to
natural isomorphism. In order to accommodate these essential examples, the
general definition of a monoidal category must ``categorify'' the monoid laws.

\begin{definition}[Monoidal category] \label{def:monoidal-category}
  A \emph{monoidal category} $(\cat{C},\otimes,I)$, sometimes called a
  \emph{tensor category}, is a category $\cat{C}$ together with a binary
  operation $\otimes$, the \emph{monoidal product}, and a fixed object
  $I \in \cat{C}$, the \emph{monoidal unit}, such that every pair of objects $x$
  and $y$ has a product object $x \otimes y$ and every pair of morphisms
  $f: x \to y$ and $g: w \to z$ has a product morphism
  $f \otimes g: x \otimes w \to y \otimes z$, subject to the \emph{interchange
    laws}: $1_{x \otimes y} = 1_x \otimes 1_y$ for all $x,y \in \cat{C}$, and
  \begin{equation*}
    (f \otimes g) \cdot (h \otimes k) = (f \cdot h) \otimes (g \cdot k)
  \end{equation*}
  for all morphisms $u \xrightarrow{f} v \xrightarrow{h} w$ and
  $x \xrightarrow{g} y \xrightarrow{k} z$. Moreover, there are natural
  isomorphisms
  \begin{enumerate}
  \item (Unitors) $\lambda_x: I \otimes x \xrightarrow{\cong} x$ and
    $\rho_x: x \otimes I \xrightarrow{\cong} x$ for any $x \in \cat{C}$, and
  \item (Associators)
    $\alpha_{x,y,z}: (x \otimes y) \otimes z \xrightarrow{\cong} x \otimes (y
    \otimes z)$ for any $x,y,z \in \cat{C}$,
  \end{enumerate}
  which must satisfy two coherence axioms, not listed here. If all the unitors
  and associators are identities, the monoidal category is said to be
  \emph{strict}.
\end{definition}

\begin{remark}[Coherence]
  The two coherence axioms, known as the ``pentagon equation'' and the
  ``triangle equation,'' equate different compositions of natural isomorphisms
  for iterated monoidal products, such as
  $w \otimes (x \otimes (y \otimes z)) \cong (((w \otimes x) \otimes y) \otimes
  z)$. According to Mac Lane's coherence theorem for monoidal categories
  \cite[\S VII.2]{maclane1998}, a fundamental result, this short list of
  coherence axioms is enough to ensure that \emph{all} composites of associators
  and unitors commute and hence any two bracketings of an iterated product of
  objects are canonically isomorphic. A consequence is that every monoidal
  category is monoidally equivalent to a strict monoidal category. In practical
  terms, we can usually pretend that all monoidal categories are strict, a
  conceit exploited by the string diagram calculus.
\end{remark}

As already noted, the category of sets is a monoidal category in two different
ways, as $(\Set,\times,1)$, where the monoidal unit $1$ is the singleton set
$\{*\}$, and as $(\Set,+,0)$, where the unit $0$ is the empty set $\emptyset$.
The category of vector spaces over a field $\K$ is also a monoidal category in
two different ways, as $(\Vect_\K,\otimes,\K)$, where the monoidal product is
the tensor product and the unit is the one-dimensional vector space $\K$, and as
$(\Vect_\K,\oplus,0)$, where the monoidal product is the direct sum and the unit
is the zero-dimensional vector space $\{0\}$. The functor category
$[G,\Vect_\K]$ of linear representations of a group $G$ is a monoidal category
under the tensor product of representations: the product of representations
$\rho: G \to \GL(V)$ and $\tau: G \to \GL(W)$ is defined pointwise as
$\phi \otimes \tau: G \to \GL(V \otimes W)$,
$g \mapsto \rho(g) \otimes \tau(g)$. Similarly, this category is a monoidal
category under the direct sum of representations. None of these monoidal
categories are strict.

A monoidal product need not resemble the set-theoretic product or coproduct. For
example, given a category $\cat{C}$, there is a category
$\End(\cat{C}) := [\cat{C},\cat{C}]$ of endofunctors of $\cat{C}$ and natural
transformations between them. As in any functor category, composition is
vertical composition of natural transformations. But the endofunctor category
$\End(\cat{C})$ is moreover a (strict) monoidal category, with the monoidal
product defined on objects by composition, $F \otimes G := F \cdot G$, and on
morphisms by horizontal composition, $\alpha \otimes \beta := \alpha * \beta$.

In the graphical language of string diagrams, a morphism $f: x \to y$ is
represented by a box labeled ``$f$'', with an incoming wire labeled ``$x$'' and
an outgoing wire labeled ``$y$'':
\begin{center}
  \input{wiring-diagrams/category-theory/box}
\end{center}
The composite $f \cdot g: x \to z$ of morphisms $f: x \to y$ and $g: y \to z$ is
represented by juxtaposition in series:
\begin{equation*}
  \input{wiring-diagrams/category-theory/compose-literal} =
  \input{wiring-diagrams/category-theory/compose}
\end{equation*}
The product $f \otimes g: x \otimes w \to y \otimes z$ of morphisms $f: x \to y$
and $g: w \to z$ is represented by juxtaposition in parallel:
\begin{equation*}
  \input{wiring-diagrams/category-theory/product-literal} =
  \input{wiring-diagrams/category-theory/product}
\end{equation*}
Identity morphisms $1_x: x \to x$ are drawn simply as wires, and the monoidal
unit $I$ is not drawn at all. As a result of these conventions, the
associativity and unitality laws, of both composition and monoidal products, are
fully implicit in the graphical syntax. The interchange laws, relating
composition and products, are also implicit.

In most monoidal categories encountered in practice, it is possible to permute
the objects in a product $x \otimes y$ to obtain the product $y \otimes x$. Such
monoidal categories are called \emph{symmetric}.

\begin{definition}[Symmetric monoidal category]
  \label{def:symmetric-monoidal-category}
  A \emph{symmetric monoidal category} is a monoidal category
  $(\cat{C},\otimes,I)$ together with natural isomorphisms
  $\Braid_{x,y}: x \otimes y \to y \otimes x$, $x,y \in \cat{C}$, called
  \emph{braidings} or \emph{symmetries} and depicted as crossed wires:
  \begin{equation*}
    \input{wiring-diagrams/category-theory/braid-literal} =
    \input{wiring-diagrams/category-theory/braid}
  \end{equation*}
  The braidings must satisfy an involutivity axiom,
  $\Braid_{x,y}^{-1} = \Braid_{y,x}$ for all $x,y \in \cat{C}$, or
  \begin{equation*}
    \input{wiring-diagrams/category-theory/braid-involution-lhs} =
    \input{wiring-diagrams/category-theory/braid-involution-rhs},
  \end{equation*}
  as well as two coherence axioms,
  \begin{equation*}
    \input{wiring-diagrams/category-theory/braid-hexagon1-lhs} =
    \input{wiring-diagrams/category-theory/braid-hexagon1-rhs}
    \quad\text{and}\quad
    \input{wiring-diagrams/category-theory/braid-hexagon2-lhs} =
    \input{wiring-diagrams/category-theory/braid-hexagon2-rhs},
  \end{equation*}
  asserting that certain braidings for product objects can be constructed out of
  the braidings for the original objects.
\end{definition}

All the monoidal categories listed above are actually symmetric monoidal
categories, with the exception of the endofunctor categories $\End(\cat{C})$.

That the braidings in symmetric monoidal category are \emph{natural}
isomorphisms means that for any two morphisms $f: x \to y$ and $g: w \to z$,
\begin{equation*}
  \input{wiring-diagrams/category-theory/braid-naturality-lhs} =
  \input{wiring-diagrams/category-theory/braid-naturality-rhs}.
\end{equation*}
Taken together, the axioms of a symmetric monoidal category imply that for every
finite set of objects $x_1,\dots,x_n$ and every permutation
$\sigma: \{1,\dots,n\} \to \{1,\dots,n\}$, there exists a canonical isomorphism
$x_1 \otimes \cdots \otimes x_n \to x_{\sigma(1)} \otimes \cdots \otimes
x_{\sigma(n)}$ constructed out of braidings, identities, associators, and
unitors. This is the content of Mac Lane's coherence theorem for symmetric
monoidal categories \cite[\S XI.1]{maclane1998}.

Small monoidal categories are most often defined by the method of generators and
relations. We are free to assume that monoidal categories presented by
generators and relations are strict, and so we shall.

\begin{example}[Theory of monoids] \label{ex:monoids}
  The \emph{theory of monoids} $\Theory{\Mon}$ is the monoidal category
  generated by one object $x$ and two morphisms $\mu: x \otimes x \to x$ and
  $\eta: I \to x$, depicted as
  \begin{equation*}
    \input{wiring-diagrams/category-theory/mmerge} \qquad\text{and}\qquad
    \input{wiring-diagrams/category-theory/create},
  \end{equation*}
  and subject to the equations of associativity and unitality,
  \begin{equation*}
    \input{wiring-diagrams/category-theory/monoid-associative-lhs} =
    \input{wiring-diagrams/category-theory/monoid-associative-rhs}
    \qquad\text{and}\qquad
    \input{wiring-diagrams/category-theory/monoid-unital-lhs} =
    \input{wiring-diagrams/category-theory/monoid-id} =
    \input{wiring-diagrams/category-theory/monoid-unital-rhs}.
  \end{equation*}
  The \emph{theory of comonoids} $\Theory{\Comon}$ is the opposite category
  $\Theory{\Mon}^\opposite$. Explicitly, it is the monoidal category generated
  by one object $x$ and two morphisms $\delta: x \to x \otimes x$ and
  $\epsilon: x \to I$, depicted as
  \begin{equation*}
    \input{wiring-diagrams/category-theory/mcopy} \qquad\text{and}\qquad
    \input{wiring-diagrams/category-theory/delete},
  \end{equation*}
  and subject to the mirror images of the monoid equations.
\end{example}

\begin{example}[Theory of commutative monoids] \label{ex:commutative-monoids}
  The \emph{theory of commutative monoids} $\Theory{\CMon}$ is the symmetric
  monoidal category with the same presentation as $\Theory{\Mon}$ but augmented
  with the commutativity equation
  \begin{equation*}
    \input{wiring-diagrams/category-theory/monoid-commutative} =
    \input{wiring-diagrams/category-theory/mmerge}.
  \end{equation*}
  The \emph{theory of commutative comonoids} $\Theory{\CComon}$ is the opposite
  category $\Theory{\CMon}^\opposite$.
\end{example}

\begin{example}[Theory of bimonoids] \label{ex:bimonoids}
  The \emph{theory of bimonoids} $\Theory{\Bimon}$ is the symmetric monoidal
  category generated by one object $x$ and four morphisms
  $\mu: x \otimes x \to x$, $\eta: I \to x$, $\delta: x \to x \otimes x$, and
  $\epsilon: x \to I$, depicted as
  \begin{equation*}
    \input{wiring-diagrams/category-theory/mmerge}, \quad
    \input{wiring-diagrams/category-theory/create}, \quad
    \input{wiring-diagrams/category-theory/mcopy}, \quad\text{and}\quad
    \input{wiring-diagrams/category-theory/delete}.
  \end{equation*}
  They are subject to the same laws of associativity and unitality as in the
  theories of monoids and comonoids, plus the \emph{bimonoid laws}:
  \begin{gather*}
    \input{wiring-diagrams/category-theory/bimonoid-copy-merge-lhs} =
    \input{wiring-diagrams/category-theory/bimonoid-copy-merge-rhs}
    \\
    \input{wiring-diagrams/category-theory/bimonoid-delete-merge-lhs} =
    \input{wiring-diagrams/category-theory/bimonoid-delete-merge-rhs}
    \qquad\qquad
    \input{wiring-diagrams/category-theory/bimonoid-copy-create-lhs} =
    \input{wiring-diagrams/category-theory/bimonoid-copy-create-rhs}
    \qquad\qquad
    \input{wiring-diagrams/category-theory/bimonoid-delete-create} = \qquad.
  \end{gather*}
  The blank in the last equation represents the identity morphism on the
  monoidal unit $I$. Together, the bimonoid laws say that comonoid operations
  are monoid homomorphisms or, equivalently, that the monoid operations are
  comonoid homomorphisms. The \emph{theory of bicommutative bimonoids}
  $\Theory{\CBimon}$ extends the presentation of $\Theory{\Bimon}$ with
  commutativity laws for both of $\mu: x \otimes x \to x$ and
  $\delta: x \to x \otimes x$, as in the previous example.
\end{example}

A strict symmetric monoidal category whose monoid of objects is freely generated
by a single object is called a \emph{PROP} (``products and permutations
category''). A strict monoidal category satisfying the same condition is called
a \emph{PRO} (``products category''). The theories of monoids and comonoids are
PROs, while the theories of commutative monoids and commutative comonoids are
PROPs. The theory of bimonoids, even in the noncommutative case, is also a PROP,
as the first bimonoid law involves a braiding. Since the free monoid on one
object is isomorphic to the natural numbers $(\N,+,0)$, the objects of a PRO or
PROP can be identified with natural numbers. The theory of monoids, for example,
then has generators $\mu: 2 \to 1$ and $\eta: 0 \to 1$.

As these examples suggest, PROs and PROPs, and small monoidal categories
generally, are often regarded as logical theories. The logical analogy is
completed by variants of functors and natural transformations that respect the
extra structure of monoidal categories.

A monoidal functor of monoidal categories ought to preserve the monoidal product
and unit, suggesting that a monoidal functor $F$ should satisfy
$F(x \otimes y) = F(x) \otimes F(y)$ and $F(I) = I$. However, in non-strict
monoidal categories, this condition is often too stringent. Consider the duality
functor $(-)^*: \Vect_\K^\opposite \to \Vect_\K$ that a sends a vector space $V$
to its dual space $V^* := \Vect_\K(V,\K)$ and a linear map $f: V \to W$ to its
transpose $f^*: W^* \to V^*$. Is this functor monoidal? The spaces
$(V \otimes W)^*$ and $V^* \otimes W^*$ are not strictly equal, but there is a
natural isomorphism $(V \otimes W)^* \cong V^* \otimes W^*$ obtained from the
unitor isomorphism $k \cong k \otimes k$. This motivates the general definition
of a monoidal functor.

\begin{definition}[Monoidal functor] \label{def:monoidal-functor}
  A \emph{(strong) monoidal functor} between monoidal categories
  $(\cat{C},\otimes_{\cat{C}},I_{\cat{C}})$ and
  $(\cat{D},\otimes_{\cat{D}},I_{\cat{D}})$ is a functor
  $F: \cat{C} \to \cat{D}$ together with natural isomorphisms
  $\Phi_{x,y}: F(x) \otimes_{\cat{D}} F(y) \to F(x \otimes_{\cat{C}} y)$, for
  $x,y \in \cat{C}$, and an isomorphism $\phi: I_{\cat{D}} \to F(I_{\cat{C}})$,
  satisfying an \emph{associativity law}
  \begin{equation*}
    \input{wiring-diagrams/category-theory/monoidal-functor-associative-lhs} =
    \input{wiring-diagrams/category-theory/monoidal-functor-associative-rhs},
  \end{equation*}
  where, in the diagrammatic notation, we suppress an associator
  $\alpha^{\cat{D}}_{F(x),F(y),F(z)}$ on the top and an associator
  $F(\alpha^{\cat{C}}_{x,y,z})$ on the bottom, and also a \emph{unitality law}
  \begin{equation*}
    \input{wiring-diagrams/category-theory/monoidal-functor-unital-lhs} =
    \input{wiring-diagrams/category-theory/monoidal-functor-id} =
    \input{wiring-diagrams/category-theory/monoidal-functor-unital-rhs},
  \end{equation*}
  where we suppress unitors $\lambda_{F(x)}^{\cat{D}}$ and
  $\rho_{F(x)}^{\cat{D}}$ on the top and unitors $F(\lambda_x^{\cat{C}})$ and
  $F(\rho_x^{\cat{C}})$ on the bottom. If all of the isomorphisms $\Phi_{x,y}$
  and $\phi$ are identities, then the monoidal functor is called \emph{strict}.

  A \emph{symmetric monoidal functor} is a monoidal functor of symmetric
  monoidal categories such that the structure morphisms $\Phi_{x,y}$ commute
  with the braidings:
  \begin{equation*}
    \begin{tikzcd}[column sep=large]
      F(x) \otimes_{\cat{D}} F(y)
        \arrow[r, "\Braid_{F(x),F(y)}^{\cat{D}}"]
        \arrow[d, swap, "\Phi_{x,y}"] &
      F(y) \otimes_{\cat{D}} F(x)
        \arrow[d, "\Phi_{y,x}"] \\
      F(x \otimes_{\cat{C}} y)
        \arrow[r, swap, "F(\Braid_{x,y}^{\cat{C}})"] &
      F(y \otimes_{\cat{C}} x)
    \end{tikzcd}.
  \end{equation*}
\end{definition}

The situation for monoidal natural transformations is similar. A natural
transformation $\alpha$ of monoidal functors ought to preserve monoidal
products, suggesting that $\alpha_{x \otimes y} = \alpha_x \otimes \alpha_y$;
however, if the monoidal functors are not strict, then the structure morphisms
must be accounted for.

\begin{definition}[Monoidal natural transformation]
  \label{def:monoidal-natural-transformation}
  Let $(F,\Phi,\phi)$ and $(G,\Gamma,\gamma)$ be parallel monoidal functors
  between monoidal categories $(\cat{C},\otimes_{\cat{C}},I_{\cat{C}})$ and
  $(\cat{D},\otimes_{\cat{D}},I_{\cat{D}})$. A \emph{monoidal natural
    transformation} between $F$ and $G$ is a natural transformation
  $\alpha: F \to G$ such that the two diagrams commute:
  \begin{equation*}
    \begin{tikzcd}[column sep=large]
      F(x) \otimes_{\cat{D}} F(y)
        \arrow[r, "\alpha_x \otimes_{\cat{D}} \alpha_y"]
        \arrow[d, swap, "\Phi_{x,y}"] &
      G(x) \otimes_{\cat{D}} G(y)
        \arrow[d, "\Gamma_{x,y}"] \\
      F(x \otimes_{\cat{C}} y)
        \arrow[r, swap, "\alpha_{x \otimes_{\cat{C}} y}"] &
      G(x \otimes_{\cat{C}} y)
    \end{tikzcd}
    \qquad\qquad
    \begin{tikzcd}[column sep=small]
      & I_{\cat{D}}
        \arrow[dl, swap, "\phi"]
        \arrow[dr, "\gamma"] & \\
      F(I_{\cat{C}})
        \arrow[rr, swap, "\alpha_{I_{\cat{C}}}"] & &
      G(I_{\cat{C}})
    \end{tikzcd}.
  \end{equation*}
  No extra condition is needed in the symmetric case: a \emph{symmetric monoidal
    natural transformation} is just a monoidal natural transformation between
  symmetric monoidal functors.
\end{definition}

Monoidal categories and monoidal functors assemble into a category $\MonCat$.
Likewise, symmetric monoidal categories and symmetric monoidal functors form a
category $\SMonCat$. The subcategories $\PRO$ and $\PROP$ of $\MonCat$ and
$\SMonCat$ have as objects PRO(P)s and as morphisms the strict (symmetric)
monoidal functors preserving the object generators. When the monoidal natural
transformations are included, all these categories even become 2-categories,
although we have not said exactly what that means.

In the dictionary of categorical logic, a \emph{monoidal theory} is a small,
strict monoidal category and a \emph{model} of a monoidal theory $\cat{C}$ is a
monoidal functor $\cat{C} \to (\Set,\times,1)$. Generalizing, a \emph{model} of
$\cat{C}$ in a monoidal category $\cat{S}$ is a monoidal functor
$\cat{C} \to \cat{S}$. When $\cat{C}$ is a PRO, its models are often called
\emph{algebras} of $\cat{C}$. We adopt the convention that if $\cat{C}$ is a
PRO, whose objects we identify with $\N$, then a model $F: \cat{C} \to \cat{S}$
always takes the left-associative form
\begin{equation*}
  F(n) = x^{\otimes n} :=
    \underbrace{((x \otimes x) \dots) \otimes x}_{n\text{ times}},
\end{equation*}
where the structure isomorphisms
$\Phi_{m,n}: x^{\otimes m} \otimes x^{\otimes n} \to x^{\otimes (m+n)}$ are the
unique coherence isomorphisms made out of associators and
unitors.\footnote{Restricting the form of the monoidal functor is not
  mathematically significant but does ensure that models of a PRO are in
  one-to-one correspondence with models as conventionally understood. The
  particular functional form chosen here, while arbitrary, agrees with that of
  \cite{fong2019supply}.} The \emph{category of models} of $\cat{C}$ in
$\cat{S}$, denoted $[\cat{C},\cat{S}]$, has models $\cat{C} \to \cat{S}$ as
objects and monoidal natural transformations as morphisms.

These definitions and conventions also apply, with obvious modifications, to
symmetric monoidal categories and PROPs.

\begin{example}[Monoid objects] \label{ex:monoid-objects}
  Many familiar algebraic structures can be reconstructed as models of the
  theory of monoids $\Theory{\Mon}$ in a suitable monoidal category $\cat{S}$.
  Such models are called \emph{monoids}, or \emph{monoid objects}, in $\cat{S}$
  \cite[Chapter 15]{street2007}. A monoid in $(\Set,\times,1)$ is just a monoid.
  A monoid in $(\Top,\times,1)$ is a topological monoid, a weakened version of a
  topological group. A monoid in $(\Ab,\otimes,\Z)$, the tensor category of
  abelian groups, is a (unital) ring, whereas a monoid in $(\CMon,\otimes,\N)$,
  the tensor category of commutative monoids, is a \emph{rig} (``ring without
  negatives''). A monoid in $(\Vect_\K,\otimes,\K)$ is an (associative, unital)
  algebra over the field $\K$. Somewhat circularly, a monoid in
  $(\Cat,\times,1)$ is a strict monoidal category. As a more exotic example, a
  monoid in $(\End(\cat{C}),\cdot,1_{\cat{C}})$ is known as a \emph{monad} on
  the category $\cat{C}$.
\end{example}

\begin{example}[Commutative monoid objects]
  \label{ex:commutative-monoid-objects}
  In all preceding examples except the last, the monoidal categories involved
  are actually symmetric, so we can consider commutative monoids in them to
  obtain commutative rings, commutative rigs, commutative algebras, and so on.
  Be warned that a \emph{commutative monoidal category}, namely a commutative
  monoid in $(\Cat,\times,1)$, is much stricter than a strict symmetric monoidal
  category, because the equations $x \otimes y = y \otimes x$ and
  $f \otimes g = g \otimes f$ must hold strictly for all objects $x,y$ and all
  morphisms $f,g$. Commutative monoidal categories are rare, although they occur
  naturally in the operational semantics of Petri nets \cite{baez2018}.
\end{example}

The setting for much statistical modeling is a real vector space or a structured
subset thereof, such as an affine subspace, a convex cone, or a convex set.
Algebraically, each kind of set admits a different kind of combination of its
elements: linear, affine, conical, or convex. A uniform treatment of the algebra
of combinations is enabled by working in a generic commutative rig $R$. In
\cref{ex:monoid-objects}, a \emph{rig}, better known as a \emph{semiring}, was
defined as a monoid in the category of commutative monoids. A rig is like a
ring, except that its elements need not have additive inverses. Two important
rigs are the real numbers $\R$ and the nonnegative real numbers $\R_+$, each
with their usual operations of addition and multiplication.

\begin{example}[Theories of linear combinations] \label{ex:linear-combinations}
  For any commutative rig $R$, the \emph{theory of $R$-linear combinations} is
  the PROP $\LinComb_R$ generated by the hom-sets
  \begin{equation*}
    \LinComb_R(x^{\otimes n},x) :=
      R^n = \{(r_i)_{i=1}^n = (r_1,\dots,r_n): r_i \in R\},
    \qquad n \geq 0,
  \end{equation*}
  subject to the equations:
  \begin{enumerate}
  \item (Distributivity) For all $n, m_1, \dots, m_n \geq 0$ and all
    $r \in R^n, s_1 \in R^{m_i}, \dots, s_n \in R^{m_n}$,
    \begin{equation*}
      ((s_{1,j})_{j=1}^{m_1} \otimes \cdots \otimes (s_{n,j})_{j=1}^{m_n}) \cdot
      (r_i)_{i=1}^n =
      (r_i \cdot s_{i,j})_{j=1,i=1}^{m_i,n} \in R^{m_1+\cdots+m_n}.
    \end{equation*}
  \item (Equivariance) For all $n \geq 0$, permutations $\sigma \in S_n$, and
    $r \in R^n$,
    \begin{equation*}
      \sigma \cdot (r_1,\dots,r_n) = (r_{\sigma(1)},\dots,r_{\sigma(n)}),
    \end{equation*}
    where we identify a permutation $\sigma: \{1,\dots,n\} \to \{1,\dots,n\}$
    with the corresponding symmetry isomorphism
    $\Braid: x^{\otimes n} \xrightarrow{\cong} x^{\otimes n}$.
  \end{enumerate}
  In particular, the singleton list $(1)$, where $1$ is the multiplicative unit
  in $R$, is the identity morphism $x \to x$, and the empty list $()$ is the
  unique morphism $I \to x$. As special cases, the \emph{theory of (real) linear
    combinations} is the PROP $\LinComb_\R$ and the \emph{theory of conical
    combinations} is the PROP $\ConeComb := \LinComb_{\R_+}$.

  The theory of $R$-linear combinations admits a presentation with fewer
  generators and relations. Take the theory of commutative monoids
  (\cref{ex:commutative-monoids}) and add a generator $r: x \to x$ for each
  $r \in R$, representing scalar multiplication by $r$. Add equations asserting
  that multiplication in $R$ is respected by composition,
  \begin{equation*}
    \input{wiring-diagrams/category-theory/combination-product-lhs} =
    \input{wiring-diagrams/category-theory/combination-product-rhs}
    \qquad\text{and}\qquad
    \input{wiring-diagrams/category-theory/combination-unit-lhs} =
    \input{wiring-diagrams/category-theory/combination-unit-rhs},
  \end{equation*}
  and that the scalar multiplications are monoid homomorphisms,
  \begin{equation*}
    \input{wiring-diagrams/category-theory/combination-mmerge-lhs} =
    \input{wiring-diagrams/category-theory/combination-mmerge-rhs}
    \qquad\text{and}\qquad
    \input{wiring-diagrams/category-theory/combination-create-lhs} =
    \input{wiring-diagrams/category-theory/combination-create-rhs}.
  \end{equation*}
  It can be shown that these equations present $\LinComb_R$ (cf.\ \cite[Theorem
  3.7]{giraudo2015}). The presentation is said to be \emph{biased} because it
  favors operations of certain arities, specifically the binary and nullary
  monoid operations and the unary scalar multiplications. In comparison with the
  original \emph{unbiased} definition, the monoid multiplication
  $\mu: x \otimes x \to x$ corresponds to the list $(1,1)$ and the monoid unit
  $\eta: I \to x$ to the empty list $()$. The scalar multiplication $r: x \to x$
  corresponds to the singleton list $(r)$.
\end{example}

\begin{example}[Theories of affine combinations] \label{ex:affine-combinations}
  For any commutative rig $R$, the \emph{theory of $R$-affine combinbations} is
  the PROP $\AffComb_R$ generated by the hom-sets
  \begin{equation*}
    \AffComb_R(x^{\otimes n},x) :=
      \{(r_i)_{i=1}^n \in R^n: r_1 + \cdots + r_n = 1\}
  \end{equation*}
  and subject to the same equations of distributivity and equivariance as in
  $\LinComb_R$. Thus, $\AffComb_R$ is a sub-PROP of $\LinComb_R$. As special
  cases, the \emph{theory of (real) affine combinations} is the PROP
  $\AffComb_\R$ and the \emph{theory of convex combinations} is the PROP
  $\ConvComb := \AffComb_{\R_+}$. At least in these cases, where $R$ equals $\R$
  or $\R_+$, the theory of $R$-affine combinations also admits a biased
  presentation via morphisms $\mu_{r,s}: x \otimes x \to x$ parameterized by
  numbers $r,s \in R$ such that $r+s=1$ or, since $s=1-r$, simply via morphisms
  $\mu_r: x \otimes x \to x$ parameterized by $r \in R$. The details are not
  given here.
\end{example}

The models of the various theories of combinations are more general than one
might first expect. A space with linear combinations, or model of $\LinComb_\R$,
is a set $X$ equipped with operations of addition $+: X \times X \to X$ and
scalar multiplication $a \cdot -: X \to X$, $a \in \R$, and a zero element
$0 \in X$, such that $(X,+,0)$ is a commutative monoid and
\begin{equation*}
  a \cdot (x + y) = ax + ay, \qquad
  a \cdot 0 = 0, \qquad
  (a \cdot b) \cdot x = a \cdot (b \cdot x), \qquad
  1 \cdot x = x,
\end{equation*}
for all $x,y \in X$ and $a,b, \in \R$. Any real vector space has linear
combinations in this sense, but vector spaces are not the only examples, because
the monoid $(X,+,0)$ need not have additive inverses and, more importantly,
because the equations
\begin{equation*}
  (a + b) \cdot x = ax + bx \qquad\text{and}\qquad
  0 \cdot x = 0
\end{equation*}
need not be satisfied.

Examples not satisfying the extra equations occur surprisingly naturally. Let
$X$ be a random variable taking values in a vector space $V$. Its first moment
vector $\E[X]$, second moment matrix $\E[X \otimes X]$, and tensors of
higher-order moments $\E[X^{\otimes n}]$, $n \geq 3$, all inhabit the space
$\Sym(V) := \bigoplus_{n \in \N} \Sym^n(V) \subseteq \bigoplus_{n \in \N}
V^{\otimes n}$ of symmetric tensors on $V$. Multiplying the random vector $X$ by
a scalar $a \in \R$ transforms its moments according to the rule
$x^{\otimes n} \mapsto (a x)^{\otimes n} = a^n x^{\otimes n}$. The corresponding
action of $\R$ on the symmetric tensor space $\Sym(V)$ models the theory of
linear combinations, yet this scalar multiplication is plainly nonlinear in $a$.
Considering the central moments or cumulants of $X$, instead of the non-central
moments, leads to the same conclusion.

Nevertheless, it is certainly important to be able to express stronger theories,
such as those of abelian groups or vector spaces. What the missing equations all
have in common is that some variable $x$ appears multiple times on one side of
an equation, as in $x + (-x) = 0$ or $(a+b) \cdot x = ax + bx$, or not at all on
one side of an equation, as in $0 \cdot x = 0$. Indeed, a central feature of
PROs and PROPs, and monoidal theories generally, is that each variable appearing
in an equation must appear exactly once on each side of that equation. Lifting
this restriction leads to algebraic theories and cartesian categories,
historically the original setting of categorical logic.

\section{Cartesian categories and algebraic theories}
\label{sec:cartesian-categories}

Most of the classical structures of abstract algebra, such as groups, rings,
modules, and associative algebras, are axiomatizable in the purely equational
form of \emph{algebraic theories}.\footnote{A notable exception is the theory of
  fields, which is not algebraic because the operation of division is undefined
  at zero.} The logic of algebraic theories was first studied in generality,
independent of any particular structure, under the name of \emph{universal
  algebra} \cite{burris1981}. Categorical logic was, in turn, born as an
amplification of universal algebra, unifying algebraic theories and their models
with concepts from category theory \cite{lawvere1963}. Under this
correspondence, algebraic theories are represented by \emph{cartesian
  categories}, which are symmetric monoidal categories whose monoidal products
enjoy a certain universal property or, equivalently, whose objects are equipped
with natural operations for copying and discarding data. In this section, the
rudiments of cartesian categories and algebraic theories are developed in a
style that is somewhat unconventional, yet is more easily mutated than the
classical formulation.

Operations for copying and discarding data can be defined succinctly as
commutative comonoids. Suppose a symmetric monoidal category
$(\cat{C},\otimes,I)$ is equipped, at every object $x \in \cat{C}$, with
morphisms $\Copy_x: x \to x \otimes x$ and $\Delete_x: x \to I$ making $x$ into
a commutative comonoid object (cf.\ \cref{ex:commutative-monoid-objects}). The
morphism $\Copy_x: x \to x \otimes x$ is interpreted as \emph{copying} or
\emph{duplication} and the morphism $\Delete_x: x \to I$ as \emph{deleting} or
\emph{discarding}, and they are depicted as
\begin{equation*}
  \input{wiring-diagrams/category-theory/supply-mcopy}
  \qquad\text{and}\qquad
  \input{wiring-diagrams/category-theory/supply-delete}.
\end{equation*}
In $(\Set,\times,1)$, these maps are, for any set $X$, the diagonal map
$\Copy_X: x \mapsto (x,x)$ and the terminal map $\Delete_X: x \mapsto *$.

The copying and deleting operations defined at each object should compatible
with the monoidal product. Thus, copying the product $x \otimes y$ should be the
same as copying $x$ and in parallel copying $y$, up to permutation of the
outputs, and deleting the product $x \otimes y$ should be the same as deleting
$x$ and deleting $y$. To be more precise, the equations
\begin{equation*}
  \input{wiring-diagrams/category-theory/supply-coherence-mcopy-lhs} =
  \input{wiring-diagrams/category-theory/supply-coherence-mcopy-rhs}
  \qquad\text{and}\qquad
  \input{wiring-diagrams/category-theory/supply-coherence-delete-lhs} =
  \input{wiring-diagrams/category-theory/supply-coherence-delete-rhs}
\end{equation*}
should hold for all objects $x,y$. Copying and deleting should also be
compatible with the monoidal unit, in that
$I \xrightarrow{\Copy_I} I \otimes I \xrightarrow{\cong} I$ and
$I \xrightarrow{\Delete_I} I$ are both equal to the identity $1_I$. When these
conditions are satisfied, the category $\cat{C}$ is said to \emph{supply
  commutative comonoids}.

The notion of a supply of commutative comonoids is usefully generalized to a
supply of an arbitrary PROP \cite{fong2019supply}.

\begin{definition}[Supply] \label{def:supply}
  Let $\cat{P}$ be a PROP, assumed to have objects $(\N,+,0)$, and let $\cat{C}$
  be a symmetric monoidal category. A \emph{supply of $\cat{P}$ in $\cat{C}$}
  consists of, for each object $x \in \cat{C}$, a strong symmetric monoidal
  functor $s_x: \cat{P} \to \cat{C}$, such that
  \begin{enumerate}[(i),nosep]
  \item $s_x(n) = x^{\otimes n}$ for each $n \in \N$,
  \item the structure isomorphism
    $(\Phi_x)_{m,n}: x^{\otimes m} \otimes x^{\otimes n} \to x^{\otimes(m+n)}$ is
    the unique coherence isomorphism for each $m,n \in \N$, and
  \item for every $x,y \in \cat{C}$ and every morphism $\mu: m \to n$ in
    $\cat{P}$, the diagrams
    \begin{equation*}
      \begin{tikzcd}[column sep=huge]
        x^{\otimes m} \otimes y^{\otimes m}
          \arrow[r, "s_x(\mu) \otimes s_y(\mu)"]
          \arrow[d, swap, "\Braid"] &
        x^{\otimes n} \otimes y^{\otimes n}
          \arrow[d, "\Braid"] \\
        (x \otimes y)^{\otimes m}
          \arrow[r, "s_{x \otimes y}(\mu)"] &
        (x \otimes y)^{\otimes n}
      \end{tikzcd}
      \qquad\qquad
      \begin{tikzcd}[column sep=small]
        & I
          \arrow[dl, swap, "\Braid"]
          \arrow[dr, "\Braid"] & \\
        I^{\otimes m}
          \arrow[rr, "s_I(\mu)"] & &
        I^{\otimes n}
      \end{tikzcd}
    \end{equation*}
    commute, where the $\Braid$'s are canonical symmetry isomorphisms.
  \end{enumerate}
\end{definition}

The theory of commutative comonoids was presented as a PROP in
\cref{ex:commutative-monoids}. Naturally, a supply of this PROP is a supply of
commutative comonoids as originally defined.

The first two axioms of a supply are merely conventions, as explained in
\cref{sec:monoidal-categories}. The essential mathematical content lies in the
third axiom, requiring that the models of $\cat{P}$ defined at each object of
$\cat{C}$ be compatible with the symmetric monoidal structure of $\cat{C}$. One
can also ask that some or all morphisms in $\cat{C}$ commute with the morphisms
supplied by $\cat{P}$, leading to the notion of supply homomorphism.

\begin{definition}[Homomorphic supply] \label{def:homomorphic-supply}
  Let $s$ be a supply of a PROP $\cat{P}$ in a symmetric monoidal category
  $\cat{C}$. A morphism $f: x \to y$ in $\cat{C}$ is an \emph{$s$-homomorphism}
  if for every morphism $\mu: m \to n$ in $\cat{P}$, the diagram
  \begin{equation*}
    \begin{tikzcd}
      x^{\otimes m} \arrow[r, "s_x(\mu)"] \arrow[d, swap, "f^{\otimes m}"] &
        x^{\otimes n} \arrow[d, "f^{\otimes n}"] \\
      y^{\otimes m} \arrow[r, swap, "s_y(\mu)"] &
        y^{\otimes n}
    \end{tikzcd}
  \end{equation*}
  commutes. If every morphism in $\cat{C}$ is an $s$-homomorphism, then $s$ is a
  \emph{homomorphic supply}.
\end{definition}

In a category $\cat{C}$ supplying commutative comonoids, a morphism $f: x \to y$
is a supply homomorphism if it commutes with copying and deleting:
\begin{equation*}
  \input{wiring-diagrams/category-theory/supply-homomorphism-mcopy-lhs} =
  \input{wiring-diagrams/category-theory/supply-homomorphism-mcopy-rhs}
  \qquad\text{and}\qquad
  \input{wiring-diagrams/category-theory/supply-homomorphism-delete-lhs} =
  \input{wiring-diagrams/category-theory/supply-homomorphism-delete-rhs}.
\end{equation*}
The supply is homomorphic if these equations hold for all morphisms $f$ in
$\cat{C}$ or, equivalently, the copying morphisms $(\Copy_x)_{x \in \cat{C}}$
are the components of a natural transformation from the identity functor
$1_{\cat{C}}$ to the diagonal functor
$\Copy_{\cat{C}}: \cat{C} \to \cat{C}, {x \mapsto x \otimes x \atop f \mapsto f
  \otimes f}$ and, similarly, the deleting morphisms
$(\Delete_x)_{x \in \cat{C}}$ are the components of a natural transformation
from the identity $1_{\cat{C}}$ to the functor
$\Delete_{\cat{C}}: \cat{C} \to \cat{C}, {x \mapsto I \atop f \mapsto 1_I}$. In
general, a supply $s$ of a PROP $\cat{P}$ in $\cat{C}$ is homomorphic if and
only if for every morphism $\mu: m \to n$ in $\cat{P}$, the morphisms
$(s_x(\mu))_{x \in \cat{C}}$ assemble into a natural transformation between the
endofunctors $(-)^{\otimes m}$ and $(-)^{\otimes n}$ of $\cat{C}$.

\begin{definition}[Cartesian category] \label{def:cartesian-category}
  A \emph{cartesian monoidal category}, or a \emph{cartesian category} for
  short, is a symmetric monoidal category that homomorphically supplies
  commutative comonoids. Dually, a \emph{cocartesian (monoidal) category} is a
  symmetric monoidal category that homomorphically supplies commutative monoids.
\end{definition}

For any supply $s$ in $\cat{C}$, the $s$-homomorphisms form a symmetric monoidal
subcategory of $\cat{C}$, denoted $\cat{C}_s$ \cite[Theorem
3.16]{fong2019supply}. It is the largest subcategory of $\cat{C}$ on which the
supply is homomorphic. When $s$ is a supply of commutative comonoids, this
subcategory is known as the \emph{cartesian center} of $\cat{C}$
\cite{hasegawa1997,selinger2010}.

Most concrete monoidal categories whose underlying product on sets is the
cartesian product are in fact cartesian categories. Besides the category
$(\Set,\times,1)$ itself, examples of such cartesian categories include the
categories $(\Top,\times,1)$ and $(\Meas,\times,1)$ of topological and
measurable spaces, with their cartesian products, and the category
$(\Vect_\K,\oplus,0)$ of vector spaces over the field $\K$, with its direct sum.
In all cases, the copying and deleting morphisms agree with those of $\Set$. The
category $\Rel$ of relations, also a symmetric monoidal category under the
cartesian product, supplies commutative comonoids, but the supply is not
homomorphic. The comonoid homomorphisms in $\Rel$ are precisely the relations
that are the graphs of functions, or, said differently, the cartesian center of
$\Rel$ is $\Set$.\footnote{Taking this observation seriously leads to Carboni
  and Walters' abstract theory of relations, the \emph{bicategory of relations}
  \cite{carboni1987a,fong2019relations,patterson2017relations}.} Another
important category that supplies commutative comonoids non-homomorphically is
the category of Markov kernels under the independent product
(\cref{ch:algebra-statistics}). Finally, the category $(\Vect_\K,\otimes,\K)$ of
vector spaces under its tensor product does not supply comonoids at all, since
the copying map $x \mapsto x \otimes x$ is nonlinear, among other difficulties.

A \emph{Lawvere theory}, also known as an \emph{algebraic theory} or a
\emph{finite-products theory}, is a small, strict cartesian category whose
monoid of objects is freely generated by a single object. Thus, Lawvere theories
play the role for cartesian categories that PROs and PROPs do for monoidal and
symmetric monoidal categories.

\begin{example}[Theory of groups] \label{ex:groups}
  The \emph{theory of groups} $\Theory{\Group}$ is the Lawvere theory presented
  by augmenting the theory of monoids (\cref{ex:monoids}) with another
  generating morphism $i: x \to x$, subject to the equations:
  \begin{equation*}
    \input{wiring-diagrams/category-theory/group-inverse-lhs} =
    \input{wiring-diagrams/category-theory/monoid-id} =
    \input{wiring-diagrams/category-theory/group-inverse-rhs}.
  \end{equation*}
  The \emph{theory of abelian groups} $\Theory{\Ab}$ is the Lawvere theory
  obtained in the same way from the theory of commutative monoids
  (\cref{ex:commutative-monoids}). In the abelian case, the inverse operation
  $i$ is often written as $-1$, reflecting the identification of abelian groups
  with $\Z$-modules (cf.\ \cref{ex:linear-spaces} below).
\end{example}

A \emph{cartesian (monoidal) functor}, or \emph{finite-product preserving
  functor}, is a symmetric monoidal functor between cartesian categories that
preserves the copying and deleting morphisms. Heuristically, this means that the
functor $F$ satisfies $F(\Copy_x) = \Copy_{F(x)}$ and
$F(\Delete_x) = \Delete_{F(x)}$ for all objects $x$, but as always the structure
isomorphisms must be accounted when the monoidal categories are not strict. For
a general supply, not necessarily of commutative comonoids, the definition is:

\begin{definition}[Preservation of supply] \label{def:supply-preservation}
  Let $\cat{P}$ be a PROP and let $s$ and $t$ be supplies of $\cat{P}$ in
  symmetric monoidal categories $\cat{C}$ and $\cat{D}$. A symmetric monoidal
  functor $(F,\Phi): \cat{C} \to \cat{D}$ \emph{preserves the supply} if for
  every object $x \in \cat{C}$ and every morphism $\mu: m \to n$ in $\cat{P}$,
  the following diagram commutes:
  \begin{equation*}
    \begin{tikzcd}[column sep=large]
      F(x)^{\otimes m} \arrow[r, "t_{F(x)}(\mu)"] \arrow[d, swap, "\Phi"] &
        F(x)^{\otimes n} \arrow[d, "\Phi"] \\
      F(x^{\otimes m}) \arrow[r, "F(s_x(\mu))"] &
        F(x^{\otimes n})
    \end{tikzcd}
  \end{equation*}
\end{definition}

Cartesian categories and cartesian functors assemble into a category $\Cart$.
Inside $\Cart$ lies the subcategory $\Lawvere$ of Lawvere theories and strict
cartesian functors preserving the object generators.

In the logical setting, a \emph{model} of a Lawvere theory $\cat{C}$ is a
cartesian functor $\cat{C} \to \Set$ and, generalizing, a \emph{model} of
$\cat{C}$ in a category $\cat{S}$ supplying commutative comonoids is a
supply-preserving functor $\cat{C} \to \cat{S}$. We follow the same conventions
on functional form as for models of PROs and PROPs. A model of a Lawvere theory
in a general category $\cat{S}$ supplying commutative comonids always lies in
the cartesian center of $\cat{S}$. This follows from the basic result that
supply-preserving functors send homomorphisms to homomorphisms \cite[Proposition
4.7]{fong2019supply}. That is, if a functor $F: (\cat{C},s) \to (\cat{D},t)$
preserves supply, then it restricts to a symmetric monoidal functor
$F_{s,t}: \cat{C}_s \to \cat{D}_t$.

\begin{example}[Group objects] \label{ex:group-objects} A \emph{group}, or
  \emph{group object}, in a cartesian category $\cat{S}$ is model in $\cat{S}$
  of the theory of groups $\Theory{\Group}$. A group in $\Set$ is just a group.
  A group in $\Top$ is a topological group, whereas a group in the category
  $\Man$ of smooth manifolds and smooth maps is a Lie group. In contrast to
  monoid objects (\cref{ex:monoid-objects}), group objects in tensor categories
  such as $(\Ab,\otimes,\Z)$ or $(\Vect_\K,\otimes,\K)$ do not make sense, as
  these categories are not cartesian. This is one reason why it is useful to
  have weak logical systems as well as strong ones: the weaker the logic in
  which a theory can be expressed, the more categories in which the theory can
  have models.
\end{example}

\begin{example}[Theory of groups, revisited]
  In \cref{ex:groups}, the theory of groups is presented following the standard
  axioms for a group, but other axiomatizations yield other presentations. For
  example, a group could be defined as a set $G$ equipped with a binary
  operation $(g,h) \mapsto g/h$ and a constant $e \in G$, satisfying the axioms
  \begin{equation*}
    g/g = e, \qquad g/e = g, \qquad (g/k) / (h/k) = g/h
  \end{equation*}
  for all $g,h,k \in G$. The definition is equivalent to the standard one via
  the assignments $g \cdot h := g / (e/h)$ and $g^{-1} := e/g$, in one
  direction, and $g/h := g \cdot h^{-1}$, in the other \cite[\S 1.3]{hall1959}.
  Let $\Theory{\Group}'$ be the Lawvere theory corresponding to these
  alternative axioms, generated by two morphisms $\delta: x \otimes x \to x$ and
  $\eta: I \to x$ and subject to three equations. Define a cartesian functor
  $\Theory{\Group} \to \Theory{\Group}'$ by
  \begin{equation*}
    \mu \mapsto
      (1_x \otimes ((\eta \otimes 1_x) \cdot \delta) \cdot \delta, \qquad
      \eta \mapsto \eta, \qquad
      i \mapsto (\eta \otimes 1_x) \cdot \delta.
  \end{equation*}
  By the equivalence of the two axiomatizations, the functor is an isomorphism
  of Lawvere theories: $\Theory{\Group} \cong \Theory{\Group}'$.

  The example highlights a beautiful aspect of categorical logic: when logical
  theories are reconstructed as algebraic structures, they become invariant to
  syntactic differences and exist independently of any particular presentation.
  Groups admit many axiomatizations, some insightful and others only curious.
  One striking axiomization consists of a single operation, the division
  operation, and a single monstrous equation \cite{higman1952,mccune1993}. No
  matter which axiomization is used, if it can be rendered as an algebraic
  theory, then it presents the same Lawvere theory, up to isomorphism.
\end{example}

As suggested at the end of \cref{sec:monoidal-categories}, the theories of
linear and affine combinations extend to Lawvere theories of vector spaces and
affine spaces.

\begin{example}[Theory of $R$-modules] \label{ex:linear-spaces}
  For any commutative rig $R$, the \emph{theory of $R$-modules} is the Lawvere
  theory $\Theory{\Mod_R}$ presented exactly as the theory of $R$-linear
  combinations $\LinComb_R$ (\cref{ex:linear-combinations}) but with additional
  laws for duplication,
  \begin{equation*}
    (\Copy_x \otimes 1_{x^{\otimes n}}) \cdot (r_1,\dots,r_{n+2})
      = (r_1+r_2, r_3, \dots, r_{n+2}),
  \end{equation*}
  and for discarding,
  \begin{equation*}
    (\Delete_x \otimes 1_{x^{\otimes n}}) \cdot (r_1,\dots,r_n)
      = (0,r_1,\dots,r_n),
    \qquad 
  \end{equation*}
  holding for all $r_i \in R$ and $n \geq 0$. Taking all the axioms together,
  the equations
  \begin{equation*}
    \Delta_{x^{\otimes n}} \cdot ((r_1,\dots,r_n) \otimes (s_1,\dots,s_n))
      \cdot (1,\dots,1) = (r_1 + s_1, \dots, r_n + s_n)
  \end{equation*}
  and
  \begin{equation*}
    \Delete_{x^{\otimes n}} \cdot (r_1,\dots,r_n) = (0,\dots,0)
  \end{equation*}
  can be derived. Alternatively, the theory $\Theory{\Mod_R}$ admits a biased
  presentation, extending that of $\LinComb_R$ with further axioms for
  duplication and discarding,
  \begin{equation*}
    \input{wiring-diagrams/category-theory/combination-sum-lhs} =
    \input{wiring-diagrams/category-theory/combination-sum-rhs}
    \qquad\text{and}\qquad
    \input{wiring-diagrams/category-theory/combination-zero-rhs} =
    \input{wiring-diagrams/category-theory/combination-zero-lhs}.
  \end{equation*}

  As special cases, the \emph{theory of (real) vector spaces} is
  $\Theory{\Vect_\R} := \Theory{\Mod_\R}$ and the \emph{theory of conical
    spaces} is $\Theory{\Cone} := \Theory{\Mod_{\R_+}}$. The nomenclature does
  not mislead: the category of models of $\Theory{\Mod_R}$ is indeed the usual
  category $\Mod_R$ of $R$-modules and $R$-linear maps; consequently, the
  category of models of $\Theory{\Vect_\R}$ is the usual category $\Vect_\R$ of
  real vector spaces and linear maps. A \emph{conical space} is, by definition,
  a model of the theory $\Theory{\Cone}$. This abstract structure is less
  familiar than a vector space but its main example, a convex cone in a real
  vector space, is ubiquitous in applied mathematics.
\end{example}

\begin{example}[Theory of $R$-affine spaces] \label{ex:affine-spaces}
  For any commutative rig $R$, the \emph{theory of $R$-affine spaces} is the
  Lawvere theory $\Theory{\Aff_R}$ obtained by restricting the presentation of
  the theory of $R$-modules to the generators $(r_i)_{i=1}^n$ such that
  $r_1 + \cdots + r_n = 1$. Thus, $\Theory{\Aff_R}$ is a sub-Lawvere theory of
  $\Theory{\Mod_R}$. Like the theory of $R$-affine combinations
  (\cref{ex:affine-combinations}), the theory of $R$-affine spaces admits a
  smaller, biased presentation in terms of morphisms
  $\mu_{r,s}: x \otimes x \to x$ indexed by numbers $r,s \in R$ such that
  $r+s = 1$.

  As special cases, the \emph{theory of (real) affine spaces} is
  $\Theory{\Aff_\R}$ and the \emph{theory of convex spaces} is
  $\Theory{\Conv} := \Theory{\Aff_{\R_+}}$. The category of models of
  $\Theory{\Aff_\R}$ is the usual category $\Aff_\R$ of affine spaces and affine
  maps. A \emph{convex space}, or model of $\Theory{\Conv}$, is an abstract
  structure whose main example is a convex set in a real vector space.
\end{example}

This introduction to categorical logic, covering the main definitions and the
examples needed later, is only the barest beginning of a fascinating subject.
The duality between syntax and semantics, namely that any Lawvere theory is
equivalent to the opposite of a certain subcategory of its category of models,
has not been discussed. This is a deep result that could not have been stated,
much less discovered, prior to the algebraization of logic. Categorical logic
has also expanded far beyond its original setting of algebraic theories into a
wide-ranging dictionary between logical systems and categorical structures
(\cref{fig:categorical-logics}). In theoretical computer science, the connection
between the simply typed lambda calculus and cartesian closed categories has
spurred a large research program. Within mathematical logic, categorical
logicians have charted a hierarchy of increasingly expressive subsystems of
first-order logic, culminating in topos theory. Rather than pursuing these more
expressive logics, this thesis will develop categorical logic in a different
direction, toward logics for probabilistic and statistical reasoning
(\cref{fig:categorical-logics-stats}).

The exposition in this chapter has been unorthodox in one respect: cartesian
categories are usually defined by a universal property, not as a homomorphic
supply of commutative comonoids. In the last part of this section, we digress to
explain how the two definitions are equivalent, a connection that is important
generally but not needed in this text.

\begin{definition}[Products and coproducts]
  In a category $\cat{C}$, a \emph{(binary) product} of a pair of objects $x$
  and $y$ is an object $x \times y$, equipped with \emph{projection} morphisms
  $\pi_x: x \times y \to x$ and $\pi_y: x \times y \to y$, such that for any
  morphisms $f: w \to x$ and $g: w \to y$, there exists a unique morphism
  $h: w \to x \times y$ making the diagram
  \begin{equation*}
    \begin{tikzcd}
      & w \arrow[d, dashed, "h", "\exists !"']
          \arrow[dl, "f"'] \arrow[dr, "g"] & \\
      x & x \times y \arrow[l, "\pi_x"] \arrow[r, "\pi_y"'] & y
    \end{tikzcd}
  \end{equation*}
  commute. An object $1$ in $\cat{C}$ is \emph{terminal} (a \emph{nullary
    product}) if for every object $x$, there exists a unique morphism $x \to 1$.

  Dually, a \emph{coproduct} of a pair of objects $x$ and $y$ is an object
  $x+y$, equipped with \emph{inclusion} morphisms $\iota_x: x \to x+y$ and
  $\iota_y: y \to x+y$, such that for any morphisms $f: x \to z$ and
  $g: y \to z$, there exists a unique morphism $h: x+y \to z$ making the diagram
  \begin{equation*}
    \begin{tikzcd}
      x \arrow[r, "\iota_x"] \arrow[dr, "f"']
        & x+y \arrow[d, dashed, near start, "h", "\exists !"']
        & y \arrow[l, "\iota_y"'] \arrow[dl, "g"] \\
      & z &
    \end{tikzcd}
  \end{equation*}
  commute. An object $0$ in $\cat{C}$ is \emph{initial} (a \emph{nullary
    coproduct}) if for every object $x$, there exists a unique morphism
  $0 \to x$.
\end{definition}

Whenever they exist, products and coproducts are unique up to canonical
isomorphism. Thus, in a mild abuse of language, it is common to speak of ``the''
product or ``the'' coproduct. In $\Set$, the product is the cartesian product
and the coproduct is the disjoint union, while the terminal object is any
singleton set and the initial object is the empty set. In general, the product
$x \times y$ in a category $\cat{C}$ classifies pairs of morphism \emph{into}
$x$ and $y$, in the sense that for every object $w$, there is a bijection of
hom-sets
\begin{equation*}
  \cat{C}(w,x) \times \cat{C}(w,y) \cong \cat{C}(w, x \times y).
\end{equation*}
Dually, the coproduct $x+y$ classifies pairs of morphisms \emph{out
  of} $x$ and $y$, in that for every object $z$, there is a bijection
\begin{equation*}
  \cat{C}(x,z) \times \cat{C}(y,z) \cong \cat{C}(x+y, z).
\end{equation*}
Although immediate from the definitions, these properties of products and
coproducts are often useful in their own right.

Classically, \emph{cartesian categories} monoidal categories whose monoidal
product is the category-theoretic product and whose monoidal unit is the
terminal object. \emph{Cocartesian categories} are monoidal categories whose
monoidal structure is the given by the coproduct. In fact, these definitions are
consistent with the previous \cref{def:cartesian-category}, an elegant result
due to Fox \cite{fox1976}.

\begin{theorem}[Fox's theorem] \label{thm:fox}
  A symmetric monoidal category $\cat{C}$ is a cartesian category, with its
  monoidal structure given by the categorical product, if and only if it has a
  homomorphic supply of commutative comonoids.
\end{theorem}

In particular, if a symmetric monoidal category has \emph{any} homomorphic
supply of commutative comonoids, then it has one uniquely, a fact not evident
from the definition.

\section{Interacting supplies in monoidal categories}
\label{sec:interacting-supplies}

There is a certain tension between the structuralist view of mathematics,
embodied by category theory, and the more anarchic approach to structure
prevalent in much of analysis and applied mathematics. In algebra, one typically
studies categories of structured objects and maps that fully preserve the
structure. So, in the category $\Group$, the morphisms are group homomorphisms,
and in the category $\Vect_\K$, the morphisms are linear maps. In the more
analytical parts of mathematics, there is often neither a single relevant class
of structured objects, nor, for a given class of objects, a single relevant
class of morphisms. It is not uncommon to speak of an affine map or a nonlinear
map between vector spaces, or of a convex-linear map between a convex cone and
an affine space (which makes sense because objects are, in particular, convex
sets). Depending on the context, the morphisms in a category of metric spaces
might reasonably be taken to include the maps that are: isometries,
nonexpansive, Lipschitz continuous, H\"older continuous, uniformly continuous,
continuous, measurable (with respect to the Borel $\sigma$-algebra), or nothing
at all.

Our approach to accommodating heterogeneous structure is twofold. To account for
heterogeneity among morphisms, supplies of PROPs in monoidal categories will
generally \emph{not} be homomorphic, thus encompassing a broad class of
morphisms while retaining the capability to assert, purely equationally, that
certain morphisms are supply homomorphisms. As for heterogeneity among objects,
the definition of supply will be extended from a single PROP to a family of
interacting PROPs. In this way, objects of different type will be able to
communicate along their maximum common substructure.

When working in the logic of cartesian categories, as in the previous section,
the theories of vector spaces, affine spaces, conical spaces, and convex spaces
are usually defined as Lawvere theories
(\cref{ex:linear-spaces,ex:affine-spaces}). In subsequent chapters, the
structure of these and other spaces is needed within the logic itself,
necessitating a level shift. So that they may be supplied in a monoidal
category, we reinterpret each Lawvere theory as a PROP by simply applying the
forgetful functor $\Lawvere \to \PROP$. In terms of presentations, the new PROP
is presented by joining the presentations of the Lawvere theory $\cat{C}$ and
the PROP $\Theory{\CComon}$, then adding equations making each generating
morphism in $\cat{C}$ into a comonoid homomorphism. As an example of this
procedure, the theory $\Theory{\Vect_\K}$ of vector spaces over a field $\K$ is
presented as a PROP in \cite[Equations 1-18]{baez2015}.

The resulting PROPs for vector spaces, affine spaces, conical spaces, and convex
spaces, as well as for commutative comonoids and bicommutative bimonoids
(\cref{ex:bimonoids}), form the objects of a category
\begin{equation*}
  \begin{tikzcd}[row sep=small]
    & \Theory{\CBimon} \arrow[r,hook] & \Theory{\Cone} \arrow[dr,hook] & \\
    \Theory{\CComon} \arrow[ur,hook] \arrow[dr,hook] & & & \Theory{\Vect_\R} \\
    & \Theory{\Conv} \arrow[uur,hook] \arrow[r,hook] &
      \Theory{\Aff_\R} \arrow[ur,hook] &
  \end{tikzcd}
\end{equation*}
whose arrows are embeddings of PROPs. This category, call it
$\TheoryLat{\Space}$, is a thin subcategory of $\PROP$. As a preorder, it has
\emph{meets} (greatest lower bounds) and a \emph{top} (greatest) object, making
it a \emph{meet-semilattice}. For example, the meet of the theories for $\Cone$
and $\Aff_\R$ is
\begin{equation*}
  \Theory{\Cone} \wedge \Theory{\Aff_\R} = \Theory{\Conv}
\end{equation*}
and the top object is $\top = \Theory{\Vect_R}$. In general, a thin subcategory
of $\PROP$ forming a meet-semilattice will be called a \emph{semilattice of
  PROPs}.

Supplies of a single PROP (\cref{def:supply,def:homomorphic-supply}) generalize
to supplies of a semilattice of PROPs. The idea is that if objects $x$ and $y$
in a symmetric monoidal category $\cat{C}$ are models of PROPs $\cat{P}_x$ and
$\cat{P}_y$ in a semilattice $\cat{L}$, then their monoidal product
$x \otimes y$ in $\cat{C}$ should be a model of the PROP having the maximum
common substructure of $\cat{P}_x$ and $\cat{P}_y$, namely the meet
$\cat{P}_x \wedge \cat{P}_y$ in $\cat{L}$. Likewise, the monoidal unit $I$ in
$\cat{C}$ should be a model of the PROP having the maximum structure, namely the
top object in $\cat{L}$.

\begin{definition}[Supply] \label{def:lattice-supply}
  A \emph{supply} of a semilattice $\cat{L}$ of PROPs in a symmetric monoidal
  category $(\cat{C}, \otimes, I)$ consists of a monoid homomorphism
  \begin{equation*}
    \cat{P}: (\Ob{\cat{C}}, \otimes, I) \to (\cat{L}, \wedge, \top), \quad
    x \mapsto \cat{P}_x,
  \end{equation*}
  and for each object $x \in \cat{C}$, a strong monoidal functor
  $s_x: \cat{P}_x \to \cat{C}$, such that
  \begin{enumerate}[(i),noitemsep]
  \item $s_x(m) = x^{\otimes m}$ for each $m \in \N$,
  \item the structure isomorphism
    $x^{\otimes m} \otimes x^{\otimes n} \to x^{\otimes (m+n)}$ is the unique
    coherence isomorphism for each $m, n \in \N$, and
  \item for every $x,y \in \cat{C}$ and every morphism $\mu: m \to n$ in
    $\cat{P}_{x \otimes y}$, the diagrams
    \begin{equation*}
      \begin{tikzcd}[column sep=9em]
        x^{\otimes m} \otimes y^{\otimes m}
          \arrow[r, "s_x(\pi_{x,y}(\mu)) \otimes s_y(\pi_{x,y}'(\mu))"]
          \arrow[d, "\Braid"'] &
        x^{\otimes n} \otimes y^{\otimes n}
          \arrow[d, "\Braid"] \\
        (x \otimes y)^{\otimes m}
          \arrow[r, "s_{x \otimes y}(\mu)"] &
        (x \otimes y)^{\otimes n}
      \end{tikzcd}
      \qquad\qquad
      \begin{tikzcd}[column sep=small]
        & I
          \arrow[dl, "\Braid"']
          \arrow[dr, "\Braid"] & \\
        I^{\otimes m}
          \arrow[rr, "s_I(\mu)"] & &
        I^{\otimes n}
      \end{tikzcd}
    \end{equation*}
    commute, where
    $\cat{P}_{x \otimes y} = \cat{P}_x \wedge \cat{P}_y \xrightarrow{\pi_{x,y}}
    \cat{P}_x$ and
    $\cat{P}_{x \otimes y} = \cat{P}_x \wedge \cat{P}_y \xrightarrow{\pi_{x,y}'}
    \cat{P}_y$ are the projections in $\cat{L}$ and the $\Braid$'s are the
    canonical symmetry isomorphisms.
  \end{enumerate}

  A morphism $f: x \to y$ in $\cat{C}$ is an \emph{$s$-homomorphism with respect
    to a PROP $\cat{Q}$}, for given PROP morphisms
  $\cat{Q} \xrightarrow{q_x} \cat{P}_x$ and
  $\cat{Q} \xrightarrow{q_y} \cat{P}_y$, if for every morphism $\mu: m \to n$ in
  $\cat{Q}$, the diagram
  \begin{equation*}
    \begin{tikzcd}[column sep=large]
      x^{\otimes m} \arrow[r, "s_x(q_x(\mu))"] \arrow[d, "f^{\otimes m}"'] &
        x^{\otimes n} \arrow[d, "f^{\otimes n}"] \\
      y^{\otimes m} \arrow[r, "s_y(q_y(\mu))"'] &
        y^{\otimes n}
    \end{tikzcd}
  \end{equation*}
  commutes. In particular, the morphism $f$ is an \emph{$s$-homomorphism} if it
  is an $s$-homomorphism with respect to $\cat{P}_x \wedge \cat{P}_y$, where
  $q_x = \pi_{x,y}$ and $q_y = \pi_{x,y}'$ are the projections. If every
  morphism in $\cat{C}$ is an $s$-homomorphism, then $s$ is a \emph{homomorphic
    supply}.
\end{definition}

When the map $\cat{P}: \Ob{\cat{C}} \to \Ob{\cat{L}}$ is constant, particularly
when the semilattice $\cat{L}$ consists of a single PROP, the original
definitions of supply and homomorphic supply are recovered.

\begin{example}[Vector space subsets] \label{ex:extrinsic-spaces}
  The \emph{category of vector space subsets}, $\VectSet_\K$, has as objects the
  pairs $(V,A)$, where $V$ is a vector space over $\K$ and $A$ is a a subset of
  $V$, and as morphisms $(V,A) \to (W,B)$ the functions $A \to B$. A symmetric
  monoidal category $(\VectSet_\K,\times,1)$ is defined by
  \begin{equation*}
    (V,A) \times (W,B) := (V \oplus W, A \times B), \qquad 1 := (0,\{0\}).
  \end{equation*}
  The category $\VectSet_\R$ supplies the semilattice $\TheoryLat{\Space}$
  by simply marking out the vector space subsets that are closed under linear
  combinations, affine combinations, conical combinations, convex combinations,
  sums, or nothing at all, thus distinguishing the linear subspaces, affine
  subspaces, convex cones, convex sets, additive monoids, and mere subsets. The
  supply is homomorphic with respect to $\Theory{\CComon}$, making $\VectSet_\R$
  into a cartesian category, but the supply is not homomorphic generally. A
  morphism in $\VectSet_\R$ is a homomorphism if it is linear, affine,
  conic-linear, convex-linear, or additive, as appropriate.
\end{example}

Vector space subsets are defined \emph{extrinsically}, with respect to an
ambient vector space, making for an easy construction of the semilattice supply.
With slightly more effort, one can construct a supply in a category of
heterogeneous objects defined \emph{intrinsically}, taking into account that, in
the absence of an ambient space, a set being a particular kind of object is an
extra structure, not a property.

\begin{example}[Vector and other spaces] \label{ex:intrinsic-spaces}
  Dual to $\TheoryLat{\Space}$, a thin subcategory of $\PROP$, is a thin
  subcategory $\SpaceCat$ of $\SMonCat$
  \begin{equation*}
    \begin{tikzcd}[row sep=small]
      & (\Cone,\oplus,0) \arrow[ddr] \arrow[r] & (\CMon,\oplus,0) \arrow[dr] \\
      (\Vect_\R,\oplus,0) \arrow[ur] \arrow[dr] & & & (\Set,\times,1) \\
      & (\Aff_\R,\times,1) \arrow[r] & (\Conv,\times,1) \arrow[ur]
    \end{tikzcd}
  \end{equation*}
  whose morphisms are forgetful functors. Importantly, all the
  products\footnote{The choice of $\oplus$ or $\times$ in the notation reflects
    whether the monoidal product is a biproduct or only a (categorical) product.
    In all cases, the underlying operation on sets is the cartesian product.} in
  these symmetric monoidal categories are compatible in the sense that the
  forgetful functors are strong symmetric monoidal functors. By duality, since
  $\TheoryLat{\Space}$ is a meet-semilattice, $\SpaceCat$ is a join-semilattice.
  Let $\vee$ and $\bot$ denote the \emph{join} (least upper bound) and
  \emph{bottom} (minimum) in $\SpaceCat$.

  Define a category $\Space$ whose objects are the disjoint union
  $\Ob{\Space} := \bigsqcup_{\cat{S} \in \SpaceCat} \Ob{\cat{S}}$ and whose
  morphisms are all functions between the underlying sets. Next, define a
  symmetric monoidal category $(\Space,\times,1)$ on objects by
  \begin{equation*}
    (\cat{S},X) \times (\cat{T},Y) :=
      (\cat{S} \vee \cat{T},\,
       U_{\cat{S},\cat{T}}(X) \otimes U_{\cat{S},\cat{T}}'(Y)),
    \qquad
    1 := (\bot, I_{\bot}),
  \end{equation*}
  where $\cat{S} \xrightarrow{U_{\cat{S},\cat{T}}} \cat{S} \vee \cat{T}$ and
  $\cat{T} \xrightarrow{U_{\cat{S},\cat{T}}'} \cat{S} \vee \cat{T}$ are the
  forgetful functors, and on morphisms by the cartesian product in $\Set$.
  Finally, using the correspondence between $\TheoryLat{\Space}$ and
  $\SpaceCat$, let the symmetric monoidal category $\Space$ supply
  $\TheoryLat{\Space}$ in the evident way. Consistent with the preceding
  example, the supply is homomorphic with respect to $\Theory{\Comon}$, so that
  $\Space$ is a cartesian category, but the supply is not homomorphic generally.
\end{example}

Generalizing a functor that preserves the supply of a single PROP
(\cref{def:supply-preservation}), a functor preserves the supply of a
semilattice of PROPs if it sends each object to another having \emph{at least}
as much as structure, in a compatible way.

\begin{definition}[Preservation of supply]
  \label{def:lattice-supply-preservation}
  Let $(\cat{P},s)$ and $(\cat{Q},t)$ be supplies of a semilattice $\cat{L}$ of
  PROPs in symmetric monoidal categories $\cat{C}$ and $\cat{D}$. A strong
  symmetric monoidal functor $(F,\Phi): \cat{C} \to \cat{D}$ \emph{preserves the
    supply} if for every object $x \in \cat{C}$, there exists a (unique) PROP
  morphism $\cat{P}_x \xrightarrow{i_x} \cat{Q}_{F(x)}$ in $\cat{L}$ and for
  every morphism $\mu: m \to n$ in $\cat{P}_x$, the diagram
  \begin{equation*}
    \begin{tikzcd}[column sep=huge]
      F(x)^{\otimes m} \arrow[r, "t_{F(x)}(i_x(\mu))"] \arrow[d, "\Phi"'] &
       F(x)^{\otimes n} \arrow[d, "\Phi"] \\
      F(x^{\otimes m}) \arrow[r, "F(s_x(\mu))"] &
       F(x^{\otimes n})
    \end{tikzcd}
  \end{equation*}
  commutes. If, moreover, $\cat{P}_x = \cat{Q}_{F(x)}$ for all $x \in \cat{C}$,
  then the supply is \emph{strictly} preserved.
\end{definition}

As an example, the embedding functor from the extrinsically-defined category
$\VectSet_\R$ to the intrinsically-defined category $\Space$, interpreting each
vector space subset as an abstract space of maximum possible structure, strictly
preserves the supply.

No further conditions need be imposed on a monoidal natural transformation
between supply preserving functors, as the components of the transformation are
always supply homomorphisms.

\begin{proposition} \label{prop:supply-natural-transformation}
  Let $\cat{C}$ and $\cat{D}$ be symmetric monoidal categories with supplies
  $(\cat{P},s)$ and $(\cat{Q},t)$ of a semilattice $\cat{L}$, and let $(F,\Phi)$
  and $(G,\Gamma)$ be supply preserving functors $\cat{C} \to \cat{D}$. Every
  component of a monoidal natural transformation $\alpha: F \to G$ is a
  $t$-homomorphism with respect to the PROP $\cat{P}_x$ and the PROP morphisms
  $\cat{P}_x \xrightarrow{i_x} \cat{Q}_{Fx}$ and
  $\cat{P}_{x} \xrightarrow{j_x} \cat{Q}_{Gx}$ in $\cat{L}$.

  In particular, if the functors $F$ and $G$ preserve the supply strictly, then
  every component of $\alpha$ is a $t$-homomorphism.
\end{proposition}
\begin{proof}
  We must show that, for any $x \in \cat{C}$, the component
  $\alpha_x: Fx \to Gx$ is a $t$-homomorphism with respect to $\cat{P}_x$,
  meaning that for every morphism $\mu: m \to n$ in $\cat{P}_x$, the diagram
  \begin{equation*}
    \begin{tikzcd}[column sep=large]
      F(x)^{\otimes m} \arrow[r, "\alpha_x^{\otimes m}"]
          \arrow[d, "t_{Fx}(i_x(\mu))"']
        & G(x)^{\otimes m} \arrow[d, "t_{Gx}(j_x(\mu))"] \\
      F(x)^{\otimes n} \arrow[r, "\alpha_x^{\otimes n}"]
        & G(x)^{\otimes n}
    \end{tikzcd}
  \end{equation*}
  commutes. Consider the composite diagram
  \begin{equation*}
    \begin{tikzcd}[column sep=large]
      F(x)^{\otimes m} \arrow[r, "\Phi"] \arrow[d, "t_{Fx}(i_x(\mu))"']
        & F(x^{\otimes m}) \arrow[r, "\alpha_{x^{\otimes m}}"]
          \arrow[d, "F(s_x(\mu))"']
        & G(x^{\otimes m}) \arrow[r, "\Gamma^{-1}"] \arrow[d, "G(s_x(\mu))"]
        & G(x)^{\otimes m} \arrow[d, "t_{Gx}(j_x(\mu))"] \\
      F(x)^{\otimes n} \arrow[r, "\Phi"]
        & F(x^{\otimes n}) \arrow[r, "\alpha_{x^{\otimes n}}"]
        & G(x^{\otimes n}) \arrow[r, "\Gamma^{-1}"]
        & G(x)^{\otimes n}.
    \end{tikzcd}
  \end{equation*}
  The middle square commutes because $\alpha$ is a natural transformation, while
  the left and right squares commute because the functors $F$ and $G$ preserve
  the supply. Moreover, since $\alpha$ is a \emph{monoidal} natural
  transformation, the top horizontal composite equals $\alpha_x^{\otimes m}$ and
  the bottom horizontal composite equals $\alpha_x^{\otimes n}$. Thus the
  commutativity of the composite diagram proves the proposition.
\end{proof}

When working with natural transformations between functors out of a category
presented by generators and relations, as we often will in
\cref{ch:algebra-statistics,ch:zoo-statistics}, it is useful to know that in
order to prove naturality, it suffices to prove it on a generating set of
morphisms. In the cases of transformations between categories or monoidal
categories, this fact must be considered well known, even if explicit statements
are not easily found in the literature.

\begin{lemma}[Natural transformations and generators]
  \label{lemma:natural-transformation-generators}
  Let $F,G: \cat{C} \to \cat{D}$ be functors out of a small category $\cat{C}$.
  A collection of morphisms $\alpha_x: Fx \to Gx$, $x \in \cat{C}$, are the
  components of a natural transformation $\alpha: F \to G$ if, for every
  $f: x \to y$ in a generating set of morphisms for $\cat{C}$, the naturality
  square
  \begin{equation*}
    \begin{tikzcd}
      Fx \arrow[r, "\alpha_x"] \arrow[d, "Ff"']
      & Gx \arrow[d, "Gf"] \\
      Fy \arrow[r, "\alpha_y"]
      & Gy
    \end{tikzcd}
  \end{equation*}
  commutes. The same statement holds when $\cat{C}$ and $\cat{D}$ are
  (symmetric) monoidal categories, $F$ and $G$ are (symmetric) monoidal
  functors, and the components $(\alpha_x)_{x \in \cat{C}}$ preserve the
  monoidal product and unit, as in \cref{def:monoidal-natural-transformation}.
  Finally, the statement still holds when $\cat{C}$ and $\cat{D}$ supply a
  semilattice $\cat{L}$ of PROPs, $F$ and $G$ are supply preserving functors,
  and the components $(\alpha_x)_{x \in \cat{C}}$ preserve the monoidal
  structure and are also supply homomorphisms, in the sense of
  \cref{prop:supply-natural-transformation}.
\end{lemma}
\begin{proof}
  First, assume that $\cat{C}$ and $\cat{D}$ are categories and $F$ and $G$ are
  functors. Every morphism in $\cat{C}$ is formed from the generating morphisms
  through composites and identities, so by structural induction, it suffices to
  show that naturality is preserved by composition and identities. When
  $f = 1_x$ is an identity, the naturality condition reduces to the trivial
  equation $\alpha_x = \alpha_x$. Given composable morphisms
  $x \xrightarrow{f} y \xrightarrow{g} z$ in $\cat{C}$, the diagram
  \begin{equation*}
    \begin{tikzcd}
      Fx \arrow[r, "Ff"] \arrow[d, "\alpha_x"']
          \arrow[rr, bend left, "F(f \cdot g)"]
        & Fy \arrow[r, "Fg"] \arrow[d, "\alpha_y"']
        & Fz \arrow[d, "\alpha_z"] \\
      Gx \arrow[r, "Gf"] \arrow[rr, bend right, "G(f \cdot g)"']
        & Gy \arrow[r, "Gg"] & Gz
    \end{tikzcd}
  \end{equation*}
  shows that if naturality holds for $f$ and $g$, then it also holds for the
  composite $f \cdot g$. This proves the first statement.

  Now suppose that $\cat{C}$ and $\cat{D}$ are (symmetric) monoidal categories
  and $(F,\Phi)$ and $(G,\Gamma)$ are (symmetric) monoidal functors. For
  simplicity, assume that $\cat{C}$ is a strict monoidal category, as all small
  monoidal categories in this text are strict. Then every morphism in $\cat{C}$
  is formed from the generators through composites, identities, monoidal
  products, and, when the monoidal categories are symmetric, braidings. To show
  that naturality is preserved by monoidal products, take any morphisms
  $x \xrightarrow{f} y$ and $w \xrightarrow{g} z$ in $\cat{C}$ and consider the
  diagram
  \begin{equation*}
    \begin{tikzcd}[column sep=large]
      F(x \otimes w) \arrow[r, "\Phi^{-1}_{x,w}"]
          \arrow[d, "\alpha_{x \otimes w}"']
        & Fx \otimes Fw \arrow[r, "Ff \otimes Fg"]
          \arrow[d, "\alpha_x \otimes \alpha_w"']
        & Fy \otimes Fz \arrow[r, "\Phi_{y,z}"]
          \arrow[d, "\alpha_y \otimes \alpha_z"]
        & F(y \otimes z) \arrow[d, "\alpha_{y \otimes z}"] \\
      G(x \otimes w) \arrow[r, "\Gamma^{-1}_{x,w}"]
        & Gx \otimes Gw \arrow[r, "Gf \otimes Gg"]
        & Gy \otimes Gz \arrow[r, "\Gamma_{y,z}"]
        & G(y \otimes z).
    \end{tikzcd}
  \end{equation*}
  If naturality holds for $f$ and $g$, then the middle square commutes, and the
  left and right squares commute by hypothesis. Thus the whole diagram commutes.
  Moreover, since $\Phi$ and $\Gamma$ are natural isomorphisms, the top
  composite is $F(f \otimes g)$ and the bottom composite is $G(f \otimes g)$.
  This proves that naturality holds for the product $f \otimes g$. The proof
  that naturality holds for the braidings
  $x \otimes y \xrightarrow{\Braid_{x,y}} y \otimes x$ in $\cat{C}$ is similar.
  Because braidings are natural isomorphisms, the middle square in the diagram
  \begin{equation*}
    \begin{tikzcd}[column sep=large]
      F(x \otimes y) \arrow[r, "\Phi^{-1}_{x,y}"]
          \arrow[d, "\alpha_{x \otimes y}"']
        & Fx \otimes Fy \arrow[r, "\Braid_{Fx,Fy}"]
          \arrow[d, "\alpha_x \otimes \alpha_y"']
        & Fy \otimes Fx \arrow[r, "\Phi_{y,x}"]
          \arrow[d, "\alpha_y \otimes \alpha_x"]
        & F(y \otimes x) \arrow[d, "\alpha_{y \otimes x}"] \\
      G(x \otimes y) \arrow[r, "\Gamma^{-1}_{x,y}"]
        & Gx \otimes Gy \arrow[r, "\Braid_{Gx,Gy}"]
        & Gy \otimes Gx \arrow[r, "\Gamma_{y,x}"]
        & G(y \otimes x),
    \end{tikzcd}
  \end{equation*}
  commutes and thus the whole diagram does. Moreover, since $F$ and $G$ are
  \emph{symmetric} monoidal functors, the top composite is $F(\Braid_{x,y})$ and
  the bottom composite is $G(\Braid_{x,y})$.
  
  Finally, suppose that the categories $(\cat{C},\cat{P},s)$ and
  $(\cat{D},\cat{Q},t)$ supply a semilattice of PROPs and that $F$ and $G$ are
  supply preserving functors. As the morphisms of $\cat{C}$ now include the
  supplied operations, we must show that naturality holds for each morphism
  $s_x(\mu): x^{\otimes m} \to x^{\otimes n}$, where $x \in \cat{C}$ and
  $\mu: m \to n$ belongs to $\cat{P}_x$. This follows from the hypothesis that
  $\alpha_x: Fx \to Gx$ is a supply homomorphism by inverting the argument in
  the proof of \cref{prop:supply-natural-transformation}.
\end{proof}

\section{Notes and references}

\paragraph{General category theory} Seventy-five years after its inception at
the hands of Eilenberg and Mac Lane \cite{eilenberg1942,eilenberg1945}, category
theory is now a large field of mathematics and is larger still when account is
taken of its interactions with other fields. General introductions to category
theory, arranged in order of increasing difficulty, are
\cite{lawvere2009,leinster2014,awodey2010,riehl2016,maclane1998,borceux1994a}.
Of these, Riehl's textbook \cite{riehl2016} is especially recommended for its
judicious choice of topics and diverse list of examples. For applications
outside of pure mathematics, the edited collection \emph{New Structures for
  Physics} \cite{coecke2011} includes introductions to category theory by
Abramsky and Tzevelekos \cite{abramsky2010}, focusing on categorical logic; Baez
and Stay \cite{baez2010}, drawing analogies between parts of physics, topology,
logic and computation; and Coecke and Paquette \cite{coecke2010}, developing
parts of categorical quantum mechanics. More recently, Spivak \cite{spivak2014}
and Fong and Spivak \cite{fong2019sketches} have published textbooks on applied
category theory.

\paragraph{Monoidal categories and string diagrams} Originally called
``categories with multiplication,'' monoidal categories were introduced
independently by B{\'e}nabou \cite{benabou1963} and Mac Lane \cite{maclane1963}.
The latter author also established the coherence theorem for monoidal
categories, subsequently improved by Kelly \cite{kelly1964}. String diagrams
originated in physics as the \emph{Penrose notation} for tensor calculus
\cite{penrose1971,penrose1984} but were first formalized and made rigorous by
Joyal and Street \cite{joyal1991,joyal1995}. Consequently, string diagrams are
sometimes called the \emph{Joyal-Street calculus}. Expositions of monoidal
categories and string diagrams include
\cite{baez2010,coecke2010,street2012,savage2018}. Selinger has written an
encyclopedic survey of the many graphical languages based on string diagrams
\cite{selinger2010}. Book-length treatments of monoidal categories, emphasizing
connections to Hopf algebras, are \cite{aguiar2010} and \cite{etingof2015}.

The notion of a \emph{supply} of a PROP in a symmetric monoidal category was
introduced recently by Fong and Spivak \cite{fong2019supply}, generalizing and
systematizing a supply of commutative comonoids. The latter notion had been in
widespread use for some time, albeit under different names, such as a
\emph{monoidal category with diagonals} \cite{selinger1999} or, in the case
where $\Copy_x: x \to x \otimes x$ and $\Delete_x: x \to I$ are natural in $x$,
\emph{uniform copying and deleting} \cite{heunen2019}. The generalization of
supply to interacting PROPs, described in \cref{sec:interacting-supplies}, is
original, although the modeling of implicit type conversions as a thin category
is well established in programming language theory (see Notes to
\cref{ch:semantic-enrichment}).

The characterization of a cartesian category as a
symmetric monoidal category with a homomorphic supply of commutative comonoids
is due to Fox \cite{fox1976}. A diagrammatic proof of this theorem appears in
the textbook by Heunen and Vicary \cite{heunen2019} and earlier in lecture notes
by the same authors \cite{heunen2013}.

\paragraph{Categorical logic} The field of categorical logic was launched by
Lawvere in a seminal PhD thesis, \emph{Functorial Semantics of Algebraic
  Theories} \cite{lawvere1963}, which introduced the functorial view of
semantics, established a connection between cartesian categories and algebraic
theories, and discovered a deep duality between syntax and semantics, now known
as Lawvere duality. The connection between cartesian closed categories and the
typed lambda calculus was subsequently established by Lambek
\cite{lambek1980,lambek1986}. Today, the standard reference on categorical logic
and topos theory is Johnstone's multi-volume treatise \cite{johnstone2002}.
Introductions to the subject include the lecture notes
\cite{awodey2019,shulman2016} and the textbooks
\cite{borceux1994b,crole1993,goldblatt1984,reyes2004}. Of these, Reyes, Reyes,
and Zolfaghari exposit the logic of $\cat{C}$-sets (\cref{sec:categories}) in a
concrete style \cite{reyes2004}. Crole gives a thorough and self-contained
treatment of algebraic theories, the simply typed lambda calculus, and their
algebraic semantics in cartesian categories and cartesian closed categories
\cite{crole1993}.

\paragraph{Linear and other spaces} Linear, affine, conical, and convex spaces
are all, to varying degrees, well-known, but despite their structural
similarities they are not often treated in a unified way. The theories of
linear, affine, conical, and convex combinations are more naturally defined as
\emph{operads} than as PROPs, after the works of Giraudo \cite{giraudo2015} and
of Leinster (unpublished but reported by Baez and Fritz \cite[\S 2]{baez2014}).
Moreover, since the theories make no reference to the rig's additive structure,
operads of combinations can be defined over any monoid $(M,\cdot,1)$, not
necessarily commutative \cite{giraudo2015}. For the sake of uniformity, we have
restricted ourselves to commutative rigs and have presented the theories of
combinations as PROPs, following the standard procedure for generating PROPs
from operads \cite[Example 60]{markl2008}.

The definition of a vector space is, of course, completely standard. The
category of vector spaces and its sibling, the category of linear relations,
have been studied as PROPs or Lawvere theories independently by several groups
\cite{baez2015,wadsley2015,bonchi2017}. Affine spaces can be defined in numerous
different but equivalent ways. Perhaps the most common is as a set equipped with
a simply transitive action by a vector space (or rather by its underlying
abelian group) \cite[Ch.\ 2]{berger1987}. Our definition is equivalent to this
one, but more closely resembles the alternative definition as an algebra of the
affine combinations monad \cite[\S 5.2]{riehl2016}. Conical spaces are also
known as \emph{semimodules} over the nonnegative real numbers and belong to the
general study of semirings (rigs) and semimodules (modules over rigs). Finally,
convex spaces, viewed as abstract structures, are somewhat obscure but go back
at least to Marshall Stone \cite{stone1949}. Capraro and Fritz compare several
different axiomatizations \cite{capraro2013}. Every convex subset of a real
vector space is convex space, yet not every convex space embeds into a real
vector space. A certain \emph{cancellation property} is a necessary and
sufficient condition for an embedding to exist \cite[Theorem
4]{stone1949,capraro2013}.

\chapter{The algebra of statistical theories and models}
\chaptermark{Algebra of statistical models}
\label{ch:algebra-statistics}

In theoretical statistics, a statistical model is formally defined to be a
parameterized family of probability distributions. For each parameter $\theta$
in a parameter space $\Omega$, the model specifies a probability distribution
$P_\theta$ supported on a common sample space $\sspace{X}$. Given an observation
$x \in \sspace{X}$, assumed to be sampled from one of the distributions
$P_{\theta_0}$, the problem of statistical inference is to determine, as
precisely as possible, the unknown parameter $\theta_0$ out of all possible
parameters $\theta \in \Omega$. So that the inference problem is not impossible,
the model is usually required to be \emph{identifiable} in that
$P_\theta \neq P_{\theta'}$ whenever $\theta \neq \theta'$. Within this modest
framework, one can already state many of the essential definitions of
theoretical statistics, such as sufficiency, minimal sufficiency, and
ancillarity, and prove classical results such as the Fisher-Neyman factorization
criterion and Basu's theorem.

Despite this, the formal definition of a statistical model is too minimalistic
to describe a large part of statistical modeling. Often the sample space is
already determined, or at least tightly constrained, before the model is even
formulated, whereas the parameter space is rarely of intrinsic interest, because
the model can always be reparametrized by an invertible transformation
$\Omega \to \Omega'$. Rather, the scientist's interest lies in how the mapping
$\theta \mapsto P_\theta$ defining the model relates the parameters to the
observed data; how the model is related to other, competing models; and how this
whole family of models is related to any relevant background scientific
theories. The general definition of a statistical model provides no guidance on
how to answer these questions, as it reduces statistical models to black boxes
for mapping parameters into probability distributions. To a practitioner not
steeped in the lore of theoretical statistics, it may not even be apparent how
the formal definition accommodates everyday statistical models like the linear
model or logistic regression.

Of course, this situation has not prevented statisticians and data scientists
from understanding the internal structure of statistical models, from comparing
competing models, or from arguing that a certain model supports or fails to
support some scientific theory. They have simply done so without the benefit of
any formal mathematical system. Our philosophy is that every element of
scientific knowledge that can be talked about at all, can be talked about
rigorously, and that doing so promotes clarity in thought, communication, and
computational representation.

This chapter develops the algebra of statistical models as a probabilistic form
of categorical logic. The formalism distinguishes between \emph{statistical
  theories}, which are purely algebraic structures, and \emph{models} of
statistical theories, which are, by a pun on the word ``model,'' also
statistical models as ordinarily understood. Statistical theories are finitary
descriptions of the structure of statistical models, amenable to machine
formalization. Morphisms of statistical theories, together with induced
morphisms between categories of models, make precise the notion of a
relationship between models. Altogether, the formalism offers a rigorous
language for describing the internal structure of, and the relationships
between, statistical models. It does not directly address relationships between
statistical models and scientific theories, much less the notorious
philosophical problem of how statistical inference can support or criticize a
scientific theory. Nevertheless, it is hoped that by building a bridge between
mathematical logic and statistics, an advance will have been made in the larger
program to understand the interlocking roles of theories and models in logic,
statistics, and science.

Apart from introductions to the category of Markov kernels and its abstraction
as a Markov category in \cref{sec:markov-kernels,sec:markov-categories}, the
content of this chapter is mostly original. A detailed account of the provenance
of this set of ideas is provided at the end of the chapter.

\section{Markov kernels in statistics}
\label{sec:markov-kernels}

A statistical model $\{P_\theta\}_{\theta \in \Omega}$ on a sample space
$\sspace{X}$, depending measurably on its parameter $\theta$, can be interpreted
as a Markov kernel $P: \Omega \to \sspace{X}$. This section introduces the
symmetric monoidal category of Markov kernels, accompanied by numerous examples
drawn from probability and statistics.

A Markov kernel is the probabilistic analogue of a function, assigning to every
point in its domain a probability distribution over its codomain. It can also be
regarded as a conditional probability distribution. Formally:

\begin{definition}[Markov kernel]
  A \emph{Markov kernel} $M: \sspace{X} \to \sspace{Y}$ from one measurable
  space $(\sspace{X}, \Sigma_{\sspace{X}})$ to another
  $(\sspace{Y}, \Sigma_{\sspace{Y}})$, also known as a \emph{probability kernel}
  or a \emph{stochastic kernel}, is a function
  $M: \sspace{X} \times \Sigma_{\sspace{Y}} \to [0,1]$ such that
  \begin{enumerate}[(i),noitemsep]
  \item for every point $x \in \sspace{X}$, the map
    $M(x; -): \Sigma_{\sspace{Y}} \to [0,1]$ is a probability measure on
    $\sspace{Y}$;
  \item for every set $B \in \Sigma_{\sspace{Y}}$, the map
    $M(-; B): \sspace{X} \to [0,1]$ is measurable.
  \end{enumerate}
  In agreement with the standard notation for conditional probability, the
  probability measure $M(x; dy)$ is often written $M(dy \given x)$.
\end{definition}

Equivalently, a Markov kernel $M: \sspace{X} \to \sspace{Y}$ is a measurable map
$\sspace{X} \to \ProbSpace(\sspace{Y})$, where $\ProbSpace(\sspace{Y})$ is the
space of all probability measures on $\sspace{Y}$ under the $\sigma$-algebra
generated by the evaluation functionals $\mu \mapsto \mu(B)$,
$B \in \Sigma_{\sspace{Y}}$ \cite[Lemma 1.40]{kallenberg2002}. From this
perspective, it is natural to denote the distribution $M(x; -)$ at
$x \in \sspace{X}$ simply as $M(x)$.

Yet another perspective is that Markov kernels are linear operators on spaces of
measures \cite[\S 3.3]{worm2010}. Let $M: \sspace{X} \to \sspace{Y}$ be a Markov
kernel. For any measure $\mu$ on $\sspace{X}$, its image under $M$ is the
measure $\mu M$ on $\sspace{Y}$ defined by
\begin{equation*}
  (\mu M)(B) := \int_{\sspace{X}} M(B \given x)\,\mu(dx),
  \qquad B \in \Sigma_{\sspace{Y}}.
\end{equation*}
With this definition, $M$ is a \emph{Markov operator}: if
$\MeasSpace(\sspace{X})$ is the space of all finite signed measures on
$\sspace{X}$, then $M$ acts as a linear map
$\MeasSpace(\sspace{X}) \to \MeasSpace(\sspace{Y})$ that preserves the total
mass, $\mu M(\sspace{Y}) = \mu(\sspace{X})$. In particular, $M$ acts as a
convex-linear map $\ProbSpace(\sspace{X}) \to \ProbSpace(\sspace{Y})$ between
spaces of probability measures.

All the parametric families of probability distributions that constitute the
basic material of statistics can be represented as Markov kernels.

\begin{example}[Normal family] \label{ex:normal-family}
  The \emph{normal} or \emph{Gaussian} family is the Markov kernel
  $\Normal = \Normal_1: \R \times \R_+ \to \R$ given by
  \begin{equation*}
    \Normal(dx \given \mu, \sigma^2) := \frac{1}{\sqrt{2\pi \sigma^2}}
      e^{-(x-\mu)^2/2\sigma^2} dx,
  \end{equation*}
  when $\sigma^2 > 0$, and equal to a point mass, $\Normal(\mu,0) = \delta_\mu$,
  when $\sigma^2 = 0$. The normal family is parametrized by mean and variance:
  if $X \sim \Normal(\mu, \sigma^2)$, then $\E(X) = \mu$ and
  $\Var(X) = \sigma^2$.

  In higher dimension $d$, let $\PSD^d$ denote the cone of $d \times d$ positive
  semi-definite, real-valued matrices. The $d$-dimensional \emph{multivariate
    normal} or \emph{Gaussian} family is the unique Markov kernel
  $\Normal_d: \R^d \times \PSD^d \to \R^d$ such that for all vectors
  $v \in \R^d$,
  \begin{equation*}
    \inner*{v,\, \Normal_d(\mu,\Sigma)} =
      \Normal_1(\inner{v, \mu}, \inner{v, \Sigma\, v})
  \end{equation*}
  \cite[\S 3.1]{mardia1979}. This is well-defined by the Cram{\'e}r-Wold theorem
  \cite[Corollary 5.5]{kallenberg2002}. An explicit formula
  \begin{equation*}
    \Normal_d(dx \given \mu, \Sigma) = \frac{1}{\sqrt{(2\pi)^d |\Sigma|}}
      \exp\left(-\frac{1}{2} \inner*{x - \mu, \Sigma^{-1} (x - \mu)} \right) dx
  \end{equation*}
  is available when the covariance matrix $\Sigma$ is nondegenerate
  (invertible), but when it is degenerate the Gaussian measure does not have a
  density with respect to Lebesgue measure.\footnote{Some authors define normal
    families only in the nondegenerate case. While this choice simplifies the
    definition, it gives up the essential property of being closed under
    arbitrary linear combinations.}
\end{example}

\begin{example}[Exponential families] \label{ex:exponential-family} A
  $d$-dimensional \emph{exponential family} is a Markov kernel
  $P: \Omega \to \sspace{Y}$ of form
  \begin{equation*}
    P(dy \given \theta) = e^{\inner{\theta,\, t(y)} - \psi(\theta)}\,\nu(dy),
  \end{equation*}
  where the \emph{canonical parameter} $\theta$ belongs to a parameter space
  $\Omega$ in $\R^d$; the \emph{sufficient statistic} or \emph{canonical
    statistic} $t(y)$ is a measurable map $t: \sspace{Y} \to \R^d$; the
  \emph{base measure} $\nu$ is a $\sigma$-finite measure on $\sspace{Y}$,
  typically having a density with respect to counting measure or Lebesgue
  measure; and the \emph{normalizing function} $\psi: \Omega \to \R$ is given by
  \begin{equation*}
    \psi(\theta) := \log
      \int_{\sspace{Y}} e^{\inner{\theta,\, t(y)}}\,\nu(dy).
  \end{equation*}
  An exponential family is often given its maximal domain of definition,
  \begin{equation*}
    \Omega = \left\{\theta \in \R^d:
      \int_{\sspace{Y}} e^{\inner{\theta,\, t(y)}}\,\nu(dy) < \infty \right\},
  \end{equation*}
  in which case $\Omega$ is called the \emph{canonical parameter space} and the
  family is called \emph{full}. The canonical parameter space is always a convex
  set in $\R^d$ and, under regularity conditions, it is also open.

  Exponential families are ubiquitous in statistics, both theoretical
  \cite{brown1986} and applied \cite{sundberg2019}. Suitably parameterized, most
  of the common families of probability distributions, such as the normal,
  exponential, gamma, beta, Bernoulli, binomial, and Poisson, are exponential
  families. Important exceptions include the uniform and Cauchy families.
\end{example}

Despite their interpretation as randomized functions, working with Markov
kernels does not preclude working with ordinary functions. Functions are
recovered as the \emph{deterministic} Markov kernels.

\begin{definition}[Deterministic kernels]
  A Markov kernel $M: \sspace{X} \to \sspace{Y}$ is \emph{deterministic} if
  $M(x)$ is a point mass for every point $x \in \sspace{X}$ or, equivalently,
  there exists a measurable map $f: \sspace{X} \to \sspace{Y}$ such that
  $M(x) = \delta_{f(x)}$ for all $x \in \sspace{X}$.
\end{definition}

A measurable map $f$ is often identified with the deterministic Markov kernel
$\delta_f$ through a mild abuse of notation. As an example, the kernel
$\Normal_d(-,0): \R^d \to \R^d$ is deterministic (in fact, it is the identity
function).

As the notation suggests, Markov kernels are composable and thus form the
morphisms of a category. Composition in this category has already been
implicitly used in \cref{ex:normal-family}, where the normal family
$\Normal_d: \R^d \times \PSD^d \to \R^d$ is composed with the linear form
$\inner{v,-}: \R^d \to \R$ to obtain a Markov kernel into $\R$.

\begin{definition}[Category of kernels] \label{def:markov-kernel-category}
  The \emph{category of Markov kernels}, denoted $\Markov$, has Polish
  measurable spaces as objects\footnote{That is, the objects are Polish spaces
    (separable, completely metrizable topological spaces), measurable under
    their Borel $\sigma$-algebras. The category of Markov kernels is just as
    easily defined to include all measurable spaces, but it will be more
    practical to adopt a regularity condition, and rule out measure-theoretic
    pathologies, at the outset.} and the Markov kernels between them as
  morphisms. The composite of a Markov kernel $M: \sspace{X} \to \sspace{Y}$
  with another kernel $N: \sspace{Y} \to \sspace{Z}$ is the kernel
  $M \cdot N: \sspace{X} \to \sspace{Z}$ given by
  \begin{equation*}
    (M \cdot N)(C \given x) := \int_{\sspace{Y}} N(C \given y) M(dy \given x),
    \qquad x \in \sspace{X}, \quad C \in \Sigma_{\sspace{Z}}.
  \end{equation*}
  The identity $1_\sspace{X}: \sspace{X} \to \sspace{X}$ is the usual identity
  map on $\sspace{X}$, construed as a Markov kernel.
\end{definition}

For a proof that composition of Markov kernels is associative, see \cite[Lemma
5.6]{cencov1982} or \cite[Proposition 3.2]{panangaden1999}.

The composition law has a natural probabilistic interpretation. For fixed
$x \in \sspace{X}$, form a joint probability distribution on random variables
$(Y,Z) \in \sspace{Y} \times \sspace{Z}$ as $M(dy \given x) N(dz \given y)$, so
that $M(x)$ is the marginal distribution of $Y$ and $N$ is the conditional
distribution of $Z$ given $Y$. Then compute the marginal distribution of $Z$ by
integrating out $Y$. Applying this procedure to every $x \in \sspace{X}$ defines
a map $\sspace{X} \to \ProbSpace(\sspace{Z})$, hence a Markov kernel
$\sspace{X} \to \sspace{Z}$. This kernel is the composite of $M$ and $N$.

Composition of Markov kernels generalizes composition of functions: for any
composable measurable maps $f$ and $g$, the deterministic kernels $\delta_f$ and
$\delta_g$ satisfy $\delta_f \cdot \delta_g = \delta_{f \cdot g}$. Thus there is
an embedding of categories $\Meas \hookrightarrow \Markov$, where $\Meas$ is the
category of Polish measurable spaces and measurable maps.

The action of Markov kernels as Markov operators is a special case of
composition. Letting $I := \{*\}$ be the singleton measurable space, a
probability measure on $\sspace{X}$ can be identified with a Markov kernel
$\mu: I \to \sspace{X}$. Its image under a kernel $M: \sspace{X} \to \sspace{Y}$
is exactly the composite $\mu \cdot M: I \to \sspace{Y}$. Moreover, if
$N: \sspace{Y} \to \sspace{Z}$ is another kernel, then
$\mu(M \cdot N) = (\mu M) N$ by associativity, showing that the composition laws
for Markov kernels and linear maps are compatible. Thus the mapping
$\sspace{X} \mapsto \MeasSpace(\sspace{X})$ extends to a functor
$\MeasSpace: \Markov \to \Vect_\R$. Similarly, there is a functor
$\ProbSpace: \Markov \to \Conv$ into the category of convex spaces. Both
functors are easily seen to be embeddings \cite[Theorem 5.2 and Lemma
5.10]{cencov1982}, making the category of Markov kernels isomorphic to at least
two different concrete categories.

Many parametric families and statistical models arise as composites of simpler
ones. Rather trivially, the composite of the normal family
$\Normal: \R \times \R_+ \to \R$ with the exponential function
$\exp: \R \to \PosR$ is the \emph{log-normal} family
$\LogNormal: \R \times \R_+ \to \PosR$, so called because its logarithm is
normally distributed. The log-normal family is a common model for effects that
accrue multiplicatively rather than additively. For fixed $n \in \N$, the
\emph{beta-binomial family} is the composite of the beta family
$\BetaDist: (\PosR)^2 \to [0,1]$ with the binomial $\Binom(n,-): [0,1] \to \N$.
As another example, the \emph{noncentral chi-squared family with one degree of
  freedom}, $\chi_1^2: \R_+ \to \R_+$, is characterized by the equation
\begin{equation*}
  \Normal(\mu,1)^2 = \chi_1^2(\mu^2)
  \qquad\text{or}\qquad
  \begin{tikzcd}
    \R \arrow{d}[swap]{\Normal(-,1)} \arrow[r, "(-)^2"]
      & \R_+ \arrow[d, "\chi_1^2"] \\
    \R \arrow[r, "(-)^2"] & \R_+
  \end{tikzcd}.
\end{equation*}
Countless relationships between parametric families are known
\cite{springer1979,leemis2008}. Some are useful in statistical modeling, others
play an important role in sampling random variables \cite{devroye1986}, still
others are merely curious. In order to express many of these relationships
algebraically, additional structure must be introduced into the category of
Markov kernels, starting with a monoidal product.

Recall that if $\mu$ and $\nu$ are probability measures on spaces
$(\sspace{X},\Sigma_{\sspace{X}})$ and $(\sspace{Y},\Sigma_{\sspace{Y}})$, then
their \emph{product measure} $\mu \otimes \nu$ on the product space
$(\sspace{X} \times \sspace{Y}, \Sigma_{\sspace{X}} \otimes
\Sigma_{\sspace{Y}})$ is defined on measurable rectangles by
\begin{equation*}
  (\mu \otimes \nu)(A \times B) := \mu(A) \nu(B), \qquad
  A \in \Sigma_{\sspace{X}}, \quad B \in \Sigma_{\sspace{Y}}.
\end{equation*}
Probabilistically, the joint distribution $(X,Y) \sim \mu \otimes \nu$ makes $X$
and $Y$ independent with marginal distributions $X \sim \mu$ and $Y \sim \nu$.
In the case of Markov kernels, taking products pointwise defines a monoidal
product.

\begin{definition}[Independent product] \label{def:independent-product}
  Define a monoidal product, the \emph{independent product}, on the category of
  Markov kernels as follows. The product $\sspace{X} \otimes \sspace{Y}$ of
  objects $\sspace{X}$ and $\sspace{Y}$ is the product space\footnote{This usage
    of ``product space'' is unambiguous because, in a Polish space, the product
    and Borel $\sigma$-algebras are compatible \cite[Lemma
    1.2]{kallenberg2002}.}

  $\sspace{X} \times \sspace{Y}$. The product
  $M \otimes N: \sspace{W} \otimes \sspace{X} \to \sspace{Y} \otimes \sspace{Z}$
  of morphisms $M: \sspace{W} \to \sspace{Y}$ and $N: \sspace{X} \to \sspace{Z}$
  is given pointwise as
  \begin{equation*}
    (M \otimes N)(w,x) := M(w) \otimes N(x), \qquad
    w \in \sspace{W}, \quad x \in \sspace{X}.
  \end{equation*}
  The monoidal unit is the singleton space $I := \{*\}$. With this definition,
  $(\Markov, \otimes, I)$ is a symmetric monoidal category, where the braidings,
  associators, and unitors are the usual maps construed as Markov kernels.
\end{definition}

On any measurable space $\sspace{X}$, copying and deleting maps are defined in
the usual way by
\begin{equation*}
  \Copy_{\sspace{X}}: \sspace{X} \to \sspace{X} \otimes \sspace{X},\
    x \mapsto (x,x), \qquad
  \Delete_{\sspace{X}}: \sspace{X} \to I,\ x \mapsto *.
\end{equation*}
Construed as Markov kernels, they equip the symmetric monoidal category
$\Markov$ with a supply of commutative comonoids
(\cref{sec:cartesian-categories}). Before characterizing the supply
homomorphisms in the next section, the symmetric monoidal structure and the
comonoid supply are illustrated by several examples. The first is a standard
construction in probability theory.

\begin{example}[Kernel product measures] \label{ex:kernel-product-measure}
  The \emph{product} of a probability measure $\mu$ on $\sspace{X}$ and a Markov
  kernel $M: \sspace{X} \to \sspace{Y}$ is a probability measure on the product
  space $\sspace{X} \times \sspace{Y}$, defined on measurable rectangles by
  \begin{equation*}
    A \times B \mapsto \int_A M(B \given x)\,\mu(dx),
    \qquad A \in \Sigma_{\sspace{X}}, \quad B \in \Sigma_{\sspace{Y}}.
  \end{equation*}
  As a morphism $I \to \sspace{X} \times \sspace{Y}$, the product measure is the
  Markov kernel $\mu \cdot \Copy_{\sspace{X}} \cdot (1_{\sspace{X}} \otimes M)$
  or, in the graphical syntax,
  \begin{equation*}
    \input{wiring-diagrams/algebra-statistics/markov-kernel-product-measure}.
  \end{equation*}
  Be warned that this measure is often denoted as $\mu \times M$ or
  $\mu \otimes M$ in standard texts \cite{klenke2013}, notations that are
  incompatible with the monoidal product in $\Markov$. On the other hand, when
  $M: \sspace{X} \to \sspace{Y}$ is the constant kernel $x \mapsto \nu$ at some
  distribution $\nu$, that is, when $M = \Delete_{\sspace{X}} \cdot \nu$, then
  the short calculation
  \begin{equation*}
    \input{wiring-diagrams/algebra-statistics/product-measure-intermediate-1} =
    \input{wiring-diagrams/algebra-statistics/product-measure-intermediate-2} =
    \input{wiring-diagrams/algebra-statistics/product-measure}
  \end{equation*}
  recovers the usual product measure $\mu \otimes \nu$ on
  $\sspace{X} \times \sspace{Y}$.
\end{example}

Many more parametric families can be realized as composites using the newly
introduced structure. Generalizing an earlier example, the \emph{noncentral
  chi-squared family with $k$ degrees of freedom}, $\chi_k^2: \R_+ \to \R_+$, is
characterized by the equation
\begin{equation*}
  \input{wiring-diagrams/algebra-statistics/noncentral-chi-squared-lhs}
  \quad=\quad
  \input{wiring-diagrams/algebra-statistics/noncentral-chi-squared-rhs},
\end{equation*}
where the unfilled circle denotes addition. Defining the chi-squared family
equationally is almost always preferable to defining it directly by its
probability density function, a complicated expression involving modified Bessel
functions. As another example \cite[\S 4.4]{holmes2018}, the \emph{negative
  binomial family} $\NegBinom: \R_+ \times (0,1) \to \N$ may be defined as a
\emph{gamma-Poisson mixture}
\begin{equation*}
  \input{wiring-diagrams/algebra-statistics/negative-binomial} :=
  \input{wiring-diagrams/algebra-statistics/gamma-poisson-mixture},
\end{equation*}
where $\GammaDist: \R_+^2 \to \R_+$ is the gamma family, parametrized by shape
and scale; $\Pois: \R_+ \to \N$ is the Poisson family,\footnote{In the gamma and
  Poisson families, the shape and scale parameter spaces are extended beyond
  their usual definitions to include zero, so that, for example,
  $\Pois(0) = \delta_0$. This convention is uncommon but not unheard of,
  especially in statistical software. As in the normal family, the purpose is to
  improve algebraic closure properties.} and the odds function is
$\pi \mapsto \pi/(1-\pi)$. When the data is underdispersed or overdispersed
under a Poisson model, the negative binomial model is a flexible alternative
allowing the variance to differ from the mean.

For continuous data, scale transformations offer a more general remedy for
underdispersion or overdispersion. The following example, the last in the
section, is more involved than the previous ones. It illustrates how the concept
of a Markov kernel allows a formal analogy to be drawn between scale
transformations in statistics and convex analysis. It also introduces the
exponential dispersion model, a family of probability distributions important
for generalized linear models (\cref{sec:glm}).

\begin{example}[Scale transformations] \label{ex:exponential-dispersion-family}
  Let $\sspace{X}$ and $\sspace{Y}$ be any convex sets in $\R^d$ closed under
  multiplication by positive scalars, such as convex cones. The \emph{scale
    transform} of a Markov kernel $P: \sspace{X} \to \sspace{Y}$ is the Markov
  kernel $\tilde P: \sspace{X} \times \PosR \to \sspace{Y}$ defined by
  \begin{equation*}
    \tilde P(\mu,\sigma) := \sigma P(\mu/\sigma)
    \qquad\text{or}\qquad
    \input{wiring-diagrams/algebra-statistics/scale-transform-notation} :=
    \input{wiring-diagrams/algebra-statistics/scale-transform}.
  \end{equation*}
  If the original family $P$ is parameterized by mean, in that
  $\E(P(\mu)) = \mu$ for all $\mu \in \sspace{X}$, then the new family
  $\tilde P$ is parameterized by mean with respect to its first argument, since
  $\E(\tilde P(\mu,\sigma)) = \sigma\, \E(P(\mu/\sigma)) = \mu$. Moreover, the
  new family has the more flexible variance
  $\Var(\tilde P(\mu,\sigma)) = \sigma^2\, V(\mu/\sigma)$, where
  $V(\mu) := \Var(P(\mu))$ is the variance function of $P$. For example, the
  scale transform of the normal location family $\Normal(-,1): \R \to \R$ is the
  normal location-scale family, parameterized by standard deviation rather than
  variance. Although it has neither mean nor variance, the Cauchy location
  family $\Cauchy(-,1): \R \to \R$ can also be scale transformed, yielding the
  Cauchy location-scale family.

  The scale transform of a Markov kernel is formally identical to the
  perspective transform of a convex function. According to a standard definition
  of convex analysis \cite{hiriart-urruty1993,combettes2018}, the \emph{scale
    transform}, or \emph{perspective}, of an arbitrary function
  $f: \sspace{X} \to \R$ is the function
  $\tilde f: \sspace{X} \times \PosR \to \R$ given by
  $\tilde f(x,\lambda) := \lambda f(x/\lambda)$. Importantly, the perspective
  $\tilde f$ is jointly convex if and only if $f$ is convex.

  The convex perspective transform in fact leads to a \emph{different} notion of
  scale transform for exponential families (\cref{ex:exponential-family}). Let
  $P: \Omega \to \sspace{Y}$ be an exponential family reduced to its sufficient
  statistic, so that $t(y)=y$. By general properties of exponential families,
  the normalizing function $\psi: \Omega \to \R$ is both a convex function and
  the cumulant generating function for the family via
  $K(u \given \theta) := \psi(\theta + u) - \psi(\theta)$. The second fact
  implies that the family's mean vector and variance-covariance matrix are given
  by the gradient and Hessian of $\psi$:
  \begin{equation*}
    \mu(\theta) := \E(P(\theta)) = \nabla \psi(\theta)
    \qquad\text{and}\qquad
    V(\theta) := \Var(P(\theta)) = \nabla^2 \psi(\theta).
  \end{equation*}
  If, for fixed $\lambda > 0$, the perspective
  $\tilde \psi(-, \lambda): \Omega \to \R$ of $\psi$ at scale $\lambda$ is the
  normalizing function for \emph{another} exponential family, say
  \begin{equation*}
    P_\lambda(dy \given \theta)
      := e^{\inner{\theta,y} - \tilde\psi(\theta,\lambda)}\, \nu_\lambda(dy)
      = e^{\inner{\theta,y} - \lambda \psi(\theta/\lambda)}\, \nu_\lambda(dy)
  \end{equation*}
  with base measure $\nu_\lambda$, then this family will have mean and variance
  functions
  \begin{gather*}
    \E(P_\lambda(\theta))
      = \nabla_\theta\, \tilde\psi(\theta,\lambda)
      = \nabla \psi(\theta)
      = \mu(\theta) \\
    \Var(P_\lambda(\theta))
      = \nabla_\theta^2\, \tilde\psi(\theta,\lambda)
      = \lambda^{-1}\, \nabla^2 \psi(\theta)
      = \lambda^{-1}\, V(\theta).
  \end{gather*}
  Taking the subset $\Lambda \subseteq \PosR$ of all viable values of $\lambda$
  (which at least includes 1) and making the change of parameter
  $\theta \mapsto \lambda\theta$, a Markov kernel
  $Q: \Omega \times \Lambda \to \sspace{Y}$ with enlarged domain is defined by
  \begin{equation*}
    Q(dy \given \theta, \lambda)
      := P_{\lambda}(dy \given \lambda\theta)
      = e^{\lambda(\inner{\theta,y} - \psi(\theta))}\, \nu_\lambda(dy).
  \end{equation*}
  The family $Q$ is called an \emph{exponential dispersion model} with
  \emph{index parameter} $\lambda$; alternatively, making another change of
  parameter $\phi = 1/\lambda$, the family
  $\tilde Q(\theta,\phi) := Q(\theta,1/\phi)$ is an exponential dispersion model
  with \emph{dispersion parameter} $\phi$ \cite{jorgensen1987,jorgensen1992}.
  The primary use of exponential dispersion models is as a component of
  generalized linear models.

  The exponential dispersion model has its roots in a simpler and ubiquitous
  construction on Markov kernels. Let $P: \Omega \to \sspace{Y}$ be a Markov
  kernel, not necessarily an exponential family, with its sample space
  $\sspace{Y}$ a convex set in $\R^d$. The \emph{sample mean} $\bar P_n$ over
  $n$ i.i.d.\ observations of $P$ is the Markov kernel
  \begin{equation*}
    \input{wiring-diagrams/algebra-statistics/markov-kernel-mean-notation} :=
    \input{wiring-diagrams/algebra-statistics/markov-kernel-mean},
  \end{equation*}
  where the coefficients $\vec{1}_n/n = (1/n,\dots,1/n)$ define a convex
  combination, the sample mean in $\sspace{Y}$. If the original kernel $P$ has
  mean function $\mu(\theta) := \E(P(\theta))$ and variance function
  $V(\theta) := \Var(P(\theta))$, then its sample mean $\bar P_n$ has the same
  mean, $\E(\bar P_n(\theta)) = \mu(\theta)$, but the reduced variance
  $\Var(\bar P_n(\theta)) = n^{-1} V(\theta)$. If, moreover, the kernel $P$ has
  moment generating function
  $M(u \given \theta) := \displaystyle \int_{\sspace{Y}} e^{\inner{u,y}} P(dy
  \given \theta)$ and cumulant generating function
  $K(u \given \theta) := \log M(u \given \theta)$, then by properties of the
  Laplace transform, the corresponding functions for $\bar P_n$ are
  $\bar M_n(u \given \theta) = M(u/n \given \theta)^n$ and
  \begin{equation*}
    \bar K_n(u \given \theta)
      = \log(M(u/n \given \theta)^n)
      = n K(u/n \given \theta)
      = \tilde K(u \given \theta, n),
  \end{equation*}
  where $\tilde K(- \given \theta, n)$ is the perspective of
  $K(- \given \theta)$ at scale $n$.

  In the case of an exponential family, the cumulant generating function
  corresponding to the perspective of $\psi$ at scale $\lambda$ is
  $K_\lambda(u \given \theta) := \tilde\psi(\theta + u, \lambda) -
  \tilde\psi(\theta,\lambda)$. Making the same change of parameter
  $\theta \mapsto \lambda\theta$ as before yields
  \begin{equation*}
    \tilde K(u \given \theta, \lambda) =
      \lambda (\psi(\theta + u/\lambda) - \psi(\theta)),
  \end{equation*}
  which is simultaneously the cumulant generating function of the exponential
  dispersion model and the perspective of
  $K(u \given \theta) = \psi(\theta + u) - \psi(\theta)$ at scale $\lambda$.
  Thus, exponential dispersion models can be seen as an analytical extension,
  specific to exponential families, of the sample mean from a natural number $n$
  to a continuous parameter $\lambda$.
\end{example}

\section{Algebraic reasoning about Markov kernels}
\label{sec:markov-categories}

The definition of the category of Markov kernels and its illustrations in
statistics have thoroughly blended algebraic and analytical reasoning. The aim
of this section is to disentangle the algebra from the analysis and axiomatize
the properties essential to statistical modeling. This will be achieved in a
qualified sense.

Consider the question of when a Markov kernel is a comonoid homomorphism. A
generic Markov kernel $M: \sspace{X} \to \sspace{Y}$ preserves deleting but not
copying, that is,
\begin{equation*}
  \input{wiring-diagrams/algebra-statistics/markov-kernel-delete-lhs} =
  \input{wiring-diagrams/algebra-statistics/markov-kernel-delete-rhs}
  \qquad\text{but \emph{not}}\qquad
  \input{wiring-diagrams/algebra-statistics/markov-kernel-mcopy-lhs} =
  \input{wiring-diagrams/algebra-statistics/markov-kernel-mcopy-rhs}.
\end{equation*}
The second equation holds exactly when, for every $x \in \sspace{X}$, the
deterministic and independent couplings of $M(x)$ with itself are equal:
\begin{equation*}
  M(x) \cdot \Delta_{\sspace{Y}} = M(x) \otimes M(x).
\end{equation*}
The following proposition characterizes the solutions to this equation over an
arbitrary probability measure $\mu$.

\begin{proposition}
  For any measurable space $\sspace{X}$, the probability measures $\mu$ on
  $\sspace{X}$ whose deterministic and independent couplings are equal,
  \begin{equation*}
    \mu \cdot \Delta_{\sspace{X}} = \mu \otimes \mu,
  \end{equation*}
  are exactly the extreme points of $\ProbSpace(\sspace{X})$.
\end{proposition}
\begin{proof}
  The condition $\mu \otimes \mu = \mu \cdot \Delta_{\sspace{X}}$ says that
  $\mu(A) \mu(B) = \mu(A \cap B)$ for all measurable sets $A$ and $B$; in
  particular, $\mu(A)^2 = \mu(A)$ for every measurable set $A$. This condition
  is equivalent to $\mu(A)$ being equal to either 0 or 1 for every $A$. Let us
  call such a probability measure $\mu$ a \emph{0-1 measure}.

  So we must show that the extreme points of $\ProbSpace(\sspace{X})$ are
  exactly the 0-1 measures on $\sspace{X}$, a well-known fact \cite[Example
  8.16]{simon2011}. Suppose that a 0-1 measure $\mu$ is expressed as a convex
  combination $\mu = t \mu_1 + (1-t) \mu_2$ for $t \in (0,1)$. Fix a measurable
  set $A$ and consider cases. If $\mu(A) = 0$, then since $\mu_1, \mu_2 \geq 0$
  and $0 < t < 1$, we must have $\mu_1(A) = \mu_2(A) = 0 = \mu(A)$. On the other
  hand, if $\mu(A) = 1$, then as $\mu_1, \mu_2 \leq 1$, we must have
  $\mu_1(A) = \mu_2(A) = 1 = \mu(A)$. Since this holds for any measurable set
  $A$, we conclude that $\mu_1 = \mu_2 = \mu$ and hence that $\mu$ is an extreme
  point. Conversely, suppose that the probability measure $\mu$ is not a 0-1
  measure. Then there exists a measurable set $B$ such that $0 < \mu(B) < 1$,
  and $\mu$ can be expressed the nontrivial convex combination
  \begin{equation*}
    \mu = \mu(B) \cdot \frac{\mu(- \cap B)}{\mu(B)} +
      (1 - \mu(B)) \cdot \frac{\mu(- \setminus B)}{\mu(\sspace{X} \setminus B)}.
  \end{equation*}
  Therefore, $\mu$ is not an extreme point of $\ProbSpace(\sspace{X})$.
\end{proof}

\begin{corollary}
  For any measurable spaces $\sspace{X}$ and $\sspace{Y}$, a Markov kernel
  $M: \sspace{X} \to \sspace{Y}$ preserves copying,
  \begin{equation*}
    \input{wiring-diagrams/algebra-statistics/markov-kernel-mcopy-lhs} =
    \input{wiring-diagrams/algebra-statistics/markov-kernel-mcopy-rhs},
  \end{equation*}
  if and only if for every $x \in \sspace{X}$, the distribution $M(x)$ is an
  extreme point of $\ProbSpace(\sspace{Y})$.
\end{corollary}

By analogy to the finite-dimensional probability simplex, one would expect that
the extreme points of any space of probability measures would be exactly the
point masses. Certainly, a point mass is a 0-1 measure, hence an extreme point.
However, a regularity condition is needed for the other direction. The following
result is classic; proofs are given in \cite[Theorem 15.9]{aliprantis2006} and
\cite[Example 8.16]{simon2011}.

\begin{theorem} \label{thm:extreme-points-prob-space}
  If $\sspace{X}$ is a Polish space, then the extreme points of
  $\ProbSpace(\sspace{X})$ are exactly the point masses $\delta_x$ for
  $x \in \sspace{X}$.
\end{theorem}

In view of the previous corollary, we immediately deduce:

\begin{corollary}[Cartesian center of $\Markov$]
  \label{cor:markov-kernel-determinism}
  For any measurable space $\sspace{X}$ and any Polish space $\sspace{Y}$, a
  Markov kernel $M: \sspace{X} \to \sspace{Y}$ preserves copying,
  \begin{equation*}
    \input{wiring-diagrams/algebra-statistics/markov-kernel-mcopy-lhs} =
    \input{wiring-diagrams/algebra-statistics/markov-kernel-mcopy-rhs},
  \end{equation*}
  if and only if it is deterministic. In particular, the comonoid homomorphisms
  in $\Markov$ are exactly the deterministic kernels, or equivalently, the
  cartesian center of $\Markov$ can be identified with the category $\Meas$ of
  measurable maps.
\end{corollary}

Determinism is thus characterized equationally inside the category of Markov
kernels, assuming the regularity conditions of
\cref{def:markov-kernel-category}. The next proposition shows that isomorphisms
in this category are in a sense trivial.

\begin{proposition}[Isomorphisms in $\Markov$]
  \label{prop:markov-kernel-isomorphims}
  Every isomorphism (invertible morphism) in $\Markov$ is deterministic. That
  is, whenever kernels $M: \sspace{X} \to \sspace{Y}$ and
  $N: \sspace{Y} \to \sspace{X}$ in $\Markov$ satisfy
  $M \cdot N = 1_{\sspace{X}}$ and $N \cdot M = 1_{\sspace{Y}}$, then both $M$
  and $N$ are deterministic.
\end{proposition}
\begin{proof}
  If $M: \sspace{X} \to \sspace{Y}$ is a isomorphism in $\Markov$, then $M$, as
  a Markov operator, is a convex-linear isomorphism
  $\ProbSpace(\sspace{X}) \to \ProbSpace(\sspace{Y})$ and hence preserves
  extreme points. But by \cref{thm:extreme-points-prob-space}, the extreme
  points are exactly the point masses. Thus, for every $x \in \sspace{X}$, there
  exists $y \in \sspace{Y}$ such that $M(x) = \delta_x M = \delta_y$, proving
  that $M$ is deterministic.
\end{proof}

This concludes a survey of the basic structural properties of the category of
Markov kernels. The most basic property of all is captured by the following
abstraction, which has been studied under different names and sometimes only
implicitly. We follow Fritz in adopting the suggestive name ``Markov category''
\cite{fritz2020}.

\begin{definition}[Markov category]
  A \emph{Markov category} is a symmetric monoidal category supplying
  commutative comonoids, such that every morphism $f: x \to y$ is a supply
  homomorphism with respect to deleting,
  \begin{equation*}
    \input{wiring-diagrams/category-theory/supply-homomorphism-delete-lhs} =
    \input{wiring-diagrams/category-theory/supply-homomorphism-delete-rhs},
  \end{equation*}
  but not necessarily with respect to copying.
\end{definition}

Cartesian categories, in the sense of \cref{def:cartesian-category}, are
evidently Markov categories, but these are not the intended examples. The
prototypical example is $\Markov$, \emph{the} category of Markov kernels.
Embedded in $\Markov$ is the full subcategory $\CAT{FinMarkov}$ of finite
measurable spaces and Markov kernels. This Markov category is equivalent to the
category of right stochastic matrices, in which the objects are natural numbers,
composition is matrix multiplication, and the monoidal product is the matrix
direct sum. As a nonprobabilistic example, the category of sets and multivalued
functions is a Markov category. It is a subcategory of $\Rel$ and its cartesian
center is $\Set$. Further examples of Markov categories may be found in \cite[\S
3-9]{fritz2020}.

In a Markov category, the monoidal unit is \emph{terminal}: for any object $x$,
there exists a unique morphism $x \to I$, namely the deleting map
$\Delete_x: x \to I$. A symmetric monoidal category whose monoidal unit is
terminal has been a called \emph{semicartesian category}, in view of
\cref{thm:fox}, or a \emph{monoidal category with projections}, because for any
objects $x$ and $y$, there are well-behaved projections
$\pi_{x,y}: x \otimes y \to x$ and $\pi_{x,y}': x \otimes y \to y$ given by
\begin{equation*}
  x \otimes y \xrightarrow{1_x \otimes \Delete_y} x \otimes I
    \xrightarrow{\cong} x
  \quad\text{and}\quad
  x \otimes y \xrightarrow{\Delete_x \otimes 1_y} I \otimes y
    \xrightarrow{\cong} y.
\end{equation*}
Every Markov category is a semicartesian category.

Several fundamental constructions on Markov kernels can be rephrased in purely
algebraic terms in a Markov category. The product of a kernel with a measure
(\cref{ex:kernel-product-measure}), and the inverse operation of
\emph{disintegrating} a product measure \cite[Theorem
1.23]{chang1997,kallenberg2017}, carry over immediately, and likewise for the
lesser known, but more general, operation of \emph{disintegrating} a Markov
kernel \cite[Theorem 1.25]{kallenberg2017}. Disintegration finds an important
statistical application in Bayesian inference. Variants of the following
definition appear as \cite[Definition 3.5]{cho2019} and \cite[Definition
11.5]{fritz2020}.

\begin{definition}[Disintegration] \label{def:disintegration}
  In a Markov category, a \emph{disintegration} of a morphism
  $f: x \to y \otimes z$ with respect to $y$ consists of a pair of morphisms
  $f_y: x \to y$ and $f_{z \given y}: x \otimes y \to z$ such that
  \begin{equation*}
    \input{wiring-diagrams/algebra-statistics/disintegration-lhs} =
    \input{wiring-diagrams/algebra-statistics/disintegration-rhs}.
  \end{equation*}
  Disintegration of the morphism $f: x \to y \otimes z$ with respect to $z$ is
  defined similarly.
\end{definition}

When a disintegration of $f$ with respect to $y$ exists, the morphism $f_y$ is
equal to the marginal $x \xrightarrow{f} y \otimes z \xrightarrow{\pi_{y,z}} y$,
as can be seen by post-composing both sides of the defining equation with the
projection $\pi_{y,z}$. The morphism $f_{z \given y}$ is generally not unique.
In the category $\Markov$ of well-behaved Markov kernels, the conditional kernel
is unique up to sets of probability zero \cite[Corollary 1.26]{kallenberg2017}.

Taking the domain $x$ to be the monoidal unit $I$ recovers the simpler notion of
disintegrating a distribution. Bayesian inference can then be formulated in any
Markov category in which the required disintegrations exist: given a sampling or
likelihood morphism $p: \theta \to x$ and a prior $\pi_0: I \to \theta$, first
integrate with respect to $\theta$ to obtain a joint distribution
$I \to \theta \otimes x$, then disintegrate with respect to $x$ to obtain a
posterior $\pi_1: x \to \theta$ and a marginal likelihood $p_x: I \to x$.

Conditional independence and exchangeability can also be formulated in any
Markov category. Equivalent notions of independence are discussed in \cite[\S
6]{cho2019} and \cite[\S 12]{fritz2020}.

\begin{definition}[Independence] \label{def:independence}
  In a Markov category, a morphism $f: x \to y_1 \otimes \cdots \otimes y_n$ has
  \emph{(conditionally) independent} components $y_1,\dots,y_n$ if there exist
  morphisms $f_i: x \to y_i$, $i=1,\dots,n$, such that
  \begin{equation*}
    f = \Copy_{x,n} \cdot (f_1 \otimes \cdots \otimes f_n),
  \end{equation*}
  where $\Copy_{x,n}: x \to x^{\otimes n}$ is the $n$-fold copying morphism. In
  this case, each $f_i$ is equal to the corresponding marginal
  $x \xrightarrow{f} y_1 \otimes \cdots \otimes y_n \xrightarrow{\pi_i} y_i$.
  When $n=2$, the defining condition appears as
  \begin{equation*}
    \input{wiring-diagrams/algebra-statistics/independence-lhs} =
    \input{wiring-diagrams/algebra-statistics/independence-rhs}.
  \end{equation*}
  If, in addition, all the $y_i$'s are equal and all the $f_i$'s are equal, then
  $f$ is said to have \emph{independently and identically distributed (i.i.d)}
  components.
\end{definition}

\begin{definition}[Exchangeability] \label{def:exchangeability}
  A morphism $f: x \to y^{\otimes n}$ in a Markov category has
  \emph{(conditionally) exchangeable} components if for all permutations
  $\sigma \in S_n$, the composite
  $x \xrightarrow{f} y^{\otimes n} \xrightarrow{\sigma} y^{\otimes n}$ is equal
  to $f$. When $n=2$, this condition reduces to the single equation
  \begin{equation*}
    \input{wiring-diagrams/algebra-statistics/exchangeability-lhs} =
    \input{wiring-diagrams/algebra-statistics/exchangeability-rhs}.
  \end{equation*}
\end{definition}

In a Markov category, as in classical probability, independence implies
exchangeability. For example, when $n=2$, one calculates that
\begin{equation*}
  \input{wiring-diagrams/algebra-statistics/exchangeability-lhs} =
  \input{wiring-diagrams/algebra-statistics/independence-exchangeability-1} =
  \input{wiring-diagrams/algebra-statistics/independence-exchangeability-2} =
  \input{wiring-diagrams/algebra-statistics/independence-exchangeability-3} =
  \input{wiring-diagrams/algebra-statistics/exchangeability-rhs}.
\end{equation*}
The second equation uses naturality of the braiding isomorphisms and the third
uses commutativity of the copying morphisms.

The characterization of deterministic Markov kernels
(\cref{cor:markov-kernel-determinism}) becomes a definition in an abstract
Markov category.

\begin{definition}[Determinism] \label{def:determinism}
  A morphism $f: x \to y$ in a Markov category is \emph{deterministic} if it
  preserves copying,
  \begin{equation*}
    \input{wiring-diagrams/category-theory/supply-homomorphism-mcopy-lhs} =
    \input{wiring-diagrams/category-theory/supply-homomorphism-mcopy-rhs}.
  \end{equation*}
  Thus, by definition, the cartesian center of a Markov category is the
  subcategory of deterministic morphisms.
\end{definition}

The subtleties surrounding determinism illustrate the gap between well-behaved
Markov kernels and morphisms in a general Markov category. In the Markov
category of \emph{all} measurable spaces and Markov kernels, the concrete
definition of determinism is not equivalent to the abstract one. Only under mild
regularity conditions, of the sort imposed on $\Markov$, are the two definitions
equivalent (\cref{cor:markov-kernel-determinism}). In a similar vein, the fact
that all isomorphisms in $\Markov$ are deterministic
(\cref{prop:markov-kernel-isomorphims}) is not true for arbitrary Markov
kernels, hence it cannot be deduced from the axioms of a Markov category. In a
generic Markov category, the most that can be said is that if two morphisms are
mutually inverse and \emph{one} of them is deterministic, then so is the other
\cite[Lemma 10.9]{fritz2020}. The definition of a Markov category should
therefore be considered minimalistic, since the Markov kernels used in
statistical applications hardly ever exhibit measure-theoretic pathologies.

\section{Linear algebraic reasoning about Markov kernels}
\label{sec:linear-algebraic-markov-categories}

The question of what further axioms, if any, should be imposed on a well-behaved
Markov category will not be answered here. But there is another, more
specialized class of structure that is indispensable to statistical modeling.
The examples of \cref{sec:markov-kernels} make liberal use of vector spaces and
other spaces, and hardly any useful statistical model can be formulated without
reference to such structure. The remainder of this section introduces Markov
categories with linear or related structure.

Recall from \cref{ch:category-theory} the theories of vector spaces, affine
spaces, conical spaces, and convex spaces
(\cref{ex:linear-spaces,ex:affine-spaces}), as well as the theories of
commutative comonoids and bicommutative bimonoids
(\cref{ex:commutative-monoids,ex:bimonoids}). In
\cref{sec:interacting-supplies}, these theories were assembled into a
subcategory of $\PROP$:
\begin{equation*}
  \begin{tikzcd}[row sep=small]
    & \Theory{\CBimon} \arrow[r,hook] & \Theory{\Cone} \arrow[dr,hook] & \\
    \Theory{\CComon} \arrow[ur,hook] \arrow[dr,hook] & & & \Theory{\Vect_\R} \\
    & \Theory{\Conv} \arrow[uur,hook] \arrow[r,hook] &
    \Theory{\Aff_\R} \arrow[ur,hook] &
  \end{tikzcd}
\end{equation*}
This thin category, denoted $\TheoryLat{\Space}$, is a meet-semilattice. A
surprisingly large part of statistical modeling can be formulated algebraically
in a Markov category supplying this semilattice of PROPs, in the sense of
\cref{def:lattice-supply}.

\begin{definition}[Linear algebraic category]
  A symmetric monoidal category supplying the semilattice $\TheoryLat{\Space}$,
  not necessarily homomorphically, is called a \emph{linear algebraic (monoidal)
    category}.
\end{definition}

In particular, every linear algebraic category supplies commutative comonoids.
The modifier ``linear algebraic'' is used predictably: a linear algebraic Markov
category is a linear algebraic category that is also a Markov category, and a
linear algebraic cartesian category is a linear algebraic category that is also
a cartesian category.

As always in categorical logic, linear algebraic categories come in the small
and in the large. Small linear algebraic Markov categories are the better part
of \emph{statistical theories}, the topic of the next section. Both of
\cref{ex:extrinsic-spaces,ex:intrinsic-spaces} are large linear algebraic
categories, albeit cartesian ones. Our primary example of a large linear
algebraic Markov category, providing the intended semantics of statistical
theories, is the following.\footnote{This extrinsic definition of $\Stat$ is
  modeled on \cref{ex:extrinsic-spaces}. An intrinsic definition in the spirit
  of \cref{ex:intrinsic-spaces} could also be given, but is omitted in the
  interest of simplicity.}

\begin{definition}[Statistical semantics]
  The \emph{category of statistical semantics}, $\Stat$, has as objects the
  pairs $(V,A)$, where $V$ is a finite-dimensional real vector space and $A$ is
  a measurable\footnote{Every $d$-dimensional real vector space $V$ has a unique
    topology making it into a Hausdorff topological vector space, and is
    isomorphic as such to the Euclidean space $\R^d$ \cite[\S 9]{treves1967}. In
    particular, $V$ is a measurable space under its Borel $\sigma$-algebra.}
  subset of $V$, and as morphisms $(V,A) \to (W,B)$ the Markov kernels
  $A \to B$. In the symmetric monoidal category $(\Stat,\otimes,I)$, the
  monoidal product is defined on objects by
  \begin{equation*}
    (V,A) \otimes (W,B) := (V \oplus W, A \times B), \qquad I := (0,\{0\})
  \end{equation*}
  and on morphisms by the independent product (\cref{def:independent-product}).
  The category $\Stat$ then becomes a linear algebraic Markov category by simply
  marking out the vector space subsets that are closed under linear, affine,
  conical, convex, or additive combinations, or that are not closed at all.
\end{definition}

The category $\Stat$ can be seen as existing inside the category of Markov
kernels via a forgetful functor $\Stat \to \Markov$, defined on objects by
$(V,A) \mapsto A$ and on morphisms by the identity.

Discussion of objects and morphisms in linear algebraic categories is simplified
by the following conventions, mostly self-explanatory. In a linear algebraic
category $\cat{C}$ with supply $(\cat{P},s)$, a \emph{vector space object}, or a
\emph{vector space in $\cat{C}$}, is an object $x \in \cat{C}$ such that
$\cat{P}_x = \Theory{\Vect_\R}$. Note that when $\cat{C}$ is not a concrete
category, a vector space object may not actually be a vector space, for it may
have no underlying set. Nevertheless, we will sometimes abuse terminology by
calling it a ``vector space.'' We similarly speak of affine space, conical
space, convex space, and additive monoid objects. A \emph{discrete object} in
$\cat{C}$ is an object $x \in \cat{C}$ with $\cat{P}_x = \Theory{\CComon}$. As
for the morphisms, a morphism $f: x \to y$ in $\cat{C}$ between vector space
objects is \emph{linear} if it is an $s$-homomorphism with respect to the theory
of linear combinations, $\Theory{\LinComb_\R}$, a sub-theory of
$\Theory{\Vect_\R}$. In other words, $f: x \to y$ is linear if
\begin{equation*}
  \input{wiring-diagrams/algebra-statistics/additivity-lhs} =
  \input{wiring-diagrams/algebra-statistics/additivity-rhs}
  \qquad\text{and}\qquad
  \input{wiring-diagrams/algebra-statistics/homogeneity-lhs} =
  \input{wiring-diagrams/algebra-statistics/homogeneity-rhs},
  \quad \forall c \in \R.
\end{equation*}
Similarly, a morphism is affine, conic-linear, convex-linear, or additive if it
is a homomorphism with respect to the relevant theory of combinations
(\cref{ex:linear-combinations,ex:affine-combinations}). Finally, borrowing the
terminology of \cite{carboni1987a}, deterministic morphisms are also called
\emph{maps}, so that, for example, linear maps are morphisms that are both
deterministic and linear.

In a perhaps surprising fact, a Markov kernel in $\Stat$ that is linear must
already be deterministic.

\begin{theorem}[Linear Markov kernels] \label{thm:linear-markov-kernels}
  Let $M: V \to W$ be a Markov kernel between finite-dimensional real vector
  spaces $V$ and $W$. Suppose that $M$ is linear, so that
  \begin{equation*}
    \input{wiring-diagrams/algebra-statistics/markov-kernel-additivity-lhs} =
    \input{wiring-diagrams/algebra-statistics/markov-kernel-additivity-rhs}
    \qquad\text{and}\qquad
    \input{wiring-diagrams/algebra-statistics/markov-kernel-homogeneity-lhs} =
    \input{wiring-diagrams/algebra-statistics/markov-kernel-homogeneity-rhs},
    \quad \forall c \in \R.
  \end{equation*}
  Then the kernel $M$ is also deterministic, hence a linear map.
\end{theorem}
\begin{proof}
  We first prove the result in the one dimensional case, using Fourier analysis.
  Suppose that $M: \R \to \R$ is a linear Markov kernel. Let $X$ and $X'$ be
  i.i.d.\ random variables with distribution $M(1)$, and let
  $\varphi := \varphi_X: \R \to \C$ be the characteristic function of
  $X \sim M(1)$, defined by
  \begin{equation*}
    \varphi(t) := \varphi_X(t) := \E[e^{itX}]
      = \int_\R e^{itx}\, M(dx \given 1).
  \end{equation*}
  By the homogeneity of $M$, we have $aX \sim M(a)$ and
  $bX' \sim M(b)$, so by additivity and then homogeneity again,
  \begin{equation*}
    aX + bX' \sim M(a+b) = M((a+b)1) \sim (a+b) X.
  \end{equation*}
  Thus, by the convolution and scaling properties of the Fourier transform, the
  characteristic function satisfies
  \begin{equation*}
    \varphi_X(at) \cdot \varphi_X(bt) = \varphi_{aX + bX'}(t)
      = \varphi_{(a+b)X}(t) = \varphi_X((a+b)t), \qquad \forall a, b, t \in \R.
  \end{equation*}
  Setting $t=1$, we obtain Cauchy's multiplicative functional equation
  \begin{equation*}
    \varphi(a) \cdot \varphi(b) = \varphi(a+b), \qquad \forall a,b \in \R,
  \end{equation*}
  in the unknown characteristic function $\varphi: \R \to \C$.

  Let us solve this equation. If, for some $t_0 \in \R$, we had
  $\varphi(t_0) = 0$, then for \emph{any} $t \in \R$, we would have
  $\varphi(t) = \varphi(t - t_0) \varphi(t_0) = 0$, so that $\varphi$ is
  identically zero. This is impossible, since all characteristic functions have
  $\varphi(0) = 1$. Thus, $\varphi: \R \to \C^*$ vanishes nowhere. Moreover,
  since all characteristic functions satisfy $\norm{\varphi}_\infty \leq 1$, we
  must have $|\varphi(t)| = 1$ everywhere, for otherwise $|\varphi(t)| < 1$
  would imply that $|\varphi(-t)| = |\varphi(t)|^{-1} > 1$. Thus,
  $\varphi: \R \to \mathbb{T}$ takes values in the unit circle
  $\mathbb{T} := \{z \in \C: |z| = 1\}$. By a lemma of probability theory
  \cite[\S XV.1, Lemma 4]{feller1971}, this already implies that $X$ is a
  concentrated at a point. Alternatively, we can observe that the functional
  equation $\varphi(a+b) = \varphi(a) \varphi(b)$ makes
  $\varphi: \R \to \mathbb{T}$ into a (uniformly continuous) character of the
  additive group of real numbers. According to a famous result of Fourier
  analysis \cite[Proposition 7.1.1]{deitmar2005}, any such function has the form
  $\varphi(t) = e^{itc}$ for some constant $c \in \R$. Inverting the Fourier
  transform, $M(1)$ is the point mass $\delta_c$ and, by homogeneity,
  $M: \R \to \R$ is the deterministic kernel $M(x) = \delta_{cx}$.

  In the general case, we may assume that $V = \R^m$ and $W = \R^n$ for some
  dimensions $m$ and $n$. Let $M: \R^m \to \R^n$ be a linear Markov kernel. For
  each $i = 1,\dots,m$ and $j = 1,\dots,n$, let $\iota_i: \R \to \R^m$ be the
  inclusion into the $i$th coordinate and let $\pi_j: \R^n \to \R$ be the
  projection onto the $j$th coordinate. Reducing to the one dimensional case,
  each composite kernel $\iota_i \cdot M \cdot \pi_j: \R \to \R$ is linear and
  therefore deterministic. Since the only couplings of point masses are point
  masses, each kernel $\iota_i \cdot M: \R \to \R^n$ is also deterministic.
  Finally, using the linearity of $M$, it follows that $M: \R^m \to \R^n$ is
  deterministic.
\end{proof}

Although linear Markov kernels are deterministic, kernels obeying closely
related properties need not be. Under its standard parametrization by mean and
variance, the normal family is additive:
\begin{equation*}
  \input{wiring-diagrams/algebra-statistics/normal-additivity-lhs} =
  \input{wiring-diagrams/algebra-statistics/normal-additivity-rhs}.
\end{equation*}
Or, stated conventionally, if $X_1 \sim \Normal(\mu_1, \sigma_1^2)$ and
$X_2 \sim \Normal(\mu_2, \sigma_2^2)$ are independent random variables, then
their sum is $X_1+X_2 \sim \Normal(\mu_1+\mu_2, \sigma_1^2+\sigma_2^2)$. The
normal family is also homogeneous with exponents one and two, in the sense that
\begin{equation*}
  \input{wiring-diagrams/algebra-statistics/normal-homogeneity-lhs} =
  \input{wiring-diagrams/algebra-statistics/normal-homogeneity-rhs},
  \qquad
  \forall c \in \R.
\end{equation*}
Equivalently, if $X \sim \Normal(\mu, \sigma^2)$ and $c \in \R$, then
$cX \sim \Normal(c\mu, c^2 \sigma^2)$. As will be shown, the two properties
actually characterize the normal family, up to linear and conic-linear
transformations of the location and scale parameters.

Such properties are best understood within the more general class of
\emph{stable distributions} \cite{feller1971,nolan2018,samorodnitsky1994}. A
probability distribution on $\R^d$ is \emph{stable} if for two independent
random vectors $X$ and $X'$ having that distribution and for every pair of
constants $a,b > 0$, there exist constants $c > 0$ and $\vec{d} \in \R^d$ such
that $aX + bX' \xeq{d} cX + \vec{d}$. The distribution is \emph{strictly stable}
if this holds for $\vec{d} = 0$, and it is \emph{symmetric stable} if it is
stable and symmetric about the origin, meaning that $-X \xeq{d} X$. A symmetric
stable distribution is strictly stable, and a strictly stable distribution is
stable.

It can be shown that unless the stable distribution is concentrated at a point,
the scalars $a$, $b$, and $c$ in the defining equation must satisfy
$c = (a^\alpha + b^\alpha)^{1/\alpha}$ for some constant $0 < \alpha \leq 2$
\cite[Theorem 2.1.2]{samorodnitsky1994}. The distribution is then called
\emph{$\alpha$-stable}, with $\alpha$ being the \emph{index of stability} or
\emph{characteristic exponent}. The 2-stable distributions are exactly the
multivariate normal distributions. All other stable distributions, for
$0 < \alpha < 2$, are heavy tailed and have infinite variance. In statistics,
linear regression with stable errors offers an alternative to ordinary
least-squares regression when the data are heavy tailed \cite{nolan2013}.

The property of strict $\alpha$-stability can be stated in any linear algebraic
Markov category.

\begin{definition}[$\alpha$-stability]
  For any $0 < \alpha \leq 2$, a morphism $g: x \to y$ in a linear algebraic
  Markov category, whose codomain $y$ is a vector space object, is
  \emph{strictly $\alpha$-stable} if
  \begin{equation*}
    \input{wiring-diagrams/algebra-statistics/strict-stable-lhs} =
    \input{wiring-diagrams/algebra-statistics/strict-stable-rhs},
    \qquad\qquad
    \parbox{1.5in}{for all $a,b,c \in \R_+$ \\
      with $a^\alpha + b^\alpha = c^\alpha$.}
  \end{equation*}
  The morphism $g: x \to y$ is \emph{symmetric $\alpha$-stable} if this equation
  holds and, in addition,
  \begin{equation*}
    \input{wiring-diagrams/algebra-statistics/symmetric-stable-lhs} =
    \input{wiring-diagrams/algebra-statistics/symmetric-stable-rhs}.
  \end{equation*}
\end{definition}

According to the definition, a Markov kernel $M: \sspace{X} \to V$ taking values
in a finite-dimensional vector space $V$ is strictly $\alpha$-stable if, at
every point $x \in \sspace{X}$, the probability distribution $M(x)$ is strictly
$\alpha$-stable, and likewise for symmetric $\alpha$-stability.

From the algebraic viewpoint, rather than directly asserting a Markov kernel to
be $\alpha$-stable, it is natural to impose a stronger set of equations
characterizing the kernel as a scale family of $\alpha$-stable distributions.
The following definition generalizes the form of homogeneity obeyed by the
centered normal family.

\begin{definition}
  For any $\alpha > 0$, a morphism $g: s \to y$ in a linear algebraic Markov
  category, whose domain $s$ is a conical space object and codomain $y$ is a
  vector space object, is \emph{positively homogeneous with exponent $\alpha$},
  or \emph{positively $\alpha$-homogeneous}, if
  \begin{equation*}
    \input{wiring-diagrams/algebra-statistics/symmetric-homogeneity-lhs} =
    \input{wiring-diagrams/algebra-statistics/positive-homogeneity-rhs}
  \end{equation*}
  for all scalars $c \in \R_+$. The morphism $g: s \to y$ is \emph{homogeneous
    with exponent $\alpha$}, or \emph{$\alpha$-homogeneous}, if this equation
  holds for all scalars $c \in \R$.
\end{definition}

This property, in conjunction with additivity, implies $\alpha$-stability.

\begin{proposition} \label{prop:stability-by-homogeneity}
  Let $g: s \to y$ be a morphism in a linear algebraic Markov category, whose
  domain $s$ is a conical space and codomain $y$ is a vector space. For any
  $0 < \alpha \leq 2$,
  \begin{enumerate}[(i),nosep]
  \item if $g$ is additive and positively $\alpha$-homogeneous, then $g$ is
    strictly $\alpha$-stable;
  \item if $g$ is additive and $\alpha$-homogeneous, then $g$ is symmetric
    $\alpha$-stable.
  \end{enumerate}
\end{proposition}
\begin{proof}
  Under the hypotheses of part (i), for any scalars $a,b,c \in \R_+$ with
  $a^\alpha + b^\alpha = c^\alpha$,
  \begin{equation*}
    \input{wiring-diagrams/algebra-statistics/stability-by-homogeneity-1} =
    \input{wiring-diagrams/algebra-statistics/stability-by-homogeneity-2} =
    \input{wiring-diagrams/algebra-statistics/stability-by-homogeneity-3} =
    \input{wiring-diagrams/algebra-statistics/stability-by-homogeneity-4} =
    \input{wiring-diagrams/algebra-statistics/stability-by-homogeneity-5}.
  \end{equation*}
  Hence $g$ is strictly $\alpha$-stable. For part (ii), taking $c=-1$ in the
  definition of $\alpha$-homogeneity implies that $g$ is symmetric. Then part
  (i) implies that $g$ is symmetric $\alpha$-stable.
\end{proof}

The location-scale families derived from symmetric stable distributions will be
now presented equationally, with the normal and Cauchy location-scale families
as important special cases. For this, it is helpful to explicitly parameterize
the symmetric stable families. In the univariate case, parameterize the
symmetric $\alpha$-stable family $\SymStable{\alpha}: \R_+ \to \R$ by letting
$\SymStable{\alpha}(\lambda)$ have characteristic function
\begin{equation*}
  \varphi(t; \alpha, \lambda) := e^{- \lambda |t|^\alpha}.
\end{equation*}
In particular, $\SymStable{\alpha}(0)$ is the point mass at zero,
$\SymStable{1}(\lambda)$ is the Cauchy distribution $\Cauchy(0,\lambda)$, and
$\SymStable{2}(\lambda)$ is the normal distribution $\Normal(0,\sigma^2)$ with
variance $\sigma^2 = 2\lambda$.

In higher dimensions, when $0 < \alpha < 2$, the symmetric $\alpha$-stable
distributions are not parameterized by a finite-dimensional vector but by an
infinite-dimensional space of measures \cite[Theorem 2.4.3]{samorodnitsky1994}.
Specifically, let $\SymMeasSpace_+^d$ be the conical space of finite, symmetric,
nonnegative measures on the unit sphere $S^{d-1}$ in $\R^d$. Parameterize the
$d$-dimensional symmetric $\alpha$-stable family
$\SymStable{\alpha}_d: \SymMeasSpace_+^d \to \R^d$ by letting
$\SymStable{\alpha}_d(\Lambda)$ have characteristic function
\begin{equation*}
  \varphi(\vec{t}; \alpha, \Lambda) := \exp\left(
    - \int_{S^{d-1}} |\inner{\vec{t},\vec{s}}|^\alpha\, \Lambda(d\vec{s})
  \right).
\end{equation*}
In one dimension, $\SymStable{\alpha}_1(\Lambda)$ recovers the univariate
distribution $\SymStable{\alpha}(\lambda)$ with parameter
$\lambda = \Lambda\{\pm 1\} = 2 \Lambda\{1\}$. When $\alpha = 2$, identify the
family $\SymStable{\alpha}_d$ with the centered, rescaled normal family
$\sqrt{2} \cdot \Normal_d(0,-): \PSD^d \to \R^d$, where the rescaling is made
for consistency with the univariate case.

\begin{theorem}[Presentation of location-scale $\alpha$-stable families]
  \label{thm:stable-family-presentation}
  Let $M: V \times K \to W$ be a Markov kernel in $\Stat$, where $V$ and $W$ are
  vector spaces and $K$ is a convex cone. Suppose that $M$ is
  additive,
  \begin{equation*}
    \input{wiring-diagrams/algebra-statistics/location-scale-additivity-lhs} =
    \input{wiring-diagrams/algebra-statistics/location-scale-additivity-rhs},
  \end{equation*}
  and also that, for some $0 < \alpha \leq 2$, $M$ is homogeneous with exponents
  1 and $\alpha$,
  \begin{equation*}
    \input{wiring-diagrams/algebra-statistics/location-scale-homogeneity-lhs} =
    \input{wiring-diagrams/algebra-statistics/location-scale-homogeneity-rhs},
    \qquad
    \forall c \in \R.
  \end{equation*}
  Then the kernel $M: V \times K \to W$ is a location-scale family derived from
  the symmetric $\alpha$-stable distributions, in the following sense. Choose
  any basis identifying $W$ with $\R^d$, where $d = \dim W$, and let $L$ be the
  convex cone $\SymMeasSpace_+^d$ when $0 < \alpha < 2$ or else $\PSD^d$ when
  $\alpha = 2$. There exist a linear map $f: V \to W$ and a conic-linear map
  $g: K \to L$ such that
  \begin{equation*}
    M(x,s) = f(x) + \SymStable{\alpha}_d(g(s)),
    \qquad x \in V, \quad s \in K.
  \end{equation*}
\end{theorem}
\begin{proof}
  Since the composite of the kernel $M: K \times V \to W$ with any linear
  isomorphism $W \xrightarrow{\cong} \R^d$ preserves the stated properties of
  $M$, we may assume that $W = \R^d$. Define the Markov kernels $f: V \to W$ and
  $P: K \to W$ by inclusion into the first and second components of $M$:
  \begin{equation*}
    \input{wiring-diagrams/algebra-statistics/location-component-lhs} :=
    \input{wiring-diagrams/algebra-statistics/location-component-rhs}
    \qquad\text{and}\qquad
    \input{wiring-diagrams/algebra-statistics/scale-component-lhs} :=
    \input{wiring-diagrams/algebra-statistics/scale-component-rhs}.
  \end{equation*}
  Since $M$ is additive, it decomposes as the sum of $f$ and $P$:
  \begin{equation*}
    \input{wiring-diagrams/algebra-statistics/location-scale-decomposition-lhs}=
    \input{wiring-diagrams/algebra-statistics/location-scale-decomposition-rhs}.
  \end{equation*}
  Furthermore, short calculations using the properties of $M$ show that $f$ is
  linear and that $P$ is additive and $\alpha$-homogeneous. Therefore, by
  \cref{thm:linear-markov-kernels}, the kernel $f: V \to W$ is deterministic,
  hence a linear map, and by \cref{prop:stability-by-homogeneity}, the kernel
  $P: K \to W$ is symmetric $\alpha$-stable. The latter statement means that
  $P(s)$ is a symmetric $\alpha$-stable distribution at every point $s \in K$,
  hence there exists a function $g: K \to L$ such that $P$ decomposes as
  $P = g \cdot \SymStable{\alpha}_d$. Moreover, by the additivity of $P$ and
  $\SymStable{\alpha}_d$,
  \begin{equation*}
    \input{wiring-diagrams/algebra-statistics/conic-stable-additivity-1} =
    \input{wiring-diagrams/algebra-statistics/conic-stable-additivity-2} =
    \input{wiring-diagrams/algebra-statistics/conic-stable-additivity-3} =
    \input{wiring-diagrams/algebra-statistics/conic-stable-additivity-4} =
    \input{wiring-diagrams/algebra-statistics/conic-stable-additivity-5}.
  \end{equation*}
  Since the family $\SymStable{\alpha}_d$ is identifiable, it follows that $g$
  is additive. Similarly, the calculation
  \begin{equation*}
    \input{wiring-diagrams/algebra-statistics/conic-stable-homogeneity-1} =
    \input{wiring-diagrams/algebra-statistics/conic-stable-homogeneity-2} =
    \input{wiring-diagrams/algebra-statistics/conic-stable-homogeneity-3} =
    \input{wiring-diagrams/algebra-statistics/conic-stable-homogeneity-4} =
    \input{wiring-diagrams/algebra-statistics/conic-stable-homogeneity-5},
    \qquad \forall c \in \R_+
  \end{equation*}
  shows that $g$ is homogeneous. Thus the function $g: K \to L$ is conic-linear,
  which completes the proof.
\end{proof}

Two special cases of the theorem should be noted. When $\alpha = 2$, any Markov
kernel $M: \R^m \times \R_+ \to \R^n$ satisfying the hypotheses has the form
$M(x, \sigma^2) = \Normal_n(Ax, \sigma^2 V)$ for some matrices
$A \in \R^{n \times m}$ and $V \in \PSD^n$. This is precisely the sampling
distribution of a weighted linear model. Also, when $\alpha = 1$, a Markov
kernel $M: \R \times \R_+ \to \R$ satisfying the hypotheses is given by
$M(x, \gamma) = \Cauchy(ax, c \gamma)$ for some scalars $a \in \R$ and
$c \in \R_+$.

Another use of theorem is to present the isotropic multivariate stable families,
up to an absolute scale. The most important case is the isotropic normal family.

\begin{corollary}[Presentation of isotropic normal family]
  \label{cor:isotropic-normal-presentation}
  In any dimension $d$, a linear algebraic Markov category $\cat{C}$ containing
  a morphism $p: y^{\otimes d} \otimes s \to y^{\otimes d}$ can be presented
  such that for any supply preserving functor $M: \cat{C} \to \Stat$ with
  $M(y) = \R$ and $M(s) = \R_+$, the Markov kernel
  $M(p): \R^d \times \R_+ \to \R^d$ is the isotropic normal family, up to an
  absolute scale. That is, there exists a constant $\sigma_0^2 \in \R_+$ such
  that
  \begin{equation*}
    M(p)(\mu, \phi) = \Normal_d(\mu, \phi \sigma_0^2 I_d),
    \qquad \mu \in \R^d, \quad \phi \in \R_+.
  \end{equation*}
\end{corollary}
\begin{proof}
  Present a linear algebraic Markov category $\cat{C}$ as follows. Introduce a
  generating morphism $q: y \otimes s \to y$, where $y$ is a vector space object
  and $s$ is a conical space object. Add equations making $q$ be additive and
  homogeneous with exponents 1 and 2. Add the further equation that
  $(1_y \otimes 0_s) \cdot q = 1_y: y \to y$. This completes the presentation of
  $\cat{C}$.

  By \cref{thm:stable-family-presentation}, if $M: \cat{C} \to \Stat$ is a
  supply preserving functor with $M(y) = \R$ and $M(s) = \R_+$, then the kernel
  $M(q): \R \times \R_+ \to \R$ has the form
  $M(q)(\mu, \phi) = \Normal(\mu, \phi \sigma_0^2)$ for some constant
  $\sigma_0^2 \in \R_+$. Thus, if the morphism
  $p: y^{\otimes d} \otimes s \to y^{\otimes d}$ is defined by
  \begin{equation*}
    \input{wiring-diagrams/algebra-statistics/normal-present},
  \end{equation*}
  then the kernel $M(p): \R^d \times \R_+ \to \R^d$ has the stated form.
\end{proof}

This presentation of the isotropic normal family relies on its stability under
linear combinations. Another possible presentation is based on a
characterization by spherical symmetry of the normal distribution, named after
James Clerk Maxwell \cite[Proposition 13.2]{kallenberg2002}. Recall that a
random vector is $Y \in \R^d$ is \emph{spherically symmetric} about the origin,
or \emph{orthogonally invariant}, if for every orthogonal matrix
$Q \in \Orth(d)$, the rotated vector $Q Y$ has the same distribution as $Y$.

\begin{proposition}[Maxwell's theorem] \label{thm:maxwell}
  In any dimension $d \geq 2$, a random vector $Y \in \R^d$ has i.i.d.\ centered
  normal distribution if and only if $Y$ is spherically symmetric and has
  independent components.
\end{proposition}

Both assumptions, spherical symmetry and independence, are crucial. For example,
if $Y \sim \Normal_d(0, \sigma^2 I_d)$, then the random vectors $Y / \norm{Y}$,
uniformly distributed on the unit sphere $S^{d-1}$, and
$Y / \sqrt{\chi_\nu^2/\nu}$, multivariate $t$-distributed with $\nu$ degrees of
freedom, are spherically symmetric but do not have independent components. Like
the stable distributions, the spherically symmetric distributions can be
characterized by the form of their characteristic functions \cite[Theorem
2.1]{fang1990}. Also, Maxwell's theorem clearly fails in one dimension, since a
random variable $Y$ can be symmetric ($Y \xeq{d} -Y)$ under many distributions
besides the normal.

Because multiplication by a fixed matrix is reducible to a composite of copies,
sums, and scalar multiplications, Maxwell's theorem can be used to present the
isotropic normal family. For example, in dimension $d=2$, the equations of
spherical symmetry are
\begin{equation*}
  \input{wiring-diagrams/algebra-statistics/normal-present-maxwell-lhs} =
  \input{wiring-diagrams/algebra-statistics/normal-present-maxwell-rhs},
  \qquad \forall
  \begin{pmatrix} a_{11} & a_{12} \\ a_{21} & a_{22} \end{pmatrix}
  \in \Orth(2).
\end{equation*}

Maxwell's theorem is only one of many characterization of the normal
distribution \cite{bryc1995}. Even in the bivariate case, the characterization
by spherical symmetry does not depend on all orthogonal plane transformations,
but on only two: rotations through angles $\pi/2$ and $\pi/4$. This surprising
fact was first proved by George P\'olya \cite{polya1923}; a contemporary proof
appears in \cite[Theorem 3.1.1]{bryc1995}.

\begin{proposition}[P\'olya's theorem] \label{thm:polya}
  If $X$ and $Y$ are i.i.d.\ random variables such that
  \begin{equation*}
    X \xeq{d} \frac{1}{\sqrt{2}}(X + Y),
  \end{equation*}
  then $X$ is centered normal.
\end{proposition}

Although perhaps less conceptually satisfying than Maxwell's theorem, P\'olya's
theorem allows the isotropic normal family to be presented using far fewer
equations. The central equation in this presentation is
\begin{equation*}
  \input{wiring-diagrams/algebra-statistics/normal-present-polya-lhs} =
  \input{wiring-diagrams/algebra-statistics/normal-present-polya-rhs}.
\end{equation*}

Probabilists and statisticians have cataloged characterizations of many
probability distributions besides the normal and other stable distributions
\cite{nagaraja2006}, and it is interesting to ask which of them may be
formulated equationally inside a linear algebraic Markov category. Nevertheless,
it is not in the spirit of this work, or of structuralist mathematics generally,
to insist that each set of axioms admit a single interpretation. As the next
section demonstrates, statistical theories having many models are equally
important, as they bring out the structural commonalities between different
models.

\section{Statistical theories, models, and their morphisms}

The central notions of a statistical theory, a model of a statistical theory,
and a morphism between models are now defined. In this and the next section,
examples are drawn mainly from simple models of discrete data.
\cref{ch:zoo-statistics} presents a selection of more complex models for
regression and other continuous data.

\begin{definition}[Statistical theory]
  A \emph{statistical theory} $(\cat{T},p)$ is a small linear algebraic Markov
  category $\cat{T}$, together with a distinguished morphism $p: \theta \to x$.
  The morphism $p$ is called the \emph{sampling morphism}, its domain $\theta$
  is the \emph{parameter space object}, and its codomain $x$ is the \emph{sample
    space object}.
\end{definition}

The first two theories we consider are trivial but play special roles in the
algebra.

\begin{example}[Initial and terminal theories]
  \label{def:initial-terminal-theories}

  The \emph{initial theory} is the statistical theory whose underlying category
  $\gen{p}$ is freely generated by two discrete objects $\theta$ and $x$ and one
  morphism $p: \theta \to x$, which is also the sampling morphism.

  The \emph{terminal theory}, or \emph{discrete theory}, is the statistical
  theory whose underlying category $\gen{\emptyset}$ is freely generated on the
  empty set of generators. Thus, the category contains exactly one object $I$,
  necessarily the monoidal unit and a vector space object, and exactly one
  morphism, necessarily the identity morphism on $I$ and the sampling morphism
  for the theory.
\end{example}

The next two examples codify \cref{def:independence,def:exchangeability} as
statistical theories.

\begin{example}[i.i.d.\ samples] \label{ex:theory-iid-samples}
  For any $n \in \N$, the \emph{theory of $n$ independent and identically
    distributed (i.i.d.)\ samples} is freely generated by one morphism
  $p_0: \theta \to x$ on discrete objects $\theta$ and $x$, and the sampling
  morphism $p: \theta \to x^{\otimes n}$ is
  \begin{equation*}
    \input{wiring-diagrams/algebra-statistics/theory-iid-sampling-lhs} :=
    \input{wiring-diagrams/algebra-statistics/theory-iid-sampling-rhs}.
  \end{equation*}
  Note that this theory has the same underlying category as the initial theory,
  up to isomorphism, but a different sampling morphism. Hence, it is a different
  statistical theory.
\end{example}

\begin{example}[Exchangeable samples] \label{ex:theory-exchangeable-samples}
  For any $n \in \N$, the \emph{theory of $n$ exchangeable samples} is generated
  by discrete objects $\theta$ and $x$ and by one morphism
  $p: \theta \to x^{\otimes n}$, also the sampling morphism, subject to the
  equations
  \begin{equation*}
    \begin{tikzcd}
      \theta \arrow[r, "p"] \arrow[dr, "p"']
        & x^{\otimes n} \arrow[d, "\sigma"] \\
        & x^{\otimes n}
    \end{tikzcd}
    \qquad\text{for all permutations $\sigma \in S_n$}.
  \end{equation*}
  Since the symmetric group is generated by the adjacent transpositions, this
  set of equations may be replaced by the much smaller set consisting of
  $p \cdot (1_{x^{\otimes k}} \otimes \Braid_{x,x} \otimes 1_{x^{\otimes
      (n-k-2)}}) = p$ for every $0 \leq k \leq n-2$.
\end{example}

Apart from the terminal theory, the statistical theories considered so far have
been highly generic, admitting many different models. The following theory,
intended for discrete data, is slightly more restrictive.

\begin{example}[i.i.d.\ counts] \label{ex:theory-iid-counts}
  For any $n \in \N$, the \emph{theory of $n$ i.i.d.\ counts} is freely
  generated by one morphism $p_0: \theta \to x$, with $\theta$ a discrete object
  and $x$ an additive monoid object, and the sampling morphism $p: \theta \to x$
  is
  \begin{equation*}
    \input{wiring-diagrams/algebra-statistics/theory-iid-counts-sampling},
  \end{equation*}
  where the ellipsis indicates
  $p_0^{\otimes n}: \theta^{\otimes n} \to x^{\otimes n}$, the $n$-fold product
  of $p_0$.
\end{example}

For the definition of a model of a statistical theory, recall the notion of a
symmetric monoidal functor that preserves the supply, not necessarily strictly
(\cref{def:lattice-supply-preservation}).

\begin{definition}[Statistical model]
  A \emph{model} of a statistical theory $(\cat{T}, \theta \xrightarrow{p} x)$,
  or for short a \emph{statistical model}, is a supply preserving functor
  $M: \cat{T} \to \Stat$. The Markov kernel $P := M(p)$ is called the
  \emph{sampling distribution} or the \emph{data distribution}, its domain
  $\Omega := M(\theta)$ is the \emph{parameter space}, and its codomain
  $\sspace{X} := M(x)$ is the \emph{sample space}.
\end{definition}

In classical statistics it is the sampling distribution
$P: \Omega \to \sspace{X}$ that would be considered the statistical model, with
no reference to a theory. For the remainder of this text, a ``statistical
model'' will be a model of a definite statistical theory unless otherwise noted.
This entails no loss of generality, as any Markov kernel
$P: \Omega \to \sspace{X}$ in $\Stat$ is the sampling distribution of a unique
model of the initial theory $\gen{p}$. The model $M: \gen{p} \to \Stat$ is
completely determined by the assignment $M(p) := P$. This is one sense in which
the initial statistical theory is initial. Of course, when the sampling
distribution has extra structure, it is more informative to view it as a model
of a richer theory. The terminal statistical theory is still more trivial,
having exactly one model: the identity morphism on the zero-dimensional vector
space.

The following models of the theory of i.i.d.\ counts are more concrete.

\begin{example}[Models of i.i.d.\ counts] \label{ex:models-iid-counts}
  Let $(\cat{T}, p)$ be the theory of $n$ i.i.d.\ counts, defined in
  \cref{ex:theory-iid-counts}. The \emph{binomial model} $M: \cat{T} \to \Stat$
  on $n$ trials assigns $M(\theta) = [0,1]$, the unit interval; $M(x) = \N$, the
  natural numbers; and $M(p_0) = \Ber: [0,1] \to \N$, the Bernoulli family,
  given by $\Ber(\pi) := \pi \delta_1 + (1-\pi) \delta_0$. By functorality, the
  sampling distribution $M(p)$ is the binomial family
  $\Binom(n;-): [0,1] \to \N$ on $n$ trials.

  Generalizing the binomial model, for any $k \geq 2$, the \emph{multinomial
    model} $M_k: \cat{T} \to \Stat$ on $k$ classes and $n$ trials assigns the
  morphism $p_0$ to be the categorical family on $k$ classes,
  \begin{equation*}
    \CatDist_k: \Delta^{k-1} \to \N^k, \quad
    (\pi_1,\dots,\pi_k) \mapsto \pi_1 \delta_{e_1} + \cdots + \pi_k \delta_{e_k}.
  \end{equation*}
  Here $\Delta^{k-1}$ is the $(k-1)$-dimensional probability simplex in $\R^k$
  and $e_1,\dots,e_k$ are the standard basis vectors. The sampling distribution
  $M_k(p)$ is then the multinomial family
  $\Multinom_k(n, -): \Delta^{k-1} \to \N^k$ on $k$ classes and $n$ trials. When
  $k=2$, the multinomial model is not identical to the binomial model; however,
  as models of the initial theory, the two statistical models are isomorphic, as
  will be seen shortly.

  The binomial and multinomial models are the most obvious models of the theory
  of i.i.d.\ counts, but they are not the only ones. A negative binomial model
  $\cat{T} \to \Stat$ assigns $p_0$ to be the geometric family
  \begin{equation*}
    \GeomDist: [0,1] \mapsto \N, \quad
    \pi \mapsto \sum_{k=0}^\infty (1-\pi)^k \pi\, \delta_k,
  \end{equation*}
  so that the sampling distribution is the reparameterized negative binomial
  $\pi \mapsto \NegBinom(n, 1-\pi)$. Yet another model, the Poisson model,
  assigns $p_0$ to be the Poisson family $\Pois: \R_+ \to \N$, under which the
  sampling distribution is the rescaled family
  $\Pois(n \cdot -): \mu \mapsto \Pois(n\mu)$.

  The binomial, multinomial, negative binomial, and Poisson models of the theory
  of i.i.d.\ counts all possess extra structure that can be described by richer
  statistical theories. The discrete object $\theta$ may in all cases be
  replaced by a convex space object. For the multinomial specifically, let
  $(\cat{T}_k,p)$ be the \emph{theory of $n$ i.i.d.\ $k$-dimensional counts},
  freely generated by a discrete (or convex space) object $\theta$, an additive
  monoid object $x$, and a single morphism $p_0: \theta \to x^{\otimes k}$. The
  sampling morphism $p$ is constructed from $p_0$ as before. The new multinomial
  model $M_k: \cat{T}_k \to \Stat$ assigns the objects $\theta$ and $x$ to be
  $\Delta^{k-1}$ and $\N$, respectively, and the morphism $p_0$ to be the
  categorical family $\CatDist_k$. The sampling distribution is again the
  multinomial family $\Multinom_k(n,-)$. However, neither the binomial nor the
  negative binomial are models of the theory, as the natural numbers are not the
  $k$-fold power of another set for any $k \geq 2$.

  A richer statistical theory for the Poisson counts model might take $\theta$
  to be a conical space object, $x$ to be an additive monoid object, and
  $p_0: \theta \to x$ to be an additive morphism,
  \begin{equation*}
    \input{wiring-diagrams/algebra-statistics/theory-iid-counts-additivity-lhs}
    =
    \input{wiring-diagrams/algebra-statistics/theory-iid-counts-additivity-rhs}.
  \end{equation*}
  This equation states the well known additivity property of the Poisson family,
  that if $X_1 \sim \Pois(\mu_1)$ and $X_2 \sim \Pois(\mu_2)$ are independent,
  then $X_1 + X_2 \sim \Pois(\mu_1 + \mu_2)$. Defining the sampling morphism
  $p: \theta \to x$ as before, the equation $p = n \cdot p_0$ is easily deduced
  within the theory itself. None of the binomial, multinomial, or negative
  binomial are models of the theory, as their parameter spaces are not convex
  cones.
\end{example}

Several statements made above are clarified by the concept of a morphism between
statistical models.

\begin{definition}[Morphism of statistical models]
  Let $M, N: \cat{T} \to \Stat$ be two models of a statistical theory
  $(\cat{T},p)$. A \emph{morphism of models}, or \emph{model homomorphism}, from
  $M$ to $N$ is a monoidal natural transformation $\alpha: M \to N$.
\end{definition}

Every statistical theory $(\cat{T},p)$ has a \emph{category of models}, denoted
$\Model{\cat{T},p}$, with models $\cat{T} \to \Stat$ as objects and model
homomorphisms as morphisms. Since the category of models does not depend on the
choice of $p$, it will also be abbreviated as $\Model{\cat{T}}$. As in any
functor category, composition is given by vertical composition of natural
transformations.

A morphism $\alpha: M \to N$ of statistical models certainly preserves the
sampling morphism $p: \theta \to x$, making the diagram
\begin{equation*}
  \begin{tikzcd}
    M(\theta) \arrow[r, "\alpha_\theta"] \arrow[d, "M(p)"']
      & N(\theta) \arrow[d, "N(p)"] \\
    M(x) \arrow[r, "\alpha_x"] & N(x)
  \end{tikzcd}
\end{equation*}
commute. When the theory is the initial theory, that is almost the only
requirement, but regardless of the theory, a model homomorphism must preserve
every morphism in it. Thus, as a general principle, the richer the statistical
theory, the fewer the morphisms between its models.

According to \cref{prop:supply-natural-transformation}, every component of a
monoidal natural transformation between supply preserving functors
$\cat{C} \to \cat{D}$ is a supply homomorphism with respect to the supply
assignments in $\cat{C}$. In the context of statistical models, this means that
for any morphism $\alpha: M \to N$ between models of a statistical theory
$(\cat{T},p)$, every component $\alpha_x: M(x) \to N(x)$ is a supply
homomorphism with respect to the supply at $x \in \cat{T}$. An especially
important consequence is:
\begin{proposition}
  Every component of a morphism of statistical models is deterministic.
\end{proposition}
\noindent
Indeed, while statistical models are inherently stochastic, there seems little
reason to think that morphisms of models should be.\footnote{Other contexts may
  call for different reasoning. In another work by the author
  \cite{patterson2019transport}, it is the models that are deterministic and the
  model homomorphisms that are stochastic, leading to a form of optimal
  transport for structured data.}

Given that statistical theories are usually presented by generators and
relations, it is useful to know that a model homomorphism is completely
determined by its components on a generating set of objects for the theory.
Moreover, in order to establish naturality, it is enough that the components be
supply homomorphisms and that the naturality condition hold on a generating set
of morphisms for the theory. This is the content of the more precisely stated
\cref{lemma:natural-transformation-generators}. We use this fact regularly and
tacitly, as in the following example.

\begin{example}[Morphisms of i.i.d.\ count models]
  Let $M$ be the binomial model and $M_k$ be the multinomial model from
  \cref{ex:models-iid-counts}, each on $n$ trials.

  The binomial model is clearly ``the same,'' in some sense, as the multinomial
  model on $k=2$ classes, and we expect this sameness to be reflected by an
  isomorphism $M \cong M_2$ of models. Define a transformation
  $\alpha: M \to M_2$ by
  \begin{equation*}
    \alpha_\theta: \begin{aligned}
      & [0,1] \to \Delta^2 \\[-2\jot] & \pi \mapsto (\pi,1-\pi)
    \end{aligned}
    \qquad\text{and}\qquad
    \alpha_x: \begin{aligned}
      & \N \to \N^2 \\[-2\jot] & m \mapsto (m,n-m).
    \end{aligned}
  \end{equation*}
  This transformation preserves the sampling morphism $p: \theta \to x$, as the
  diagram
  \begin{equation*}
    \begin{tikzcd}
      {[0,1]} \arrow[r, "\alpha_\theta"] \arrow[d, "{\Binom(n,-)}"']
        & \Delta^2 \arrow[d, "{\Multinom_k(n,-)}"] \\
      \N \arrow[r, "\alpha_x"] & \N^2
    \end{tikzcd}
  \end{equation*}
  commutes, and hence $\alpha: M \to M_2$ is a model homomorphism with respect
  to the initial statistical theory. Another transformation $\beta: M_2 \to M$,
  where both $\beta_\theta: \Delta^2 \to [0,1]$ and $\beta_x: \N^2 \to \N$ are
  projections onto the first coordinate, also preserves the sampling morphism
  and is mutually inverse to $\alpha$. Thus, $\alpha: M \cong M_2$ is a model
  isomorphism, again with respect to the initial theory.

  The qualification about the initial theory is subtle but important. The map
  $\alpha_x: \N \to \N^2$ is not additive and thus cannot be the component of
  model homomorphism $\alpha: M \to M_2$ when $M$ and $M_2$ are regarded as
  models of the theory of $n$ i.i.d.\ counts, as originally intended. Indeed,
  the theory of $n$ i.i.d.\ counts implicitly contains, through its underlying
  category, not just the sum of $n$ trials but a sum of $n'$ trials for
  \emph{every} number $n' \in \N$, yet the map $\alpha_x: m \mapsto (m,n-m)$
  depends on the fixed number $n$. The transformation $\alpha$ does not preserve
  the full structure of the theory.

  For any $k>1$, the multinomial model $M_k$ of the theory of $n$ i.i.d.\ counts
  has nontrivial symmetries, which manifest as model automorphisms. Given a
  permutation $\sigma \in S_k$, define the transformation
  $\alpha(\sigma): M_k \to M_k$ by
  \begin{equation*}
    \alpha(\sigma)_\theta: \begin{aligned}
      \Delta^k &\to \Delta^k \\[-2\jot]
      (\pi_1,\dots,\pi_k) &\mapsto (\pi_{\sigma(1)},\dots,\pi_{\sigma(k)}),
    \end{aligned}
    \qquad
    \alpha(\sigma)_x: \begin{aligned}
      \N^k &\to \N^k \\[-2\jot]
      (m_1,\dots,m_k) &\mapsto (m_{\sigma(1)},\dots,m_{\sigma(k)}).
    \end{aligned}
  \end{equation*}
  Then $\alpha(\sigma)_\theta$ is a convex-linear map, $\alpha(\sigma)_x$ is an
  additive map, and $\alpha(\sigma)$ preserves the generating morphism
  $p_0: \theta \to x$, since the diagram
  \begin{equation*}
    \begin{tikzcd}
      \Delta^k \arrow[r, "\alpha(\sigma)_\theta"] \arrow[d, "\CatDist_k"']
        & \Delta^k \arrow[d, "\CatDist_k"] \\
      \N^k \arrow[r, "\alpha(\sigma)_x"] & \N^k
    \end{tikzcd}
  \end{equation*}
  commutes. Thus $\alpha(\sigma)$ is an endomorphism of the multinomial model
  $M_k$. It is also invertible, with inverse $\alpha(\sigma^{-1})$, making it a
  model automorphism. One can further check that every automorphism of $M_k$ has
  this form and that the map $\alpha: S_k \to \Aut(M_k)$ is even a group
  isomorphism. In summary, the automorphism group of the multinomial model $M_k$
  of the theory of $n$ i.i.d.\ counts is isomorphic to the symmetric group
  $S_k$, confirming the intuition that the multinomial model is invariant under
  relabeling of the classes.
\end{example}

Statistical theories and models are \emph{Bayesian} when they are accompanied by
a prior.

\begin{definition}[Bayesian theories and models] \label{def:bayesian-theory}
  A \emph{Bayesian (statistical) theory} $(\cat{T},p,\pi)$ is a statistical
  theory $(\cat{T}, \theta \xrightarrow{p} x)$ together with a distinguished
  morphism $I \xrightarrow{\pi} \theta$, the \emph{prior morphism}.

  A \emph{model} of a Bayesian theory is a model $M: \cat{T} \to \Stat$ of the
  underlying statistical theory. The probability distribution $M(\pi)$ is called
  the \emph{prior distribution} and the distribution $M(\pi \cdot p)$ is called
  the \emph{marginal distribution} or the \emph{prior predictive distribution}.
\end{definition}

Morphisms of Bayesian models, and the category of models of a Bayesian theory,
are those of the underlying statistical models and theory. In practice, however,
extending a ``frequentist'' statistical theory $(\cat{T},p)$ to a Bayesian one
typically requires the category $\cat{T}$ to be enlarged with another morphism,
representing the prior, and this changes the class of models and their
morphisms.

\section{Morphisms of statistical theories and model migration}

Although morphisms of statistical models are useful for formalizing what it
means for two models to be isomorphic, or for defining the group of symmetries
of a model, it is arguably the morphisms of statistical \emph{theories} that are
more important, as they enable relationships to be stated between models of
different theories. Morphisms of statistical theories come in several variants.
The simplest are the strict morphisms.

\begin{definition}[Strict theory morphisms] \label{def:strict-theory-morphism}
  A \emph{(strict) morphism} from one statistical theory $(\cat{T},p)$ to
  another $(\cat{T}',p')$ is a supply preserving functor
  $F: \cat{T} \to \cat{T}'$ that strictly preserves the sampling morphism,
  satisfying $F(p) = p'$.
\end{definition}

Statistical theories and theory morphisms form a category, with composition and
identities defined as usual for functors. The initial and terminal theories
(\cref{def:initial-terminal-theories}) derive their names from the following
fact.

\begin{proposition}
  In the category of statistical theories, the initial statistical theory is an
  initial object and the terminal theory is a terminal object.
\end{proposition}
\begin{proof}
  For any statistical theory $(\cat{T},p)$, a theory morphism
  $F: (\gen{p},p) \to (\cat{T},p)$ must assign $F(p) = p$, and since $\gen{p}$
  is generated by $p$, this assignment uniquely determines a supply preserving
  functor $F: \gen{p} \to \cat{T}$. In the other direction, a theory morphism
  $(\cat{T},p) \to (\gen{\emptyset}, 1_I)$ necessarily assigns every object in
  $\cat{T}$ to the unique object $I$ in $\gen{\emptyset}$ and every morphism in
  $\cat{T}$, including $p$, to the unique morphism $1_I$ in $\gen{\emptyset}$.
\end{proof}

Theory morphisms commonly represent inclusions of one theory into another that
is larger or richer. The universal morphisms out of the initial theory are
extreme examples of such morphisms. The subsequent examples are more concrete
but still very simple.

\begin{example}[i.i.d.\ and exchangeable samples]
  Let $(\cat{T}_{\mathrm{iid}}, p)$ be the theory of $n$ i.i.d.\ samples from
  \cref{ex:theory-iid-samples}, and let $(\cat{T}_{\mathrm{ex}}, p)$ be the
  theory of $n$ exchangeable samples from \cref{ex:theory-exchangeable-samples}.
  Define a supply preserving functor
  $F: \cat{T}_{\mathrm{ex}} \to \cat{T}_{\mathrm{iid}}$ as the identity on
  objects and by $F(p) := p_0 \cdot \Copy_{x,n} = p$ on morphisms. The functor
  $F$ is well-defined because, as was seen in \cref{sec:markov-categories},
  being independent and identically distributed implies being exchangeable. By
  construction, $F$ preserves the sampling morphism and is thus a morphism
  $(\cat{T}_{\mathrm{ex}}, p) \to (\cat{T}_{\mathrm{iid}}, p)$ of statistical
  theories.
\end{example}

\begin{example}[i.i.d.\ counts]
  Let $(\cat{T},p)$ be the theory of $n$ i.i.d.\ counts and $(\cat{T}_k,p)$ be
  the theory of $n$ i.i.d.\ $k$-dimensional counts from
  \cref{ex:theory-iid-counts,ex:models-iid-counts}. The supply preserving
  functor $F: \cat{T} \to \cat{T}_k$ determined by $F(\theta) := \theta$,
  $F(x) := x^{\otimes k}$, and $F(p_0) := p_0$ preserves the sampling morphism
  ($F(p) = p$) and is therefore a theory morphism
  $(\cat{T},p) \to (\cat{T}_k,p)$. This morphism expresses the mundane idea that
  a $k$-dimensional object $x^{\otimes k}$ can be regarded as a basic object
  $x'$ by forgetting about its division into $k$ components.
\end{example}

It may seem backwards that every i.i.d.\ probability model is exchangeable,
while the theory morphism goes from the theory of exchangeable samples to the
theory of i.i.d.\ samples, or that every $k$-dimensional counts model is a
counts model, while the theory morphism goes from the theory of counts to the
theory of $k$-dimensional counts. But this is no accident. The directionality of
the functor reflects a contravariance that is universal to categorical logic.
Namely, every morphism between statistical theories induces a morphism between
the corresponding categories of statistical models, but \emph{going in the
  opposite direction}.

\begin{definition}[Pullback model migration]
  Let $F: (\cat{T},p) \to (\cat{T}',p')$ be a morphism of statistical theories.
  The \emph{pullback functor} $F^*: \Model{\cat{T}'} \to \Model{\cat{T}}$ from
  the category of models of $\cat{T}'$ to the category of models of $\cat{T}$ is
  defined on objects by pre-composition,
  \begin{equation*}
    (M: \cat{T}' \to \Stat) \quad\xmapsto{F^*}\quad (FM: \cat{T} \to \Stat),
  \end{equation*}
  and on morphisms by pre-whiskering,
  \begin{equation*}
    (\alpha: M \to N) \quad\xmapsto{F^*}\quad (F\alpha: FM \to FN),
  \end{equation*}
  where the transformation $F\alpha$ has components
  $(F\alpha)_x := \alpha_{Fx}: M(F(x)) \to N(F(x))$.
\end{definition}

The pullback construction is a recurring theme in categorical logic and its
applications. In the context of relational databases, Spivak has called $F^*$
the ``pullback data migration functor'' induced by a database schema translation
$F$ \cite{spivak2012}. By analogy, we call the functor $F^*$ induced by a
morphism $F$ of statistical theories a \emph{pullback model migration functor}.
Note that the operation of taking pullbacks is itself contravariantly
functorial. That is, $(F \cdot G)^* = G^* \cdot F^*$ whenever $F$ and $G$ are
composable functors, and also $1_{\cat{T}}^* = 1_{\Model{\cat{T}}}$.

Returning to the examples, for any model $M$ of the theory
$\cat{T}_{\mathrm{iid}}$ of $n$ i.i.d.\ samples, applying the pullback functor
$F^*: \Model{\cat{T}_{\mathrm{iid}}} \to \Model{\cat{T}_{\mathrm{ex}}}$ yields a
model $F^*(M)$ of the theory $\cat{T}_{\mathrm{ex}}$ of $n$ exchangeable
samples, as expected. Likewise for the theories of $n$ i.i.d.\ counts.
Similarly, the universal morphisms \emph{out of} the initial statistical theory
induce forgetful functors \emph{into the} initial theory's category of models.
This formalizes the earlier observation that every statistical model can be
regarded as a model of the initial theory.

Morphisms of statistical theories, as defined so far, cannot express certain
relationships that one would like to formalize. Letting $(\cat{T},p_n)$ be the
theory of $n$ i.i.d.\ samples from \cref{ex:theory-iid-samples}, one would
expect that for every pair of numbers $m \leq n$, there would be a theory
morphism $(\cat{T},p_m) \to (\cat{T},p_n)$, embedding an i.i.d\ sample of size
$m$ as a sub-sample of an i.i.d.\ sample of size $n$. The underlying functor
$\cat{T} \to \cat{T}$ should be the identity, but then there cannot be a theory
morphism because $p_m \neq p_n$ whenever $m \neq n$. This situation, where
statistical theories have sampling morphisms with varying parameter or sample
spaces, occurs commonly when a family of theories is indexed by parameter
dimensionality or sample size. Such families are accommodated by relaxing the
definition of a theory morphism.

\begin{definition}[Lax and colax theory morphisms]
  \label{def:lax-theory-morphism}
  A \emph{lax morphism} from one statistical theory
  $(\cat{T},\, \theta \xrightarrow{p} x)$ to another
  $(\cat{T}',\, \theta' \xrightarrow{p'} x')$ consists of a supply preserving
  functor $F: \cat{T} \to \cat{T}'$, together with morphisms
  $f_0: F\theta \to \theta'$ and $f_1: Fx \to x'$ in $\cat{T}'$, such that the
  diagram
  \begin{equation*}
    \begin{tikzcd}
      F\theta \arrow[r, "Fp"] \arrow[d, "f_0"'] & Fx \arrow[d, "f_1"] \\
      \theta' \arrow[r, "p'"] & x'
    \end{tikzcd}
  \end{equation*}
  commutes. Dually, a \emph{colax theory morphism} is a supply preserving
  functor $F: \cat{T} \to \cat{T}'$, together with morphisms
  $f_0: \theta' \to F\theta$ and $f_1: x' \to Fx$ in $\cat{T}'$, such that the
  diagram
  \begin{equation*}
    \begin{tikzcd}
      \theta' \arrow[r, "p'"] \arrow[d, "f_0"'] & x' \arrow[d, "f_1"] \\
      F\theta \arrow[r, "Fp"] & Fx
    \end{tikzcd}
  \end{equation*}
  commutes. When both components $f_0$ and $f_1$ of a lax or colax morphism of
  theories are isomorphisms, the theory morphism is called \emph{strong}. When
  moreover both components are identities, the theory morphism is \emph{strict},
  recovering the original \cref{def:strict-theory-morphism}.
\end{definition}

\begin{example}[i.i.d.\ samples of different sizes]
  \label{ex:morphism-iid-sample-sizes}
  For any numbers $m \leq n$, a lax morphism
  $(1_\cat{T}, 1_\theta, \pi): (\cat{T},p_n) \to (\cat{T},p_m)$ between the
  theories of $n$ and $m$ i.i.d.\ samples is defined by taking
  $\pi: x^{\otimes n} \to x^{\otimes m}$ to be a projection morphism, such as
  the projection $\pi_{m,n-m}$, which discards the last $n-m$ components of
  $x^{\otimes n}$, or the projection $\pi_{n-m,m}'$, which discards the first
  $n-m$ components of $x^{\otimes n}$. In the former case, the laxness condition
  is the equality
  \begin{equation*}
    \input{wiring-diagrams/algebra-statistics/theory-iid-sampling-projection-lhs}
    \quad=\quad
    \input{wiring-diagrams/algebra-statistics/theory-iid-sampling-projection-rhs}.
  \end{equation*}
  A colax morphism $(1_\cat{T}, 1_\theta, \pi): (\cat{T},p_m) \to (\cat{T},p_n)$
  is defined in exactly the same way (but note the reversed direction). The lax
  morphism can be interpreted as a projection of statistical theories and the
  colax morphism as an inclusion of theories.

  On the other hand, for any positive numbers $m \neq n$, there are no lax or
  colax morphisms between the theories of $m$ and $n$ i.i.d.\ counts. Due to the
  aggregation happening in the sampling morphisms, it is not possible to project
  or include from one sample to another of differing size.
\end{example}

\begin{example}[Exchangeable samples of different sizes]
  \label{ex:morphism-exchangeable-sample-sizes}
  Let $(\cat{T}_n, p_n)$ be the theory of $n$ exchangeable samples and for
  numbers $m \leq n$, let $\pi: x^{\otimes n} \to x^{\otimes m}$ be any
  projection morphism, as in the previous example. Define a supply preserving
  functor $F: \cat{T}_m \to \cat{T}_n$ as the identity on objects and by
  $F(p_m) := p_n \cdot \pi$ on morphisms. The functor is well-defined because
  exchangeability is preserved under projection. Thus, by construction,
  $(F,1_\theta,\pi): (\cat{T}_m, p_m) \to (\cat{T}_n, p_n)$ is a colax morphism
  of statistical theories.

  In contrast to the case of i.i.d.\ samples, there is no lax morphism
  $(\cat{T}_n, p_n) \to (\cat{T}_m, p_m)$ going in the opposite direction, since
  the generating morphism $p_n$ in $\cat{T}_n$ has no evident image in
  $\cat{T}_m$. This situation, where a ``smaller'' theory has a colax morphism
  including it in a ``larger'' theory but has no lax morphism in the opposite
  direction, seems to be the more common one.
\end{example}

Lax morphisms of statistical theories, like strict ones, are composable. The
\emph{composite} of lax morphisms
\begin{equation*}
  (\cat{T}, \theta \xrightarrow{p} x)
  \xrightarrow{(F, f_0, f_1)}
  (\cat{T}', \theta' \xrightarrow{p'} x')
  \xrightarrow{(G, g_0, g_1)}
  (\cat{T}'', \theta'' \xrightarrow{p''} x'')
\end{equation*}
is the lax morphism $(F \cdot G,\, Gf_0 \cdot g_0, Gf_1 \cdot g_1)$, where the
laxness condition is verified by the pasting of commutative squares
\begin{equation*}
  \begin{tikzcd}[column sep=large]
    G(F\theta) \arrow[r, "G(Fp)"] \arrow[d, "G f_0"']
      & G(Fx) \arrow[d, "G f_1"] \\
    G\theta' \arrow[r, "Gp'"] \arrow[d, "g_0"']
      & Gx' \arrow[d, "g_1"] \\
    \theta'' \arrow[r, "p''"]
      & x''.
  \end{tikzcd}
\end{equation*}
The composite of colax morphisms is defined dually. In both cases, the identity
morphism on $(\cat{T}, \theta \xrightarrow{p} x)$ is the (strict) morphism
$(1_{\cat{T}}, 1_\theta, 1_x)$. Thus, there is a category of statistical
theories and lax morphisms and also of statistical theories and colax
morphisms.\footnote{For the reader familiar with 2-category theory, we note that
  the category of statistical theories and (co)lax theory morphisms is the
  (co)lax coslice 2-category of the 2-category of small linear algebraic Markov
  categories, supply preserving functors, and monoidal natural transformations,
  under the underlying category of the initial statistical theory. This
  observation leads to the correct definition of a 2-morphism between
  statistical theory morphisms, not given in the main text.} Both categories
contain the strong theory morphisms as a subcategory, and the strict theory
morphisms as a subcategory of those.

With respect to lax or colax morphisms, or even strong morphisms, the terminal
theory is still terminal but the initial theory is only \emph{weakly} initial.
That is, every statistical theory has a morphism out of the initial theory but
this morphism need not be unique. Suppose, for example, that $F$ is the (unique)
strict morphism from the initial theory into the theory
$(\cat{T}, \theta \xrightarrow{p} x^{\otimes n})$ of $n$ exchangeable samples.
If $\Braid: x^{\otimes n} \to x^{\otimes n}$ is any symmetry isomorphism, then
$(F, 1_\theta, \Braid)$ is a strong theory morphism.

A large part of statistics is about \emph{hypothesis tests}, or formal tests
that the data conforms to a sub-model of a larger statistical model
\cite{lehmann2005}. Hypothesis testing is classically formulated by partitioning
the parameter space $\Omega$ of a statistical model $P: \Omega \to \sspace{X}$
into subsets $\Omega_0$ and $\Omega_1$, so that
$\Omega = \Omega_0 \sqcup \Omega_1$. One of these sets, say $\Omega_0$, is
designated as the class of \emph{null} parameters, and the other set $\Omega_1$
as the class of \emph{alternative} parameters. The problem is then to test the
\emph{null hypothesis} $H_0: \theta \in \Omega_0$ against the \emph{alternative
  hypothesis} $H_1: \theta \in \Omega_1$, where $\theta$ is the true parameter
according to the full probability model.

The standard formalism is misleading in treating null and alternative hypotheses
symmetrically. Practically speaking, the null and alternative are nearly always
treated asymmetrically during both model specification and statistical
inference. The null hypothesis is understood to be a meaningful sub-model of the
full model and the alternative is ``everything else.'' In the extreme case of a
\emph{point null hypothesis}, the null $\Omega_0 = \{\theta_0\}$ consists of a
single point. In general, the null hypothesis $\Omega_0$ is a subset of $\Omega$
possessing special structure of logical or scientific interest, while the
alternative $\Omega_1 = \Omega \setminus \Omega_0$ fails to possess this
structure.

From the viewpoint of categorical logic, a statistical hypothesis is better
understood as a morphism from a statistical theory $(\cat{T},p)$, representing
the full model, to another theory $(\cat{T}_0,p_0)$, representing the sub-model.
For example, to define a point null hypothesis for a given theory
$(\cat{T}, \theta \xrightarrow{p} x)$, let $\cat{T}_0$ be generated by $\cat{T}$
together with a map $I \xrightarrow{\theta_0} \theta$. Then a model $M$ of the
theory $(\cat{T}_0, I \xrightarrow{\theta_0 \cdot p} x)$ is a model of the
original theory, plus a chosen parameter $M(\theta_0)$ in the parameter space
$M(\theta)$. The null hypothesis is represented by the colax morphism
\begin{equation*}
  (\iota, \theta_0, 1_x): (\cat{T}, p) \to (\cat{T}_0, \theta_0 p),
\end{equation*}
where $\iota: \cat{T} \hookrightarrow \cat{T}_0$ is the inclusion functor. True
point null hypotheses occur rarely, as the parameter of interest, say a mean, is
usually accompanied by an unknown nuisance parameter, such as a variance.

More substantial examples are hypotheses of independence or homogeneity in
two-way contingency tables.

\begin{example}[Independence in contingency tables]
  Contingency tables are among the simplest models of discrete data in
  widespread practical use \cite[Chapter 2]{agresti2019}. A \emph{two-way
    contingency table} studies the relationship between two discrete random
  variables $X_1$ and $X_2$. If each $X_i$ takes on $k_i$ distinct values, say
  $X_i \in \{1,\dots,k_i\}$, then a dataset of i.i.d.\ samples of pairs
  $(X_1,X_2)$ is summarized by a $k_1 \times k_2$ matrix, where the
  $(j_1,j_2)$-th entry is the number of observations with $X_1 = j_1$ and
  $X_2 = j_2$. In the common case where both $X_1$ and $X_2$ are binary valued,
  the contingency table is called a \emph{$2 \times 2$ table}. Higher order
  contingency tables, involving three or more discrete random variables, are
  defined similarly.

  The \emph{theory of a two-way contingency table with total count $n$} is
  generated by a convex space object $\theta$, two discrete objects $x_1$ and
  $x_2$, an additive monoid object $x$, a map $m_x: x_1 \otimes x_2 \to x$, and
  a morphism $p_0: \theta \to x_1 \otimes x_2$. The sampling morphism
  $p: \theta \to x$ is
  \begin{equation*}
    \input{wiring-diagrams/algebra-statistics/theory-contingency-full-sampling}.
  \end{equation*}
  Denote this theory by $(\cat{T},p)$. In an intended model $M$, each $M(x_i)$
  is the set $\sspace{X}_i := \{1,\dots,k_i\}$ for some number $k_i$;
  $M(\theta)$ is the convex set of probability vectors in
  $\R^{\sspace{X}_1 \times \sspace{X}_2}$, isomorphic to the standard
  $(k_1k_2 - 1)$-simplex in $\R^{k_1 k_2}$; $M(p_0)$ is the categorical family
  on $\sspace{X}_1 \times \sspace{X}_2$, given by
  \begin{equation*}
    M(p_0): \big( \pi_{j_1,j_2} \big)_{j_1,j_2=1}^{k_1,k_2} \mapsto
      \sum_{j_1,j_2=1}^{k_1,k_2} \pi_{j_1,j_2}\, \delta_{(j_1,j_2)};
  \end{equation*}
  $M(x)$ is the additive monoid $\N^{\sspace{X}_1 \times \sspace{X}_2}$,
  isomorphic to $\N^{k_1 k_2}$; and $M(m_x)$ is the map
  $(j_1, j_2) \mapsto e_{j_1,j_2}$ sending each element of
  $\sspace{X}_1 \times \sspace{X}_2$ to the corresponding basis element of
  $\N^{\sspace{X}_1 \times \sspace{X}_2}$. Consequently, the sampling
  distribution $M(p)$ is the full multinomial family, for $n$ observations, on
  $\sspace{X}_1 \times \sspace{X}_2$.\footnote{Alternatively, in the theory, one
    might let each $x_i$ be an additive monoid and, in the model $M$, take
    $M(x_i)$ to be $\N^{k_i}$ and $M(x)$ to be the tensor product
    $\N^{k_1} \otimes \N^{k_2} \cong \N^{k_1 k_2}$. Then $M(m_x)$ should be the
    biadditive map $\N^{k_1} \times \N^{k_2} \to \N^{k_1} \otimes N^{k_2}$ given
    by the universal property of the tensor product. An advantage of this
    approach is that the biadditivity of $m_x$ can be axiomatized within the
    theory.}

  The null hypothesis that the random variables $X_1$ and $X_2$ are
  \emph{independent} is classically stated as
  \begin{equation*}
    H_0: \pi_{j_1,j_2} = \pi_{j_1,+} \pi_{+,j_2},\ \forall j_1, j_2,
  \end{equation*}
  where $\pi_{j_1,+}$ and $\pi_{+,j_2}$ are the marginal distributions. As a
  statistical theory, let $\cat{T}_{\mathrm{ind}}$ be presented as $\cat{T}$,
  plus two convex space objects $\theta_1$ and $\theta$, two morphisms
  $p_{0,1}: \theta_1 \to x_1$ and $p_{0,2}: \theta_2 \to x_2$, and a
  convex-bilinear\footnote{That is, the map is convex-linear in each argument,
    with the other argument held fixed.} map
  $m_\theta: \theta_1 \otimes \theta_2 \to \theta$, subject to the equation
  \begin{equation*}
    \input{wiring-diagrams/algebra-statistics/theory-contingency-independence-lhs} =
    \input{wiring-diagrams/algebra-statistics/theory-contingency-independence-rhs}.
  \end{equation*}
  The sampling morphism $p_{\mathrm{ind}}: \theta_1 \otimes \theta_2 \to x$ is
  \begin{equation*}
    \input{wiring-diagrams/algebra-statistics/theory-contingency-independent-sampling}.
  \end{equation*}
  The intended model $M$ extends the previous one by taking each $M(\theta_i)$
  to be the probability simplex $\Delta^{k_i-1}$, each $M(p_{0,i})$ to be the
  categorical family of on $\sspace{X}_i$, and $M(\theta)$ to be the
  convex-bilinear map $(\pi_1, \pi_2) \mapsto \pi_1 \otimes \pi_2$. The asserted
  equation is indeed satisfied under this interpretation, so $M$ is a
  well-defined model.

  The null hypothesis of independence is represented by the colax theory
  morphism
  $(\iota, m_\theta, 1_x): (\cat{T}, p) \to (\cat{T}_{\mathrm{ind}},
  p_{\mathrm{ind}})$, where
  $\iota: \cat{T} \hookrightarrow \cat{T}_{\mathrm{ind}}$ is the evident
  inclusion functor. The colaxness conditions holds due to the asserted equation
  in $\cat{T}_{\mathrm{ind}}$ and the fact that $m_\theta$ is deterministic and
  hence distributes over copies.
\end{example}

\begin{example}[Homogeneity in contingency tables]
  Contingency tables arise from sampling schemes besides the full multinomial
  family on the joint distribution. When one discrete variable
  $X \in \{1,\dots,k\}$ is regarded as explanatory and another discrete variable
  $Y \in \{1,\dots,\ell\}$ as a response, it is common to assign $X$ a fixed
  value $i$ and then sample $Y$ conditionally on $X=i$. If $n_i$ samples are
  taken at each level $X=i$, then the resulting $k \times l$ contingency table
  will have fixed row totals $n_1, \dots, n_k$. This is a \emph{conditional
    multinomial} or, when $\ell = 2$, a \emph{conditional binomial}, sampling
  scheme. For example, in a randomized experiment, $X$ might be the assignment
  of a new drug or a placebo and $Y$ the response to treatment (success or
  failure), resulting in a $2 \times 2$ table.

  In the \emph{theory of a two-way contingency table with row counts
    $n_1,\dots,n_k$}, denoted $(\cat{T},p)$, the category $\cat{T}$ is freely
  generated by a discrete (or convex space) object $\theta$, an additive monoid
  object $y$, and a morphism $p_0: \theta \to y^{\otimes \ell}$. The sampling
  morphism
  $p: \theta^{\otimes k} \to (y^{\otimes \ell})^{\otimes k} \cong y^{\otimes
    k\ell}$ is
  \begin{equation*}
    \input{wiring-diagrams/algebra-statistics/theory-contingency-conditional-sampling}.
  \end{equation*}
  For each $1 \leq i \leq k$, this theory colaxly includes the theory of $n_i$
  i.i.d.\ $\ell$-dimensional counts (\cref{ex:models-iid-counts}). The intended
  model $M$ takes $M(\theta) = \Delta^{\ell-1}$, $M(x) = \N$, and $M(p_0)$ to be
  the categorical family $\CatDist(\ell,-): \Delta^{\ell-1} \to \N^\ell$. The
  sampling morphism is then the independent product of $k$ multinomial families.

  The null hypothesis of \emph{homogeneous} conditional distributions would
  traditionally be stated as
  \begin{equation*}
    H_0: \pi_{j|1} = \pi_{j|2} = \cdots = \pi_{j | k},\
    \forall j = 1, \dots, \ell,
  \end{equation*}
  where $\pi_{-|i}$ is the conditional distribution of $Y$ given $X=i$.
  Algebraically, homogeneity corresponds to the reduced sampling morphism
  $p_{\mathrm{hom}}$ given by
  $\theta \xrightarrow{\Delta_{\theta,k}} \theta^{\otimes k} \xrightarrow{p}
  y^{\otimes k \ell}$ or
  \begin{equation*}
    \input{wiring-diagrams/algebra-statistics/theory-contingency-homogeneous-sampling}.
  \end{equation*}
  As a lax morphism $(\cat{T}, p_{\mathrm{hom}}) \to (\cat{T}, p)$, or a colax
  morphism $(\cat{T}, p) \to (\cat{T}, p_{\mathrm{hom}})$, the null hypothesis
  is simply $(1_{\cat{T}}, \Copy_{\theta,k}, 1_{y^{\otimes k \ell}})$.
\end{example}

Contingency tables illustrate the important lesson that the presentation of data
in a particular format generally says little or nothing about the experimental
design, the sampling scheme, or what would be an appropriate statistical model.
Full multinomial and conditional multinomial sampling of two discrete random
variables both yield data in the form of a two-way contingency table, yet the
sampling schemes are very different, as reflected by their different statistical
theories and models. Nor are these the only possible sampling schemes; another
is \emph{Poisson sampling}, under which not even the total count of the table is
restricted. Likewise, independence and homogeneity are both hypotheses of ``no
association'' between variables, but they correspond to different statistical
theories and models. Statistical theories thus serve an important purpose in
making precise and explicit the background information that cannot be discerned
from a display of the data.

Discussion of null and alternative hypotheses notwithstanding, this work takes
no stance on the proper role of formal hypothesis testing in science. The
philosophy implicit in the algebraic approach, insofar as it has one, is that a
null hypothesis is just a morphism of statistical theories, not essentially
different than any other morphism. Thus, null hypotheses enjoy no special
logical status in the larger web of relationships between statistical theories.
But even if this is true, that does not imply anything about the role of formal
methods in theory and model selection generally.

The notion of a morphism between statistical theories extends straightforwardly
to Bayesian theories, defined at the end of the previous section
(\cref{def:bayesian-theory}).

\begin{definition}[Bayesian theory morphisms]
  A \emph{lax morphism} from one Bayesian theory
  $(\cat{T}, \theta \xrightarrow{p} x, I \xrightarrow{\pi} \theta)$ to another
  $(\cat{T}', p', \pi')$ is a lax morphism
  $(F, f_0, f_1): (\cat{T}, p) \to (\cat{T}', p')$ between the underlying
  statistical theories such that the diagram
  \begin{equation*}
    \begin{tikzcd}
      I \arrow[r, "F \pi"] \arrow[dr, "\pi'"']
        & F\theta \arrow[d, "f_0"] \\
        & \theta'
    \end{tikzcd}
  \end{equation*}
  commutes. Dually, a \emph{colax morphism} from $(\cat{T}, p, \pi)$ to
  $(\cat{T}', p', \pi')$ is a colax morphism
  $(F, f_0, f_1): (\cat{T}, p) \to (\cat{T}', p')$ such that the diagram
  \begin{equation*}
    \begin{tikzcd}
      I \arrow[r, "\pi'"] \arrow[dr, "F \pi"']
        & \theta' \arrow[d, "f_0"] \\
        & F\theta
    \end{tikzcd}
  \end{equation*}
  commutes. A lax or colax morphism of Bayesian theories is \emph{strong} or
  \emph{strict} if the underlying morphism of statistical theories is.
\end{definition}

In the category of Bayesian statistical theories and (co)lax morphisms,
composition and identities are those of the (co)lax morphisms between the
underlying statistical theories. Thus, by construction, there is a forgetful
functor from the category of Bayesian theories to the category of statistical
theories, which discards the prior. Another forgetful functor performs
marginalization, taking a Bayesian theory $(\cat{T}, p, \pi)$ to the statistical
theory $(\cat{T}, \pi \cdot p)$ and a Bayesian theory morphism
$(F,f_0,f_1): (\cat{T}, p, \pi) \to (\cat{T}', p', \pi')$ to the theory morphism
$(F,1_I,f_1): (\cat{T}, \pi \cdot p) \to (\cat{T}', \pi' \cdot p')$. The laxness
condition is verified by the commutative diagram
\begin{equation*}
  \begin{tikzcd}
    I \arrow[r, "F\pi"] \arrow[dr, "\pi'"']
      & F\theta \arrow[r, "Fp"] \arrow[d, "f_0"] & Fx \arrow[d, "f_1"] \\
      & \theta' \arrow[r, "p'"'] & x'
  \end{tikzcd}
  \qquad\leadsto\qquad
  \begin{tikzcd}
    I \arrow[r, "F(\pi \cdot p)"] \arrow[dr, "\pi' \cdot p'"']
      & Fx \arrow[d, "f_1"] \\
      & x'.
  \end{tikzcd}
\end{equation*}
The colaxness condition is dual.

\section{Notes and references}

\paragraph{Models in theoretical statistics}

Beginning with his 1939 paper \cite{wald1939} and culminating in his 1950 book
\cite{wald1950}, Abraham Wald introduced statistical decision theory as a
general framework for theoretical statistics, encompassing estimation and
hypothesis testing. Inspired by von Neumann’s game theory, the theory of
statistical decisions formalizes a statistical model as a parameterized family
$\{P_\theta\}_{\theta \in \Omega}$ of probability distributions on a sample
space $\sspace{X}$, as sketched at the beginning of the chapter; adds to this a
space $\sspace{A}$ of actions and a loss function
$L: \Omega \times \sspace{A} \to \R$, yielding a decision-theoretic problem; and
then defines criteria for a decision rule $d: \sspace{X} \to \sspace{A}$ to be
optimal or admissible with respect to the loss. The books by Lehmann et al
\cite{lehmann1998,lehmann2005} are now the standard references on statistical
decision theory, whereas the texts by Ferguson \cite{ferguson1967} and Berger
\cite{berger1985} are more introductory.

Although not included in the formalism of this chapter, the decision-theoretic
elements of theoretical statistics are also compositional in nature. Given a
sampling distribution $P: \Omega \to \sspace{X}$, a loss function
$L: \Omega \times \sspace{A} \to \R$, and a possibly randomized decision rule
$d: \sspace{X} \to \sspace{A}$, the composite Markov kernel
\begin{equation*}
  \input{wiring-diagrams/algebra-statistics/risk}
\end{equation*}
gives the distribution of the loss $L(\theta, d(X))$ under data
$X \sim P_\theta$ at every parameter $\theta \in \Omega$. Its expectation
\begin{equation*}
  R(\theta,d) := \E_\theta[L(\theta, d(X))]
\end{equation*}
is the \emph{risk}, the central quantity of statistical decision theory. The
statistical theories introduced here could conceivably be extended to include
morphisms for the loss function or even a preferred decision rule.

\paragraph{Statistical models as Markov kernels}

Markov kernels are a standard topic in advanced textbooks on probability theory,
such as by Kallenberg \cite{kallenberg2002} or Klenke \cite{klenke2013}. Another
book by Kallenberg is a comprehensive reference on the closely related topic of
\emph{random measures} \cite{kallenberg2017}.

The interpretation of a statistical model $\{P_\theta\}_{\theta \in \Omega}$ as
a Markov kernel $P: \Omega \to \sspace{X}$, and likewise for a possibly
randomized decision rule $d: \sspace{X} \to \sspace{A}$, is so natural that it
cannot easily be separated from the origin of statistical decision theory. The
first explicitly compositional (category-theoretic) study of statistical models
as Markov kernels was conducted by N.\ N.\ \v{C}encov, originally in Russian
\cite{cencov1965,cencov1972} and eventually translated into English
\cite{cencov1978,cencov1982}. Through this work \v{C}encov also made early
contributions to the differential-geometric study of statistical models, known
today as \emph{information geometry}. These two strands can be separated: our
algebraic study of statistics involves no differential geometry, and later work
on \v{C}encov's characterization of the Fisher information metric has eschewed
the language of category theory \cite{campbell1986,lebanon2004}.

In his master's thesis \cite{fong2012}, Fong develops an elegant algebraic
perspective on directed graphical models, also known as \emph{Bayesian
  networks}. Given a \emph{causal structure} in the form of a directed acyclic
graph, Fong builds a small symmetric monoidal category called a \emph{causal
  theory}. Functors out of this category into $\Markov$ are joint probability
distributions compatible with the causal structure \cite[Theorem 4.5]{fong2012}.
In relation to the present work, causal theories comprise a special class of
statistical theories that are freely generated by discrete objects
$v_1,\dots,v_n$ and morphisms representing the \emph{causal mechanisms}, and
that have sampling morphisms of form $p: I \to v_1 \otimes \cdots \otimes v_n$.
Thus, Fong gives a recipe for constructing a whole class of interesting
statistical theories. Since Bayesian networks often have unknown numerical
parameters that must be estimated from the data, it would be natural to extend
this formalism to sampling morphisms
$p: \theta \to v_1 \otimes \cdots \otimes v_n$ having nontrivial parameter space
objects $\theta$.

Remarkably, Fritz has recently demonstrated that sufficiency, ancillarity,
completeness, and minimal sufficiency may be defined, and versions of the
Neyman-Fisher factorization theorem, Basu's theorem, and Bahadur's theorem
proved, in the purely synthetic setting of a Markov category \cite{fritz2020}.
All of these belong to the classic definitions and abstract results of
statistical decision theory.

\paragraph{Markov kernels in categorical probability}

The earliest category-theoretic study of Markov kernels, outside of statistics
and independently from \v{C}encov, was made by William Lawvere in an unpublished
appendix to a 1962 grant proposal. Later, Giry extended and published this work
\cite{giry1982}, defining what is now called the \emph{Giry monad} on the
category of measurable spaces and maps. The Kleisli category of this monad is
the category of Markov kernels (without topological restrictions). Although
monads do not figure explicitly in this thesis, the Giry monad and other
\emph{probability monads} are now among the best-studied aspects of categorical
probability \cite{perrone2018}. Theoretical computer scientists have used the
category of Markov kernels to reason about probabilistic systems and programs,
with early works by Blute, Desharnais, Edalat, and Panangaden
\cite{blute1997,panangaden1999}.

As a synthetic setting for probability and statistics, Markov categories have
been studied, under various names, by Fong, Cho and Jacobs, Fritz and others
\cite{fong2012,cho2019,fritz2020}. Further references can be found in Fritz's
survey of previous work \cite{fritz2020}. The idea of characterizing the
``functions'' or ``maps'' in a symmetric monoidal category as the morphisms that
preserve copying and deleting goes back at least as far as the
\emph{bicategories of relations} and \emph{bicategories of partial maps} of
Carboni and Walters \cite{carboni1987a,carboni1987b}. Although probability and
linear algebra have each been studied separately from the categorical viewpoint,
their conjunction in the form of a linear-algebraic Markov category is original.
This seems to be the natural synthetic setting for a large part of everyday
statistical modeling, as demonstrated in \cref{ch:zoo-statistics}.

The characterizations in this chapter of the deterministic morphisms
(\cref{cor:markov-kernel-determinism}) and the isomorphisms
(\cref{prop:markov-kernel-isomorphims}) in a category of well-behaved Markov
kernels are well known \cite{belavkin2013,fong2012}.

\paragraph{McCullagh on statistical models}

McCullagh's paper on ``What is a statistical model?'' is a rare example of a
category-theoretic treatment of statistical models by a professional
statistician \cite{mccullagh2002}. The paper aims to formalize the prospect,
regarded as necessary for meaningful statistical inference, of extending the
parameter and sample spaces of a statistical model to include additional
observational units, such as new subjects or future points in time. In
comparison with this work, McCullagh's central commutative diagram
\begin{equation*}
  \begin{tikzcd}
    \mathcal{P}(\mathcal{S})
      & \Theta_{\Omega} \arrow[l, "P_\psi"'] \\
    \mathcal{P}(\mathcal{S}') \arrow[u, "\varphi_d^\dagger"]
      & \Theta_{\Omega'} \arrow[l, "P_{\psi'}"'] \arrow[u, "\varphi_c^*"'],
  \end{tikzcd}
\end{equation*}
reproduced from \cite[Equation 1]{mccullagh2002}, bears a strong formal
similarity to the commutative diagram
\begin{equation*}
  \begin{tikzcd}
    \theta' \arrow[r, "p'"] \arrow[d, "f_0"'] & x' \arrow[d, "f_1"] \\
    F\theta \arrow[r, "Fp"] & Fx
  \end{tikzcd}
  \qquad\text{or equivalently}\qquad
  \begin{tikzcd}
    Fx & F\theta \arrow[l, "Fp"'] \\
    x' \arrow[u, "f_1"] & \theta' \arrow[l, "p'"'] \arrow[u, "f_0"']
  \end{tikzcd}
\end{equation*}
obeyed by a colax morphism $(F,f_0,f_1): (\cat{T},p) \to (\cat{T}',p')$ of
statistical theories (\cref{def:lax-theory-morphism}). In other respects the two
mathematical formalisms are very different. In particular, the notion of
statistical theory is not present in McCullagh's work, and thus no separation is
made between the algebraic and analytical aspects of a statistical model. We
suspect that many of the coherency conditions between statistical models
proposed by McCullagh could be formulated as families of statistical theories
connected by colax theory morphisms.
\cref{ex:morphism-iid-sample-sizes,ex:morphism-exchangeable-sample-sizes} can be
seen as simple examples of such conditions, but a full account is beyond the
scope of this work.

\chapter{A zoo of statistical theories and models}
\chaptermark{Zoo of statistical models}
\label{ch:zoo-statistics}

Statistical theories, models, and their morphisms have so far been illustrated
mainly through simple models of discrete data, such as the binomial and
multinomial models and models for contingency tables. With the basic formalism
established, this chapter builds on the previous one by presenting a wider range
of statistical theories and models, primarily for regression and classification.
The statistical methods progress in complexity from linear models for
categorical or continuous predictors; to Bayesian, hierarchical, and mixed
linear models; and finally to generalized linear models. We do not aspire to
encyclopedic coverage, much less to a complete treatment, of the great variety
of statistical models devised by statisticians, computer scientists, and domain
scientists over hundreds of years. Instead, we aim to show how the algebraic
view of statistical models usefully formalizes and brings out the relations
between some of the most essential models in statistics.

The linear model is the point of departure for a large part of statistics and
machine learning, and so it is for this chapter. In its standard formulation, a
normal \emph{linear model} is any statistical model of form
\begin{equation*}
  y \sim \Normal_n(X \beta, \sigma^2 I_n),
\end{equation*}
where $X \in \R^{n \times p}$ is a fixed but arbitrary design matrix,
$y \in \R^n$ is the response vector, $\beta \in \R^p$ is an unknown vector of
coefficients, and $\sigma^2 \geq 0$ is an unknown variance.\footnote{The linear
  model admits a mild generalization to multivariate data, where the vectors
  $y \in \R^n$ and $\beta \in \R^p$ are replaced by matrices
  $Y \in \R^{n \times q}$ and $B \in \R^{p \times q}$.} However, depending on
the form of the design matrix, the model will be analyzed and interpreted in
different ways, say as an analysis of variance or as a linear regression. So,
although it admits a unified mathematical treatment, the linear model should be
understood not as single statistical method but as a family of closely related
statistical methods. Each method has its own statistical theory, in the sense of
\cref{ch:algebra-statistics}, and the relations between these theories are
formalized by theory morphisms.

The family of normal distributions plays a central role in linear modeling. In
the previous chapter, the isotropic normal family has already been presented as
a statistical theory, up to an absolute scale
(\cref{cor:isotropic-normal-presentation}). From the algebraic perspective,
however, it is natural to ask only for a family of distributions that is
additive and homogeneous with exponents one and two
(\cref{thm:stable-family-presentation}). Such families will be called
\emph{linear-quadratic}. Any linear quadratic family of probability
distributions is a location-scale normal family, up to linear transformations of
the location and scale parameters.

\section{Linear models with discrete predictors}
\label{sec:lm-discrete}

Linear models with discrete predictors encompass such statistical methods as the
one-sample test, the two-sample test, analysis of variance (ANOVA), and multiway
ANOVA. The first of these models, the one-sample model, can seen as a degenerate
case of a linear model with a single discrete predictor taking a single value.

\paragraph{One-sample normal model}

The simplest of all linear models is the univariate, one-sample normal model
\begin{equation*}
  y_i \xsim{\mathrm{iid}} \Normal(\mu, \sigma^2), \quad i=1,\dots,n,
\end{equation*}
with parameters $\mu \in \R$ and $\sigma^2 \geq 0$. In the statistical
\emph{theory of one normal sample (of size $n$)}, denoted $(\cat{S}, p_n)$, the
category $\cat{S}$ is presented by vector space objects $\mu$ and $y$, a conical
space object $\sigma^2$, and a linear-quadratic morphism
$q: \mu \otimes \sigma^2 \to y$, meaning that
\begin{equation*}
  \input{wiring-diagrams/zoo-statistics/normal-linear-lhs} =
  \input{wiring-diagrams/zoo-statistics/normal-linear-rhs}
  \qquad \forall \vec{a} \in \R^2,
\end{equation*}
where $\vec{a} := (a_1, a_2)$ is a linear combination and
$\vec{a}^2 := (a_1^2, a_2^2)$ is a conic combination. Also, note that
``$\sigma^2$'' is only a suggestive symbol; the object $\sigma^2$ is not squared
in any sense. The theory's sampling morphism
$p_n: \mu \otimes \sigma^2 \to y^{\otimes n}$ is the morphism
$\Copy_{\mu \otimes \sigma^2} \cdot q^{\otimes n}$ or
\begin{equation*}
  \input{wiring-diagrams/zoo-statistics/one-sample-normal-sampling}.
\end{equation*}
The intended univariate model $M$ takes $M(\mu)$ and $M(y)$ to be the real
numbers, $M(\sigma^2)$ to be the nonnegative real numbers, and $M(q)$ to be the
normal family $\Normal: \R \times \R_+ \to \R$.

The symmetries of this model are simply the dilations (changes of scale), along
with reflection across the origin (multiplication by $-1$).

\begin{proposition} \label{prop:univariate-normal-morphisms}
  The endomorphisms $\alpha: M \to M$ of the univariate normal model are
  isomorphic, as a monoid, to the multiplicative monoid of real numbers $\R$.

  In particular, the automorphism group of $M$ is isomorphic to
  $\R^* := \{a \in \R: a \neq 0\}$.
\end{proposition}
\begin{proof}
  By \cref{lemma:natural-transformation-generators}, a morphism
  $\alpha: M \to M$ consists of linear maps $\alpha_\mu, \alpha_y: \R \to \R$
  (scalars) and a conic-linear map $\alpha_{\sigma^2}: \R_+ \to \R$ (a
  nonnegative scalar) making the diagram
  \begin{equation*}
    \begin{tikzcd}
      \R \times \R_+ \arrow[r, "\Normal"]
        \arrow[d, "\alpha_\mu \times \alpha_{\sigma^2}"']
        & \R \arrow[d, "\alpha_y"] \\
      \R \times \R_+ \arrow[r, "\Normal"]
        & \R
    \end{tikzcd}
  \end{equation*}
  commute. That is, for every $\mu \in \R$ and $\sigma^2 \geq 0$,
  \begin{equation*}
    \Normal(\alpha_\mu\, \mu, \alpha_{\sigma^2}\, \sigma^2) =
    \alpha_y\, \Normal(\mu, \sigma^2),
  \end{equation*}
  which happens if and only if $\alpha_\mu = \alpha_y = a$ and
  $\alpha_{\sigma^2} = a^2$ for some scalar $a \in \R$. Moreover, if morphisms
  $\alpha,\beta: M \to M$ correspond to scalars $a,b \in \R$, then their
  composite $\alpha \cdot \beta$ corresponds to $ab$, and the identity morphism
  $1_M$ corresponds to the scalar $1$.
\end{proof}

The theory of one normal sample also has multivariate models. For any dimension
$d$, define the model $M_d$ that maps the objects $\mu$ and $y$ to $\R^d$ and
$\sigma^2$ to $\PSD^d$, and the morphism $q$ to the $d$-dimensional normal
family $\Normal_d$. Then the sampling distribution under $M_d$ can be written as
\begin{equation*}
  y_i \xsim{\mathrm{iid}} \Normal_d(\mu, \Sigma), \quad i=1,\dots,n,
\end{equation*}
with parameters $\mu \in \R^d$ and $\Sigma \in \PSD^d$. Generalizing the
univariate case, the symmetries of the multivariate model are arbitrary
invertible linear transformations.

\begin{proposition} \label{prop:multivariate-normal-morphisms}
  The category of multivariate normal models and model homomorphisms is
  isomorphic to the category $\Mat_\R$ of real matrices.

  In particular, the automorphism group of the $d$-dimensional model $M_d$ is
  isomorphic to the general linear group $\GL(d,\R)$.
\end{proposition}
\begin{proof}
  Arguing exactly as before, a morphism $\alpha: M_d \to M_{d'}$ consists of
  linear maps $\alpha_\mu, \alpha_y: \R^d \to \R^{d'}$ and a conic-linear map
  $\alpha_{\sigma^2}: \PSD^d \to \PSD^{d'}$ making the diagram
  \begin{equation*}
    \begin{tikzcd}
      \R^d \times \PSD^d \arrow[r, "\Normal_d"]
        \arrow[d, "\alpha_\mu \times \alpha_{\sigma^2}"']
        & \R^d \arrow[d, "\alpha_y"] \\
      \R^{d'} \times \PSD^{d'} \arrow[r, "\Normal_{d'}"]
        & \R^{d'}
    \end{tikzcd}
  \end{equation*}
  commute. That is, for every $\mu \in \R^d$ and $\Sigma \in \PSD^d$,
  \begin{equation*}
    \Normal_{d'}(\alpha_\mu\, \mu, \alpha_{\sigma^2} \Sigma) =
    \alpha_y\, \Normal_d(\mu, \Sigma),
  \end{equation*}
  which happens if and only if there exists a matrix $A \in \R^{d' \times d}$
  such that $\alpha_\mu = \alpha_y = A$ and $\alpha_{\sigma^2}$ is the map
  $\Sigma \mapsto A \Sigma A^\top$. The correspondence between $\alpha$ and $A$
  is clearly functorial.
\end{proof}

Another model of the theory, the isotropic multivariate normal model, will be
given later, but for now let us return to the theory itself. The theory
$(\cat{S}, p_n)$ of one normal sample of size $n$ is a refinement of the theory
$(\cat{T}, p_n)$ of $n$ i.i.d.\ samples, defined in
\cref{ex:theory-iid-samples}. More precisely, the functor
$F: \cat{T} \to \cat{S}$ sending the objects $\theta$ and $x$ to
$\mu \otimes \sigma^2$ and $y$, respectively, and the morphism $p_0$ to $q$
defines a strict morphism of theories $F: (\cat{T}, p_n) \to (\cat{S}, p_n)$.
Also, the null hypothesis
\begin{equation*}
  H_0: \mu = 0
\end{equation*}
of zero mean, as might be tested by a one-sample $t$-test, is represented by the
colax theory morphism $(1_{\cat{S}}, 0_\mu \otimes 1_{\sigma^2}, 1_y)$ from the
one sample theory $(\cat{S}, p_n)$ to a reduced theory $(\cat{S}, p_{n,0})$,
whose sampling morphism $p_{n,0} :=(0_\mu \otimes 1_{\sigma^2}) \cdot p_n$ is
\begin{equation*}
  \input{wiring-diagrams/zoo-statistics/one-sample-normal-null-lhs} =
  \input{wiring-diagrams/zoo-statistics/one-sample-normal-null-rhs}.
\end{equation*}

Regarding the normal distribution as a location-scale family, one might expect
that the morphisms of the normal model would include location transformations.
These are absent because model homomorphisms, being natural transformations,
must preserve all structure of the theory, including the linear structure of
vector space objects. However, there is a slightly weaker theory that takes the
mean and response objects to be affine spaces. Let $\cat{S}^{\aff}$ be presented
by affine space objects $\mu$ and $y$, a conical space object $\sigma^2$, and an
\emph{affine-quadratic} morphism $q: \mu \otimes \sigma^2 \to y$, meaning that
\begin{equation*}
  \input{wiring-diagrams/zoo-statistics/normal-affine-lhs} =
  \input{wiring-diagrams/zoo-statistics/normal-affine-rhs},
\end{equation*}
where the affine combination $\vec{t} := (t,1-t)$ and conical combination
$\vec{t}^2 := (t^2, (1-t)^2)$ range over all $t \in \R$. In the statistical
theory $(\cat{S}^{\aff}, p_n)$, the sampling morphism $p_n$ is defined exactly
as before.

Having defined this new theory, the inclusion functor
$\aff: \cat{S}^{\aff} \hookrightarrow \cat{S}$ is a strict theory morphism
$(\cat{S}^{\aff}, p_n) \to (\cat{S}, p_n)$ and the model migration functor
$\aff^*: \Model{\cat{S}} \to \Model{\cat{S}^{\aff}}$ interprets the multivariate
models $M_d$ of $\cat{S}$ as models of $\cat{S}^{\aff}$. Arguing as in
\cref{prop:univariate-normal-morphisms,prop:multivariate-normal-morphisms}, a
model homomorphism $\alpha: \aff^*(M_d) \to \aff^*(M_{d'})$ is seen to consist
of, for any matrix $A \in \R^{d \times d'}$ and vector $b \in \R^{d'}$, the
affine transformation $\alpha_\mu = \alpha_y: y \mapsto Ay + b$ and the
conic-linear transformation $\alpha_{\sigma^2}: \Sigma \mapsto A \Sigma A^\top$.
In particular, the automorphism group of the model $\aff^*(M_d)$ is isomorphic
to the general affine group $\AffGroup(d, \R)$. Besides models inherited from
the linear theory $(\cat{S},p_n)$, the affine theory $(\cat{S}^{\aff}, p_n)$
also has normal models whose mean and sample spaces are affine subspaces of
$\R^d$ not containing the origin. These are not models of the original theory.

In the affine theory $(\cat{S}^{\aff}, p_n)$, the zero morphism
$0_\mu: I \to \mu$ does not exist and the null hypothesis $H_0: \mu = 0$ of zero
mean cannot be stated. Indeed, the origin is not preserved by affine
transformations, so the hypothesis of zero mean would not be preserved under
model homomorphism. To introduce an affine point null hypothesis
$H_0: \mu = \mu_0$, present $\cat{S}^{\aff}_0$ as $\cat{S}^{\aff}$ together with
an (affine) map $\mu_0: I \to \mu$ and form the statistical theory
$(\cat{S}^{\aff}_0, p_{n,0})$ with sampling morphism $p_{n,0}$ equal to
$(\mu_0 \otimes 1_{\sigma^2}) \cdot p_n$ or
\begin{equation*}
  \input{wiring-diagrams/zoo-statistics/one-sample-normal-general-null}.
\end{equation*}
The relationships between the theories and hypotheses for the one-sample normal
model are summarized by the commutative diagram
\begin{equation*}
  \begin{tikzcd}[column sep=huge]
    (\cat{S}^{\aff}, p_n)
      \arrow[r, "{(\iota_0, \mu_0 \otimes 1_{\sigma^2}, 1_x)}"]
      \arrow[d, "\aff"']
      & (\cat{S}^{\aff}_0, p_{n,0}) \arrow[d, "\aff_0"] \\
    (\cat{S}, p_n)
      \arrow[r, "{(1_{\cat{S}}, 0_\mu \otimes 1_{\sigma^2}, 1_x)}"]
      & (\cat{S}, p_{n,0}),
  \end{tikzcd}
\end{equation*}
where the functor $\aff_0: \cat{S}^{\aff}_0 \to \cat{S}$ extends the inclusion
functor $\aff: \cat{S}^{\aff} \hookrightarrow \cat{S}$ by assigning
$\mu_0 \mapsto 0_\mu$. Most of the statistical theories presented below using
vector space objects also have affine analogues, but in the interest of brevity
usually only the linear setting is treated.

\paragraph{Two-sample normal model}

After the one-sample model, the next simplest normal model is the univariate,
two-sample model, consisting of two independent samples
\begin{equation*}
  y_{i,j} \xsim{\mathrm{ind}} \Normal(\mu_i, \sigma_i^2), \quad
  i=1,2, \quad j=1,\dots,n_i
\end{equation*}
of sizes $n_1$ and $n_2$, with parameters $\mu_1, \mu_2 \in \R$ and
$\sigma_1^2 ,\sigma_2^2 \geq 0$. As a statistical theory, the \emph{theory of
  two normal samples (of sizes $n_1$ and $n_2$)} has the same underlying
category $\cat{S}$ as the theory of one normal sample. Abbreviating
$\vec{n} := (n_1,n_2)$ and $n := n_1 + n_2$, its sampling morphism
$p_{\vec{n}}: \mu^{\otimes 2} \otimes (\sigma^2)^{\otimes 2} \to y^{\otimes n}$
is defined by
\begin{equation*}
  \input{wiring-diagrams/zoo-statistics/two-sample-normal-sampling}.
\end{equation*}
Choosing this sampling morphism instead of the product
$p_{n_1} \otimes p_{n_2}: (\mu \otimes \sigma^2)^{\otimes 2} \to y^{\otimes n}$
is purely conventional, as the two statistical theories are isomorphic via the
strong theory isomorphism
\begin{equation*}
  (1_{\cat{S}}, 1_\mu \otimes \Braid_{\mu,\sigma^2} \otimes 1_{\sigma^2},
   1_{y^{\otimes n}}):
  (\cat{S}, p_{\vec{n}}) \cong (\cat{S}, p_{n_1}\otimes p_{n_2}).
\end{equation*}
The theory of two normal samples contains two separate copies of the theory of
one normal sample, via colax theory morphisms
\begin{equation*}
  (\cat{S}, p_{n_i}) \xrightarrow{(1_{\cat{S}}, \pi_i^0, \pi_i^1)}
  (\cat{S}, p_{n_1} \otimes p_{n_2}) \cong (\cat{S}, p_{\vec{n}}),
  \quad i=1,2,
\end{equation*}
whose components
$(\mu \otimes \sigma^2)^{\otimes 2} \xrightarrow{\pi_i^0} \mu \otimes \sigma^2$
and $y^{\otimes n} \xrightarrow{\pi_i^1} y^{\otimes n_i}$ are the evident
projections on the parameter and sample space objects. Because it has the same
underlying category, the two-sample theory has the same models and the same
model homomorphisms as the one-sample theory, including the univariate and
multivariate models defined above.

The two-sample normal model is, not strictly speaking, a linear model, due to
the unequal variances of the two groups.\footnote{If the ratio of the two
  variances is known, the model is a \emph{weighted} linear model, often fit by
  the method of weighted least squares \cite[\S 3.10]{seber2003}. A statistical
  theory for this model introduces two further morphism generators, conic-linear
  maps $v_1, v_2: \sigma^2 \to \sigma^2$, and a suitable sampling morphism.} The
\emph{theory of two homoscedastic normal samples} has the same underlying
category $\cat{S}$ and the sampling morphism
$q_{\vec{n}}: \mu^{\otimes 2} \otimes \sigma^2 \to y^{\otimes n}$ given by
\begin{equation*}
  \input{wiring-diagrams/zoo-statistics/two-sample-normal-equal-var}.
\end{equation*}
The assumption of equal variances, $\sigma_1^2 = \sigma_2^2$, is represented by
the colax theory morphism
$(1_{\cat{S}}, 1_{\mu^{\otimes 2}} \otimes \Copy_{\sigma^2}, 1_{y^{\otimes n}}):
(\cat{S}, p_{\vec{n}}) \to (\cat{S}, q_{\vec{n}}).$ The further hypothesis
\begin{equation*}
  H_0: \mu_1 = \mu_2
\end{equation*}
of equal means, as might be tested by a two-sample $t$-test, is represented by
the colax morphism
$(1_{\cat{S}}, \Copy_\mu \otimes 1_{\sigma^2}, 1_{y^{\otimes n}}): (\cat{S},
q_{\vec{n}}) \to (\cat{S}, p_n)$, where the reduced sampling morphism
$p_n: \mu \otimes \sigma^2 \to y^{\otimes n}$ is
\begin{equation*}
  \input{wiring-diagrams/zoo-statistics/two-sample-normal-equal-mean-var-lhs} =
  \input{wiring-diagrams/zoo-statistics/two-sample-normal-equal-mean-var-rhs}.
\end{equation*}
Thus, under the null hypothesis of equal means, we recover the sampling morphism
of the theory of one normal sample of size $n = n_1 + n_2$.

The theories of one or two normal samples generalize easily to the theory of $k$
normal samples, for any number $k$. Models of this theory are linear models with
a single discrete predictor taking $k$ distinct values. The global null
hypothesis $H_0: \mu_1 = \cdots = \mu_k$ of all means being equal is classically
tested by an \emph{analysis of variance}, or \emph{ANOVA}. Extensions of this
model to two or more discrete predictors are known as \emph{two-way} or
\emph{multiway} ANOVAs. A large family of statistical theories ensues.

\paragraph{Normal means model}

The normal means model, also known as the normal sequence model, is widely
studied in theoretical statistics and signal processing as a simple model that
already exhibits important and generic features of high-dimensional estimation
\cite{johnstone2019}. Classically, the model is written as
\begin{equation*}
  y_i \xsim{\mathrm{ind}} \Normal(\mu_i, \sigma^2), \quad i = 1,\dots,n,
\end{equation*}
or, equivalently, as
\begin{equation*}
  y \sim \Normal_n(\mu, \sigma^2 I_n).
\end{equation*}
Although the two expressions are mathematically equivalent, they suggest
different interpretations, leading to different statistical theories.

In the first case, we think of observing $n$ independent normal variables with
equal variance but possibly different means, say from a signal measured at
discrete times indexed by $i=1,\dots,n$. The \emph{theory of a normal sequence
  of length $n$} has the underlying category $\cat{S}$ and the sampling morphism
$q_n: \mu^{\otimes n} \otimes \sigma^2 \to y^{\otimes n}$ given by
\begin{equation*}
  \input{wiring-diagrams/zoo-statistics/normal-means-sampling}.
\end{equation*}
The intended model $M: \cat{S} \to \Stat$ is the univariate normal model defined
earlier, assigning $\mu$ and $y$ to the real numbers, $\sigma^2$ to the
nonnegative real numbers, and $q: \mu \otimes \sigma^2 \to y$ to the normal
family $\Normal: \R \times \R_+ \to \R$. The sampling distribution
$M(q_n): \R^n \times \R_+ \to \R^n$ is then the isotropic $n$-dimensional normal
family.

Alternatively, consider a single observation of an isotropic normal vector in
$n$-dimensional space. The corresponding statistical theory $(\cat{S},q)$ is
that of one normal sample of size 1, and the intended model
$M_n^{\mathrm{iso}}: \cat{S} \to \Stat$ assigns $\mu$ and $y$ to $\R^n$,
$\sigma^2$ to $\R_+$, and $q: \mu \otimes \sigma^2 \to y$ to the isotropic
$n$-dimensional normal family.

The sampling distributions $M(q_n)$ and $M_n^{\mathrm{iso}}(q)$ may be the same,
but the statistical theories $(\cat{S}, q_n)$ and $(\cat{S}, q)$, the
statistical models $M$ and $M_n^{\mathrm{iso}}$, and the model homomorphisms are
all importantly distinct. In \cref{prop:univariate-normal-morphisms}, the
symmetries of $M$ were seen to be rescalings and reflection, so that
$\Aut(M) \cong \R^*$. The second model has a much larger group of symmetries.

\begin{proposition} \label{prop:isotropic-normal-morphisms}
  The category of isotropic multivariate normal models and model homomorphisms
  is isomorphic to the category of real matrices $A$ such that $A A^\top$ is
  proportional to the identity.

  In particular, the automomorphism group of the isotropic $d$-dimensional model
  $M_d^{\mathrm{iso}}$ is isomorphic to the conformal orthogonal group $\CO(d)$,
  the direct product of the orthogonal group $\Orth(d)$ with the group of
  dilations $\R^*_+ := \{a \in \R: a > 0\}$.
\end{proposition}
\begin{proof}
  Arguing as in \cref{prop:multivariate-normal-morphisms}, a morphism
  $\alpha: M_d^{\mathrm{iso}} \to M_{d'}^{\mathrm{iso}}$ consists of linear maps
  $\alpha_\mu, \alpha_y: \R^d \to \R^{d'}$ and a conic-linear map
  $\alpha_{\sigma^2}: \R_+ \to \R_+$ (a nonnegative scalar) such that
  \begin{equation*}
    \Normal_{d'}(\alpha_\mu\, \mu, \alpha_{\sigma^2} \sigma^2 I_{d'}) =
    \alpha_y\, \Normal_d(\mu, \sigma^2 I_d)
  \end{equation*}
  for every $\mu \in \R^d$ and $\sigma^2 \geq 0$. This happens if and only if
  there exists a matrix $A \in \R^{d' \times d}$ such that
  $\alpha_\mu = \alpha_y = A$ and $A A^\top = \alpha_{\sigma^2} I_{d'}$. Thus,
  $A$ has the form $A = aU$, where $\alpha_{\sigma^2} = a^2$ and $Q$ is a
  semi-orthogonal matrix satisfying $Q Q^\top = I_{d'}$.
\end{proof}

The normal means model highlights the principle, implicit throughout the whole
development, that choosing a statistical theory amounts to deciding what
structure is essential to the problem at hand and that this choice determines
what are the allowed models and model homomorphisms. If the components
$y_1,\dots,y_n$ are regarded as distinct observational units, say observations
at particular times or of particular subjects, then one ought to choose the
theory of a normal sequence $(\cat{S}, q_n)$. The observational units are
explicitly recorded in the sampling morphism
$q_n: \mu^{\otimes n} \otimes \sigma^2 \to y^{\otimes n}$ and the model
homomorphisms are not permitted to mix components, which would destroy the
interpretation of the data. On the other hand, if the components $y_1,\dots,y_n$
are expressed in an arbitrary coordinate system, say from a measurement in
three-dimensional space ($n=3$) with fixed origin, then the alternate theory
$(\cat{S}, q)$ may be preferable, as the model homomorphisms allow for arbitrary
rotations of the coordinate system.

Due to the existence of theory morphisms, choosing one statistical theory does
not preclude consideration of other theories. The supply preserving functor
$F_n: \cat{S} \to \cat{S}$ mapping the objects $\mu$, $\sigma^2$, and $y$ to
$\mu^{\otimes n}$, $\sigma^2$, and $y^{\otimes n}$ and the morphism $q$ to $q_n$
is a strict morphism of theories $F_n: (\cat{S}, q) \to (\cat{S}, q_n)$. Thus,
the theory of a normal sequence $(\cat{S}, q_n)$ can be seen as a specialization
of the theory $(\cat{S}, q)$. The model migration functor
$F_n^*: \Model{\cat{S}} \to \Model{\cat{S}}$ recovers the $n$-dimensional
isotropic model $F_n^*(M) = M_n^{\mathrm{iso}}$ from the univariate model $M$.

\section{Linear models with general design}
\label{sec:lm-general}

The $k$-sample normal model and normal means model are both special cases of the
linear regression model, obtained by suitable choices of the design matrix. In
its general form, the linear regression model has several formulations as a
statistical theory, differing according to whether one, both, or neither of the
dimensions of the design matrix are made explicit. Following tradition, we
denote the number of observations by $n$ and the number of predictors by $p$, so
that the design matrix is an $n \times p$ real matrix.

\paragraph{Linear models}

The weakest of the theories, the \emph{theory of a linear model}, has underlying
category $\cat{LM}$ presented by vector space objects $\beta$, $\mu$, and $y$, a
conical space object $\sigma^2$, a linear map $X: \beta \to \mu$, and a
linear-quadratic morphism $q: \mu \otimes \sigma^2 \to y$. The sampling morphism
$p: \beta \otimes \sigma^2 \to y$ is
\begin{equation*}
  \input{wiring-diagrams/zoo-statistics/lm-sampling}.
\end{equation*}
The intended models $M: \cat{LM} \to \Stat$ take $M(\beta)$ to be $\R^p$, for
some dimension $p$; both $M(\mu)$ and $M(y)$ to be $\R^n$, for some dimension
$n$; $M(\sigma^2)$ to be the nonnegative reals; $M(X)$ to be any matrix
$X_M \in \R^{n \times p}$; and $M(q)$ to be the $n$-dimensional isotropic normal
family. The sampling distribution $M(p): \R^p \times \R_+ \to \R^n$ is then
\begin{equation*}
  y \sim \Normal_n(X_M \beta, \sigma^2 I_n)
\end{equation*}
with parameters $\beta \in \R^p$ and $\sigma^2 \geq 0$. Unlike in the previous
section, the theory does not have a single preferred model, but a whole of
family of models with different design matrices.

\begin{proposition} \label{prop:lm-morphisms}
  A morphism $\alpha: M \to M'$ between linear models $M$ and $M'$ is uniquely
  determined by matrices $A \in \R^{n' \times n}$ and $B \in \R^{p' \times p}$
  such that $A A^\top$ is proportional to the identity matrix and $A$ and $B$
  intertwine the design matrices:
  \begin{equation*}
    A X_M = X_{M'} B.
  \end{equation*}
  In particular, an isomorphism $\alpha: M \cong M'$ can exist only if $n=n'$
  and $p=p'$ and is then uniquely determined by matrices $A \in \CO(n)$ and
  $B \in \GL(p,\R)$ exhibiting the design matrices as equivalent:
  \begin{equation*}
    X_{M'} = A X_M B^{-1}.
  \end{equation*}
\end{proposition}
\begin{proof}
  A morphism $\alpha: M \to M'$ consists of linear maps
  $\alpha_\beta: \R^p \to \R^{p'}$ and $\alpha_\mu, \alpha_y: \R^n \to \R^{n'}$
  and a conic-linear map $\alpha_{\sigma^2}: \R_+ \to \R_+$ making the diagrams
  \begin{equation*}
    \begin{tikzcd}
      \R^p \arrow[r, "X_M"] \arrow[d, "\alpha_\beta"']
        & \R^n \arrow[d, "\alpha_\mu"] \\
      \R^{p'} \arrow[r, "X_{M'}"]
        & \R^{n'}
    \end{tikzcd}
    \qquad\qquad
    \begin{tikzcd}
      \R^n \times \R_+ \arrow[r, "\Normal_n^{\mathrm{iso}}"]
        \arrow[d, "\alpha_\mu \times \alpha_{\sigma^2}"']
        & \R^n \arrow[d, "\alpha_y"] \\
      \R^{n'} \times \R_+ \arrow[r, "\Normal_{n'}^{\mathrm{iso}}"]
        & \R^{n'}
    \end{tikzcd}
  \end{equation*}
  commute. From the proof of \cref{prop:isotropic-normal-morphisms}, we know
  that the second diagram amounts to having $\alpha_\mu = \alpha_y = A$ for some
  matrix $A \in \R^{n' \times n}$ such that
  $A A^\top = \alpha_{\sigma^2} I_{n'}$. Defining the matrix
  $B := \alpha_\beta \in \R^{p' \times p}$, the first diagram becomes the
  equation $A X_M = X_{M'} B$. The morphism $\alpha: M \to M'$ is invertible if
  and only if its components are, which happens if and only if $A$ and $B$ are
  invertible.
\end{proof}

So, isomorphic linear models have equivalent design matrices, or design matrices
of the same rank. But the change of basis matrix $A$ must be a conformal
orthogonal matrix, adding many further constraints. It is therefore \emph{not}
the case that any two linear models of equal size and equal rank are isomorphic.
Even a non-invertible morphism from $M$ to $M'$ is significantly constrained.
Except in the degenerate case that $A = 0$, the condition that $A A^\top$ is
proportional to the identity $I_{n'}$ implies that $n' \leq n$. A general
morphism of linear models can thus be interpreted as a kind of projection from a
larger model onto a smaller one. In particular, $X_{M'}$ may be obtained by
selecting a subset of rows from $X_M$, so that $A$ is the corresponding
selection matrix in $\{0,1\}^{n' \times n}$ and $B$ is the identity matrix
$I_p$.

The affine version of the linear model has a larger group of symmetries
encompassing location transformations. Let the \emph{theory of an affine
  model},\footnote{The term ``affine model'' is not standard in statistics, as
  the linear and affine aspects of regression are not usually distinguished.}
$(\cat{AM},p)$, be presented as the theory of a linear model, except that
$\beta$, $\mu$, and $y$ are affine space objects, $X: \beta \to \mu$ is an
affine map, and $q: \mu \otimes \sigma^2 \to y$ is an affine-quadratic morphism.
A model $M: \cat{AM} \to \Stat$ then allows the design $X_M: \R^p \to \R^n$ to
be an affine map, so that the regression function has a fixed, generally nonzero
intercept. Such models occur rarely; in practice, the intercept is usually
unknown and is fitted to the data by including a constant column in the design
matrix. So, consider instead the affine models obtained from linear models by
applying the model migration functor
$\aff^*: \Model{\cat{LM}} \to \Model{\cat{AM}}$, induced by the inclusion
functor $\aff: \cat{AM} \hookrightarrow \cat{LM}$.

\begin{proposition} \label{prop:lm-affine-morphisms}
  For any linear models $M$ and $M'$, a morphism
  $\alpha: \aff^*(M) \to \aff^*(M')$ is uniquely determined by a matrices
  $A \in \R^{n \times n'}$ and $B \in \R^{p' \times p}$ and vectors
  $b \in \R^{n'}$ and $c \in \R^{p'}$ such that $A A^\top$ is proportional to the
  identity matrix and the design matrices are intertwined as
  \begin{equation*}
    A X_M = X_{M'} B \qquad\text{and}\qquad b = X_{M'} c.
  \end{equation*}
\end{proposition}
\begin{proof}
  The morphism $\alpha$ consists of affine maps $\alpha_\beta: \R^p \to \R^{p'}$
  and $\alpha_\mu, \alpha_y: \R^n \to \R^{n'}$ and a conic-linear map
  $\alpha_{\sigma^2}: \R_+ \to \R_+$ making the two diagrams of
  \cref{prop:lm-morphisms} commute. The second diagram, concerning the isotropic
  normal family, implies that if $\alpha_y$ is the affine map
  $y \mapsto Ay + b$, then $\alpha_\mu = \alpha_y$ and
  $AA^\top = \alpha_{\sigma^2} I_{n'}$. If $\alpha_\beta$ is the affine map
  $\beta \mapsto B \beta + c$, then the first diagram is the equation
  \begin{equation*}
    A X_M \beta + b = X_{M'} (B \beta + c),
  \end{equation*}
  which holds for all $\beta \in \R^p$. Setting $\beta = 0$ shows that
  $b = X_{M'} c$. The equation that remains is $A X_M = X_{M'} B$.
\end{proof}

Ordinary least squares (OLS) linear regression, the most basic method of fitting
a linear model, is equivariant under model isomorphism. It is even ``laxly''
equivariant under model homomorphism.

\begin{theorem}[Equivariance of linear regression]
  \label{thm:linear-regression-equivariance}
  Let $\alpha: \aff^*(M) \to \aff^*(M')$ be a morphism of linear models $M$ and
  $M'$ with design matrices $X_M \in \R^{n \times p}$ and
  $X_{M'} \in \R^{n' \times p'}$. For any data $y \in \R^n$ and parameters
  $\beta \in \R^p$, the transformed data $y' := \alpha_y(y) := Ay + b$ and
  transformed parameters $\beta' := \alpha_\beta(\beta) := B \beta + c$ satisfy
  \begin{equation*}
    \norm{X_{M'} \beta' - y'} \leq a \norm{X_M \beta - y},
  \end{equation*}
  where $a := \sqrt{\alpha_{\sigma^2}}$. In particular, if
  $\alpha: \aff^*(M) \cong \aff^*(M')$ is an isomorphism, then it sends
  minimizers $\hat\beta \in \argmin_{\beta \in \R^p} \norm{X_M \beta - y}$ to
  minimizers
  $\hat\beta' \in \argmin_{\beta' \in \R^{p'}} \norm{X_{M'} \beta' - y'}$.
\end{theorem}
\begin{proof}
  Using the relations $A X_M = X_{M'} B$ and $b = X_{M'} c$, calculate that
  \begin{equation*}
    X_{M'} \beta' - y'
      = X_{M'} (B \beta + c) - (Ay + b)
      = A(X_M \beta - y)
  \end{equation*}
  Decompose the matrix $A$ as $a Q$, where $a^2 = \alpha_{\sigma^2}$ and
  $Q Q^\top = I_{n'}$. Since $P_Q := Q^\top Q \in \R^{n \times n}$ is an
  orthogonal projection, it follows that
  $X_{M'} \beta' - y' = a Q(X_M \beta -y)$ and
  \begin{equation*}
    \norm{X_{M'} \beta' - y'}
      = a \norm{Q (X_M \beta - y)}
      = a \norm{P_Q (X_M \beta - y)}
      \leq a \norm{X_M \beta - y},
  \end{equation*}
  proving the desired inequality. If, moreover, the homomorphism $\alpha$ is an
  isomorphism, then $y$ and $y'$, and also $\beta$ and $\beta'$, are in
  one-to-one correspondence through $\alpha_y$ and $\alpha_\beta$. Applying the
  inequality to the inverse isomorphism $\alpha^{-1}$ yields
  $\norm{X_M \beta - y} \leq a^{-1} \norm{X_{M'} \beta' - y'}$, so that the
  least-squares objectives are proportional,
  \begin{equation*}
    \norm{X_{M'} \beta' - y'} = a \norm{X_M \beta - y},
  \end{equation*}
  with nonzero constant of proportionality $a$. Thus,
  $\alpha: \aff^*(M) \cong \aff^*(M')$ establishes a one-to-one correspondence
  between least-squares solutions under the designs $X_M$ and $X_{M'}$.
\end{proof}

\paragraph{Linear models with $n$ observations}

Another formulation of the linear model makes the number of observations
explicit in the theory, as in most theories from the previous section. The
\emph{theory of a linear model on $n$ observations} has underlying category
$\cat{LM}_n$ presented by vector space objects $\beta$, $\mu$, and $y$, a
conical space object $\sigma^2$, linear maps $X_1,\dots,X_n: \beta \to \mu$, and
a linear-quadratic morphism $q: \mu \otimes \sigma^2 \to y$. The sampling
morphism $p_n: \beta \otimes \sigma^2 \to y^{\otimes n}$ is
\begin{equation*}
  \input{wiring-diagrams/zoo-statistics/lm-n-sampling}.
\end{equation*}
An intended model $M: \cat{LM}_n \to \Stat$ assigns $M(\beta) = \R^p$, for some
dimension $p$; $M(\mu) = M(y) = \R$ and $M(\sigma^2) = \R_+$;
$M(X_1),\dots,M(X_n)$ to be linear functionals
$X_{M,1}, \dots, X_{M,n}: \R^p \to \R$, identified with row vectors in
$\R^{1 \times p}$; and $M(q)$ to be the univariate normal family. The sampling
morphism $M(p_n): \R^p \times \R_+ \to \R^n$ can then be written as
\begin{equation*}
  y_i \xsim{\mathrm{ind}} \Normal(X_{M,i}\, \beta, \sigma^2), \quad i=1,\dots,n.
\end{equation*}

The analogue of \cref{prop:lm-morphisms}, with a similar proof, is:

\begin{proposition} \label{prop:lm-n-morphisms}
  A morphism $\alpha: M \to M'$ between linear models $M$ and $M'$ on $n$
  observations is uniquely determined by a scalar $a \in \R$ and a matrix
  $B \in \R^{p' \times p}$ such that
  \begin{equation*}
    a X_{M,i} = X_{M',i}\, B, \quad \forall i=1,\dots,n.
  \end{equation*}
  In particular, an isomorphism $\alpha: M \cong M'$ can exist only if $p=p'$
  and is then uniquely determined by a scalar $a \in \R^*$ and a matrix
  $B \in \GL(p,\R)$ such that
  \begin{equation*}
    X_{M',i} = a X_{M,i}\, B^{-1}, \quad \forall i=1,\dots,n.
  \end{equation*}
\end{proposition}

Under the theory of a linear model on $n$ observations, two linear models are
isomorphic if and only if their feature vectors on each observation are
\emph{simultaneously} equivalent as matrices. This is stronger than isomorphism
under the theory of a linear model. Indeed, forming the $n \times p$ design
matrices $X_M$ and $X_{M'}$ by stacking row vectors, we have
$X_{M'} = A X_M B^{-1}$, where $A := a I_n$ is a conformal orthogonal matrix.

The theories of one normal sample, two normal samples, and a normal sequence are
all specializations of the theory of a linear model on $n$ observations. In the
first case, define the supply preserving functor $F: \cat{LM}_n \to \cat{S}$
sending $\beta$ to $\mu$, all of $X_1,\dots,X_n$ to the identity $1_\mu$, and
preserving the other generators.
\begin{equation*}
  F: \input{wiring-diagrams/zoo-statistics/lm-n-design}
  \mapsto \input{wiring-diagrams/zoo-statistics/mu-id}, \quad
  i = 1,\dots,n.
\end{equation*}
Then $F: (\cat{LM}_n, p_n) \to (\cat{S}, p_n)$ is a strict morphism from the
theory of a linear model on $n$ observations into the theory of a normal sample
of size $n$. Next, for any numbers $\vec{n} = (n_1, n_2)$ with $n_1 + n_2 = n$,
define the functor $G: \cat{LM}_n \to \cat{S}$ sending $\beta$ to
$\mu^{\otimes 2}$, all of $X_1,\dots,X_n$ to the projection $\pi_{\mu,\mu}$, the
rest of $X_{n_1+1},\dots,X_n$ to the other projection $\pi_{\mu,\mu}'$, and
preserving the other generators.
\begin{equation*}
  G: \quad \begin{dcases}
    \input{wiring-diagrams/zoo-statistics/lm-n-design}
    \mapsto \input{wiring-diagrams/zoo-statistics/mu-proj1} \quad
    & i = 1, \dots, n_1 \\[2\jot]
    \input{wiring-diagrams/zoo-statistics/lm-n-design}
    \mapsto \input{wiring-diagrams/zoo-statistics/mu-proj2} \quad
    & i = n_1+1, \dots, n_1+n_2.
  \end{dcases}
\end{equation*}
Then $G: (\cat{LM}_n, p_n) \to (\cat{S}, q_{\vec{n}})$ is a morphism into the
theory of two homoscedastic normal samples of sizes $n_1$ and $n_2$. Finally,
the functor $H: \cat{LM}_n \to \cat{S}$ sending $\beta$ to $\mu^{\otimes n}$,
each $X_i$ to the $i$th projection $\mu^{\otimes n} \to \mu$, and preserving the
other generators is a morphism $H: (\cat{LM}_n, p_n) \to (\cat{S}, q_n)$ into
the theory of a normal sequence of length $n$. The model migration functors
$F^*, G^*, H^*: \Model{\cat{S}} \to \Model{\cat{LM}_n}$ yield the expected
linear models when applied to the univariate normal model.

\paragraph{Linear model with $p$ predictors}

The next formulation of the linear model reverses the convention of the previous
one, making the predictors explicit in the theory but suppressing the individual
observations. The \emph{theory of a linear model with $p$ predictors} has
underlying category $\cat{LM}_p$ presented by vector space objects
$\beta_1,\dots,\beta_p$, $\mu$, and $y$, a conical space object $\sigma^2$,
linear maps $X_j: \beta_j \to \mu$ for $j=1,\dots,p$, and a linear-quadratic
morphism $q: \mu \otimes \sigma^2 \to y$. The sampling morphism
$q_p: \beta_1 \otimes \cdots \otimes \beta_p \otimes \sigma^2 \to y$ is
\begin{equation*}
  \input{wiring-diagrams/zoo-statistics/lm-p-sampling}.
\end{equation*}
The choice of a theory with distinct objects $\beta_1,\dots,\beta_p$, rather
than with a single object $\beta$ and a sampling morphism
$\beta^{\otimes p} \otimes \sigma^2 \to y$, is not of great significance,
although it does have a few consequences. One is that a model is allowed to
assign parameter subspaces of different dimensionalities to different
$\beta_j$'s, which can be useful when parameters occur in groups of different
sizes.

Nevertheless, the default models $M: \cat{LM}_p \to \Stat$ assign all of
$\beta_1, \dots, \beta_p$ to be $\R$, both $\mu$ and $y$ to be $\R^n$, for some
number $n$; $\sigma^2$ to be $\R_+$; each of $X_1,\dots,X_p$ to be linear maps
$X_{M,1}, \dots, X_{M,p}: \R \to \R^n$, identified with vectors in $\R^n$; and
$q$ to be the $n$-dimensional isotropic normal family. The sampling morphism
$M(q_p)$ can then be written as
\begin{equation*}
  y \sim \Normal_n(
    X_{M,1}\, \beta_1 + \cdots + X_{M,p}\, \beta_p,\ \sigma^2 I_n).
\end{equation*}
For these models, the analogue of \cref{prop:lm-morphisms,prop:lm-n-morphisms}
is:

\begin{proposition} \label{prop:lm-p-morphisms}
  A morphism $\alpha: M \to M'$ between linear models $M$ and $M'$ with $p$
  predictors is uniquely determined by a matrix $A \in \R^{n' \times n}$ and
  scalars $b_1, \dots, b_p \in \R$ such that $A A^\top$ is proportional to the
  identity matrix and
  \begin{equation*}
    A X_{M,j} = b_j X_{M',j}, \quad \forall j=1,\dots,p.
  \end{equation*}
  In particular, an isomorphism $\alpha: M \cong M'$ can exist only if $n=n'$
  and is then uniquely determined by a matrix $A \in \CO(n)$ and scalars
  $b_1, \dots, b_j \in \R^*$ such that
  \begin{equation*}
    X_{M',j} = b_j^{-1} A X_{M,j}, \quad \forall j=1,\dots,p.
  \end{equation*}
\end{proposition}
\begin{proof}
  A morphism $\alpha: M \to M'$ consists of linear maps
  $\alpha_{\beta_1}, \dots, \alpha_{\beta_p}: \R \to \R$, linear maps
  $\alpha_\mu, \alpha_y: \R^n \to \R^{n'}$, and a conic-linear map
  $\alpha_{\sigma^2}: \R_+ \to \R_+$ making the diagrams
  \begin{equation*}
    \begin{tikzcd}
      \R \arrow[r, "X_{M,j}"] \arrow[d, "\alpha_{\beta_j}"']
        & \R^n \arrow[d, "\alpha_\mu"] \\
      \R \arrow[r, "X_{M',j}"]
        & \R^{n'}
    \end{tikzcd}
    \qquad\qquad
    \begin{tikzcd}
      \R^n \times \R_+ \arrow[r, "\Normal_n^{\mathrm{iso}}"]
        \arrow[d, "\alpha_\mu \times \alpha_{\sigma^2}"']
        & \R^n \arrow[d, "\alpha_y"] \\
      \R^{n'} \times \R_+ \arrow[r, "\Normal_{n'}^{\mathrm{iso}}"]
        & \R^{n'}
    \end{tikzcd}
  \end{equation*}
  commute for all $j=1,\dots,p$. By the now familiar argument, the last diagram
  is equivalent to having $\alpha_\mu = \alpha_y = A$ for some matrix
  $A \in \R^{n' \times n}$ such that $A A^\top = \alpha_{\sigma^2} I_{n'}$.
  Setting $b_j := \alpha_{\beta_j}$, the other diagrams are the equations
  $A X_{M,j} = b_j X_{M',j}$ for $j=1,\dots,p$.
\end{proof}

By defining the theory of a linear model with $p$ predictors to have distinct
objects $\beta_1, \dots, \beta_p$, a model homomorphism is allowed to rescale
each predictor individually. Such transformations are useful when the predictors
have different dimensions of measurement, say for converting one predictor from
feet to meters and another from pounds to kilograms. Similarly to the linear
model on $n$ observations, isomorphism under the theory of a linear model with
$p$ predictions is stronger than under the general theory of a linear model.
Indeed, forming the $n \times p$ design matrices $X_M$ and $X_{M'}$ by stacking
column vectors, we have $X_{M'} = A X_M B^{-1}$, where
$B := \diag(b_1,\dots,b_p)$ is an invertible diagonal matrix.

The global null hypothesis
\begin{equation*}
  H_0: \beta_1 = \cdots = \beta_p = 0,
\end{equation*}
as might be tested by an $F$-test, is represented by the colax morphism
$(1_{\cat{LM}_p}, 0_\beta \otimes 1_{\sigma^2}, 1_y): (\cat{LM}_p, q_p) \to
(\cat{LM}_p, q_{p,0})$, where the object $\beta$ is shorthand for
$\beta_1 \otimes \cdots \otimes \beta_p$ and the reduced sampling morphism
$q_{p,0}: \sigma^2 \to y$ is
\begin{equation*}
  \input{wiring-diagrams/zoo-statistics/lm-p-global-null-lhs} =
  \input{wiring-diagrams/zoo-statistics/lm-p-global-null-rhs}.
\end{equation*}
The individual null hypotheses $H_{0,j}: \beta_j = 0$ for $j=1,\dots,p$, as
might be tested by marginal $t$-tests, are expressed similarly.

As another relation in this vein, consider enlarging a linear model with $p$
predictors to a linear model with an additional $k$ predictors. The inclusion
functor $\iota: \cat{LM}_p \hookrightarrow \cat{LM}_{p+k}$ defines a colax
theory morphism
$(\iota, 1_{\beta_{1:p}} \otimes 0_{\beta_{p+1:p+k}} \otimes 1_{\sigma^2}, 1_y):
(\cat{LM}_p, q_p) \to (\cat{LM}_{p+k}, q_{p+k})$, where we use the shorthand
$\beta_{i:i+j} := \beta_i \otimes \cdots \otimes \beta_{i+j}$ for any numbers
$i,j$. The colaxness equation for this morphism
\begin{equation*}
  \input{wiring-diagrams/zoo-statistics/lm-p-ext-null} =
  \input{wiring-diagrams/zoo-statistics/lm-p-sampling}
\end{equation*}
generalizes that of the global null hypothesis, which is the case where $p=0$.
The model migration functor
$\iota^*: \Model{\cat{LM}_{p+k}} \to \Model{\cat{LM}_p}$ sends a linear model
with $p+k$ predictors to a linear model with $p$ predictors by setting the last
$k$ coefficients to zero.

\paragraph{Linear model with $n$ observations and $p$ predictors}

The last formulation of the linear model considered here jointly refines the two
previous theories by making both the observations and the predictors explicit in
the theory. The \emph{theory of a linear model on $n$ observations and $p$
  predictors} has underlying category $\cat{LM}_{n,p}$ presented by vector space
objects $\beta_1,\dots,\beta_p$, $\mu$, and $y$, a conical space object
$\sigma^2$, linear maps $X_{i,j}: \beta_j \to \mu$ for $i=1,\dots,n$ and
$j=1,\dots,p$, and a linear-quadratic morphism $q: \mu \otimes \sigma^2 \to y$.
The sampling morphism
$q_{n,p}: \beta_1 \otimes \cdots \otimes \beta_p \otimes \sigma^2 \to y^{\otimes
  n}$ is
\begin{equation*}
  \input{wiring-diagrams/zoo-statistics/lm-np-sampling}.
\end{equation*}
The intended models $M: \cat{LM}_{n,p} \to \Stat$ assign all of
$\beta_1,\dots,\beta_p$, $\mu$, and $y$ to be $\R$, the object $\sigma^2$ to be
$\R_+$, and the morphism $q$ to be the univariate normal family. The morphisms
$X_{i,j}$ are then arbitrary scalars $X_{M,i,j} \in \R$, and the sampling
distribution $M(q_{n,p})$ can be written as
\begin{equation*}
  y_i \xsim{\mathrm{ind}} \Normal(
    X_{M,i,1}\,\beta_1 + \cdots + X_{M,i,p}\,\beta_p, \sigma^2), \quad
  i = 1,\dots,n.
\end{equation*}

Following \cref{prop:lm-n-morphisms,prop:lm-p-morphisms}, it is easily shown
that:

\begin{proposition} \label{prop:lm-np-morphisms}
  A morphism $\alpha: M \to M'$ between linear models $M$ and $M'$ with $n$
  observations and $p$ predictors consists of scalars
  $a := \alpha_\mu = \alpha_y$ and $b_j := \alpha_{\beta_j}$, $j=1,\dots,p$,
  where $\alpha_{\sigma^2} = a^2$, such that
  \begin{equation*}
    a X_{M,i,j} = b_j X_{M',i,j}, \quad \forall i=1,\dots,n,\ j=1,\dots,p.
  \end{equation*}
  The morphism $\alpha$ is an isomorphism if and only if all of $a$,
  $b_1, \dots, b_p$ are nonzero.
\end{proposition}

So, under the theory with $n$ and $p$ fixed, two linear models are isomorphic if
and only if each pair of corresponding columns in the design matrices
$X_M := (X_{M,i,j})_{i,j}$ and $X_{M'} := (X_{M',i,j})_{i,j}$ are proportional.
This condition is quite strong, certainly stronger than isomorphism under any of
the previous theories. But this theory of a linear model is still not the most
explicit possible. For any fixed matrix $X \in \R^{n \times p}$, a fully
specified \emph{theory of a linear model with design matrix $X$} takes the
category $\cat{S}$ from the previous section and directly encodes matrix
multiplication by $X$ into the sampling morphism via the scalar multiplications
$X_{i,j}: \mu \to \mu$. The theories of $k$ normal samples or of a normal
sequence can be seen as arising this way. However, outside of these special
cases, this form of the theory of a linear model seems to be too explicit to be
practical.

\paragraph{Relations between theories}

The four theories of the linear model are related to each other by refinement of
the design matrix, by dividing the full matrix into rows or columns and then
dividing the rows or columns into their individual components. The relationships
are formalized by a commutative diagram of strict theory morphisms.
\begin{equation*}
  \begin{tikzcd}[column sep=-2em,row sep=small]
    & \underset{\text{general LM}}{(\cat{LM}, p)}
      \arrow[dl, "F_n"']
      \arrow[dr, "G_p" ] & \\
    \underset{\text{LM with $n$ observations}}{(\cat{LM}_n, p_n)}
      \arrow[dr, "F_{n,p}"']
    & & \underset{\text{LM with $p$ predictors}}{(\cat{LM}_p, q_p)}
      \arrow[dl, "G_{n,p}"] \\
    & \underset{\text{LM with $n$ observations and $p$ predictors}}%
      {(\cat{LM}_{n,p}, q_{n,p})} &
  \end{tikzcd}
\end{equation*}
In presenting these functors, we adopt the convention that any generator not
explicitly mapped is preserved in the sense of being mapped to the corresponding
generator with the same name. The supply preserving functor
$F_n: \cat{LM} \to \cat{LM}_n$ sends $\mu$ to $\mu^{\otimes n}$ and $y$ to
$y^{\otimes n}$; divides the design morphism by rows,
\begin{equation*}
  F_n: \input{wiring-diagrams/zoo-statistics/lm-design}
  \mapsto \input{wiring-diagrams/zoo-statistics/lm-design-to-n};
\end{equation*}
and divides the morphism $q$ accordingly
\begin{equation*}
  F_n: \input{wiring-diagrams/zoo-statistics/normal}
  \mapsto
  \input{wiring-diagrams/zoo-statistics/normal-means-sampling}.
\end{equation*}
Note that the right-hand side is the sampling morphism of the theory of a normal
sequence of length $n$. Similarly, the functor $G_p: \cat{LM} \to \cat{LM}_p$
sends $\beta$ to $\beta_1 \otimes \cdots \otimes \beta_p$ and divides the design
matrix by columns,
\begin{equation*}
  G_p: \input{wiring-diagrams/zoo-statistics/lm-design}
  \mapsto \input{wiring-diagrams/zoo-statistics/lm-design-to-p}.
\end{equation*}
This defines morphisms comprising the upper legs of the commutative diagrams.

As for the lower legs, the functor $F_{n,p}: \cat{LM}_n \to \cat{LM}_{n,p}$
sends $\beta$ to $\beta_1 \otimes \cdots \otimes \beta_p$ and divides each row
of the design matrix into its components,
\begin{equation*}
  F_{n,p}: \input{wiring-diagrams/zoo-statistics/lm-n-design}
  \mapsto \input{wiring-diagrams/zoo-statistics/lm-n-design-to-np},
  \quad i=1,\dots,n.
\end{equation*}
Finally, the functor $G_{n,p}: \cat{LM}_p \to \cat{LM}_{n,p}$ sends $\mu$ to
$\mu^{\otimes n}$ and $y$ to $y^{\otimes n}$; divides each column of the design
matrix into its components,
\begin{equation*}
  G_{n,p}: \input{wiring-diagrams/zoo-statistics/lm-p-design}
  \mapsto \input{wiring-diagrams/zoo-statistics/lm-p-design-to-np},
  \quad j=1,\dots,n;
\end{equation*}
and acts on the morphism $q$ in the same way as the functor
$F_n: \cat{LM} \to \cat{LM}_n$.

\section{Bayesian and hierarchical linear models}
\label{sec:hierarchical-lm}

The algebraic view of statistical modeling formalizes the everyday practice of
building complex statistical models out of simpler ones. By way of illustration,
this section constructs hierarchical linear models from the standard linear
model. Linear models with hierarchical structure go by many names, including but
not limited to \emph{hierarchical models}, \emph{multilevel models},
\emph{random coefficient models}, \emph{random-effects models}, and \emph{mixed}
or \emph{mixed-effects models}. As the proliferation of names suggests, there
many kinds of hierarchical models and perspectives on them. This section treats
a few basic hierarchical extensions of the linear models in
\cref{sec:lm-discrete,sec:lm-general}. Although their interpretations and
statistical inference differ, frequentist hierarchical models are also
structurally similar to Bayesian models, and so we begin with an example of the
latter.

\paragraph{Bayesian one-sample models} Among the simplest Bayesian models of
continuous data are those of a normal sample with unknown mean and variance, the
Bayesian version of one-sample normal model from \cref{sec:lm-discrete}. We
present two versions of a Bayesian theory of one normal sample, the first making
the priors on the mean and variance independent and the second giving them a
hierarchical structure.

Define the first Bayesian theory $(\cat{T}_{\mathrm{ind}}, p_n, \pi)$ of one
normal sample of size $n$ as follows. Let $\cat{T}_{\mathrm{ind}}$ be generated
by vector objects $y$, $\mu$, and $\mu_0$; conical space objects $\sigma^2$,
$\sigma_0^2$, and $\tau_0^2$; a discrete object $\nu_0$; linear-quadratic
morphisms $q: \mu \otimes \sigma^2 \to \mu$ and
$\pi_\mu: \mu_0 \otimes \tau_0^2 \to \mu$; a morphism
$\pi_{\sigma^2}: \nu_0 \otimes \sigma_0^2 \to \sigma^2$, homogeneous in its
second argument; and hyperparameters $\tilde \mu_0$, $\tilde \sigma_0^2$,
$\tilde \tau_0^2$, and $\tilde \nu_0$. Here the \emph{hyperparameter}
$\tilde \theta$ corresponding to an object $\theta$ is decorated with a tilde
and is a deterministic generalized element of that object, that is, a map
$\tilde \theta: I \to \theta$. The sampling morphism
$p_n: \mu \otimes \sigma^2 \to y^{\otimes n}$ of the Bayesian theory is that of
the theory of one normal sample of size $n$, defined in \cref{sec:lm-discrete}.
Finally, the prior morphism $\pi: I \to \mu \otimes \sigma^2$ is the independent
product
\begin{equation*}
  \input{wiring-diagrams/zoo-statistics/one-sample-normal-prior}.
\end{equation*}
The marginal, or prior predictive, morphism $\pi \cdot p_n: I \to y^{\otimes n}$
is then
\begin{equation*}
  \input{wiring-diagrams/zoo-statistics/one-sample-normal-prior-predictive}.
\end{equation*}

The standard univariate models $M: \cat{T}_{\mathrm{ind}} \to \Stat$ of this
Bayesian theory assign the objects $y$, $\mu$, and $\mu_0$ to $\R$; the objects
$\sigma^2$, $\sigma_0^2$, and $\tau_0^2$ to $\R_+$; the object $\nu_0$ to
$\PosR = (0,\infty)$; the morphisms $q$ and $\pi_\mu$ to the univariate normal
family; and the morphism $\pi_{\sigma^2}$ to the reparameterized inverse-gamma
family
\begin{equation*}
  M(\pi_{\sigma^2}): \PosR \times \R_+ \to \R_+,\
    (\nu_0, \sigma_0^2) \mapsto
      \InvGammaDist\left(\frac{\nu_0}{2}, \frac{\nu_0 \sigma_0^2}{2}\right).
\end{equation*}
Here the \emph{inverse-gamma distribution} $\InvGammaDist(\alpha,\beta)$ is the
distribution of the reciprocal of a $\GammaDist(\alpha,\beta)$ random variable,
where the gamma distribution is parameterized by shape $\alpha$ and rate
$\beta$. The reparameterized family $M(\pi_{\sigma^2})$ is sometimes called the
\emph{scaled inverse chi-squared family} \cite[\S 2.6]{gelman2013}. Each
Bayesian model $M$ also includes a choice of hyperparameters
$\mu_{0,M} := M(\tilde\mu_0) \in \R$; $\sigma_{0,M}^2, \tau_{0,M}^2 \in \R_+$;
and $\nu_{0,M} \in \R^*_+$.

In classical notation, this Bayesian model is specified by a list of assertions
about the conditional distributions of random variables:
\begin{align*}
  \mu &\sim \Normal(\mu_0, \tau_0^2) \\
  \sigma^2 &\sim \InvGammaDist(\nu_0/2, \nu_0 \sigma_0^2/2) \\
  y_1, \dots, y_n \given \mu, \sigma^2
    &\xsim{\mathrm{iid}} \Normal(\mu, \sigma^2).
\end{align*}
The independence of $\mu$ and $\sigma^2$ under the prior is implicit in the
notation, as is the distinction between parameters and hyperparameters.

Unless their hyperparameters bear certain relations to each other, two Bayesian
models $M$ and $M'$ will not be related by any nontrivial model homomorphisms
$M \to M'$. A generic Bayesian model likewise has no nontrivial automorphisms.
This situation reflects a common criticism of invariance principles in
statistics, namely that the prior information we usually possess nullifies
geometrical symmetries of location or scale. On the other hand, if
noninformative priors are desired, then invariance principles provide a way of
generating them \cite[Chapter 9]{robert2007}.

Bayesian theories tend to have a greater number of plausible models than the
corresponding frequentist theories, as the choice of priors and hyperparameters
is not canonical. The inverse-gamma family is, for reasons of analytical
convenience, the classic choice of prior for the top-level scale parameter in a
Bayesian model, but it has been argued that the half $t$-family makes for a
better default \cite{gelman2006b,polson2012}. A model $M$ of the Bayesian theory
above would then assign $M(\nu_0) = \{1,2,\dots\}$, $M(\sigma_0^2) = \R_+$, and
$M(\pi_{\sigma^2})$ to be the \emph{half $t$-family}, that is, the distribution
of the absolute value of a centered $t$-random variable, having $\nu_0$ degrees
of freedom and scale parameter $\sigma_0^2$. Restricting to $\nu_0 = 1$ degree
of freedom yields the popular \emph{half Cauchy family}.

The theory $(\cat{T}_{\mathrm{ind}}, p_n, \pi)$ is the most obvious Bayesian
theory of a normal sample, but it is not the most standard. The traditional
reason for this is that under the independent prior
$\pi_\mu \otimes \pi_{\sigma^2}$ for $\mu$ and $\sigma^2$, the inverse-gamma
model for $\sigma^2$ is not a conjugate prior, but only conditionally conjugate
\cite[\S 6.1]{hoff2009}. A conjugate prior is obtained by making the prior
variance $\tau_0^2$ of the mean of $y$ proportional to the variance of $y$
\cite[\S 5.3]{hoff2009}.

Specifically, define a new Bayesian theory $(\cat{T}, p_n, \pi)$ by presenting
$\cat{T}_{\mathrm{ind}}$ as $\cat{T}$, except that hyperparameter
$\tilde \tau_0^2: I \to \tau_0^2$ is replaced by a conic-linear map
$\tilde \kappa_0: \sigma^2 \to \tau_0^2$. The sampling morphism
$p_n: \mu \otimes \sigma^2 \to y^{\otimes n}$ is again that of the theory of one
normal sample of size $n$, but the prior morphism
$\pi: I \to \mu \otimes \sigma^2$ is now
\begin{equation*}
  \input{wiring-diagrams/zoo-statistics/one-sample-normal-prior-coupled}.
\end{equation*}
The standard univariate model $M: \cat{T} \to \Stat$ is defined as before. In
conventional notation, the Bayesian model is:
\begin{align*}
  \sigma^2 &\sim \InvGammaDist(\nu_0/2, \nu_0 \sigma_0^2/2) \\
  \mu \given \sigma^2 &\sim \Normal(\mu_0, \kappa_0 \sigma^2) \\
  y_1, \dots, y_n \given \mu, \sigma^2
    &\xsim{\mathrm{iid}} \Normal(\mu, \sigma^2).
\end{align*}

Both Bayesian statistical theories extend the statistical theory
$(\cat{S}, p_n)$ of one normal sample defined in \cref{sec:lm-discrete}.
Formally, the inclusion functors $\iota: \cat{S} \hookrightarrow \cat{T}$ and
$\iota_{\mathrm{ind}}: \cat{S} \hookrightarrow \cat{T}_{\mathrm{ind}}$ define
strict theory morphisms
\begin{equation*}
  \begin{tikzcd}
    (\cat{T}, p_n)
    & (\cat{S}, p_n) \arrow[l, "\iota"']
      \arrow[r, "\iota_{\mathrm{ind}}"]
    & (\cat{T}_{\mathrm{ind}}, p_n)
  \end{tikzcd}
\end{equation*}
into the underlying statistical theories of the Bayesian theories. The model
migration functors $\iota^*: \Model{\cat{T}} \to \Model{\cat{S}}$ and
$\iota_{\mathrm{ind}}^*: \Model{\cat{T}_{\mathrm{ind}}} \to \Model{\cat{S}}$
then return the sampling distributions of the Bayesian models. There are also
colax theory morphisms
\begin{equation*}
  \begin{tikzcd}[column sep=huge]
    (\cat{T}, \pi \cdot p_n)
    & (\cat{S}, p_n) \arrow[l, "{(\iota, \pi, 1_{y^{\otimes n}})}"']
      \arrow[r, "{(\iota_{\mathrm{ind}}, \pi, 1_{y^{\otimes n}})}"]
    & (\cat{T}_{\mathrm{ind}}, \pi \cdot p_n)
  \end{tikzcd}
\end{equation*}
into the marginalized Bayesian theories.

\paragraph{$k$-sample normal model with random effects}

Generalizing the one and two sample theories from \cref{sec:lm-discrete}, let
$(\cat{S}, \mu^{\otimes k} \otimes \sigma^2 \xrightarrow{q_{\vec{n}}} y^{\otimes
  n})$ be the theory of $k$ homoscedastic normal samples of sizes
$\vec{n} := (n_1,\dots,n_k)$, with total sample size $n := n_1 + \cdots + n_k$.
In a hierarchical model, the $k$ groups do not comprise a fixed class of groups;
rather, they are regarded as being sampled from a larger population of groups.
Consider, for example, sampling $k$ specific schools from some population of
schools and then, from each school $i = 1,\dots,k$, sampling $n_i$ students.

To describe this two-level sampling scheme, the statistical \emph{theory of $k$
  normal samples (of sizes $n_1, \dots, n_k$) with random effects} takes the
underlying category $\cat{S}^{(2)}$ generated by $\cat{S}$ together with a
vector space object $\mu_0$, a conical space $\sigma_0^2$, and another
linear-quadratic morphism $r: \mu_0 \otimes \sigma_0^2 \to \mu$. The sampling
morphism
$p_{\vec{n}}^{(2)}: \mu_0 \otimes \sigma_0^2 \otimes \sigma^2 \to y^{\otimes n}$
is
\begin{equation*}
  \input{wiring-diagrams/zoo-statistics/k-sample-hierarchical-sampling}.
\end{equation*}
The intended univariate model $M$ takes $M(\mu)$, $M(\mu_0)$, and $M(y)$ to be
the real numbers, $M(\sigma^2)$ and $M(\sigma_0^2)$ to be the nonnegative reals,
and $M(q)$ and $M(r)$ to be the normal family.

The random effects model is conventionally written in various styles depending
on how the observations are indexed. Using the flattened indexing scheme
suggested by the diagram, the model is
\begin{align*}
  \mu_i &\sim \Normal(\mu_0, \sigma_0^2), \quad i=1,\dots,k \\
  y_j &\sim \Normal(\mu_{i[j]}, \sigma^2), \quad j=1,\dots,n,
\end{align*}
where the map $i[-]: \{1,\dots,n\} \to \{1,\dots,k\}$ assigns observations to
groups. Alternatively, indexing the observations at two levels, the model is
\begin{align*}
  \mu_i &\sim \Normal(\mu_0, \sigma_0^2), \quad i=1,\dots,k \\
  y_{ij} &\sim \Normal(\mu_i, \sigma^2), \quad i=1,\dots,k, \quad j=1,\dots,n_i.
\end{align*}
As the definition of the statistical theory shows, under the algebraic approach
to statistical modeling, two-level and higher-level hierarchical models can be
constructed recursively without recourse to index manipulation.

Despite having more structure, the normal model with random effects has exactly
the same symmetries as the normal model.

\begin{proposition} \label{prop:hierarchical-normal-morphisms}
  The endomorphisms of the univariate normal random effects model $M$ are
  isomorphic to the multiplicative monoid $\R$. In particular, the
  automomorphism group of $M$ is isomorphic to $\R^*$.
\end{proposition}
\begin{proof}
  A morphism $\alpha: M \to M$ consists of scalars
  $\alpha_\mu, \alpha_{\mu_0}, \alpha_y \in \R$ and
  $\alpha_{\sigma^2}, \alpha_{\sigma_0^2} \in \R_+$. As in
  \cref{prop:univariate-normal-morphisms}, naturality with respect to
  $q: \mu \otimes \sigma^2 \to \mu$ is equivalent to
  $a := \alpha_\mu = \alpha_y$ and $\alpha_{\sigma^2} = a^2$. Naturality with
  respect to $r: \mu_0 \otimes \sigma_0^2 \to \mu$ then gives that
  $\alpha_{\mu_0} = \alpha_\mu = a$ and $\alpha_{\sigma_0^2} = a^2$.
\end{proof}

The theory $(\cat{S}^{(2)}, p_{\vec{n}})$ of $k$ normal samples with random
effects can be seen as a composite of the theory $(\cat{S}, p_k)$ of one normal
sample of size $k$ with the theory $(\cat{S}, p_{\vec{n}})$ of $k$ homoscedastic
normal samples. Specifically, if $\gen{\mu}$ is generated by a vector space
object $\mu$, then the diagram
\begin{equation*}
  \begin{tikzcd}[row sep=small, column sep=small]
    & \gen{\mu} \arrow[dl, hook] \arrow[dr, hook] & \\
    \cat{S} \arrow[dr, tail, "F^{(1)}"'] & &
      \cat{S} \arrow[dl, tail, "F^{(2)}"] \\
    & \cat{S}^{(2)} &
  \end{tikzcd}
\end{equation*}
is a pushout of linear-algebraic Markov categories, where the ``level-one''
embedding functor $F^{(1)}: \cat{S} \rightarrowtail \cat{S}^{(2)}$ maps the
morphism $q: \mu \otimes \sigma^2 \to y$ in $\cat{S}$ to its counterpart $q$ in
$\cat{S}^{(2)}$ and the ``level-two'' functor
$F^{(2)}: \cat{S} \rightarrowtail \cat{S}^{(2)}$ maps $q$ in $\cat{S}$ to the
morphism $r: \mu_0 \otimes \sigma_0^2 \to \mu$ in $\cat{S}^{(2)}$ and acts
accordingly on objects. A more general approach to composing statistical
theories is suggested in \cref{sec:future-work}; however, a careful development
of the compositionality of statistical theories and models is beyond the scope
of this text.

\paragraph{Linear mixed models}

A \emph{linear model with mixed effects}, or for short a \emph{linear mixed
  model}, is a linear model that combines fixed predictors, as in an ordinary
linear model, with random effects. The fixed and random effects each have their
own design matrices, so that a linear mixed model with $n$ observations, $p$
fixed effects, and $q$ random effects is specified by an $n \times p$ design
matrix $X$ and an $n \times q$ design matrix $Z$. Like the theory of a linear
model in \cref{sec:lm-general}, the theory of a linear mixed model admits many
variations, depending on which of $n$, $p$, and $q$ are made explicit. Only some
of the more general theories are presented.

The statistical \emph{theory of a linear mixed model} has underlying category
$\cat{LMM}$ presented by vector space objects $\beta$, $\mu$, $b$, and $y$;
conical space objects $\sigma_b^2$ and $\sigma^2$; linear maps
$X: \beta \to \mu$ and $Z: b \to \mu$; a linear-quadratic morphism
$q: \mu \otimes \sigma^2 \to y$; and a quadratic morphism $r: \sigma_b^2 \to b$.
The sampling morphism $p: \beta \otimes \sigma_b^2 \otimes \sigma^2 \to y$ is
\begin{equation*}
  \input{wiring-diagrams/zoo-statistics/lm-mixed-sampling}.
\end{equation*}
The intended models $M: \cat{LMM} \to \Stat$ take $M(\beta)$ and $M(b)$ to be
$\R^p$ and $\R^q$, for some dimensions $p$ and $q$; both $M(\mu)$ and $M(y)$ to
be $\R^n$, for some dimension $n$; $M(\sigma^2) = \R_+$ and $M(q)$ to be the
$n$-dimensional isotropic normal family; and $M(\sigma_b^2) = \PSD^q$ and $M(r)$
to be the $q$-dimensional centered normal family. The linear maps $M(X)$ and
$M(Z)$ are arbitrary matrices $X_M \in \R^{n \times p}$ and
$Z_M \in \R^{n \times q}$. Under this model, the sampling distribution
$M(p): \R^p \times \PSD^q \times \R_+ \to \R^n$ is
\begin{align*}
  b &\sim \Normal_q(0, \Sigma_b) \\
  y &\sim \Normal_n(X\beta + Zb, \sigma^2 I_n),
\end{align*}
with parameters $\beta \in \R^p$, $\Sigma_b \succeq 0$, and $\sigma^2 \geq 0$.

The symmetries of the linear model, described in \cref{prop:lm-morphisms},
generalize to:

\begin{proposition} \label{prop:lmm-morphisms}
  A morphism $\alpha: M \to M'$ between linear mixed models $M$ and $M'$ is
  uniquely determined by matrices $A \in \R^{n' \times n}$,
  $B \in \R^{p' \times p}$, and $C \in \R^{q' \times q}$ such that $A A^\top$ is
  proportional to the identity and both pairs of design matrices are
  intertwined:
  \begin{equation*}
    A X_M = X_{M'} B \qquad\text{and}\qquad A Z_M = Z_{M'} C.
  \end{equation*}
  In particular, an isomorphism $\alpha: M \cong M'$ can exist only if all three
  dimensions are equal and is then uniquely determined by matrices
  $A \in \CO(n)$, $B \in \GL(p,\R)$, and $C \in \GL(q,\R)$ exhibiting both pairs
  of design matrices as equivalent:
  \begin{equation*}
    X_{M'} = A X_M B^{-1} \qquad\text{and}\qquad Z_{M'} = A Z_M C^{-1}.
  \end{equation*}
\end{proposition}

A linear mixed model with no random effects is just a linear model. To formalize
this relationship, define the supply preserving functor
$F: \cat{LMM} \twoheadrightarrow \cat{LM}$ that sends the objects $\sigma_b^2$
and $b$ to the monoidal unit $I$, thus forcing $r: \sigma_b^2 \to b$ to the
trivial morphism $1_I$; sends $Z: b \to \mu$ to the zero map $0_\mu: I \to \mu$;
and preserves all other generators. Then $F: (\cat{LMM}, p) \to (\cat{LM}, p)$
is a strict morphism of statistical theories. The model migration functor
$F^*: \Mod(\cat{LM}) \to \Mod(\cat{LMM})$ transforms linear models into linear
mixed models with $q=0$, hence having no random effects, and sends the model
homomorphisms in \cref{prop:lm-morphisms} to those in \cref{prop:lmm-morphisms}
having as $C$ the degenerate $0 \times 0$ matrix.

Estimating the unconstrained covariance matrix $\Sigma_b$ can be problematic
when $q$ is large. In this case, the covariance is often restricted to a
function $\Psi$ of some parameter $\theta$ of smaller dimension, so that
$\Sigma_b = \Psi(\theta)$. As a statistical theory $(\cat{LMM}', p')$, let
$\cat{LMM}'$ be generated by $\cat{LMM}$ together with an object $\theta$ and a
map $\psi: \theta \to \sigma_b^2$, and let
$p': \beta \otimes \theta \otimes \sigma^2 \to y^{\otimes n}$ be the evident
sampling morphism. If $\iota: \cat{LMM} \hookrightarrow \cat{LMM}'$ is the
inclusion functor, then the colax theory morphism
\begin{equation*}
  (\iota, 1_\beta \otimes \psi \otimes 1_{\sigma^2}, 1_{y^{\otimes n}}):
  (\cat{LMM}, p) \to (\cat{LMM}', p')
\end{equation*}
represents the restriction of the unconstrained model to the constrained one. In
the other direction, let $P: \cat{LMM}' \twoheadrightarrow \cat{LMM}$ be the
functor sending $\theta$ to $\sigma_b^2$ and $\psi$ to $1_{\sigma_b^2}$ and
preserving the other generators. Then the strict theory morphism
$P: (\cat{LMM}', p') \to (\cat{LMM}, p)$ recovers the unconstrained models as
models ``constrained'' by the identity function.

For future reference, consider yet another variant $(\cat{LMM}_n, p_n)$, the
\emph{theory of a linear mixed model on $n$ observations}. The category
$\cat{LMM}_n$ is presented by vector space objects $\beta$, $\mu$, $b$, and $y$;
conical space objects $\sigma_b^2$ and $\sigma^2$; linear maps
$X_1,\dots,X_n: \beta \to \mu$ and $Z_1,\dots,X_n: b \to \mu$; a
linear-quadratic morphism $q: \mu \otimes \sigma^2 \to y$; and a quadratic
morphism $r: \sigma_b^2 \to b$. The sampling morphism
$p_n: \beta \otimes \sigma_b^2 \otimes \sigma^2 \to y^{\otimes n}$ is
\begin{equation*}
  \input{wiring-diagrams/zoo-statistics/lmm-n-sampling}.
\end{equation*}
Like the theory morphism $F: (\cat{LMM}, p) \to (\cat{LM}, p)$ just defined, a
projection functor $F_n: \cat{LMM}_n \twoheadrightarrow \cat{LM}_n$ sends $b$
and $\sigma_b^2$ to the unit $I$ and all of $Z_1,\dots,Z_n$ to the zero map
$0_\mu$, thereby defining a strict theory morphism
$F_n: (\cat{LMM}_n, p_n) \to (\cat{LM}_n, p_n)$. There is also a strict theory
morphism $(\cat{LMM}, p) \to (\cat{LMM}_n, p_n)$ that splits $X$ and $Z$ into
blocks $X_1,\dots,X_n$ and $Z_1,\dots,Z_n$, analogously to the morphism
$(\cat{LM}, p) \to (\cat{LM}_n, p_n)$ from \cref{sec:lm-general}.

\section{Generalized linear models} \label{sec:glm}

Generalized linear models (GLMs) extend the linear model from the normal family
to other exponential families, thus allowing discrete as well as continuous
responses, while retaining the linear dependence on the parameters that is the
hallmark of the linear model. The theory of a generalized linear model is most
easily formulated when the number of observations $n$ is fixed, since the link
function is computed pointwise across the $n$ dimensions. The number of
parameters $p$ can also be made explicit in the theory, as in
\cref{sec:lm-general}, although we omit that formulation. We present the theory
of a generalized linear model both with and without a dispersion parameter, as
well as a few other variations.

\paragraph{Generalized linear models without a dispersion parameter}

The \emph{theory of a generalized linear model on $n$ observations without a
  dispersion parameter} has underlying category $\cat{GLM0}_n$ presented by
vector space objects $\beta$ and $\eta$, a convex space object $\mu$, and a
discrete object $y$; a map $g: \mu \to \eta$, the \emph{link morphism}, and a
map $h: \eta \to \mu$, the \emph{mean morphism}; linear maps
$X_1, \dots, X_n: \beta \to \eta$, constituting the \emph{design}; and a
morphism $q: \mu \to y$, subject to the condition that $g$ and $h$ be mutually
inverse:
\begin{equation*}
  \input{wiring-diagrams/zoo-statistics/glm-link-inverse} =
  \input{wiring-diagrams/zoo-statistics/mu-id}
  \qquad\text{and}\qquad
  \input{wiring-diagrams/zoo-statistics/glm-mean-inverse} =
  \input{wiring-diagrams/zoo-statistics/eta-id}.
\end{equation*}
The sampling morphism $p_n: \beta \to y^{\otimes n}$ of the theory is
$\Copy_{\beta,n} \cdot (X_1 h q \otimes \cdots \otimes X_n h q)$, or
\begin{equation*}
  \input{wiring-diagrams/zoo-statistics/glm-n-no-dispersion-sampling}.
\end{equation*}

One kind of \emph{generalized linear model} $M: \cat{GLM0}_n \to \Stat$ assigns
$M(\beta)$ to be $\R^p$, for some dimension $p$; $M(\eta)$ to be $\R$; and
$M(q): M(\mu) \to M(y)$ to be a one-dimensional exponential family parameterized
by its mean, which belongs to an interval $M(\mu)$ of the real line
(\cref{ex:exponential-family}). The function $g_M := M(g): \R \to M(\mu)$ is the
\emph{link function} and its inverse $h_M := M(h): M(\mu) \to \R$ is the
\emph{mean function}. When the link function parameterizes the exponential
family by its canonical parameter, it is called the \emph{canonical link
  function}. The model $M$ also specifies $M(X_1),\dots,M(X_n)$ to be any linear
functionals $X_{M,1}, \dots, X_{M,n}: \R^p \to \R$, identified with rows of the
design matrix $X_M \in \R^{n \times p}$. The sampling distribution
$M(p): \R^p \to \R^n$ is then
\begin{equation*}
  y_i \xsim{\mathrm{ind}} M(q)(g_M^{-1}(X_{M,i}\, \beta)), \quad i=1,\dots,n,
\end{equation*}
with parameter $\beta \in \R^p$.

\begin{example}[Binary logistic regression]
  \label{ex:binary-logistic-regression}
  Under the \emph{binary logistic regression model}, $M(q)$ is the Bernoulli
  family $\Ber: \pi \mapsto \pi \delta_1 + (1-\pi) \delta_0$, restricted to the
  open unit interval $M(\mu) = (0,1)$ and taking binary values in
  $M(y) = \{0,1\}$. The link function is then an invertible map
  $g_M: (0,1) \to \R$. The canonical link is the \emph{logit} or \emph{log-odds}
  function
  \begin{equation*}
    \logit(\pi) := \log\left(\frac{\pi}{1-\pi}\right),
  \end{equation*}
  with corresponding mean function the \emph{logistic sigmoid}
  \begin{equation*}
    \logistic(x) := \logit^{-1}(x) = \frac{e^x}{1+e^x} = \frac{1}{1+e^{-x}}.
  \end{equation*}
  Thus, under the canonical link, the sampling distribution is
  \begin{equation*}
    y_i \sim \Ber\left(\frac{1}{1+e^{-X_{M,i}\beta}}\right), \quad i=1,\dots,n.
  \end{equation*}
  Any inverse of a continuous, strictly increasing CDF can serve as the link
  function for the logistic model, motivated by the supposition that the
  observed response is a truncation of a continuous latent variable \cite[\S
  5.5]{agresti2019}. Common examples include the \emph{probit}, or \emph{inverse
    normal}, and the \emph{complementary log-log} links,
  \begin{equation*}
    \probit(\pi) := \Phi^{-1}(\pi) \qquad\text{and}\qquad
    \operatorname{cloglog}(\pi) := \log(-\log(1-\pi)).
  \end{equation*}
\end{example}

The logistic regression model has fewer symmetries than the linear regression
model and their existence depends on the form of the link function.

\begin{proposition}
  Suppose that binary logistic regression models $M$ and $M'$ share a link
  function that is the inverse CDF of a continuous random variable symmetric
  about zero. Then a model homomorphism $\alpha: M \to M'$ is uniquely
  determined by a permutation $\tau: \{0,1\} \to \{0,1\}$ and a matrix
  $B \in \R^{p' \times p}$ such that
  \begin{equation*}
    \sgn(\tau)\, X_{M,i} = X_{M',i}\, B, \quad \forall i=1,\dots,n.
  \end{equation*}
  In particular, $\alpha$ is an isomorphism $M \cong M'$ if and only if $p = p'$
  and $B \in \GL(p,\R)$.
\end{proposition}
\begin{proof}
  A morphism $\alpha: M \to M'$ consists of a matrix
  $B := \alpha_\beta \in \R^{p' \times p}$, a scalar $a := \alpha_\eta \in \R$,
  a convex-linear map $\alpha_\mu: (0,1) \to (0,1)$, and a function
  $\tau := \alpha_y: \{0,1\} \to \{0,1\}$ obeying the naturality conditions for
  the morphisms $g$, $h$, $X_1,\dots,X_n$, and $q$. As in
  \cref{prop:lm-n-morphisms}, naturality for $X_1,\dots,X_n$ are the equations
  $a X_{M,i} = X_{M',i} B$ for $i=1,\dots,n$. The naturality equation for
  $q: \mu \to y$ requires that $\Ber(\alpha_\mu \pi) = \tau \Ber(\pi)$ for all
  $\pi \in (0,1)$. Since the probabilities 0 and 1 are excluded, the function
  $\tau$ must be a permutation, namely the identity or the transposition
  $(0\ 1)$, in which cases $\alpha_\mu$ is the identity or the map
  $\pi \mapsto 1-\pi$. Finally, the naturality conditions for the link and mean
  functions are $a \cdot g_M(\pi) = g(\alpha_\mu \pi)$ and
  $\alpha_\mu(h_M(x)) = h_M(ax)$ for all $\pi \in (0,1)$ and $x \in \R$. When
  $\tau$ and $\alpha_\mu$ are identities, it easily follows that
  $a = 1 = \sgn(\tau)$, so assume that the other case holds.

  So, we must show that $a = -1 = \sgn(\tau)$, or equivalently that
  \begin{equation*}
    -g_M(\pi) = g_M(1-\pi) \qquad\text{and}\qquad 1 - h_M(x) = h_M(-x)
  \end{equation*}
  for all $\pi \in (0,1)$ and $x \in \R$. By assumption, there is a continuous
  symmetric random variable $Z$ such that $h_M(x) = \Prob(Z \leq x)$. Now,
  compute:
  \begin{equation*}
    h_M(-x) = \Prob(Z \leq -x) = \Prob(-Z \leq -x) = \Prob(Z \geq x)
      = 1 - \Prob(Z \leq x) = 1 - h_M(x).
  \end{equation*}
  The other equation then follows by substituting $x = g_M(\pi)$ and using that
  $g_M = h_M^{-1}$.
\end{proof}

The logit and probit links satisfy the symmetry assumption, being derived from
the centered normal and logistic distributions. The complementary log-log link
does not satisfy the assumption and logistic regression models using it are not
symmetric under permutation of the labels. However, a logistic regression model
using the complementary log-log link, $\pi \mapsto \log(-\log(1-\pi))$, is
isomorphic to a \emph{different} logistic regression model using the
\emph{log-log} link, $\pi \mapsto -\log(-\log(\pi))$, via permutation of the
labels.

\begin{example}[Poisson regression] \label{ex:poisson-regression}
  Under the \emph{Poisson regression model} $M \in \Model{\cat{GLM0}_n}$, the
  kernel $M(\mu \xrightarrow{q} y)$ is the Poisson family
  $\PosR \xrightarrow{\Pois} \N$. The canonical link function is the logarithm,
  $g_M(\lambda) := \log(\lambda)$, under which the sampling distribution is
  \begin{equation*}
    y_i \xsim{\mathrm{ind}} \Pois\left(e^{X_{M,i} \beta}\right), \quad
    i=1,\dots,n.
  \end{equation*}
  The Poisson regression model has no symmetries beyond the usual
  reparameterizations of the design matrix.
\end{example}

\begin{example}[Binomial logistic regression]
  \label{ex:binomial-logistic-regression}
  A minor variant of binary logistic regression, \emph{binomial logistic
    regression}, makes each observation a binomial random variable with its own
  sample size. To accommodate this model, the theory $(\cat{GLM0}_n, p_n)$ must
  be extended to a slightly more general theory $(\cat{GLM0}_n', p_n')$, in
  which the morphism $q: \mu \to y$ is replaced by morphisms
  $q_1,\dots,q_n: \mu \to y$ and the sampling morphism
  $p_n: \beta \to y^{\otimes n}$ is changed accordingly to
  $p_n' := \Copy_{\beta,n} \cdot (X_1 h q_1 \otimes \cdots \otimes X_n h q_n)$.
  Under the binomial logistic model $M$ of this theory, each $M(q_i)$ is the
  binomial family $\Binom(m_i,-): (0,1) \to \N$ for a fixed sample size
  $m_i \in \N$. Note that the numbers $m_1,\dots,m_n$ are part of the data of
  the model, not of the theory.

  The functor $\cat{GLM0}_n' \twoheadrightarrow \cat{GLM0}_n$ sending each
  morphism $q_i$ to $q$ projects the extended theory $(\cat{GLM0}_n', p_n')$
  onto the original one $(\cat{GLM}_n, p_n)$. When applied to logistic
  regression models, the model migration functor
  $\Model{\cat{GLM0}_n} \to \Model{\cat{GLM0}_n'}$ interprets a binary logistic
  regression as a binomial logistic regression in which each observation has a
  sample size of one.
\end{example}

Binary and binomial logistic regression generalize from two classes to $k$
classes as \emph{categorical} and \emph{multinomial logistic regression}. As
generalized linear models, they are based on $(k-1)$-dimensional, rather than
one-dimensional, exponential families and they have distinct parameters
$\beta_1, \dots, \beta_{k-1} \in \R^p$ for all but one of the classes. The
theories of a generalized linear model are easily extended to accommodate
$d$-dimensional exponential families, resulting in statistical theories
$(\cat{GLM0}_{n,d}, p_{n,d})$ and $(\cat{GLM0}_{n,d}', p_{n,d})$ for any
dimensions $n$ and $d$. The details are omitted.

\paragraph{Generalized linear models with a dispersion parameter}

Because the distribution, and hence the variance, of a binary-valued random
variable is determined by its mean, the binary logistic regression model does
not have an additional dispersion parameter. But many other generalized linear
models do, and are described by a larger statistical theory.

The \emph{theory of generalized linear model on $n$ observations} has underlying
category $\cat{GLM}_n$ presented by vector spaces $\beta$ and $\eta$, a convex
space $\mu$, a conical space $\phi$, a discrete object $y$, mutually inverse
maps $g: \mu \to \eta$ and $h: \eta \to \mu$, linear maps
$X_1,\dots,X_n: \beta \to \eta$, and a morphism $q: \mu \otimes \phi \to y$. The
sampling morphism $p_n: \beta \otimes \phi \to y^{\otimes n}$ is
\begin{equation*}
  \input{wiring-diagrams/zoo-statistics/glm-n-sampling}
\end{equation*}
The functor $\cat{GLM}_n \twoheadrightarrow \cat{GLM0}_n$ sending the object
$\phi$ to the monoidal unit $I$ and preserving the other generators projects the
new theory $(\cat{GLM}_n, p_n)$ onto the theory $(\cat{GLM0}_n, p_n)$ without a
dispersion parameter. In the other direction, the functor
$H: \cat{GLM0}_n \rightarrowtail \cat{GLM}_n$ sending the morphism
$\mu \xrightarrow{q} y$ to the composite
$\mu \xrightarrow{1_\mu \otimes 0_\phi} \mu \otimes \phi \xrightarrow{q} y$
defines a colax theory morphism
$(H, 1_\beta \otimes 0_\phi, 1_y^{\otimes n}): (\cat{GLM0}_n, p_n) \to
(\cat{GLM}_n, p_n)$. The model migrations induced by this embedding are
sometimes interesting, sometimes not.

A \emph{generalized linear model with a dispersion parameter} is a model
$M: \cat{GLM}_n \to \Stat$ that assigns $M(\beta)$ to be $\R^p$, for some
dimension $p$; $M(\eta) = \R$ and $M(\phi) = \R_+$; and
$M(q): M(\mu) \times \R_+ \to M(y)$ to be a one-dimensional exponential
dispersion family (\cref{ex:exponential-dispersion-family}) or discrete
exponential dispersion family,\footnote{A \emph{discrete exponential dispersion
    model} is an exponential dispersion model rescaled to have support on the
  integers and so be suitable for discrete data. See \cite[\S
  2.4]{jorgensen1987} or \cite[\S 3]{jorgensen1992}.} parameterized by its mean
and by a nonnegative dispersion parameter. The link function, mean function, and
design matrix are defined exactly as before. The model migration functor
$\Model{\cat{GLM0}_n} \to \Model{\cat{GLM}_n}$ induced by the projection
$\cat{GLM}_n \twoheadrightarrow \cat{GLM0}_n$ recovers the generalized linear
models without a dispersion parameter as models of the new theory.

As the name suggests, generalized linear models do indeed generalize linear
models. Define a supply preserving functor
$G_n: \cat{GLM}_n \twoheadrightarrow \cat{LM}_n$ that sends both $\mu$ and
$\eta$ to $\mu$; sends the link and mean morphisms $g$ and $h$ to the identity
$1_\mu$,
\begin{equation*}
  G_n: \input{wiring-diagrams/zoo-statistics/glm-link},
    \input{wiring-diagrams/zoo-statistics/glm-mean}
  \mapsto \input{wiring-diagrams/zoo-statistics/mu-id};
\end{equation*}
sends $\phi$ to $\sigma^2$; and preserves the other generators. Then the model
migration functor $G_n^*: \Model{\cat{LM}_n} \to \Model{\cat{GLM}_n}$ interprets
a linear model as a generalized linear model based on the normal family, which
is a one-dimensional exponential dispersion family under its standard
parameterization by mean and variance.\footnote{Suitably parameterized, the
  normal distribution is also a two-dimensional exponential family, but this is
  a special property of the normal, gamma, and inverse Gaussian families, not
  generally true of one-dimensional exponential dispersion families \cite[\S
  9.4]{sundberg2019}} In this case, applying the further model migration
$\Model{\cat{GLM}_n} \to \Model{\cat{GLM0}_n}$ yields a sampling distribution
with zero variance, that is, a fully deterministic statistical model.

Another example, now for discrete data, is:

\begin{example}[Negative binomial regression]
  Preceding \cref{ex:exponential-dispersion-family}, the negative binomial
  distribution was constructed as a Gamma-Poisson mixture, yielding a more
  flexible alternative to the Poisson distribution for underdispersed or
  overdispersed count data. Negative binomial regression is the corresponding
  extension of Poisson regression from \cref{ex:poisson-regression}. The
  negative binomial family reparameterized by the inverse odds,
  $(\kappa, \theta) \mapsto \NegBinom(\kappa, \theta/(\theta+1))$, or
  equivalently the mixture $\GammaDist \cdot \Pois$, is a discrete exponential
  dispersion family with mean $\mu := \kappa \theta$ and variance
  $\kappa\theta(\theta + 1) = \mu(1 + \mu/\kappa)$.

  In the \emph{negative binomial model} $M \in \Model{\cat{GLM}_n}$, the kernel
  $M(q): \PosR \times \R_+ \to \N$ is this family further reparameterized by its
  mean $\mu$ and dispersion parameter $\phi := 1/\kappa$. For consistency with
  Poisson regression, the standard link function is the logarithm,
  $g_M(\mu) = \log(\mu)$, although it is no longer canonical. The sampling
  distribution is then
  \begin{equation*}
    y_i \xsim{\mathrm{ind}} \NegBinom\left(
      \frac{1}{\phi},\, \frac{\phi}{\phi + e^{-X_{M,i} \beta}}\right), \quad
    i = 1,\dots,n.
  \end{equation*}
  In the limit that $\kappa \to \infty$ or equivalently $\phi \to 0$, with the
  mean $\mu$ held constant, the negative binomial family reduces to the Poisson.
  Consequently, the model migration function
  $\Model{\cat{GLM}_n} \to \Model{\cat{GLM0}_n}$ recovers the Poisson regression
  model from the negative binomial regression model.
\end{example}

\paragraph{Generalized linear mixed models}

Generalized linear models and linear mixed models have a common further
generalization in generalized linear mixed models (GLMMs). The \emph{theory of a
  generalized linear mixed model on $n$ observations}, denoted
$(\cat{GLMM}_n, p_n)$, is presented by vector spaces $\beta$, $b$, and $\eta$; a
convex space $\mu$, conical spaces $\sigma_b^2$ and $\phi$, a discrete object
$y$, mutually inverse maps $g: \mu \to \eta$ and $h: \eta \to \mu$, linear maps
$X_1,\dots,X_n: \beta \to \eta$ and $Z_1,\dots,Z_n: b \to \eta$, a quadratic
morphism $r: \sigma_b^2 \to b$, and a morphism $q: \mu \otimes \phi \to y$. The
sampling morphism $p_n: \beta \otimes \sigma_b^2 \otimes \phi \to y^{\otimes n}$
is
\begin{equation*}
  \input{wiring-diagrams/zoo-statistics/glmm-n-sampling}.
\end{equation*}
A similar theory $(\cat{GLMM0}_n, p_n)$ can be defined for GLMMs without a
dispersion parameter.

A \emph{generalized linear mixed model with a dispersion parameter} is a model
$M \in \Model{\cat{GLMM}_n}$ that assigns $M(\beta)$ and $M(b)$ to be $\R^p$ and
$\R^q$, for some dimensions $p$ and $q$; $M(\eta)$ to be $\R$;
$M(\sigma_b^2) = \PSD^q$ and $M(r)$ to be the $q$-dimensional centered normal
family; and $M(\phi) = \R_+$ and $M(q)$ to be a one-dimensional exponential
dispersion family. The sampling distribution of the model is
\begin{align*}
  b &\sim \Normal_q(0, \Sigma_b) \\
  y_i &\xsim{\mathrm{ind}} M(q)(
    g_M^{-1}(X_{M,i}\, \beta + Z_{M,i}\, b),\, \phi), \quad i=1,\dots,n,
\end{align*}
with parameters $\beta \in \R^p$, $\Sigma_b \succeq 0$, and $\phi \geq 0$.

The relationships between the theories of linear models, linear mixed models,
GLMs, and GLMMs on $n$ observations are summarized by the commutative
diagram:\footnote{For comparison, these relationships are stated in
  conventional, informal style in \cite[Table 1.4]{stroup2012}.}
\begin{equation*}
  \begin{tikzcd}[column sep=0,row sep=small]
    & \underset{\text{generalized linear mixed model}}{(\cat{GLMM}_n, p_n)}
      \arrow[dl, "H_n"']
      \arrow[dr, "K_n" ] & \\
    \underset{\text{linear mixed model}}{(\cat{LMM}_n, p_n)}
      \arrow[dr, "F_n"']
      & & \underset{\text{generalized linear model}}{(\cat{GLM}_n, p_n)}
      \arrow[dl, "G_n"] \\
    & \underset{\text{linear model}}{(\cat{LM}_n, p_n)} &
  \end{tikzcd}
\end{equation*}
The bottom legs of the diagram have already been defined in this section and the
previous one. The top legs are defined similarly. Thus, the functor
$H_n: \cat{GLMM}_n \twoheadrightarrow \cat{LMM}_n$ sends both $\mu$ and $\eta$
to $\mu$ and the link and mean morphisms to the identity $1_\mu$, while the
functor $K_n: \cat{GLMM}_n \twoheadrightarrow \cat{GLM}_n$ sends both $b$ and
$\sigma_b^2$ to the unit $I$, hence $r$ to the trivial morphism $1_I$, and sends
all of $Z_1,\dots,Z_n$ to the zero map $0_\eta: I \to \eta$. The diagram then
commutes, and the induced model migrations functors behave as expected.

\section{Notes and references}

\paragraph{Models in applied statistics}

All of the statistical models treated by this chapter are mainstays of applied
statistics. Davison \cite{davison2003} and Efron and Hastie \cite{efron2016}
give broad surveys of statistical modeling and inference, describing most of the
models here and much else besides. Linear models are a topic in nearly every
introductory text on statistics or machine learning. A thorough theoretical
treatment is given by Seber and Lee \cite{seber2003} and a more geometric
perspective by J{\o}rgensen \cite{jorgensen1993}. Textbooks on hierarchical
linear models and linear mixed models include \cite{gelman2006a} and
\cite{pinheiro2000}. Generalized linear models were invented by Nelder and
Wedderburn to unify a number of commonly used statistical models
\cite{nelder1972}. The standard reference is by McCullagh and Nelder
\cite{mccullagh1989}. Exponential families and exponential dispersion families,
the probabilistic underpinning of generalized linear models, are described from
a theoretical perspective in \cite{brown1986,jorgensen1997} and a pragmatic one
in \cite{sundberg2019}. A textbook on generalized linear mixed models and their
special cases is by Stroup \cite{stroup2012}.

The selection of Bayesian statistical theories and models here is regrettably
limited. Textbooks on applied Bayesian statistics, such as by Hoff or by Gelman
et al \cite{hoff2009,gelman2013}, present a much wider range of Bayesian
models.

\paragraph{Symmetry and statistical models}

Invariance and equivariance under symmetry are well explored in statistical
decision theory as a principle for a constraining the acceptable estimators so
as to single out a unique optimal estimator
\cite{ferguson1967,lehmann1998,lehmann2005,eaton1989,wijsman1990}. Subtle
differences between authors notwithstanding,\footnote{Compare the definitions of
  an invariant family of probability distributions in \cite[Definition 3.1 and
  Theorem 3.1]{eaton1989}, \cite[Definition 4.1.1]{ferguson1967}, and
  \cite[Definition 3.2.1]{lehmann1998}.} the generally agreed-upon paradigm is
to
\begin{enumerate}[nosep]
\item define a group action on the sample space of the probability model,
\item use identifiability of the model to transport the group action from the
  sample space to the parameter space,
\item use ``identifiability'' or ``discriminant-ability'' of the loss function
  to transport the group action from the parameter space to the action space,
\item and finally define invariant loss functions and equivariant decision rules
  using these three group actions.
\end{enumerate}
This convoluted logic is needed because the classical works do not possess the
concept of a natural transformation, which formalizes what it means for a family
of transformations of a multi-sorted mathematical structure to be compatible
with each other.

The homomorphisms of a statistical model, that is, of a model of a statistical
theory, offer a concept of symmetry that is superior to the classical one on
several grounds. It does not assume that the model is identifiable. It is not
restricted to automorphisms of a single model, but applies equally to
isomorphisms of different models and even to non-invertible morphisms between
models. It ensures that the transformations are compatible with the full
structure of the model, not just the sample and parameter spaces. And above all
else, it affirms the basic principle of algebra and logic that the homomorphisms
of a mathematical structure are not extra data to be added arbitrarily, but are
determined by the structure's axiomatization. Formulating a statistical model as
a model of a statistical theory completely determines its symmetries. If a
smaller or larger group of symmetries is desired, they may be obtained by
passing to a stronger or weaker statistical theory via theory morphisms and
their model migrations.

In the case where the linear model is identifiable (the design matrix is full
rank) and the model homomorphism is an isomorphism, the equivariance of OLS
linear regression (\cref{thm:linear-regression-equivariance}) is well known.
Terms such as ``regression equivariance,'' ``scale equivariance,'' and ``affine
equivariance'' appear commonly in the literature on robust regression \cite[\S
3.4]{rousseeuw1987}. However, it is difficult to find a statement of
equivariance that would correspond to a general, non-invertible morphism of
linear models, since such transformations are not usually considered.

\chapter{Computer program analysis of data science code}
\label{ch:program-analysis}

Computer programs, though precise enough to be unambiguously executed, are
primarily written by humans in order to be intelligible to humans. Whether
expressed in their native textual format, or in the parsed format of abstract
syntax trees, computer programs are for many purposes nearly as unintelligible
to machines as is natural language text. In particular, machines do not readily
create, introspect, or manipulate code in intelligent ways.

The difficulty of computer program understanding has several sources, some more
fundamental than others. As a practical matter, the programming languages in
common use tend to be complex. Popular languages offer many features and
conveniences, such as classes, special operators, and flavors of ``syntactic
sugar,'' that are appreciated by programmers but go considerably beyond the
simple models of computation studied by theoretical computer scientists. This is
certainly true of Python and, to a lesser extent, R. Moreover, no matter how
simple or complex the language, important properties of a program are rarely
immediate from its syntactic presentation but must be inferred through
\emph{computer program analysis}. Such properties include control flow and data
flow, opportunities for parallelism, and guarantees about correctness and
termination. Foundational results, like the famous unsolvability of halting
problem, ensure that many program properties are generally impossible to decide
algorithmically, without executing the program. Finally, when viewed as human
artifacts, computer programs embody abstract concepts and domain knowledge and
so possess a meaning going beyond what is normally understood as program
semantics in programming language theory. Intelligent inspection and
manipulation of a program requires some understanding of its meaning in this
more nebulous sense.

This chapter concerns the first two problems, of managing programming language
complexity and inferring program properties, particularly data flow. It
describes the design and implementation of a software system for transforming
programs in the Python and R languages into data-flow diagrams, introduced as
\emph{raw flow graphs} in \cref{sec:introduction-semantic-enrichment}. While
assuming nothing specifically about the subject matter of the code, the system
comprises the first half of a method for the semantic analysis of data science
code. The second half, the topic of the next chapter, introduces concepts
specific to data science.

In a departure from the rest of the dissertation, the program analysis
methodology is not presented in a formal or mathematically rigorous style. This
appears to be an inevitable feature of the domain. Like many practical
programming languages, the Python language is too complex and ambiguously
specified to admit a useful formal semantics.\footnote{In his master's thesis,
  Guth gives a large, but incomplete, operational semantics of Python 3.3
  \cite{guth2013}.} The R language is simpler than Python, but still quite
complex compared to idealized mathematical models of
computation.\footnote{Operational semantics of ``core R,'' a proper subset of
  the R language, is given in \cite{morandat2012}.} No attempt is made to
formalize the process of abstraction and simplification by which Python and R
are reduced to an idealized programming language. Indeed, the main purpose of
this process to create a computational representation that is easy reason
precisely about and manipulate.

\section{Recording flow graphs using program analysis}
\label{sec:program-analysis-general}

Raw flow graphs, or in this chapter simply \emph{flow graphs}, model the flow of
data during the execution of a computer program. Several small examples have
been given in
\cref{lst:toy-kmeans-scipy,lst:toy-kmeans-sklearn,lst:toy-kmeans-r} and
\cref{fig:toy-kmeans-scipy,fig:toy-kmeans-sklearn,fig:toy-kmeans-r} in the
Introduction. Another such example, of fitting a linear regression model using
Scikit-learn, is shown in \cref{lst:toy-regression-sklearn} and
\cref{fig:toy-regression-sklearn-raw}. The correspondence between the source
listing and the flow graph should be mostly self-evident.\footnote{The names
  \verb|read_csv| in the code and \verb|_make_parser_function.<locals>.parser_f|
  in the flow graph do not match because the \verb|read_csv| and
  \verb|read_table| functions in Pandas are dynamically generated from an
  internal higher-order function called \verb|_make_parser_function|. The
  obscure name in the flow graph is the function's ``true name'' according to
  the Python interpreter.}

\begin{listing}
  \begin{tcolorbox}
    \inputminted{python}{code/regression-sklearn-metrics.py}
  \end{tcolorbox}
  \caption{Linear regression in Python, using Pandas and Scikit-learn}
  \label{lst:toy-regression-sklearn}
\end{listing}

\begin{figure}
  \centering
  \input{wiring-diagrams/program-analysis/regression-sklearn-metrics.raw}
  \caption{Raw flow graph for \cref{lst:toy-regression-sklearn}}
  \label{fig:toy-regression-sklearn-raw}
\end{figure}

Flow graphs are \emph{wiring diagrams}, also known as \emph{string diagrams}.
Formally, they represent morphisms in a cartesian closed category, a concept
partially introduced in \cref{sec:cartesian-categories} and elaborated on in the
next chapter. From a programmer's perspective, flow graphs are an idealized
model of typed functional programming, in which functions can be defined and
composed in sequence and in parallel. Data is permitted to be duplicated or
discarded.\footnote{In wiring diagrams, the duplication and discarding
  operations are represented implicitly, as ports having multiple or zero
  incident wires, rather than explicitly, as the filled circles used in
  \cref{ch:category-theory,ch:algebra-statistics,ch:zoo-statistics}.} The type
system has basic or primitive types, product types, a unit or singleton type,
and function types.

The wires in a flow graph represent typed objects and the boxes represent
function calls in the target programming language. The types can be primitive,
such as \verb|float| or \verb|str| in Python and \verb|numeric| or
\verb|character| in R, or externally defined classes, such as NumPy's core
datatype \verb|ndarray| or R's linear model type \verb|lm|. As for the boxes,
the phrase ``function call'' is used liberally to encompass essentially any
user-invoked computational action. So, in a typical object-oriented programming
language like Python, the ``functions'' include not just standalone functions,
but also static methods, instance methods, object attribute accessors, container
indexing, and special unary and binary operators. This understanding of
functions is more closely aligned with the R programming language, with its
slogan that ``everything that happens in R is a function call''
\cite{chambers2016}.

The flow graph is constructed by two stages of program analysis, the first
static and the second dynamic. Program analysis can in general be static or
dynamic or both. A \emph{static} program analysis inspects or transforms the
program without executing it, while a \emph{dynamic} program analysis involves
executing the program. In our system, the static phase transforms the original
program to emit events whenever user-level code performs certain actions, and
then the dynamic phase executes the modified program, assembling the flow graph
based on the observed events. Because both phases are aided by tight integration
with the language interpreter, the Python program analysis system is implemented
in Python and the R system in R. Nevertheless, the two systems are
architecturally similar.

Let us now consider the static and dynamic phases in greater detail.

\paragraph{Static program analysis}

The static phase of program analysis augments the source code with additional
instructions to record whenever
\begin{enumerate}[(i),nosep]
\item a function call is about to begin,
\item a function call is returned from,
\item a variable is accessed,
\item a variable in the local scope is assigned, or
\item a variable in the local scope is deleted.\footnote{Deletion is performed
    by the \verb|del| keyword in Python and the \verb|rm| function in R. Since
    these languages have garbage collectors, explicit deletion is rarely
    encountered in practice.}
\end{enumerate}
The code is transformed by compiling it to an \emph{abstract syntax tree (AST)}
using the language's built-in parser, walking the syntax tree to find
occurrences of function calls and variable uses, and modifying those nodes to
make special callbacks. In Python, the source code is transformed in a single
pass. The resulting AST can then be pretty-printed as new source code or, more
usefully, simply executed directly. In R, the code transformations are made
on-the-fly as the programs executes, which is possible due to R's unique
combination of lazy evaluation and dynamic metaprogramming.

Consider, for concreteness, the rewriting of a function call. In Python, the
function call \verb|f(x,y)| is transformed into the expression
\begin{minted}[autogobble,xleftmargin=0.5in]{python}
  __tracer.trace_return(__tracer.trace_function(f)(
      __tracer.trace_argument(x),
      __tracer.trace_argument(y)
  ))
\end{minted}
where \verb|__tracer| refers to a hidden object that will be injected into the
code's global namespace at runtime. Each of the tracer's methods returns its
argument unmodified after suitable processing. Thus, due to Python's evaluation
order for function calls, the transformed code will
\begin{enumerate}[nosep]
\item evaluate the name \verb|f|, yielding a function or other callable object,
  then call \verb|trace_function| with this value,
\item evaluate the name \verb|x|, then call \verb|trace_argument| with this
  value,
\item evaluate the name \verb|y|, then call \verb|trace_argument| with this
  value, and finally
\item call the object \verb|f| with the objects \verb|x| and \verb|y| as
  arguments, then call \verb|trace_return| with the return value.
\end{enumerate}
The effect is to evaluate the expression \verb|f(x,y)| while capturing enough
information to reconstruct what function was called, what were its arguments,
and what value it returned.

One might wonder why a simpler transformation is not performed instead. For
example, the function call \verb|f(x,y)| might be rewritten as
\verb|trace_call(f,[x,y])| or \verb|trace_call(f,x,y)|, where the function
\verb|trace_call| would evaluate \verb|f| with arguments \verb|x| and \verb|y|,
record all relevant information about the function call, and return the computed
value. However, this can have unwanted effects when the function \verb|f| being
evaluated is not \emph{referentially transparent}, that is, when it depends not
just on its inputs but on some property of its calling environment. A correct
code transformation preserves the call stack of the program by evaluating the
function \verb|f| in its original context, whereas the simpler transformations
suggested here do not. Obscure though it may seem, this issue does arise in
practice. For instance, the Python package Patsy \cite{patsy2018}, which
provides R-style model formulas for the popular package Statsmodels
\cite{statsmodels2010}, inspects the call stack to circumvent Python's lack of
metaprogramming facilities. It will break if the call stack is not preserved.

The prime directive of any program transformation, whether for optimization or
introspection, is to not alter the observed behavior of the program. This can be
surprisingly subtle, as the case above illustrates. Another difficult case,
non-standard evaluation in R, will be considered later.

\paragraph{Dynamic program analysis}

The dynamic phase of program analysis executes the transformed program,
gradually constructing the flow graph as callbacks occur. The sequence of
function calls, defining the boxes of the diagram, is immediate from the
callbacks, but the data flow between function calls, defining the wires of the
diagram, must also be recorded. Data is generally passed between functions
either directly, from one function call to another by function composition, or
indirectly, though variable assignments and accesses. Extra bookkeeping is
needed to track this data flow.

The program analysis system tracks the flow of data both across the call stack
and within a single level of the call stack. An empty wiring diagram is created
when the program begins. Whenever a function is about to be called, if the
function is defined within the program being analyzed, then a new wiring diagram
is created, pushed onto a stack maintained by the system, and the recording
process is restarted recursively. Otherwise, the function must be defined by an
external library or by the language itself and is treated as an \emph{atomic},
or indecomposable, computation. In either case, a new box is added to the outer
wiring diagram,\footnote{As a technical note, in Python, the new box can be
  created and connected \emph{before} the function is called, whereas in R this
  must happen \emph{after} the function returns, because the lazy evaluation of
  function arguments defers their availability for introspection.} the box being
a nested wiring diagram in the first case and an atomic box in the second. The
hierarchy of the wiring diagram thus mirrors the call stack of the program
throughout its execution. Once a new box has been added, the box's inputs are
wired either to outputs of previously created boxes or to inputs of the outer
box, according to the provenance of the arguments to the function call.

At each level of the call stack, the provenance of the objects in scope is
maintained by associating program events and objects with their sources in the
flow graph. For local variables, a lookup table maps each variable name to an
output port of the box representing the function call that created or most
recently mutated the variable's value. The five callbacks made by the
transformed program are then handled by:
\begin{enumerate}[(i),nosep]
\item when a function call is about to begin, create a new box for the function
  call, as described above;
\item when a function call is returned from, pass the output ports of its box as
  sources to the enclosing expression;
\item when a variable is accessed, look up variable's source in the lookup table
  and pass the source to the enclosing expression;
\item when a variable is assigned, add or replace that entry in the lookup table
  with a new source;
\item when a variable is deleted, remove that entry from the
  lookup table.
\end{enumerate}
The passing of extra information in cases (ii) and (iii) can be performed
statically, by boxing the values that pass between two statically transformed
expressions, or dynamically, by pushing them onto a queue to be emptied later.

This concludes a technical overview of the program analysis systems for Python
and R. Language-specific features and challenges are discussed in subsequent
sections, but let us first consider the general tradeoffs made by our
methodology.

\paragraph{Static versus dynamic analysis}

Although it has a static phase, the program analysis system is essentially
dynamic, for several reasons. Static analysis, especially about type inference,
is challenging for dynamically typed languages like Python and R. Moreover, in
applications to data science, it is of interest to capture values computed
during the course of the program’s execution, such as parameter estimates or the
selected hyperparameters. For this dynamic analysis is indispensable. In
general, dynamic analysis is easier to implement that static analysis, as it
skirts the computationally intractable or even undecidable problems that static
analysis easily produces \cite{landi1992}.

Of course, a disadvantage of dynamic analysis is the necessity of running the
program. Crucially, our system needs not just the code itself, but its input
data and runtime environment. These are all requirements of scientific
reproducibility, so in principle they ought to be satisfied. In practice they
are often neglected. Even assuming that all the resources are available, a
significant investment in infrastructure and curation would be needed to process
data analyses at a large scale. Thus, for metascientific purposes, it is
worthwhile to pursue purely static program analyses in addition to dynamic ones,
acknowledging the trade-offs involved in both cases.

\paragraph{User versus library code}

The program analysis system transforms and records only the code that is
directly presented to the system as input. For now, the input program is assumed
to be a single script or Jupyter notebook, as is common in data
analysis.\footnote{That said, complex data analyses may easily span multiple
  scripts or notebooks connected by intermediate files and databases. A more
  complete system would track the provenance of data both within a single
  program, as our system does, and between programs, as done by conventional
  data provenance tools \cite{simmhan2005a}.} Any code that is executed
indirectly through library or system functions is thus not recorded. This
restriction accords with the basic assumption, first stated in
\cref{ch:introduction}, that library functions are semantically meaningful
computational units.

Limiting the scope of the program analysis also reduces its runtime overhead.
For a typical data analysis script, the vast majority of its running time is
spent inside imported Python or R functions or inside the underlying C and
Fortran routine. These computations are not affected by the program analysis
system.

\paragraph{Programs versus program executions}

The presentation so far has not carefully distinguished between a computer
program and a particular \emph{execution} of that program. In general, a program
having nontrivial control flow, through conditional branching (\verb|if|
statements), looping (\verb|for| and \verb|while| statements), or recursion,
will execute a different sequence of instructions depending on the input to the
program. At present, our system ignores control flow, capturing only the
function calls made during the execution of the program, not those that might
have been made on different input data. Thus, faced with a conditional branch,
the system records only the executed branch, and faced with a loop, the system
records the unrolled sequence of iterations.

This limitation, while important, is not as severe as might initially be
supposed. Most programs, such as a graphical application or a web server, are
expected to run on many different user inputs and behave differently each time.
In contrast, data analyses are usually created for and attached to specific
datasets. As a requirement of scientific reproducibility, we expect that if the
program is run repeatedly on the same data, it will produce the same result each
time. It is therefore reasonable to ask what actually did happen on given data,
rather than what might have happened on different data. Even so, unrolling a
loop with many iterations causes a blow-up in the size of the flow graph, which
is highly impractical. Addressing this problem is left to future work.

\section{Flow graphs for Python}
\label{sec:program-analysis-python}

A clean and simple syntax belies the surprising complexity of the Python
programming language. Advertised as a multi-paradigm language, it supports
imperative, functional, and object-oriented programming styles and includes
numerous advanced features, such as generators, list comprehensions, and, in
recent versions, coroutines. Its combination of features has made Python into a
popular general-purpose programming language but causes difficulties when
translating Python code into a simpler model of computation, like that of flow
graphs. Our program analysis system resolves some but not all of these
difficulties.

\begin{listing}
  \begin{tcolorbox}
    \inputminted[fontsize=\footnotesize]{python}{code/plot_svm_regression.py}
  \end{tcolorbox}
  \caption{Support vector regression (SVR) in Scikit-learn. The code is
    reproduced without changes from Scikit-learn's official suite of examples.}
  \label{lst:plot-svm-regression}
\end{listing}

\begin{figure}
  \centering
  \resizebox{\textwidth}{!}{%
    \input{wiring-diagrams/program-analysis/plot_svm_regression.raw}}
  \caption{Raw flow graph for \cref{lst:plot-svm-regression}}
  \label{fig:plot-svm-regression-raw}
\end{figure}

\paragraph{Homogenizing the syntax}

Many instructions in Python, including access and assignment of object
attributes, indexed access and assignment of containers, and special unary and
binary operators, can be interpreted as function calls yet are not officially
regarded as such by the Python grammar. During the static phase of program
analysis, such instructions are reduced to function calls so that they may be
recorded by the same methods. Conveniently for this purpose, the standard
library module \verb|operator| provides functional aliases for nearly every
special operator in the Python language.

For example, the script in \cref{lst:plot-svm-regression}, reproduced verbatim
from the Scikit-learn documentation, uses the multiplication (\verb|*|),
subtraction (\verb|-|), array indexing and slicing (\verb|[::]|), and in-place
addition (\verb|+=|) operators. In the corresponding flow graph of
\cref{fig:plot-svm-regression-raw}, these instructions appear as calls to the
\verb|mul|, \verb|sub|, \verb|getitem|, and \verb|iadd| functions from the
\verb|operator| module, as well as the built-in function \verb|slice|.

\paragraph{Object-oriented programming}

Although the foundational NumPy and SciPy packages have mainly procedural
interfaces, many Python packages for data science, such as Pandas and
Scikit-learn, are written in an object-oriented style. Some packages, such as
Matplotlib and Statsmodels, provide both procedural and object-oriented
interfaces. In comparison, flow graphs are typed, even allowing a form of
subtyping through implicit conversion (\cref{sec:concepts-as-category}), but
they are not object-oriented in the sense of class-based programming, where
objects inherit attributes and methods through classes.

To accommodate class-based objects, the program analysis system reinterprets
calls of instance methods as calls of functions having an extra first argument
for the object instance. This protocol actually agrees with Python's own syntax
for defining an instance method, where the instance object is represented by a
first argument conventionally called \verb|self|. For example, the \verb|fit|
method of a supervised model in Scikit-learn, seen in
\cref{lst:plot-svm-regression}, has the signature \verb|fit(self, X, y)|, which
appears as a box with three input ports in \cref{fig:plot-svm-regression-raw}.
Object attribute accesses and assignments are translated into \verb|getattr| and
\verb|setattr| calls, using the method described previously, but for readability
are displayed as boxes labeled by the name of the attribute.

\paragraph{Side effects}

The most severe difficulties in analyzing Python code arise from the mismatch
between the programming models of Python, which permits mutation, global state,
and other side effects, and of flow graphs, which is purely functional. To these
problems only partial solutions are available. When it is known that a function
or method mutates its inputs, it may be reinterpreted as a non-mutating function
with an additional output for every mutated input. For example, the \verb|fit|
methods in Scikit-learn return a fitted model, which is in fact the original,
unfitted model having been mutated. Regardless of whether the return value is
consumed (compare \cref{lst:toy-regression-sklearn,lst:plot-svm-regression}),
function calls to \verb|fit| become boxes with an output for the fitted model.

That is simple enough; the trouble is knowing \emph{when} a function is
mutating. Functions corresponding to certain special operations, such as
\verb|setattr| and \verb|setitem| for attribute and indexed assignment, can
always be assumed to be mutating. When they are available, function annotations
can also be used to manually mark a function as mutating (\cref{sec:dso}). But
in general the system cannot determine whether an arbitrary function is
mutating, since the mutation can occur in library code that is not statically
analyzed. Moreover, mutations can be \emph{implicit}. For example, if a column
in a data frame is mutated through a reference to that column, then the
containing data frame should be regarded as mutated as well. Implicit mutations
are currently not addressed.

Mutation of global state is a still worse problem for the functional paradigm.
As an example, \cref{lst:plot-svm-regression} generates a plot using
Matplotlib's imperative interface, which maintains global state for the active
plot. The calls made to \verb|xlabel|, \verb|ylabel|, \verb|title|,
\verb|legend|, and \verb|show| all then appear as isolated boxes in the flow
graph of \cref{fig:plot-svm-regression-raw}. From the functional perspective,
this sequence of functional calls should be reinterpreted as a chain explicitly
passing a plot object from one call to the next. However, adopting this
convention in the flow graph would require a more drastic rewriting of the
observed program events than any considered here.

\paragraph{Weak references}

For objects with a unique identity, weak references offer a dynamic alternative
to static program transformation for tracking the data flow of objects. A
\emph{weak reference} is a reference, or pointer, to an object that does not
increase the object's reference count and hence does not prevent the object from
being garbage collected. The program analysis system may therefore freely use
weak references without causing memory leaks. Most objects in Python can be
weakly referenced. Important exceptions are the primitive objects, such as
booleans, integers, and strings, and the built-in containers for lists, tuples,
and dictionaries.

Wherever possible, the Python program analysis system maintains weak references
to the objects returned by recorded function calls, allowing the objects to be
reidentified if they appear as arguments to future function calls. This method
is complementary to the default, static approach.\footnote{The original
  implementation of the Python program analysis, described in
  \cite{patterson2017ibmjrd}, relied exclusively on weak references; however,
  this is plainly inadequate.} It is not applicable to all objects, yet where it
is applicable, it provides a strong guarantee of object identity, immune to any
dynamic trickery that might fool a static analysis.

\section{Flow graphs for R}
\label{sec:program-analysis-r}

The R programming language is, in most respects, far simpler than Python.
Despite appearances, R is quite close to being a functional language. Nearly all
of its standard functions and data structures have \emph{copy-on-modify}
semantics, meaning that ostensibly mutating operations in fact copy the data
before modifying it. The only mutable data structure in base R is the
environment type, which does not figure explicitly in most data analysis
scripts. The default and most commonly used object system, known as S3, is based
on generic functions, not classes, making it a natural fit for flow graphs. The
R language also has a very simple abstract syntax, with expressions, or abstract
syntax trees, composed of only four types of nodes:\footnote{In R, all control
  flow constructs, variable assignment operators, and unary and binary operators
  are special types of functions.} scalar constants, names, function calls, and
``pairlists.'' The first three node types are accommodated by the general
methodology of \cref{sec:program-analysis-general}, while the last is a special
type for function arguments that does not appear in ordinary code.

Before discussing a few distinctive aspects of the R language, let us examine a
data analysis written in R that is more realistic than the preceding examples. A
recent DREAM Challenge \cite{saez-rodriguez2016} asked how well clinical and
genetic covariates predict patient response to anti-TNF treatment for rheumatoid
arthritis \cite{sieberts2016}. Each team of analysts was instructed to submit
two models, one using only genetic covariates and the other using any
combination of clinical and genetic covariates.

An analysis submitted by a top-ranking team \cite{kramer2014} is displayed in
\cref{lst:dream-ra}, having been lightly modified for portability. Its flow
graph is shown in \cref{fig:dream-ra-raw}. The analysts fit two predictive
models, the first using both genetic and clinical covariates and the second
using only clinical covariates. Both models use the Cubist regression algorithm
\cite[\S 8.7]{kuhn2013}, a variant of random forests based on M5 regression
model trees \cite{wang1997}. Because the genetic data is high-dimensional, the
first model is fit using the subset of the genetic covariates selected by VIF
regression, a variable selection algorithm \cite{lin2011}. The linear regression
model created by VIF regression is used only for variable selection, not for
prediction. If this code is found to be opaque, the semantic flow graph
displayed in \cref{fig:dream-ra-semantic} of the following chapter can also be
consulted.

\begin{listing}
  \begin{tcolorbox}
    \inputminted[fontsize=\footnotesize]{r}{code/dream-ra.R}
  \end{tcolorbox}
  \caption{R code for two models from the Rheumatoid Arthritis DREAM Challenge
    \cite{sieberts2016}. The code is reproduced without essential changes from
    \cite{kramer2014}.}
  \label{lst:dream-ra}
\end{listing}

\begin{figure}
  \centering
  \resizebox{\textwidth}{!}{%
    \input{wiring-diagrams/program-analysis/dream-ra.raw}}
  \caption{Raw flow graph for \cref{lst:dream-ra}}
  \label{fig:dream-ra-raw}
\end{figure}

Although R conforms fairly neatly to the formalism of flow graphs, the language
has two distinctive features, lazy evaluation and non-standard evaluation, that
must be accounted for when transforming programs.

\paragraph{Laziness} The R language has \emph{call-by-need} semantics, under
which arguments to function calls are evaluated lazily. That is, when a function
is called, it does not receive previously computed values as arguments but
rather \emph{promises} to compute these values. A promise object is evaluated
only when it is accessed for the first time within the body of the function. In
particular, if a promise object is never accessed, then it is never evaluated.
Also, the order in which a function's arguments are evaluated depends on the
body of the function.

Since R functions can have side effects, a valid program transformation
generally cannot force promises to be evaluated earlier than they would
otherwise be. The R program analysis system meets this requirement through
non-standard evaluation, rewriting the code in the promise before the function
is called, without evaluating the promise. The evaluation order is thus exactly
the same as in the original program, which is not necessarily that given for
Python programs in \cref{sec:program-analysis-general}.

\paragraph{Non-standard evaluation}

Although R may be the world's most widely used lazy functional programming
language, most end-users of R do not explicitly rely on laziness in their code
and may not even be aware that R is a lazy language \cite{goel2019}. One
consequence of laziness in R is to make possible a dynamic form of
metaprogramming. In \emph{non-standard evaluation}, the interpretation of a
function argument having a fixed syntactical form can depend on the internal
logic of the function. Thus, in general, the semantics of R code cannot be known
statically. In practice, non-standard evaluation is often only a convenience,
used to abbreviate code and capture variable names for plot labels. Such uses
have no significant impact on the program semantics. However, the tidyverse
suite of packages makes much more extensive use of non-standard evaluation,
effectively implementing its own domain-specific language for data processing
within R \cite{wickham2019}.

Non-standard evaluation is the most serious obstacle to the program analysis of
R code because the meaning of an expression can depend arbitrarily on the
context in which it is evaluated. There appears to be no simple and generic
solution to this problem. At present, our program analysis system handles, on a
case-by-case basis, certain recurring patterns of non-standard evaluation in
base R. However, no attempt has been made to systematically accommodate all the
major uses of non-standard evaluation in the R ecosystem.

\section{Notes and references}

An early version of the Python program analysis system is described in previous
work by the author and collaborators \cite{patterson2017ibmjrd}. Significant
design and implementation changes have since been made to improve the robustness
of the system. Most importantly, the use of Python's system trace function
(\verb|sys.settrace|), which cannot trace built-in functions or C extensions,
has been eliminated in favor of more difficult to implement, but more reliable,
static program transformations
(\cref{sec:program-analysis-general,sec:program-analysis-python}). Preliminary
support for the R programming language has also been added
(\cref{sec:program-analysis-r}). The R program analysis system was announced,
but not described in any detail, in a more recent publication
\cite{patterson2018}.

\paragraph{Program analysis}

A general reference on computer program analysis, with a more theoretical
orientation than this chapter, is the text by Nielson et al \cite{nielson1999}.
As the book's contents suggest, a large part of the literature on program
analysis is about static analysis, with the main application being to optimizing
compilers \cite{aho2006}. A recent survey of dynamic program analysis is
\cite{gosain2015}. Most applications of dynamic analysis are to debugging and
testing, performance profiling, and security analysis.

Data flow analyses are also classified as \emph{intraprocedural} (within a
single procedure) or \emph{interprocedural} (between procedures in a larger
program) \cites[Chapter 2]{nielson1999}{khedker2009}. The program analysis of
this chapter is therefore a dynamic, interprocedural data flow analysis. A
common intermediate representation in an interprocedural data flow analysis is a
\emph{call graph} (or \emph{call multigraph}), where the vertices correspond to
functions and there is a directed edge from vertex $u$ to $v$ for each call of
the function $v$ from within $u$. Call graphs do not contain the information
about data dependency that a flow graph does. Also, a flow graph is not,
strictly speaking, a graph but a wiring diagram.

By building a custom version of the R interpreter, the authors of
\cite{morandat2012} implemented the \emph{TrackeR} system to capture detailed
execution traces of R programs. The dynamic program analysis here is
comparatively light-weight, being based on the standard R interpreter but
running transformed R programs.

\paragraph{Data provenance} The aim of recording the steps of a data analysis is
shared by the field of data provenance. The \emph{provenance} of a data resource
includes its origin and the process of transformation by which it was derived
\cite{simmhan2005a}. Insofar as the survey \cite{simmhan2005b} is
representative, the main difference between our system and the typical data
provenance system seems to be granularity. In data provenance, the finest
granularity of data resource is often files or database records, whereas our
system operates on a single file and traces arbitrary program objects. For
example, the StarFlow system targets data analyses written in Python but
operates at the level of scripts \cite{angelino2010}. The program analysis here
would be naturally complemented by a file-level provenance system to describe
data provenance at multiple levels of granularity.

\chapter{Semantic enrichment of data science code}
\label{ch:semantic-enrichment}

Compared to the programs from which they are derived, the raw flow graphs
generated by the program analysis system are more readily inspected, reasoned
about, and manipulated by machines. However, each flow graph is still expressed
in the vernacular of a specific programming language and set of libraries. This
state of affairs is unsatisfactory because the computations made by data
analyses are rarely inseparable from the idioms of particular languages and
packages; rather, the data analysis and its software dependencies are thought to
\emph{instantiate}, in a concrete, computable form, the abstract mathematical,
statistical, and algorithmic concepts of data science. The detailed way in which
this instantiation happens is in many respects arbitrary, depending on
conventions adopted by the authors and the community at large.

Semantic enrichment aims to bring out the abstract concepts that underlie data
science code, in a manner that is transparent to machines and insensitive to
implementation details. Specifically, semantic enrichment is the process of
transforming the raw flow graphs from the previous chapter into \emph{semantic
  flow graphs}, whose types and functions belong to the controlled vocabulary of
the Data Science Ontology. This process was outlined in the Introduction, where
the semantic flow graph corresponding to all three of
\cref{fig:toy-kmeans-scipy,fig:toy-kmeans-sklearn,fig:toy-kmeans-r} was shown in
\cref{fig:toy-kmeans-semantic}. Revisiting another small example,
\cref{fig:toy-regression-sklearn-semantic} below shows the semantic flow graph
derived from the previous chapter's raw flow graph of
\cref{fig:toy-regression-sklearn-raw}. Another, more realistic example will be
shown later.

\begin{figure}
  \centering
  \input{wiring-diagrams/semantic-enrichment/regression-sklearn-metrics.semantic}
  \caption{Semantic flow graph for \cref{lst:toy-regression-sklearn} and
    \cref{fig:toy-regression-sklearn-raw}}
  \label{fig:toy-regression-sklearn-semantic}
\end{figure}

In this chapter, the Data Science Ontology and its use in semantic enrichment
are introduced and informally described. The ontology and the semantic
enrichment process are then formalized mathematically. The concepts in the
ontology form a cartesian closed category with implicit conversion, while the
ontology's code annotations partially define a functor between two such
categories. Finally, the raw and semantic flow graphs are formalized as
morphisms in categories of elements.

\section{The Data Science Ontology}
\label{sec:dso}

The Data Science Ontology is a nascent knowledge base about statistics, machine
learning, and data processing. It aims to support automated reasoning about data
science software.

The ontology is comprised of concepts and annotations. \emph{Concepts} catalog
the abstract entities of data science, such as data tables and statistical
models, as well as the processes that manipulate them, such as loading data from
a file or fitting a model to data. Reflecting the intuitive distinction between
``things'' and ``processes,'' concepts bifurcate into two kinds: types and
functions. The terminology agrees with that of functional programming. Thus, a
\emph{type} represents a kind or species of thing in the domain of data science.
A \emph{function} is a functional relation or mapping from an input type, the
\emph{domain}, to an output type, the \emph{codomain}. In this terminology, the
concepts of a data table and of a statistical model are types, whereas the
concept of fitting a predictive model is a function that maps an unfitted
predictive model, together with predictors and response data, to a fitted
predictive model.

As a modeling assumption, software packages for data science, such as Pandas and
Scikit-learn, are regarded as instantiating the concepts of the ontology.
\emph{Annotations} say how this instantiation occurs by mapping types and
functions in software packages onto type and function concepts in the ontology.
To avoid confusion between levels of abstraction, we call the former
``concrete'' and the latter ``abstract.'' Thus, a type annotation maps a
concrete type---a primitive type or user-defined class in Python or R---onto an
abstract type---a type concept. Likewise, a function annotation maps a concrete
function onto an abstract function. As in \cref{ch:program-analysis}, concrete
function are construed in the broadest possible sense to include any programming
language construct that ``does something'': standalone functions, instance
methods, attribute getters and setters, and so on.

The division of the ontology into concepts and annotations on the one hand, and
into types and functions on the other, leads to a two-way classification.
\cref{tbl:ontology-classification} lists several examples of each of the four
combinations, drawn from the Data Science Ontology.

\begin{table}
  \centering
  \caption{Example concepts and annotations from the Data Science Ontology}
  \label{tbl:ontology-classification}
  \begin{tabular}{lp{2in}p{2in}}
    & \textbf{Concept} & \textbf{Annotation} \\ \toprule
    \textbf{Type} &
      data table &
      pandas data frame \\ \cmidrule{2-3}
    & statistical model &
      scikit-learn estimator \\ \cmidrule{2-3}
    \textbf{Function} &
      reading a tabular data file &
      \texttt{read\_csv} function in pandas \\ \cmidrule{2-3}
    & fitting a statistical model to data &
      \texttt{fit} method of scikit-learn estimators \\ \bottomrule
  \end{tabular}
\end{table}

\begin{figure}
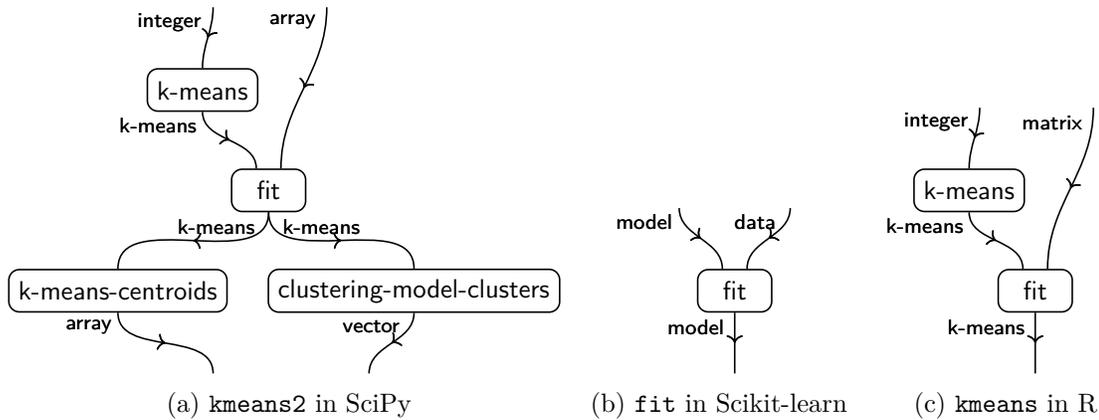

  \begin{subfigure}[b]{0.5\textwidth}
    \centering
    \input{wiring-diagrams/semantic-enrichment/annotation-python-scipy-kmeans2}
    \caption{\texttt{kmeans2} in SciPy}
    \label{fig:annotation-python-scipy-kmeans}
  \end{subfigure}%
  \begin{subfigure}[b]{0.25\textwidth}
    \centering
    \input{wiring-diagrams/semantic-enrichment/annotation-python-sklearn-fit}
    \caption{\texttt{fit} in Scikit-learn}
    \label{fig:annotation-python-sklearn-fit}
  \end{subfigure}%
  \begin{subfigure}[b]{0.25\textwidth}
    \centering
    \input{wiring-diagrams/semantic-enrichment/annotation-r-stats-kmeans}
    \caption{\texttt{kmeans} in R}
    \label{fig:annotation-r-stats-kmeans}
  \end{subfigure}
  \caption{Selected function annotations from the Data Science Ontology: (a)
    \texttt{kmeans2} function in SciPy (see \cref{lst:toy-kmeans-scipy});
    \texttt{fit} method of \texttt{BaseEstimator} class in Scikit-learn (see
    \cref{lst:toy-kmeans-sklearn}); (c) \texttt{kmeans} function in R's built-in
    \texttt{stats} package (see \cref{lst:toy-kmeans-r}).}
  \label{fig:annotations-kmeans}
\end{figure}

Significant modeling flexibility is needed to faithfully translate the widely
varying interfaces of statistical software into a single set of concepts.
\cref{lst:toy-kmeans-scipy,lst:toy-kmeans-sklearn,lst:toy-kmeans-r} show, for
example, that the concept of $k$-means clustering can be instantiated in
software in many different ways. To accommodate this diversity, function
annotations may map a single concrete function onto an arbitrary abstract
``program'' comprised of function concepts. Three function annotations related
to the fitting of $k$-means clustering models are shown in
\cref{fig:annotations-kmeans}.

An \emph{ontology language} specifies what kind of abstract ``program'' is
allowed to appear in a function annotation. It can be helpful to think of the
ontology language as a minimalistic, typed, functional programming language. The
analogy is imperfect because the ontology language is simpler than any commonly
used programming language, being designed for knowledge representation rather
than actual computing.

The ontology language has a minimalistic type system, supporting product and
unit types as well as a simple form of subtyping. A \emph{basic type} is a type
that cannot be decomposed into simpler types. Basic types must be explicitly
defined. All other types are \emph{composite}. The \emph{product} of two types
$X$ and $Y$ is another type $X \times Y$. It has the usual meaning: an element
of type $X \times Y$ is an element of type $X$ \emph{and} an element of type
$Y$, in that order. Products of three or more types are defined similarly.
Product types are similar to record types in conventional programming languages,
such as a tuple or a named tuple in Python. There is also a \emph{unit type} $1$
inhabited by a single element. It is analogous to the \verb|NoneType| type in
Python (whose sole inhabitant is \verb|None|) and the \verb|NULL| type in R
(whose sole inhabitant is also called \verb|NULL|).

A type can be declared a \emph{subtype} of one or more other types. To a first
approximation, subtyping establishes an ``is-a'' relationship between types. In
the Data Science Ontology, matrices are a subtype of both arrays (being arrays
of rank 2) and data tables (being tables whose columns all have the same data
type). As this example illustrates, subtyping in the ontology differs from
inheritance in a typical object-oriented programming language. Subtyping should
instead be understood through \emph{implicit conversion}, also known as
\emph{coercion} \cite{reynolds1980}. The idea is that if a type $X$ is a subtype
of $X'$, then there is a canonical way to convert elements of type $X$ into
elements of type $X'$. Elaborating the example, a matrix simply \emph{is} an
array (of rank 2), hence can be trivially converted into an array. A matrix is
not strictly speaking a data table but can be converted into one (of homogeneous
data type) by assigning numerical names to the columns.

A function $f: X \to Y$ in the ontology language has an input type $X$, its
\emph{domain}, and an output type $Y$, its \emph{codomain}. Like types,
functions are either basic or composite. The two basic ways of constructing
composite functions are composition and products, corresponding to the vertical
and horizontal directions in wiring diagrams. The \emph{composite} of a function
$f: X \to Y$ with $g: Y \to Z$ is a new function $f \cdot g: X \to Z$, with the
usual meaning. Algorithmically speaking, $f \cdot g$ computes \emph{in series}:
first $f$ and then $g$. The \emph{product} of functions $f: X \to Y$ and
$g: W \to Z$ is another function $f \times g: X \times W \to Y \times Z$.
Algorithmically, $f \times g$ computes $f$ and $g$ \emph{in parallel}, taking
the inputs, and returning the outputs, of both $f$ and $g$. The language also
contains special functions for permuting tuples of data, duplicating data, and
discarding data; for a more precise account, see
\cref{sec:cartesian-categories}.

Besides serving as the ``is-a'' relation ubiquitous in knowledge representation
systems, the subtype relation for objects enables \emph{ad hoc polymorphism} for
functions. The type restrictions in function composition are relaxed to allow
implicit conversion, namely, to compose a function $f: X \to Y$ with
$g: Y' \to Z$, it is not required that $Y$ equals $Y'$, but only that $Y$ be a
subtype of $Y'$. Operationally, to compute $f \cdot g$, one first computes $f$,
then coerces the result from type $Y$ to $Y'$, and finally computes $g$.
Diagrammatically, a wire connecting two boxes has valid types if and only if the
source port's type is a subtype of the target port's type. Thus implicit
conversions truly are implicit in the graphical syntax.

The ontology language also supports ``is-a'' relations between functions, called
\emph{subfunctions} in analogy to subtypes. In the Data Science Ontology, the
function concept \textsf{read-tabular-file} of reading a table from a tabular
file is a subfunction of the function concept \textsf{read-data} of reading data
from a generic data source. The meaning of this statement is as follows. The
domain of \textsf{read-tabular-file}, a tabular file, is a subtype of the domain
of \textsf{read-data}, a generic data source. The codomain of
\textsf{read-tabular-file}, a table, is a subtype of the codomain of
\textsf{read-data}, generic data. Now consider two possible computational paths
that take a tabular file and return generic data. We could apply
\textsf{read-tabular-file}, then coerce the resulting table to generic data.
Alternatively, we could coerce the tabular file to a generic data source, then
apply \textsf{read-data}. The subfunction relation asserts that these two
computations are equivalent. The definition of the subfunction relation for
general functions $f: X \to Y$ and $f': X' \to Y'$ is analogous and is stated
formally in \cref{sec:concepts-as-category}.

\section{Semantic enrichment algorithm}
\label{sec:semantic-enrichment}

The semantic enrichment algorithm, transforming raw flow graphs into semantic
flow graphs, proceeds in two independent stages, one expansionary and the other
contractionary. The expansion stage makes essential use of code annotations in
the ontology.

\paragraph{Expansion}

In the expansion stage, the annotated parts of the raw flow graph are replaced
by their abstract definitions. Each annotated box---that is, each box referring
to a concrete function annotated by the ontology---is replaced by the
corresponding abstract function. Likewise, the concrete type of each annotated
port is replaced by the corresponding abstract type. This stage of the algorithm
is ``expansionary'' since a function annotation's definition can be an arbitrary
program in the ontology language. In other words, a single box in the raw flow
graph can become an arbitrarily large subdiagram in the semantic flow graph.

The expansion procedure is \emph{functorial}, to use the jargon of category
theory. Informally, this means two things. First, notice that concrete types are
effectively annotated twice, explicitly by type annotations and implicitly by
the domain and codomain types in function annotations. Functorality requires
that these abstract types be compatible, ensuring the logical consistency of
type and function annotations. Second, expansion preserves the structure of the
ontology language, including composition and products. The expansion of a wiring
diagram is thus completely determined by its action on individual boxes (basic
functions). Functorality is a modeling decision that greatly simplifies the
semantic enrichment algorithm, at the expense of imposing restrictions on how
the raw flow graph can be transformed.

\paragraph{Contraction}

It is practically infeasible to annotate every reusable unit of data science
source code. Even if the Data Science Ontology were to grow significantly, most
real-world data analyses would use concrete types and functions without
annotations. This unannotated code has unknown semantics, so properly speaking
it does not belong in the semantic flow graph. However, it usually cannot be
deleted without altering the data flow of the wiring diagram.

As a compromise, in the contraction stage, the unannotated parts of the raw flow
graph are simplified to the extent possible. All references to unannotated types
and functions are removed, leaving behind unlabeled wires and boxes.
Semantically, the unlabeled wires are interpreted as arbitrary ``unknown'' types
and the unlabeled boxes as arbitrary ``unknown'' functions (which could have
known domain and codomain types). The diagram is then simplified by
\emph{encapsulating} unlabeled boxes. Specifically, every maximal connected
subdiagram of unlabeled boxes is encapsulated by a single unlabeled box. The
interpretation is that any composition of unknown functions is just another
unknown function. This stage is ``contractionary'' because it can only decrease
the number of boxes in the diagram.

\paragraph{Examples revisited}

To reprise the small example from the Introduction, semantic enrichment
transforms all three raw flow graphs of
\cref{fig:toy-kmeans-scipy,fig:toy-kmeans-sklearn,fig:toy-kmeans-r} into the
semantic flow graph of \cref{fig:toy-kmeans-semantic}. Expansions related to
$k$-means clustering occur in all three programs.

In the first Python program (\cref{lst:toy-kmeans-scipy} and
\cref{fig:toy-kmeans-scipy}), the \verb|kmeans2| function from SciPy expands to
a compound function that creates a $k$-means clustering model, fits it to data,
and extracts its cluster assignments and centroids, as specified by the
annotation in \cref{fig:annotation-python-scipy-kmeans}. Note that the abstract
$k$-means clustering model does \emph{not} correspond to any concrete object in
the original program. This design pattern is used throughout the ontology to
cope with functions that are not object-oriented with respect to models.

By contrast, the second Python program is written in object-oriented style
(\cref{lst:toy-kmeans-sklearn} and \cref{fig:toy-kmeans-sklearn}). The
\verb|KMeans| class from Scikit-learn expands to an abstract type
\textsf{k-means}. The \verb|fit| method of this class is not annotated in the
Data Science Ontology. However, the \verb|fit| method of its superclass
\verb|BaseEstimator| \emph{is} annotated
(\cref{fig:annotation-python-sklearn-fit}), so the expansion is performed using
this annotation. In general, subtyping and polymorphism are indispensable for
annotating object-oriented libraries parsimoniously.

The R program is intermediate between these two styles (\cref{lst:toy-kmeans-r}
and \cref{fig:toy-kmeans-r}). The \verb|kmeans| function, annotated in
\cref{fig:annotation-r-stats-kmeans}, directly receives the data and the number
of clusters, but returns an object of class \verb|kmeans|. The cluster
assignments and centroids are slots of this object, annotated separately. This
design pattern is typical in R, due to its informal type system.

Contractions also occur in the three programs. In the first Python program, the
only unannotated box is NumPy's \verb|delete| function. Contracting this box
does not reduce the size of the wiring diagram. A contraction involving multiple
boxes occurs in the second Python program. The subdiagram consisting of the
Pandas method \verb|NDFrame.drop| composed with the attribute accessor
\verb|values| is encapsulated into a single unlabeled box.

\begin{figure}
  \centering
  \resizebox{\textwidth}{!}{%
    \input{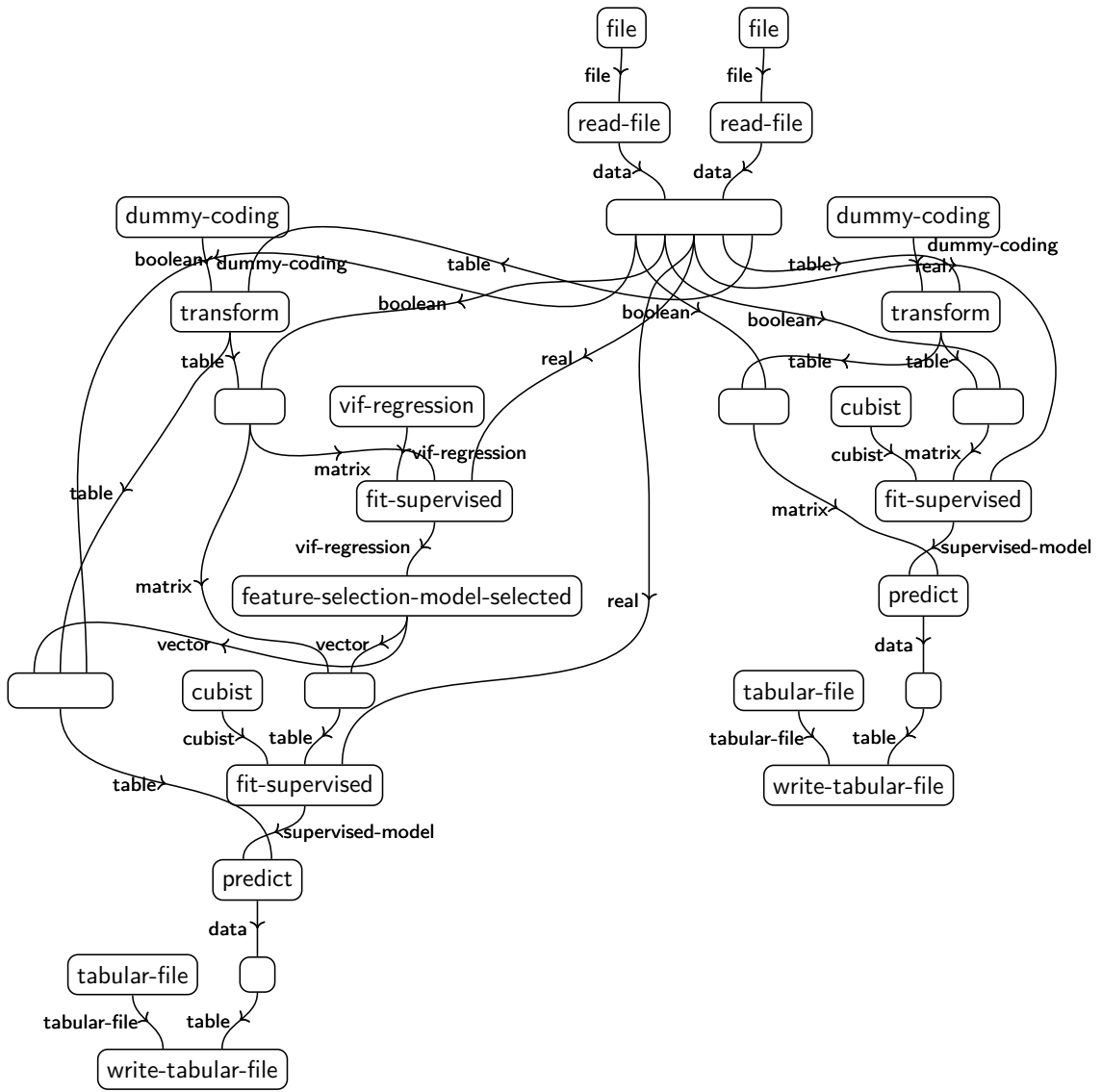}}
  \caption{Semantic flow graph for data analysis from Rheumatoid Arthritis DREAM
    Challenge (\cref{lst:dream-ra} and \cref{fig:dream-ra-raw})}
  \label{fig:dream-ra-semantic}
\end{figure}

As a more realistic example, recall the data analysis for the Rheumatoid
Arthritis DREAM Challenge presented in the previous chapter. The R code was
displayed in \cref{lst:dream-ra} and the raw flow graph in
\cref{fig:dream-ra-raw}. The semantic flow graph is now shown in
\cref{fig:dream-ra-semantic}. The two models fit by the analysis, one involving
a feature selection step and the other not, appear as the left and right
branches of the semantic flow graph. As in the $k$-means clustering example,
most of the unlabeled nodes in \cref{fig:dream-ra-semantic}, including the wide
node at the top, refer to code for data preprocessing or transformation. It is a
commonplace among data scientists that such ``data munging'' is a crucial aspect
of data analysis. There is no fundamental obstacle to representing its
semantics; it so happens that the relevant portion of the Data Science Ontology
has not yet been developed.

\section{Concepts as a category}
\label{sec:concepts-as-category}

The Data Science Ontology is formalized as a finitely presented category of a
certain kind. The type and function concepts in the ontology are, respectively,
the objects and morphisms that generate the category. Abstract programs
expressed in the language of concepts are arbitrary morphisms in the category,
constructed from the object and morphism generators through composition,
monoidal products, and other operations. In this section, we develop a
categorical structure suitable for the Data Science Ontology by augmenting
cartesian closed categories with a form of subtyping based on implicit
conversion. Ultimately, an ontology is defined to be a finite presentation of a
cartesian closed category with implicit conversion.

Cartesian closed categories are perhaps the simplest model of typed, functional
computing. Cartesian categories, reviewed in \cref{sec:cartesian-categories},
have a type system with product and unit types, while closed categories, defined
below, also have function types. In our experience, augmenting the type system
with some form of polymorphism is a practical necessity, for the sake of both
knowledge representation and parsimonious annotation of code. Our aim is not to
give a sophisticated account of polymorphism but to define the minimal
practically useful system. The following approach to polymorphism is adapted
from the works of Goguen and Reynolds \cite{goguen1978,reynolds1980}.

\begin{definition}[Implicit conversion]
  A \emph{category with implicit conversion} is a category $\cat C$ with a
  distinguished subcategory $\cat{C}_0$ that is wide as a subcategory but thin
  as a category. That is, $\cat{C}_0$ contains all the objects of $\cat{C}$ but
  contains at most one morphism between any two objects. If there exists a
  morphism $X \to X'$ in $\cat{C}_0$, we write $X \leq X'$ and say that $X$ is a
  \emph{subtype} of $X'$. The morphism $X \to X'$ itself is called an
  \emph{implicit conversion} or \emph{coercion}.\footnote{To be consistent in
    the usage of categorical and programming terminology, one might wish to say
    that $X$ is a \emph{subobject} of $X'$. However, the term ``subobject''
    already has an established meaning in categorical logic that is related to,
    but different than, the usage of ``subtype'' here.}
\end{definition}

The informal interpretation of subtyping and implicit conversion was explained
in \cref{sec:dso}. One subtle point should be noted: even when types are
interpreted as sets, implicit conversions are not necessarily set inclusions or
even injections. In the motivating example, matrices are a subtype of data
tables, yet the set of matrices is \emph{not} a subset of the set of data tables
under any plausible formalization of a data table. The implicit conversion
function must add names to the columns of the matrix, among other possible
obstructions. Hence, to the slogan that ``types are not sets'' \cite{morris1973}
it may be added that ``subtypes are not subsets.''

Mathematically speaking, the subtype relation defines a preorder on the objects
of $\cat C$. Thus, every type $X$ is a subtype of itself. If $X$ is a subtype of
$X'$ and $X'$ a subtype of $X''$, then $X$ is a subtype of $X''$. The
corresponding implicit conversions are given by identities and by composition,
respectively. In what follows, there is no mathematical obstruction to allowing
the conversions $\cat{C}_0$ to form an arbitrary category, not necessarily a
preorder. That would, however, defeat the purpose: conversions would have to be
disambiguated by \emph{name} and hence would cease to be implicit.

When a category $\cat C$ with implicit conversion is also a monoidal category,
the implicit conversions should be compatible with the monoidal product.

\begin{definition}[Monoidal implicit conversion]
  A \emph{monoidal category with implicit conversion} is a monoidal category
  $(\cat{C},\otimes,I)$ whose underlying category $\cat{C}$ has implicit
  conversions $\cat{C}_0$, such that $\cat{C}_0$ forms a \emph{monoidal}
  subcategory of $\cat{C}$.

  A symmetric monoidal category has implicit conversions if its underlying
  monoidal category does. Likewise, a cartesian category has implicit
  conversions if its underlying monoidal category $(\cat{C},\times,1)$ does.
\end{definition}

The definition requires that subtyping be compatible with product types, in the
sense that if $X \leq X'$ and $Y \leq Y'$, then $X \times Y \leq X' \times Y'$,
where the corresponding implicit conversion is given by the product of
morphisms. The subtype relation thus makes $\cat{C}_0$ into a \emph{monoidal
  preorder}.

\begin{remark}[Implicit versus explicit conversion]
  If $\cat{C}$ is a non-strict cartesian category, then the requirement that the
  implicit conversions form a monoidal subcategory implies that the associators
  and unitors belong to $\cat{C}_0$. Therefore,
  $(X \times Y) \times Z \approx X \times (Y \times Z)$ and
  $X \times 1 \approx X \approx 1 \times X$ for all objects $X,Y,Z \in \cat{C}$,
  where we write $X \approx Y$ for the equivalence relation that $X \leq Y$ and
  $Y \leq X$. It is natural to require that the associators and unitors be
  implicit conversions, since one does not typically wish to distinguish between
  different bracketings of products. In particular, the implicit conversions
  $\cat{C}_0$ do \emph{not} form a partial order when the monoidal category
  $\cat{C}$ is not strict (and need not form a partial order even when $\cat{C}$
  is strict).

  It may seem strange that $\cat{C}_0$ is not asked to inherit the cartesian or
  even the symmetric monoidal structure of $\cat{C}$. However, this leads to
  unwanted implicit conversions and to strictification of the original category.
  Namely, if $\cat{C}_0$ is a symmetric monoidal subcategory of $\cat C$, then
  the braidings $\Braid_{X,Y}: X \times Y \to Y \times X$ in $\cat C$ must
  satisfy $\Braid_{X,X} = 1_{X \times X}$. This is false under the set-theoretic
  interpretation and implies that $\cat{C}$ is symmetrically monoidally
  equivalent to a strictly commutative symmetric monoidal category
  \cite{kim2016}. The braidings must therefore be treated as \emph{explicit}
  conversions.
\end{remark}

Because the notion of subtyping is operationalized by the implicit conversions,
it can be extended from objects to morphisms through naturality squares.

\begin{definition}[Submorphism]
  Let $\cat C$ be a category with implicit conversion. A morphism $f$ in
  $\cat C$ is a \emph{submorphism} (or \emph{subfunction}) of another morphism
  $f'$, written $f \leq f'$, if in the arrow category $\cat C^\to$ there exists
  a (unique) morphism $f \to f'$ whose components are implicit conversions.

  Explicitly, if $f: X \to Y$ and $f': X' \to Y'$ are morphisms in $\cat C$,
  with $X \leq X'$ and $Y \leq Y'$, then $f \leq f'$ if and only if the diagram
  commutes:
  \begin{equation*}
    \begin{tikzcd}
      X \arrow[r, "f"] \arrow[d, "\leq"'] & Y \arrow[d, "\leq"] \\
      X' \arrow[r, "f'"] & Y'
    \end{tikzcd}
  \end{equation*}
\end{definition}

Again, see \cref{sec:dso} for the informal interpretation and examples of this
notion. Just as subtypes define a preorder on the objects of $\cat C$,
submorphisms define a preorder on the morphisms of $\cat C$. Moreover,
submorphisms respect the compositional structure of $\cat C$. They are closed
under identities, i.e., $1_X \leq 1_{X'}$ whenever $X \leq X'$, and under
composition, i.e., if $f \leq f'$ and $g \leq g'$ are composable, then
$fg \leq f'g'$. All these statements are easy to prove. To illustrate,
transitivity and closure under composition are proved by pasting commutative
squares vertically and horizontally:
\begin{equation*}
  \begin{tikzcd}
    X \arrow[r, "f"] \arrow[d, "\leq"'] & Y \arrow[d, "\leq"] \\
    X' \arrow[r, "f'"] \arrow[d, "\leq"'] & Y' \arrow[d, "\leq"] \\
    X'' \arrow[r, "f''"] & Y''
  \end{tikzcd}
  \qquad\qquad
  \begin{tikzcd}
    X \arrow[r, "f"] \arrow[d, "\leq"']
      & Y \arrow[r, "g"] \arrow[d, "\leq"']
      & Z \arrow[d, "\leq"] \\
    X' \arrow[r, "f'"] & Y' \arrow[r, "g'"] & Z'
  \end{tikzcd}
\end{equation*}
When $\cat C$ is a \emph{cartesian} category with implicit conversion,
submorphisms are also closed under products: if $f \leq f'$ and $g \leq g'$,
then $f \times g \leq f' \times g'$, because, by functorality, monoidal products
preserve commutative diagrams.

As an aside, we note that this structure is succinctly captured by the notion of
a monoidal double category \cite{bruni2002,shulman2010,hansen2019}.

\begin{proposition}
  A (monoidal) category $\cat{C}$ with implicit conversions forms a (monoidal)
  double category, in which the objects are the objects of $\cat{C}$, the
  horizontal 1-morphisms are the morphisms of $\cat{C}$, the vertical
  1-morphisms are the implicit conversions, and the 2-morphisms are the
  submorphisms.
\end{proposition}

\noindent Compared to a general monoidal double category, the one here is quite
simple, being thin with respect to both vertical 1-morphisms and 2-morphisms.

In a functional programming language, functions are treated as data and as such
they can be passed as arguments to other functions. The type for functions
$X \to Y$ is the \emph{function type}, or \emph{exponential type}, $Y^X$.
Algebraically, the existence of function types corresponds to the property of a
monoidal category of being \emph{closed}.

\begin{definition}[Closed monoidal category]
  A symmetric monoidal category $(\cat{C}, \otimes, I)$ is \emph{closed} if for
  every object $X \in \cat{C}$, the product functor
  $- \otimes X: \cat{C} \to \cat{C}$ has a right adjoint, denoted
  $[X,-]: \cat{C} \to \cat{C}$. That is, for all objects $W,X,Y \in \cat{C}$,
  there is a bijection of hom-sets
  \begin{equation*}
    \cat{C}(W \otimes X, Y) \cong \cat{C}(W, [X,Y]),
  \end{equation*}
  natural in $W$ and $Y$. The object $[X,Y]$ is called an \emph{internal hom}.

  A \emph{cartesian closed category} is a cartesian category that is also
  closed. In this case, the internal hom $[X,Y]$ is called an \emph{exponential
    object} and is often denoted $Y^X$.
\end{definition}

Taking $W$ to be the monoidal unit, the adjunction establishes a one-to-one
correspondence between morphisms $f: X \to Y$ in $\cat{C}$ and generalized
elements of type $[X,Y]$. More generally, each morphism $f: W \times X \to Y$
corresponds to a unique morphism $\lambda(f): W \to [X,Y]$, the \emph{currying}
of $f$. The inverse transformation associates each morphism $g: W \to [X,Y]$
with its \emph{uncurrying} $\lambda^{-1}(g): W \times X \to Y$.

Let $\cat{C}$ be any closed monoidal category. By the parameterized adjunction
theorem \cites[Theorem IV.7.3]{maclane1998}[Proposition 4.3.6]{riehl2016}, the
family of internal hom functors $[X,-]: \cat{C} \to \cat{C}$, parameterized by
objects $X \in \cat{C}$, assembles into a bifunctor
$[-,-]: \cat{C}^\opposite \times \cat{C} \to \cat{C}$ such that the bijection of
hom-sets $\cat{C}(W \otimes X, Y) \cong \cat{C}(W, [X,Y])$ is natural not just
in $W$ and $Y$ but also in $X$. When the category $\cat{C}$ has implicit
conversions, the conversions should be preserved by this bifunctor.

\begin{definition}[Closed implicit conversion]
  A \emph{closed monoidal category with implicit conversion} is a symmetric
  monoidal category with implicit conversion that is also closed, in such a way
  that the internal hom bifunctor
  $[-,-]: \cat{C}^\opposite \times \cat{C} \to \cat{C}$ restricts to a bifunctor
  $\cat{C}_0^\opposite \times \cat{C}_0 \to \cat{C}_0$ of implicit conversions.
  Equivalently, $[X',Y] \leq [X,Y']$ whenever $X \leq X'$ and $Y \leq Y'$.

  A cartesian closed category has implicit conversions if its underlying closed
  monoidal category does. We denote by $\CCC_\leq$ the category whose objects
  are the (small) cartesian closed categories with implicit conversion and whose
  morphisms are the cartesian closed functors that preserve implicit
  conversions. When no confusion will result, the morphisms are simply called
  ``functors.''
\end{definition}

The subtype relation $[X',Y] \leq [X,Y']$, where $X \leq X'$ and $Y \leq Y'$,
has the interpretation that any morphism $f: X' \to Y$ can be implicitly
converted to a morphism $\tilde f: X \to Y'$ by ``restricting the domain'' and
``expanding the codomain'':
\begin{equation*}
  \begin{tikzcd}
    X' \arrow[r, "f"] & Y \arrow[d, "\leq"] \\
    X \arrow[r, dashed, "\tilde f"] \arrow[u, "\leq"] & Y'.
  \end{tikzcd}
\end{equation*}
Subtypes of function types should not be confused with subfunctions. The former
is a relation between the objects and the latter between morphisms; moreover,
the former is contravariant with respect to the domain object, while the latter
is covariant.

The submorphism relation is, however, compatible with currying and uncurrying,
as the following proposition shows.

\begin{proposition}
  Let $(\cat{C}, \otimes, I)$ be a closed monoidal category with implicit
  conversions and let $W \leq W'$, $X \leq X'$, and $Y \leq Y'$ be objects of
  $\cat{C}$. For any morphisms $f: W \otimes X \to Y$ and
  $f': W' \times X \to Y'$, one has $f \leq f'$ if and only if
  $\lambda(f) \leq \lambda(f')$. Moreover, for any morphisms
  $g: W \otimes X \to Y$ and $g': W \otimes X' \to Y$, one has $g \leq g'$ if
  and only if $\lambda(g') \leq \lambda(g)$.
\end{proposition}
\begin{proof}
  The two statements are the equivalences
  \begin{equation*}
    \begin{tikzcd}
      W \otimes X \arrow[r, "f"] \arrow[d, "\leq"'] & Y \arrow[d, "\leq"] \\
      W' \otimes X \arrow[r, "f'"] & Y'
    \end{tikzcd}
    \qquad\leftrightsquigarrow\qquad
    \begin{tikzcd}
      W \arrow[r, "\lambda f"] \arrow[d, "\leq"'] & {[X,Y]} \arrow[d, "\leq"] \\
      W' \arrow[r, "\lambda f'"] & {[X,Y']}
    \end{tikzcd}
  \end{equation*}
  and
  \begin{equation*}
    \begin{tikzcd}
      W \times X \arrow[r, "g"] \arrow[d, "\leq"'] & Y \arrow[d, equal] \\
      W \times X' \arrow[r, "g'"] & Y
    \end{tikzcd}
    \qquad\leftrightsquigarrow\qquad
    \begin{tikzcd}
      W \arrow[r, "\lambda g"] \arrow[d, equal] & {[X,Y]} \\
      W \arrow[r, "\lambda g'"] & {[X',Y]}  \arrow[u, "\leq"']
    \end{tikzcd},
  \end{equation*}
  meaning that the left square commutes if and only if the right square does.
  The first equivalence follows from the naturality of the adjunction with
  respect to $W$ and $Y$ \cite[Lemma 4.1.3]{riehl2016}, while the second
  equivalence follows from the naturality of the adjunction with respect to $X$
  \cite[Theorem IV.7.3]{maclane1998}.
\end{proof}

With these preliminaries on implicit conversion, an ontology is now defined to
be nothing other than a \emph{finitely presented} cartesian closed category with
implicit conversion.

\begin{definition}[Ontology]
  An \emph{ontology} is a cartesian closed category with implicit conversion,
  given by a finite presentation. That is, it is the cartesian closed category
  with implicit conversion generated by finite sets of:
  \begin{itemize}[nosep]
  \item \emph{basic types}, or \emph{object generators}, $X$
  \item \emph{basic functions}, or \emph{morphism generators}, $f: X \to Y$,
    where $X$ and $Y$ are objects
  \item \emph{basic subtypes}, or \emph{subtype generators}, $X \leq X'$, where
    $X$ and $X'$ are objects
  \item \emph{basic subfunctions}, or \emph{submorphism generators},
    $f \leq f'$, where $f: X \to Y$ and $f': X' \to Y'$ are morphisms satisfying
    $X \leq X'$ and $Y \leq Y'$
  \item \emph{function equations}, or \emph{morphism equations}, $f=g$, where
    $f,g: X \to Y$ are morphisms with equal domains and codomains.
  \end{itemize}
  If the set of morphism equations is empty, the ontology is called \emph{free}
  or \emph{freely generated}.
\end{definition}

Strictly speaking, a finite presentation of a category is not the same as the
category it presents. The former is a finitary object that can be represented
on, and manipulated by, a machine. The latter is an algebraic structure of
infinite size, convenient for mathematical reasoning. However, we will abuse
terminology by calling both finitely presented categories, and particular
presentations thereof, ``ontologies.''

At the time of this writing, the Data Science Ontology is freely generated.
Inference in a freely generated ontology is straightforward. Assuming that the
generating subtypes and subfunctions are between basic objects and basic
functions, deciding the subtype or subfunction relations amounts to computing a
reflexive transitive closure. Deciding equality of objects is trivial. Deciding
equality of morphisms is the \emph{word problem} in a free cartesian closed
category. In the cartesian case, this problem can be solved by the congruence
closure algorithm for term graphs \cite[\S 4.4]{baader1999}. In the future, the
Data Science Ontology may include knowledge in the form of morphism equations,
creating a need for new inference procedures. If arbitrary morphism equations
are allowed, the word problem becomes undecidable.

\section{Annotations as a functor}
\label{sec:annotations-as-functor}

If the concepts form a category, then surely the annotations ought to assemble
into a functor. Let $\cat C$ be a cartesian closed category with implicit
conversion, viewed as the concepts of an ontology. Suppose $\cat L$ is another
such category, modeling a programming language and a collection of modules
written in that language. The annotations ought to define a functor
$F: \cat L \to \cat C$, saying how to translate programs in $\cat L$ into
programs in $\cat C$.

This tidy story does not quite survive contact with reality. A fairly small set
of formal concepts cannot be to exhaust the supply of informal concepts found in
real-world programs. Consequently, any ``functor'' $F: \cat L \to \cat C$
annotating $\cat L$ must be \emph{partial}, in a sense that must be made
precise. There will be both objects and morphisms in $\cat L$ on which $F$
cannot be defined, because the category $\cat C$ is not rich enough to fully
interpret $\cat L$.

Before turning to partial functors, consider the simpler case of partial
functions. In accordance with mathematical custom, the pre-theoretical idea of
``partial function'' can be reduced to the standard notion of total function.
There are two common ways to do this, the first based on pointed sets and the
second on spans. They are equivalent as far as sets and functions are concerned
but suggest different generalizations to categories and functors.

The category of pointed sets leads to one viewpoint on partiality, popular in
programming language theory. Given a set $X$, let $X_\bot := X \sqcup \{\bot\}$
be the set $X$ with a freely adjoined base point $\bot$. A \emph{partial
  function} from $X$ to $Y$ is then a function $f: X_\bot \to Y_\bot$ preserving
the base point, so that $f(\bot) = \bot$. The function $f$ is regarded as
``undefined'' on the points $x \in X$ with $f(x) = \bot$. This notion of
partiality can be transported from sets to categories using enriched category
theory. Categories enriched in pointed sets, where each hom-set has a base
morphism $\bot$, have been proposed as a qualitative model of incomplete
information \cite{marsden2016}. Such categories make partiality an
all-or-nothing affair, because their composition laws satisfy
$\bot \cdot f = f \cdot \bot = \bot$ for all morphisms $f$. That is far too
stringent. If this composition law were adopted, the semantic flow graphs would
rarely be anything besides the trivial morphism $\bot$.

A partial function can also be defined as a special kind of span of total
functions. On this view, a \emph{partial function} from $X$ to $Y$ is a span in
$\Set$
\begin{equation*}
  \begin{tikzcd}[column sep=small, row sep=small]
    & J \arrow[dl, tail, "\iota"'] \arrow[dr, "f"] & \\
    X & & Y
  \end{tikzcd}
\end{equation*}
whose left leg $\iota: J \to X$ is monic (injective). The partial function's
domain of definition is $J$, regarded as a subset of $X$. Although we shall not
need it here, we note that partial functions, and partial morphisms generally,
can be composed by taking pullbacks whenever they exist \cite[\S
5.5]{borceux1994c}.

The span above can be interpreted as \emph{partially} defining a function $f$ on
$X$, via a set of equations indexed by $J$:
\begin{equation*}
  f(x_j) := y_j, \qquad j \in J.
\end{equation*}
It is then natural to ask: what is the most general way to define a \emph{total}
function on $X$ obeying these equations? The answer is given by the pushout in
$\Set$:
\begin{equation*}
  \begin{tikzcd}
    & J \arrow[dl, tail, "\iota"'] \arrow[dr, "f"]
      \arrow[dd, phantom, very near end, "\ucorner"] & \\
    X \arrow[dr, "f_*"'] & & Y \arrow[dl, tail, "\iota_*"] \\
    & Y_* &
  \end{tikzcd}
\end{equation*}
Because $\iota: J \to X$ is monic, so is $\iota_*: Y \to Y_*$, and $Y$ can be
regarded as a subset of $Y_*$. The commutativity of the diagram says that $f_*$
satisfies the set of equations indexed by $J$. The \emph{universal property}
defining the pushout says that any other function $f': X \to Y'$ satisfying the
equations factors uniquely through $f_*$, meaning that there exists a unique
function $g: Y_* \to Y'$ making the diagram commute:
\begin{equation*}
   \begin{tikzcd}
    & J \arrow[dl, "\iota"'] \arrow[dr, "f"] & \\
    X \arrow[r, "f_*"] \arrow[dr, "f'"'] & Y_* \arrow[d, dashed, "g"] &
      Y \arrow[l, "\iota_*"] \arrow[dl, "\iota'"] \\
    & Y' &
  \end{tikzcd}
\end{equation*}
The codomain of the function $f_*: X \to Y_*$ consists of $Y$ plus a ``formal
image'' $f(x)$ for each element $x$ on which $f$ is undefined. Contrast this
with the codomain of a function $X \to Y_\bot$, which consists of $Y$ plus a
single element $\bot$ representing \emph{all} the undefined values.

This viewpoint on partiality generalizes effortlessly from $\Set$ to any
category with pushouts. Annotations can be defined as a span in $\CCC_\leq$
\begin{equation*}
   \begin{tikzcd}[column sep=small, row sep=small]
    & \cat{J} \arrow[dl, tail, "\iota"'] \arrow[dr, "F"] & \\
    \cat{L} & & \cat{C}
  \end{tikzcd}
\end{equation*}
whose left leg $\iota: \cat{J} \to \cat{L}$ is monic. We then form the pushout
in $\CCC_\leq$:
\begin{equation*}
  \begin{tikzcd}
    & \cat{J} \arrow[dl, tail, "\iota"'] \arrow[dr, "F"]
      \arrow[dd, phantom, very near end, "\ucorner"] & \\
    \cat{L} \arrow[dr, "F_*"'] & & \cat{C} \arrow[dl, "\iota_*"] \\
    & \cat{C}_* &
  \end{tikzcd}
\end{equation*}
Given a morphism $f$ in $\cat L$, which represents a concrete program, its image
$F_*(f)$ in $\cat{C}_*$ is a partial translation of the program into the
language defined by the ontology's concepts.

The universal property of the pushout in $\CCC_\leq$, stated above in the case
of $\Set$, gives an appealing intuitive interpretation to program translation.
The category $\cat C$ is not rich enough to fully translate $\cat L$ via a
functor $\cat L \to \cat C$. As a modeling assumption, we suppose that $\cat C$
has some ``completion'' $\overline{\cat C}$ for which a full translation
$\overline F: \cat L \to \overline{\cat C}$ \emph{is} possible. We do not know
$\overline{\cat C}$, or at the very least we cannot feasibly write it down.
However, if we take the pushout functor $F_*: \cat{L} \to \cat{C}_*$, we can at
least guarantee that, no matter what the complete translation $\overline{F}$ is,
it will factor through $F_*$. Thus $F_*$ defines the most general possible
translation, given the available information.

The properties of partial functions largely carry over to partial functors, with
one important exception: the ``inclusion'' functor
$\iota_*: \cat{C} \to \cat{C}_*$ need not be monic, even though
$\iota: \cat{J} \to \cat{L}$ is. Closely related is the fact that $\CCC_\leq$
(like its cousins $\Cat$, $\Cart$, and $\CCC$, but unlike $\Set$) does not
satisfy the \emph{amalgamation property} \cite{macdonald2009}. To see how
$\iota_*$ can fail to be monic, suppose that the equation $f_1 \cdot f_2 = f_3$
holds in $\cat{L}$ and that the defining equations include $F(f_i) := g_i$ for
$i=1,2,3$. Then, by the functorality of $F_*$, we must have
$g_1 \cdot g_2 = g_3$ in $\cat{C}_*$, even if $g_1 \cdot g_2 \neq g_3$ in
$\cat{C}$. Thus the existence of $F_*$ can force equations between morphisms in
$\cat{C}_*$ that do not hold in $\cat{C}$.

When the categories in question are finitely presented, the pushout functor also
admits a finitary, equational presentation, suitable for computer algebra. Just
as an ontology is defined to be a finitely presented category, an ontology with
annotations is defined to be a finitely presented functor.

\begin{definition}[Annotations]
  An \emph{ontology with annotations} is a functor between cartesian closed
  categories with implicit conversion, defined by a finite presentation.
  Explicitly, it is generated by:
  \begin{itemize}[nosep]
  \item a finite presentation of a category $\cat C$ in $\CCC_\leq$, the
    ontology category;
  \item a finite presentation of a category $\cat L$ in $\CCC_\leq$, the
    programming language category; and
  \item a finite set of equations partially defining a functor $F$ from $\cat L$
    to $\cat C$.
  \end{itemize}
\end{definition}

The equations partially defining the functor $F$ may be indexed by a category
$\cat{J}$, in which case they take the form
\begin{equation*}
  F(X_j) := Y_j \qquad\text{where}\qquad
    X_j \in \cat L, \quad Y_j \in \cat C,
\end{equation*}
for each $j \in \cat{J}$, and
\begin{equation*}
  F(f_k) := g_k \qquad\text{where}\qquad
    f_k \in \cat L(X_j, X_{j'}), \quad g_k \in \cat C(Y_j, Y_{j'}),
\end{equation*}
for each $j,j' \in \cat{J}$ and $k \in \cat{J}(j,j')$. The equations present a
span
$\cat{L} \overset{\iota}\leftarrowtail \cat{J} \overset{F}\rightarrow \cat{C}$
whose left leg is monic and the functor generated by the equations is the
pushout functor $F_*: \cat{L} \to \cat{C}_*$ described above.

\begin{remark}[Formalizing presentations]
  The two definitions involving finite presentations can be made completely
  formal using generalized algebraic theories \cite{cartmell1978,cartmell1986}.
  There is a generalized algebraic theory of cartesian closed categories with
  implicit conversion, whose category of models is $\CCC_\leq$, and a theory of
  functors between them, whose category of models is the arrow category
  $\CCC_\leq^\to$. Cartmell gives as simpler examples the theory of categories,
  with models $\Cat$, and the theory of functors, with models $\Cat^\to$
  \cite{cartmell1986}. Lambek and Scott give an equational theory of $\cat{CCC}$
  \cite[\S I.3]{lambek1986}. Any category of models of a generalized algebraic
  theory is cocomplete and admits free models defined by finite presentations.
\end{remark}

\section{Flow graphs and categories of elements}
\label{sec:categories-of-elements}

To a first approximation, the raw and semantic flow graphs are morphisms in the
categories $\cat L$ and $\cat C_*$, respectively. The expansion stage of the
semantic enrichment algorithm simply applies the annotation functor
$F_*: \cat L \to \cat C_*$ to a morphism in $\cat L$. The contraction stage,
purely syntactical, groups together morphisms in $\cat C_*$ that are not images
of $\cat C$ under the inclusion functor $\iota_*: \cat C \to \cat C_*$.

To complete the formalization of semantic enrichment, the observed elements in
the raw and semantic flow graphs must be accounted for. Flow graphs capture not
only the types and functions comprising a program, but also the values computed
by the program. In category theory, values can be bundled together with objects
and morphisms using a device called the \emph{category of elements}. The raw and
semantic flow graphs are formalized as morphisms in categories of elements.

The objects and morphisms in the ontology category $\cat C$ can, at least in
principle, be interpreted as sets and functions. Extending the definition in
\cref{sec:cartesian-categories}, a set-theoretic \emph{model} of $\cat C$ is a
cartesian closed functor $M_{\cat C}: \cat C \to \Set$. In programming language
terms, $M_{\cat C}$ is a \emph{denotational semantics} for $\cat C$. Suppose the
concrete language $\cat L$ also has a model $M_{\cat L}: \cat L \to \Set$.
Assuming that the equations partially defining the annotation functor hold in
the models of $\cat{C}$ and $\cat{L}$, the diagram
\begin{equation*}
  \begin{tikzcd}
    & \cat{J} \arrow[dl, tail, "\iota"'] \arrow[dr, "F"] & \\
    \cat{L} \arrow[dr, "M_{\cat L}"'] & & \cat{C} \arrow[dl, "M_{\cat C}"] \\
    & \Set &
  \end{tikzcd}
\end{equation*}
commutes. By the universal property of the annotation functor $F_*$, there
exists a unique model $M_{\cat C_*}: \cat C_* \to \Set$ making the diagram
commute:
\begin{equation*}
  \begin{tikzcd}
    \cat{L} \arrow[r, "F_*"] \arrow[dr, "M_{\cat L}"']
      & \cat{C}_* \arrow[d, dashed, pos=0.33, "M_{\cat{C}_*}"]
      & \cat{C} \arrow[l, "\iota_*"'] \arrow[dl, "M_{\cat C}"] \\
    & \Set &
  \end{tikzcd}
\end{equation*}

Each of these three interpretations yields a category of elements, also known as
a ``Grothendieck construction'' \cites[\S 12.2]{barr1990}[\S 2.4]{riehl2016}.
\begin{definition}
  The \emph{category of elements} of a cartesian functor $M: \cat C \to \Set$
  has as objects, the pairs $(X,x)$, where $X \in \cat C$ and $x \in M(X)$, and
  as morphisms $(X,x) \to (Y,y)$, the morphisms $f: X \to Y$ in $\cat C$
  satisfying $M(f)(x) = y$.
\end{definition}

The category of elements of a cartesian functor $M: \cat C \to \Set$ is itself a
cartesian category. Composition and identities are inherited from $\cat C$.
Products are defined on objects by
\begin{equation*}
  (X,x) \times (Y,y) := (X \times Y, (x,y))
\end{equation*}
and on morphisms exactly as in $\cat C$, and the unit object is $(1,*)$, where
$*$ is the single element of $M(1) = \{*\}$. The braidings and the supply of
commutative comonoids are also inherited from $\cat C$, the latter as
\begin{equation*}
  \Copy_{(X,x)}: (X,x) \to (X \times X, (x,x)), \qquad
  \Delete_{(X,x)}: (X,x) \to (1,*).
\end{equation*}
The category of elements of a cartesian \emph{closed} functor
$M: \cat{C} \to \Set$ also contains objects of the form $(Y^X, g)$, where $g$ is
any function $M(X) \to M(Y)$. However, the category of elements is generally not
closed, because for any pair of elements $x \in M(X)$ and $y \in M(Y)$ there can
be many functions $g \in M(Y^X)$ such that $g(x) = y$.

A \emph{raw flow graph} is finally defined to be a morphism in the category of
elements of $M_{\cat L}$. Likewise, a \emph{semantic flow graph} is a morphism
in the category of elements of $M_{\cat C_*}$. Note that the models of $\cat L$,
$\cat C$, and $\cat C_*$ are conceptual devices; we do not actually construct a
denotational semantics for the language $\cat L$ or the ontology $\cat C$.
Instead, the program analysis system observes a \emph{single} computation and
produces a \emph{single} morphism $f$ in the category of elements of
$M_{\cat L}$. By the construction of the model $M_{\cat C_*}$, applying the
annotation functor $F_*: \cat L \to \cat C_*$ to this morphism $f$ yields a
morphism $F_*(f)$ belonging to the category of elements of $M_{\cat C_*}$.

In summary, semantic enrichment amounts to applying the annotation functor in
the category of elements. The expansion stage simply computes the functor, while
the contraction stage is entirely syntactical. Contraction computes a new
syntactic expression for the expanded morphism by grouping together morphisms
that are not images under the inclusion $\iota_*: \cat{C} \to \cat{C}_*$ of the
ontology category.

\section{Notes and references}

Most of the content of this chapter recapitulates previous work by the author
and collaborators, principally from \cite{patterson2018} and its abridgment in
\cite{patterson2018demo}. The discussion of implicit conversion is here extended
from cartesian categories to cartesian closed categories.

\paragraph{Knowledge representation and program analysis}

The raw flow graphs forming the input to semantic enrichment are products of
computer program analysis, as described in the previous chapter. The history of
artificial intelligence is replete with interactions between knowledge
representation and program analysis. In the late 1980s and early 1990s,
automated planning and ruled-based expert systems figured in ``knowledge-based
program analysis'' \cite{johnson1985,harandi1990,biggerstaff1994}. Other early
systems were based on description logic \cite{devanbu1991,welty2007} and graph
parsing \cite{wills1992}. These systems were designed to help software
developers maintain large codebases (exceeding, say, a million lines of code) in
specialized industrial domains like telecommunications. In data science, the
code tends to be much shorter, the control flow more linear, and the underlying
concepts better defined. Our methodology for combining program analysis with
knowledge representation is accordingly quite different from those of the older
literature.

\paragraph{Ontologies for data science}

There already exist several ontologies and schemas related to data science, such
as the Statistics Ontology (STATO), an OWL ontology about basic statistics
\cite{gonzalez-beltran2016}; the Data Mining OPtimization Ontology (DMOP), an
OWL ontology for the data mining process \cite{keet2015}; the Predictive
Modeling Markup Language (PMML), an XML schema for data mining models
\cite{guazzelli2009}; and ML-Schema, a schema for data mining and machine
learning workflows developed by a W3C community group \cite{publio2018}. The
Data Science Ontology distinguishes itself from previous efforts by attempting
to systematically code for data analysis. This does not appear to be a design
criterion of the existing standards. For example, in STATO, concepts
representing statistical methods may have designated inputs and outputs, but
they are too imprecisely specified to map onto code, among other obstacles. PMML
is a purely static format, designed for serializing fitted models. To
successfully model computer programs, special attention must be paid to the
algebraic and logical structure of programs. For this reason, the ontology
language of the Data Science Ontology is based on cartesian closed categories,
whereas the description logic of OWL more closely resembles bicategories of
relations \cite{patterson2017relations}.

\paragraph{Ontology languages and programming languages}

The cartesian closed categories forming the basis of the ontology language are
deeply connected to the simply typed lambda calculus (see Notes to
\cref{ch:category-theory}). The latter system is a fundamental model of
functional computing. Implicit conversion, comprising the other half of the
ontology language, is a form of \emph{ad hoc polymorphism} well known in
programming language theory. Building on Goguen's work on order-sorted algebras
\cite{goguen1978,goguen1992}, Reynolds gave a categorical treatment of implicit
conversion for an Algol-like language \cite{reynolds1980}. A brief textbook
account of Reynolds' work appears in \cite[\S 3.2]{pierce1991}. The approach to
implicit conversion in this chapter differs from Reynold's in being entirely
algebraic, not dependent on any particular syntax or type theory. The algebraic
formulation leads to shorter definitions and an easy treatment of implicit
conversion for product and function types.

\chapter{Conclusion}
\label{ch:conclusion}

Two contributions towards the digitization and systematization of data analysis
have figured in this dissertation. In
\cref{ch:algebra-statistics,ch:zoo-statistics}, the notion of a statistical
model, first formalized by the statistical decision theorists of the early
twentieth century, is supplemented by that of a \emph{statistical theory}.
Statistical models are then reinterpreted as models of statistical theories,
where the word ``model'' assumes its meaning in mathematical logic. Some
consequences and many examples of this change in perspective are presented.
Statistical models possess a notion of model homomorphism that clarifies and
generalizes the classical account of symmetry in theoretical statistics.
Morphisms between statistical theories formalize commonly occurring
relationships in statistics, such as generalization and specialization of
statistical methods, null hypotheses and other model containments, and
extensions of the parameter and sample spaces. Furthermore, morphisms of
statistical theories induce model migration functors between the corresponding
categories of statistical models.

In the second part of the dissertation, comprising
\cref{ch:program-analysis,ch:semantic-enrichment}, a software system for
building semantics models of data analyses is designed and implemented. Data
analyses in the form of Python or R scripts are subjected to computer program
analysis, yielding a record of the data flow during the data analysis. This raw
flow graph, expressed in the vernacular of a particular programming language and
set of libraries, is then partially translated into a semantic flow graph,
expressed in a controlled vocabulary that is independent of any particular
programming language or library. Semantic enrichment is enabled by the Data
Science Ontology, a nascent ontology about statistics, machine learning, and
computing on data. The ontology language and the semantic enrichment process are
formalized using category-theoretic methods.

\section{Limitations and future work}
\label{sec:future-work}

Despite their placement in the text, the Data Science Ontology and associated
software predate the development of the algebra of statistical theories and
models. The ontology therefore does not incorporate structural information about
statistical models and their relationships. For example, in the ontology, the
concept of a linear model is related to that of a generalized linear model by an
``is-a'' morphism (an implicit conversion), but not in the detailed way
specified by the statistical theory morphism $\cat{GLM}_n \to \cat{LM}_n$ from
\cref{sec:glm}. The inability of the ontology language to express such
relationships was, in fact, an early impetus for the author to develop the
formalism of statistical theories.

A worthwhile future project would bring the two threads back together. Concepts
in the Data Science Ontology representing statistical models would be equipped
with statistical theories and the links between them extended to theory
morphisms. Observed data and fitted model parameters, already collected by the
program analysis software, would be embedded in models of statistical theories.
To do this accurately, parsers must be written for domain-specific modeling
languages such as the ``model formulas'' in the R language or the Patsy package
\cite{chambers1993a,zeileis2010,patsy2018}. Statistical theories could even be
taken as the mathematical foundation for a probabilistic programming language,
realized by code generators and parsers for existing languages like Stan
\cite{carpenter2017}. An ontology augmented with statistical theories might also
serve as a pedagogical resource for statistical modeling, of a very different
style than a conventional textbook.

Further investments in software and knowledge engineering are needed to
transition the program analysis software and Data Science Ontology from research
prototypes to production systems. Limitations of the program analysis systems
for Python and R have been discussed in \cref{ch:program-analysis}. Most of them
could be overcome through additional engineering effort. The more pressing
question is how to scale the ontology's concepts and annotations so as to cover
a reasonably broad class of statistical methods and software. The prospects for
automating the codification of concepts seem dim, but a combination of natural
language processing and static program analysis might plausibly allow the
annotation of library code to be automated, at least partially. That would be a
significant advance, as it is the annotation of code that is ultimately the most
burdensome.

Both mathematically and statistically, the investigation undertaken here into
the algebra of statistical theories and models is the only barest beginning of a
structuralist conception of statistics. Statisticians typically distinguish
between (1) the specification of a statistical model, (2) the method of
estimating the model's parameters, and (3) the algorithm for computing the
estimator. Only the first of these is addressed by the formalism of statistical
theories and models. Moreover, within this division, the examples of
\cref{ch:zoo-statistics} are selected from among the most fundamental of
statistical models. The formalism should be tested against a wider range of
statistical models, which may reveal the need for additional structure within
statistical theories. Just as, from the pluralist standpoint of categorical
logic, there is not a single kind of logical theory, there is no reason to
expect there to be a single kind of statistical theory.

Besides the introduction of randomness, the most essential difference between
logical and statistical models is the concept of fitting a statistical model to
data, which seems to have no counterpart in mathematical logic. Indeed, in
statistics, the selection of estimators and algorithms is no less important than
the specification of models. In Bayesian statistics, there is only one method of
fitting a model---Bayesian inference---but in frequentist statistics, a single
model may be fit by many different methods. A linear model, for instance, may be
fit by ordinary least squares (the maximum likelihood estimator under i.i.d.\
normal errors), but also by ridge regression ($\ell_2$ regularization), the
lasso ($\ell_1$ regularization), the elastic net (a mixture of $\ell_1$ and
$\ell_2$ regularization), and least absolute deviations ($\ell_1$ objective),
among countless other methods. As this example illustrates, frequentist
estimators, when they cannot be written in closed form, are usually defined as
solutions to optimization problems. A fuller account of statistical models would
likely make contact with convex analysis and mathematical optimization.

Furthermore, both Bayesian and frequentist statistics depend upon efficient
algorithms for fitting models. The distinction between estimators and algorithms
is especially important when the algorithm is not guaranteed to converge to the
intended estimator, as often happens in modern high-dimensional statistics and
machine learning. Thus, statistical computing constitutes another broad
direction for extending the algebra of statistical theories and models.

The development of statistical theories in
\cref{ch:algebra-statistics,ch:zoo-statistics} has emphasized examples over
theorems, leaving many natural mathematical questions unanswered. Little has
been said about the algebraic properties of the 2-category of statistical
theories, theory morphisms (strict, lax, or colax), and transformations of
theory morphisms or about the properties of categories of models of statistical
theories. The conditions under which model migration functors have left or right
adjoints should also be determined.

Of both mathematical and practical interest is a formal way of composing
statistical theories and models, the absence of which has been felt throughout
\cref{ch:zoo-statistics} but especially in \cref{sec:hierarchical-lm} on
hierarchical models. Here is one possible approach to composing theories. For
simplicity, restrict attention to statistical theories $(\cat{T},p)$ where
$\cat{T}$ is not just a strict symmetric monoidal category but a colored PROP,
meaning that its monoid of objects is freely generated.\footnote{All statistical
  theories presented in this text indeed have colored PROPs as underlying
  categories.} A statistical theory $(\cat{T},\, \theta \xrightarrow{p} x)$ with
supply $\cat{P}$, where $\theta = \bigotimes_{i=1}^m \theta_i$ and
$x = \bigotimes_{i=1}^n x_i$ are products of object generators
$\theta_1,\dots,\theta_m$ and $x_1,\dots,x_n$, will have domain
$(\cat{P}_{\theta_1},\dots,\cat{P}_{\theta_m})$ and codomain
$(\cat{P}_{x_1},\dots,\cat{P}_{x_n})$. Composition is defined by
\begin{equation*}
  (\cat{T},\, \theta \xrightarrow{p} x) \cdot
  (\cat{S},\, x' \xrightarrow{q} y)
    := (\cat{T} +_x \cat{S},\, \theta \xrightarrow{p \cdot q} y),
\end{equation*}
where $\cat{T} +_x \cat{S}$ denotes the pushout of $\cat{T}$ and $\cat{S}$
identifying $x_i = x_i'$ for every $i=1,\dots,n$. More simply, a monoidal
product is defined by
\begin{equation*}
  (\cat{T},\, \theta \xrightarrow{p} x) \otimes
  (\cat{S},\, \phi \xrightarrow{q} y)
    := (\cat{T} + \cat{S},\,
        \theta \otimes \phi \xrightarrow{p \otimes q} x \otimes y),
\end{equation*}
using the coproduct $\cat{T} + \cat{S}$ of $\cat{T}$ and $\cat{S}$, and the
monoidal unit is the terminal theory from \cref{def:initial-terminal-theories}.
Unless one passes to isomorphism classes of theories, both the composition and
product will be non-strict. A rigorous construction of the resulting
higher-categorical structure is left to future work.

\section{Outlook: statistics and the scientific method}
\label{sec:outlook}

Wherever there is a large gap between the scientific method, as commonly
understood and practiced, and our best theoretical account of the method, there
is an opportunity to advance science by fitting the theory and practice more
closely to each other. Properly executed, the mathematical and statistical
reconstruction of science is not an exercise in empty formalism. It is a means
of improving the efficacy of science by eliminating errors resulting from
imprecise thinking and opening new ways of understanding the world, grounded in
new methodology. At its most successful, ideas that were once the exclusive
province of philosophy are transformed into actionable scientific methodology.
For example, this process is currently underway in the field of causal
inference, where the ancient idea of causality is now being operationalized by
statistical and computational methods.

An enormous gap presently exists between scientific knowledge, as it is
conceived by scientists and philosophers of science, and the conception of
scientific knowledge implicit in mathematical statistics. According to a
statistical paradigm established in the early twentieth century, scientific
inference is performed by formally stating null hypotheses within statistical
models and then testing them against observed data using statistical decision
procedures. But is this truly how science advances---one rejected null
hypothesis, one purported falsification at a time? Taking the paradigm
completely literally would suggest that scientific knowledge is nothing more or
less than the sum total of all rejected null hypotheses, a view that hardly any
scientist or statistician could seriously entertain. So, without denying that
hypothesis testing has valid uses, such as in screening to identify promising
future studies, it seems plain that scientific knowledge has a far more
intricate structure than a na\"ive interpretation of statistical hypothesis
testing would suggest.

The chief value of statistics lies in the construction of statistical models
that usefully, if imperfectly, explain and predict natural phenomena. However,
the statistical models chosen for a particular experiment or study do not exist
in a vacuum; they are motivated by, sometimes even directly derived from, a
larger body of scientific knowledge. It is no simple matter to say exactly what
this knowledge consists of, but it surely involves scientific theories and
models, as well as experimental designs and models of experiments. Statistical
theories and models, making direct contact with experimental data, sit at the
bottom of a hierarchy of increasingly abstract and general scientific theories
and models. The ultimate aim of science is not just to adequately model a single
experiment, with its specific set of experimental conditions, but to map out the
range of conditions under which a general theory is empirically adequate. When
this generalizability cannot be achieved, science is useless, without
explanatory or predictive power. Thus, a statistics that is well matched to the
aims of science would formalize the propagation of statistical inference up the
hierarchy of scientific theories and models.

While the idea of a hierarchy of scientific models has a long history in the
philosophy of science, going back at least to Patrick Suppes \cite{suppes1966},
it has had no discernible impact on statistical methodology. One possible
explanation for this is that implementing the idea in practice would require
effective computational representations of scientific and statistical models, as
well as of experimental designs and models of experiments. Such considerations
only lend further support to the arguments made in the Introduction for
digitizing science. However, on a more basic level, the question of \emph{how}
statistical models connect to and support scientific theories is still too
poorly understood to translate into statistical methodology. Making rigorous
sense of the network of theories and models in science and statistics is an
inherently interdisciplinary project, which ought to involve mathematicians,
statisticians, philosophers, computer scientists, and domain scientists from
across the natural and social sciences. Statistics has always justified itself
through its service to science by rigorizing the scientific method. A future
statistics, more strongly connected to \emph{all} the elements of scientific
knowledge, would better serve this essential purpose.

\printbibliography[heading=bibintoc]

\end{document}